\documentclass[leqno,oneside,a4paper,final]{book}
\pdfoutput=1
\usepackage{./Latexstyle/koensstyle}

\author{}
\title{\Huge{Some computations with the $\sF$-homotopy limit spectral sequence}}
\date{}

\fxsetup{status=final}
\begin{document}

\frontmatter
\begin{titlepage}
	\centering
	{\Huge Some computations \\ with the \\
    $\sF$-homotopy limit spectral sequence \par}
	\vspace{3cm}
	{\sc 
    DISSERTATION ZUR ERLANGUNG DES DOKTORGRADES\\
DER NATURWISSENSCHAFTEN (DR.\ RER.\ NAT.)\\
DER FAKULT\"AT F\"UR MATHEMATIK\\
DER UNIVERSIT\"AT REGENSBURG\\
    \par}
	\vspace{3cm}
    {vorgelegt von \\
        Koenraad van Woerden \\
    aus Groningen, die Niederlande \\
im Jahr 2017}
\end{titlepage}
\null
\vfill
{
    \thispagestyle{empty}
\begin{table}[H]
    \raggedright
    \begin{tabular}{ll}
    Promotionsgesuch eingereicht am: &  27.~April 2017 \\
                                     &  \\
    Die Arbeit wurde angeleitet von: &  Dr.~Justin Noel  \\
     & \\
    Pr\"ufungsausschuss:  & \\
      Vorsitzender: & Prof.~Dr.~Harald Garcke \\
      1.~Gutachter: & Dr.~Justin Noel \\
      2.~Gutachter: & Prof.~Dr.~Niko Naumann \\
      weiterer Pr\"ufer: & Prof.~Dr.~Ulrich Bunke
    \end{tabular}
\end{table}
}
\clearpage
\tableofcontents
\mainmatter
\chapter{Introduction}
\section{Nilpotence in Borel equivariant cohomology}
For a finite group $G$, consider a cohomology class $u \in H^*(BG; \Ffield_p)$ that restricts to zero on all elementary abelian $p$-subgroups (groups of the form $(\Integers/p)^{\times l}$). It is a theorem of Quillen that $u$ is nilpotent. In fact, Quillen showed
\begin{Theorem}[{{\cite[Thm.\ 7.1]{quillen71spectrum}}}]
    \label{thm:quillenfp}
    The map
\begin{equation}
    \widetilde{\res} \colon H^*(BG;\Ffield_p) \to \lim_{E\subset G \mathrm{\ el.\ ab.\ }p\text{-}\mathrm{gp}.} H^*(BE;\Ffield_p)
    \label{eq:quillenfp}
\end{equation}
where the maps in the indexing category for the limit are given by restricting along subgroups and conjugation, is a uniform $\calF_p$-isomorphism (\cite[Thm.\ 7.1]{quillen71spectrum}, \cite[p. 556]{quillen71spectrum}), which means that there is an $n$ such that
\begin{enumerate}
    \item Every $u \in \ker(\widetilde{\res})$ satisfies $u^n = 0$.
    \item Every $v \in \lim_E H^*(BE; \Ffield_p) \setminus \Im(\widetilde{\res})$ satisfies $v^{p^n} \in \Im(\widetilde{\res})$. 
\end{enumerate}
\end{Theorem}
This result in group cohomology is a consequence of the following landmark theorem by Quillen in Borel equivariant cohomology, which is the main theorem of the same paper. For a $G$-space $X$, we write $H_G^*(X; \Ffield_p)$ for the mod-$p$ Borel equivariant cohomology of $X$.
\begin{Theorem}[{{\cite[Thm.\ 6.2]{quillen71spectrum}}}]
    \label{thm:quillenborfp}
    For $X$ any paracompact $G$-space of finite cohomological dimension, the map
    \begin{equation}
      \widetilde{\res} \colon H_G^*(X;\Ffield_p) \to \lim_{E\subset G \text{ el.\ ab.\ $p$-gp.}} H_E^*(X;\Ffield_p)
\end{equation}
is a uniform $\calF_p$-isomorphism.
\end{Theorem}
These results led to a lot of structural results in group cohomology. Quillen himself immediately deduced
\begin{Corollary}[{{\cite[Cor.\ 7.8]{quillen71spectrum}}}]
    The Krull dimension of $H^*(BG;\Ffield_p)$ equals the rank of the maximal elementary abelian $p$-subgroup of $G$.
\end{Corollary}
Results directly building on \cref{thm:quillenborfp} include a theorem of Duflot on the depth of $H^*(BG;\Ffield_p)$ (\cite[Thm. 1]{duflot81}), a theorem on the complexity of $kG$-modules by Alperin-Evens (\cite{alperinevens81}), Benson's description of the image of the transfer map (\cite[Thm.\ 1.1]{benson93}), and a theorem on the depth of group cohomology rings by Carlson (\cite[Thm.\ 2.3]{carlson95}). These results, and Quillen's original result all indicate the importance of the elementary abelian $p$-subgroups in Borel equivariant $\Ffield_p$-cohomology in general, and group cohomology with $\Ffield_p$-coefficients in particular.

It is natural to ask what one can say about the $n$ in \cref{thm:quillenfp}. One approach to this question is in \cite{mnn}, where the map (\ref{eq:quillenfp}) is realized as the edge homomorphism of a homotopy limit spectral sequence, called the $\sF$-homotopy limit spectral sequence,
\begin{equation}
    E_2^{s,t} = \sideset{}{^s} \lim_{\sO(G)_{\sF}^{\op}} H_H^t(X; \Ffield_p) \Rightarrow H_G^*(X;\Ffield_p)
    \label{intro:fholimss1}
\end{equation}
converging strongly to the target (\cite[Prop. \ 2.24]{mnn}), where $X$ can be any $G$-space, and $\sF$ any family of subgroups of $G$ which contains at least the family $\sE_{(p)}$ of elementary abelian $p$-subgroups. The indexing category is the subcategory of the orbit category $\sO(G)$ spanned by the orbits $G/H$ with $H$ in the family $\sF$. The key property of the $\sF$-homotopy limit spectral sequence is that it collapses at a finite page with a horizontal vanishing line (\cite[Thm.\ 2.25]{mnn}). 

This implies that every computation with the $\sF$-homotopy limit spectral sequence is a finite one. Moreover, in many concrete situations we can establish a bound on the height of the horizontal vanishing line, and a bound on which page it will appear. Besides implying when the $\sF$-homotopy limit spectral sequence will have collapsed, it can also be used to deduce differentials. This was illustrated in the computation of the cohomology of the quaternion group of order 8 in \cite[Ex.\ 5.18]{mnn}. The first prinicipal goal of this thesis is to demonstrate the computational utility of the $\sF$-homotopy limit spectral sequence and the use of the vanishing line in determining differentials. This will be done by using it to compute the cohomology of all 2-groups up to order 16. 

Varying $X$ over all $G$-spectra, this horizontal vanishing line turns out to have a uniform bound in height (\cite[Prop.\ 2.26]{mnn}). The minimal upper bound of this height is one of the equivalent definitions of the $\sE_{(p)}$-exponent $\exp_{\sE_{(p)}}\underline{H\Ffield_p}_G$. In practice one can often determine this $\sE_{(p)}$-exponent, and this leads to a quantified version of \cref{thm:quillenfp}, because one has $n \leq \exp_{\sE_{(p)}}\underline{H\Ffield_p}_G$ (\cite[Thm.\ 3.24, Rem.\ 3.26]{mnn}). The identification of $\sE_{(2)}$-exponents for 2-groups is the second principal goal of this thesis, and leads to the main theorem:
\begin{Theorem}
    Let $G$ be a finite group with a 2-Sylow subgroup of order $\leq 16$, let $X$ be any $G$-space, and let $I$ be the kernel of
    \begin{equation}
            \widetilde{\res} \colon H_G^*(X;\Ffield_2) \to \lim_{E\subset G \text{ el.\ ab.\ $2$-gp.}} H_E^*(X;\Ffield_2).
    \end{equation}
    Then $I^4 = 0$. Moreover, if $u$ is not in the image of $\widetilde{\res}$, then $u^{8}$ is.
\end{Theorem}
This follows from combining \cite[Thm.\ 3.24]{mnn} and \cref{lem:psylowhfpexp} with the upper bounds on the exponents obtained in this thesis and summarized in \cref{table:groupsleq16}. For a description of the groups appearing in the table we refer to \cref{sec:classif2gps}. The first column lists the groups of order $\leq 16$, the second column the exponent or an interval in which the exponent lies, and the last column gives a forward reference for the claim.
\begin{table}[H]
    \centering
     \begin{tabular}{r|r|r}
         $G$ & $\exp_{\sE_{(2)}} \underline{H\Ffield_2}_G$ & Reference \\
         $e$ & 1 & \cref{prop:expelab} \\
         \hline
         $C_2$ & 1 & \cref{prop:expelab} \\
         \hline
         $C_2 \times C_2$ & 1 & \cref{prop:expelab} \\
         $C_4$ & 2 & \cref{prop:expelab} \\
        $C_2^{\times 3}$ & 1 & \cref{prop:expelab} \\
        $C_2 \times C_4$ & 2 & \cref{prop:expelab} \\
        $C_8$ & 2 & \cref{prop:expelab} \\
        $D_8$ & 2 & \cref{prop:d2nexp} \\
        $Q_8$ & 4 & \cite[Ex.\ 5.18]{mnn} \\
        \hline
        $C_2^{\times 4}$ & 1 & \cref{prop:expelab}  \\
        $C_2^{\times 2} \times C_4$ & 2 & \cref{prop:expelab}  \\
        $C_4 \times C_4$ & 3 & \cref{prop:expelab} \\
        $C_8 \times C_2$ & 2 & \cref{prop:expelab}  \\
        $C_{16}$ & 2 & \cref{prop:expelab}   \\
        $D_{16}$ & 2 & \cref{prop:d2nexp}  \\
        $Q_{16}$ & 4 & \cref{prop:q2ne2exp}  \\
        $ SD_{16} =  C_8 \overset{3}{\rtimes} C_2$ & $[3,4]$ & \cref{prop:sd16e2exp}  \\
        $ M_{16} = C_8 \overset{5}{\rtimes} C_2$ & $[3,4]$ & \cref{prop:m16e2exp}  \\
        $D_8 \ast C_4$ & $4$ & \cref{d8cpc4:e2exp} \\
        $C_4 \rtimes C_4$ & $[3,4]$ & \cref{prop:c4sdc4e2exp} \\
        $(C_4 \times C_2) \overset{\psi_5}{\rtimes} C_2$ & 2 & \cref{prop:c4xc2sdc2exp} \\
        $Q_8 \times C_2$ & 4 & \cref{prop:q8xc2e2exp} \\
        $D_8 \times C_2$ & 2 & \cref{prop:d8xc2e2exp} \\
       \end{tabular}
    \caption{The $\sE_{(2)}$-exponents of the groups of order $\leq 16$.}
    \label{table:groupsleq16}
\end{table}
For specific 2-Sylow subgroups of order $\leq 16$ it is possible to improve the previous theorem by using the upper bound on the relevant exponent.

The upper bounds in \cref{table:groupsleq16} have been obtained by either using the projective bundle theorem as in \cite[Ex.\ 5.18]{mnn} or by using Euler classes and group cohomology computation with (\ref{intro:fholimss1}) as in \cite[Ex.\ 5.22]{mnn}. We also refer to \cref{projbund:prop1} for the former and \cref{cor:hfpeulerclass} for the latter method.

The lower bounds on the exponents have been obtained from the computations of group cohomology with (\ref{intro:fholimss1}). These then serve a threefold purpose: firstly, they illustrate the computational power of the $\sF$-homotopy limit spectral sequence, secondly they give us lower bounds on $\sE_{(2)}$-exponents, and lastly they can give us upper bounds on exponents. Not all our computations of group cohomology with (\ref{intro:fholimss1}) are done with the family $\sE_{(2)}$, sometimes we chose a family which was strictly larger, if it lead us to a more manageable computation.

The structure of this thesis is as follows. In \cref{ch:fhtyss} we collect some results from \cite{mnn} about the $\sF$-homotopy limit spectral sequence and $\sF$-exponents.  Some lemmas that can aid in the computation of exponents will be proven in \cref{ch:explems}. In \cref{ch:gpcohfhtylimss}, which is the main chapter, we determine upper and lower bounds on $\exp_{\sE_{(2)}} \underline{H\Ffield_2}$ for various small 2-groups, and use the $\sF$-homotopy limit spectral sequence to compute group cohomology. 

In \cref{ch:kthy}, which is mostly unrelated to the previous chapters, we consider complex equivariant $K$-theory $KU$, which also admits $\sF$-homotopy limit spectral sequences, but one needs to replace the family $\sE_{(2)}$ by the family $\sC$ to get a horizontal vanishing line on a finite page (\cite[Prop.\ 5.6]{mnn}). We prove lower bounds on $\exp_{\sC} KU$ for $E$-equivariant $K$ theory, where $E$ is an elementary abelian $2$-group. For primes $p \neq 2$ we conjecture a lower bound, which we verify for several small cases using Sage in \cref{ch:sage}.

\section{Acknowledgements}
I would like to thank my supervisor Justin Noel in general for teaching me practically everything that I know about stable equivariant homotopy theory, and in particular for the help with the statement and proofs of the exponent lemmas in \cref{sec:exponentlemmas}.

I would like to thank Akhil Mathew for a detailed explanation of the use of the projective bundle theorem in \cite[Ex.\ 5.18]{mnn}.

This thesis was written at SFB 1085 Higher Invariants, and I would like to thank all my colleagues for the great working atmosphere.

I would like to thank Kim and Christoph for all the nice discussions we had in their office, whose door was always open, for mathematics or otherwise.

I would like to thank Johannes for making more coffee than anyone else in the SFB, inventing analog divisibility testers, and helping me with the formalities of handing in a thesis.

I would like to thank Oriol for patiently answering all my questions about stable homotopy theory, and discovering and taking me to all the hidden cultural events in Regensburg.

I would like to thank Florent for initiating a shared effort to learn German, giving history lessons and bringing pastry to the office.

Lastly I would like to thank my parents and my brother, for their interest, support, and advice.

\section{Notations and conventions}
Throughout, $G$ will denote a finite group.

Let $G$ be a finite group, and $g_1,\ldots,g_n \in G$ such that $\overline{g_1},\ldots,\overline{g_n} \in G^{\ab}$ are a basis for $G^{\ab} \otimes \Ffield_p$ as an $\Ffield_p$-vector space. We then denote the cohomology classes in $H^1(BG;\Ffield_p)$ dual to $\overline{g_i} \in G^{\ab}$ by $\delta_{\overline{g_i}}$.

For $G$ a finite group, we denote by $\rmS_G$ the $\infty$-category obtained from the simplicial model category of topological $G$-spaces with the level model structure. We denote by $\Sp_G$ the $\infty$-category obtained from the symmetric monoidal simplicial model category of orthogonal $G$-spectra with the positive stable model structure (see \cite[Def.~5.2]{mandellmay}).

The cyclic group of order $n$ we denote by $C_n$.

\chapter{The $\sF$-homotopy limit spectral sequence}
\label{ch:fhtyss}
\section{Introduction}
In this chapter we recall the $\sF$-homotopy limit spectral sequence as constructed in \cite{mnn}, as well as some of the properties of this spectral sequence. For details of the constructions and for proofs of the statements, we refer to \cite{mnn} and \cite{mnnnd}.

This chapter is organized as follows.  We first recall the notion of a family of subgroups and the orbit category, which we will then use to state the $\sF$-homotopy limit spectral sequence. One of its key properties is the fact that it has a horizontal vanishing line at a finite page. The height of this vanishing line and the page on which it appears can be bounded by the $\sF$-exponent, which we then describe next. Finally we show that for specific families $\sF$, the $\sF$-homotopy spectral sequence (for $X = \pt$) is isomorphic to a Lyndon-Hochschild-Serre spectral sequence in group cohomology.

\section{Families of subgroups}
The $\sF$-homotopy limit spectral sequence depends on a family of subgroups $\sF$. We recall this notion and describe some examples that will feature in this thesis.
\begin{Definition}
    A family $\sF$ of subgroups of $G$ is a non-empty set of subgroups closed under taking subgroups and conjugation.
\end{Definition}
\begin{Examples}
    The following are examples:
    \begin{enumerate}[1.]
        \item The elementary abelian $p$-subgroups (groups isomorphic to $C_p^{\times n}$) with $p$ a prime. This family is denoted $\sE_{(p)}$.
        \item The family $\sC$ of cyclic subgroups.
        \item The family $\sA$ of abelian subgroups.
        \item The family $\All$ of all subgroups.
        \item For a family $\sF$ of subgroups of $G$, and $N$ any normal subgroup of $G$, let $\sF_N$ be the subset of groups in $\sF$ contained in $N$. Then $\sF_N$ is a family of subgroups of $N$, and also of $G$.
    \end{enumerate}
\end{Examples}

\section{Orbit categories}
The $\sF$-homotopy limit spectral sequence has an $E_2$-term the derived limit out of the opposite of a subcategory of the orbit category, of which we now recall the definition.
\begin{Definition}
    For $G$ a finite group, the orbit category $\sO(G)$ is the subcategory of the category of $G$-sets spanned by the transitive $G$-sets $G/H$, where $H$ ranges over the subgroups of $G$. For $\sF$ a family of subgroups of $G$, the category $\sO(G)_{\sF}$ is defined to be the full subcategory of $\sO(G)$ spanned by the orbits $G/H$ with the $H$ ranging over the subgroups in $\sF$.
\end{Definition}

\section{The $\sF$-homotopy limit spectral sequence}
We now state the spectral sequence which is the main object of study of this thesis. We will use the following notation.
\begin{Notation}
    We denote by $H_G^*(X; \Ffield_2)$ the mod-2 Borel equivariant cohomology of a $G$-space $X$, which is by definition the ordinary cohomology $H^*(X \times_G EG; \Ffield_2)$ of the homotopy orbits of $X$. The $G$-spectrum representing this cohomology theory we denote by $\underline{H\Ffield_2}_G$, or by $\underline{H\Ffield_2}$ if the group is clear from the context.
    \label{not:borel}
\end{Notation}
\begin{Theorem}[{{\cite[Def.\ 2.23, Thm.\ 2.25, Cor.\ 3.6, Prop.\ 5.16]{mnn}}}]
    Let $G$ be a finite group, let $X$ be a $G$-spectrum, and let $\sF$ be family of subgroups containing $\sE_{(2)}$. Then there is a spectral sequence, called the $\sF$-homotopy limit spectral sequence, given by
    \begin{equation}
        E_2^{s,t}(X) = \sideset{}{^s}\lim_{G/H \in \sO(G)_{\sF}^{\op}}H_H^t(X; \Ffield_2) \Rightarrow H_G^{s+t}(X;\Ffield_2).
            \label{eq:fholimss}
    \end{equation}
    This spectral sequence converges strongly to its target. Moreover, it collapses at a finite page with a horizontal vanishing line: there are $N, l \geq 1$ such that $E_{N+1}^{s,*}(X) = 0$ for all $s \geq l$.
    \label{thm:fholimss}
\end{Theorem}
In \cref{ch:gpcohfhtylimss}, we will illustrate the computational use of the $\sF$-homotopy limit spectral sequence by using it when $X = \pt$, is a point, to compute the cohomology of all 2-groups up to order 16.
\begin{Remark}
    We abbreviate $E_2^{s,t} = E_2^{s,t}(\pt)$ for the $\sF$-homotopy limit spectral sequence for $X = \pt$.
\end{Remark}

The $\sF$-homotopy limit spectral sequence can be constructed for any generalized equivariant cohomology theory $E$ (\cite[Def.\ 2.23]{mnn}). The case considered in \cref{thm:fholimss} is $E = \underline{H\Ffield_2}_G$. For different choices of $E$ a statement about a vanishing line analogous to the one in \cref{thm:fholimss} also holds, but one needs to change the family $\sE_{(2)}$ to the derived defect base of $E$ (\cite[remark after def.\ 1.3, Thm.\ 2.25]{mnn}). We refer to \cite{mnn} for the details.

The vanishing line mentioned in \cref{thm:fholimss} and the page on which it appears turn out to have an upper bound uniform in the spectrum $X$. To describe it, we introduce the notion of $\sF$-nilpotence in the next section.

\section{$\sF$-nilpotence}
We recall the notion of $\sF$-nilpotence from \cite{mnn}, for which we first need to recall the following space.
\begin{Definition}
    Let $\sF$ be a family of subgroups of $G$. Then the \textbf{universal $\sF$-space} $E\sF$ is the $G$-space given by the homotopy colimit
    \begin{equation}
        E\sF = \hocolim_{\sO(G)} G/H.
    \end{equation}
\end{Definition}
The space $E\sF$ is, up to $G$-equivalence, characterized by (\cite[Def.\ II.2.10]{lms})
\begin{equation}
    E\sF^H \simeq \begin{cases}
        \varnothing & \text{if $H \not\in \sF$,} \\
        \pt & \text{if $H \in \sF$.}
    \end{cases}
\end{equation}
Using the space $E\sF$, we can give one of the equivalent definitions of $\sF$-nilpotence.
\begin{Definition}[{{cf. \cite[Def.\ 1.3]{mnn}}}]
    Let $M$ be a $G$-spectrum. Then $M$ is said to be $\sF$-nilpotent if there is an $n$ such that $M$ is a retract of $F(\sk_{n-1}E\sF_+,M)$. The minimal $n \geq 0$ for which this holds is called the $\sF$-exponent of $M$, and denoted $\exp_{\sF}M$. 
\end{Definition}
Being $\sF$-nilpotent is a strong conditiont, it implies for example the following.
\begin{Proposition}[{{\cite[Prop.\ 2.8]{mnn}}}]
    Let $M$ be $\sF$-nilpotent. Then $M$ is $\sF$-complete and $\sF$-colocal, that is, the $\sF$-completion map $M \to M(E\sF_+,M)$ and the colocalization map $M \wedge E\sF_+ \to M$ are weak equivalences.
\end{Proposition}
The primary case of interested for us is $M = \underline{H\Ffield_2}_G$.
\begin{Proposition}[{{\cite[Prop.\ 5.16]{mnn}}}]
    For $G$ any group, the \\ $G$-spectrum $\underline{H\Ffield_2}_G$ is $\sE_{(2)}$-nilpotent.
\end{Proposition}
\begin{Remark}
    In fact, \cite[Prop.\ 5.16]{mnn} shows that $\sE_{(2)}$ is the minimal family $\sF$ such that $\underline{H\Ffield_2}_G$ is $\sF$-nilpotent. In the terminology of \cite{mnn} we have that $\sE_{(2)}$ is the derived defect base of $\underline{H\Ffield_2}_G$.
\end{Remark}
A trivial example is the following.
\begin{Proposition}
    Every $G$-spectrum $M$ is $\All$-nilpotent with $\exp_{\All} M \leq 1$, and $\exp_{\All} M = 0$ if and only if $M$ is contractible.
    \label{exam:allnilp}
    \label{prop:allexp}
    \label{lem:expallsubgroups}
\end{Proposition}
\begin{Proof}
    We have $E\All \simeq \pt$, and every spectrum $M$ is a retract of $F(S,M) \simeq M$. 

    A spectrum $M$ has $\All$-exponent 0 if and only if it is a retract of $F(\pt, M) \simeq \pt$, which happens if and only if $M$ is contractible.
\end{Proof}
The following two propositions are immediate with some alternative definitions of $\sF$-nilpotence, but not with the one we have given.
\begin{Proposition}[{{\cite[Prop. 6.39]{mnnnd}}}]
    A $G$-spectrum $M$ is $\sF$-nilpotent and $\sG$-nilpotent if and only if $M$ is $\sF \cap \sG$-nilpotent.
    \label{prop:intersectnilp}
\end{Proposition}
\begin{Proposition}[{{\cite[Def.\ 1.3]{mnn}}}]
    If $M$ is $\sF$-nilpotent and $\sG \supset \sF$, then $M$ is $\sG$-nilpotent.
    \label{prop:supfamnilp}
\end{Proposition}
Combining \cref{exam:allnilp}, \cref{prop:intersectnilp} and \cref{prop:supfamnilp} shows that every $G$-spectrum $M$ has a minimal family $\sF$ such that $M$ is $\sF$-nilpotent (cf.\ the remark after \cite[Def. 1.3]{mnn}). 
\begin{Definition}[{{\cite[remark after Def.\ 1.3]{mnn}}}]
    This minimal family $\sF$ such that a $G$-spectrum $M$ is $\sF$-nilpotent is called the \textbf{derived defect base} of $M$.
\end{Definition}
The following is the main case of interest for us.
\begin{Proposition}[{{\cite[Prop.\ 5.16]{mnn}}}]
    For $G$ any finite group, the derived defect base of $\underline{H\Ffield_2}_G$ is $\sE_{(2)}$.
    \label{prop:derdefbaseborel}
\end{Proposition}
We can now state the uniform upper bound in the height of the horizontal vanishing of the $\sF$-homotopy limit spectral sequence and the page on which it appears.
\begin{Proposition}[{{\cite[Prop.\ 2.26, Rem.\ 2.27]{mnn}}}]
    \label{expprop}
    Let $G$ be a finite group, and $\sF \supset \sE_{(2)}$ a family of subgroups, $M$ an $\sF$-nilpotent $G$-spectrum. Then the following integers equal:
    \begin{enumerate}[(1)]
        \item The $\sF$-exponent of $M$.
        \item The minimal $N$ such that for all $G$-spectra $X$, the $\sF$-homotopy limit spectral sequence $E_*^{*,*}(X)$ admits a vanishing line of height $N$ on the $N+1$-page: $E_{N+1}^{s,*} = E_\infty^{s,*} = 0$ for all $s \geq N$.
        \item The minimal $n$ such that the canonical map \\ $F(E\sF_+,M) \simeq M \to F(\sk_{n-1}E\sF_+, M)$ admits a section.
        \item The minimal $n'$ such that there is an $(n-1)$-dimensional CW-complex $X$ with isotropy in $\sF$ such that $M$ is a retract of $F(X_+, M)$.
        \item The minimal $m$ such that the canonical map $\sk_{m-1}E\sF \wedge M \to M$ admits a retraction. \label{expprop:smash2}
        \item The minimal $m'$ such that there is an $(m'-1)$-dimensional CW-complex $X$ with isotropy in $\sF$ such that $M$ is a retract of $X_+ \wedge M$. \label{expprop:smash}
    \end{enumerate}
    Moreover, if $M'$ is any $G$-spectrum, then the existence of an integer for $M'$ as in any one of the items from (2) to (6) implies that $M'$ is $\sF$-nilpotent.
\end{Proposition}
An immediate consequence is that exponents are invariant under suspension:
\begin{Corollary}
    Let $M$ be a $G$-spectrum, $\sF$ a family and $s$ an integer. Then $M$ is $\sF$-nilpotent with $\sF$-exponent $n$ if and only if $\Sigma^s M$ is $\sF$-nilpotent with $\sF$-exponent $n$.
    \label{cor:expsuspinv}
\end{Corollary}
\begin{Proof}
    This follows from suspending the retraction sequences from \cref{expprop} (\ref{expprop:smash2}). 
\end{Proof}
We end this section by recalling from \cite{mnn} some properties of exponents that will be used in the next chapter to prove lemmas about exponents.
\begin{Proposition}
    Let $H \in \sF$. Then $G/H_+$ is $\sF$-nilpotent with $\exp_{\sF} G/H_+ = 1$.
    \label{prop:fingsetexp}
    \label{lem:expgsets}
\end{Proposition}
\begin{Proof}
    Because $F(G/H_+, G/H_+) \simeq G/H_+ \wedge G/H_+$ by self duality of $G$-sets in $\Sp_G$ (\cite[Cor.\ 6.3]{lms}), the result follows the suspending the following composite of $G$-spaces, which is the composite of the diagonal followed by the projection
    \begin{equation}
        G/H \to G/H \times G/H \to G/H
    \end{equation}
    to get a retraction
    \begin{equation}
        G/H_+ \to G/H_+ \wedge G/H_+ \to G/H_+
    \end{equation}
\end{Proof}
\begin{Proposition}[{{\cite[Cor.\ 4.15]{mnnnd}}}]
    If $M$ is an $\sF$-nilpotent spectrum and $X$ is any $G$-spectrum, then $F(X,M)$ is $\sF$-nilpotent with \\ $\exp_{\sF} F(X,M) \leq \exp_{\sF} M$.
\end{Proposition}
\fxwarning{Fix overfull boxes}
\begin{Proposition}
    If $N$ is an $\sF$-nilpotent $G$-spectrum and $M$ is any $G$-spectrum then $M \wedge N$ is $\sF$-nilpotent with $\exp_{\sF}M \wedge N \leq \exp_{\sF} N = n$.
    \label{prop:smashexp}
\end{Proposition}
\begin{Proof}
    Because $N$ is $\sF$-nilpotent there is, by \cref{expprop} (\ref{expprop:smash}), an $(n-1)$-dimensional CW-complex with isotropy in $\sF$ and a retraction
    \begin{equation}
        N \to X_+ \wedge N \to N,
    \end{equation}
    Smashing this with $M$ gives a retraction
    \begin{equation}
        M \wedge N \to X_+ \wedge M \wedge N  \to M \wedge N,
    \end{equation}
    and another application of \cref{expprop} (\ref{expprop:smash}) gives the result.
\end{Proof}
\begin{Proposition}[{{\cite[Prop.\ 4.9]{mnnnd}}}]
    \begin{enumerate}[1.]
        \item If $M$ is a retract of an $\sF$-nilpotent spectrum $N$, then $M$ is $\sF$-nilpotent and $\exp_{\sF}M \leq \exp_{\sF} N$.
        \item If $M'$ and $M''$ are $\sF$-nilpotent and $M' \to M \to M''$ is a cofiber sequence then $M$ is $\sF$-nilpotent and $\exp_{\sF} M \leq \exp_{\sF} M' + \exp_{\sF} M''$.
    \end{enumerate}
    \label{prop:retrexp}
    \label{lem:expretracts}
\end{Proposition}
\begin{Proposition}
    Let $M_{\alpha}$ be a set of $\sF$-nilpotent spectra with $\sF$-exponents bounded uniformly by $n$. Then $\bigvee_{\alpha} M_\alpha$ is $\sF$-nilpotent with $\sF$-exponent $\leq n$.
    \label{prop:wedgeexp}
\end{Proposition}
\begin{Proof}
    By \cref{expprop} (\ref{expprop:smash2}) there are retractions
    \begin{equation}
        M_\alpha \to \sk_{n-1}E\sF_+ \wedge M_\alpha \to M_\alpha.
    \end{equation}
    Wedge them together and use that smashing commutes with arbitrary colimits to obtain a retraction
    \begin{equation}
        \bigvee_\alpha M_\alpha \to \sk_{n-1}E\sF_+ \wedge \bigvee_{\alpha}M_\alpha \to \bigvee_{\alpha}M_\alpha
    \end{equation}
    and apply \cref{expprop} (\ref{expprop:smash2}) again.
    \end{Proof}
\begin{Proposition}
    Let $X$ be an $(n-1)$-dimensional $G$-CW-complex with isotropy in $\sF$. Then $X_+$ is $\sF$-nilpotent and $\exp_{\sF} X_+ \leq n$.
\end{Proposition}
\begin{Proof}
    By induction on $n$. For $n=1$, $X$ is of the form
    \begin{equation}
    \bigvee_{H \in \sF} \left(G/H \wedge \left(\bigvee_{\alpha} S^0\right)\right), 
    \end{equation} 
    which is $\sF$-nilpotent of exponent $1$ by \cref{prop:smashexp}, \cref{prop:wedgeexp} and \cref{prop:fingsetexp}.

    Assume we have the result for $n-1$, and let $X$ be an $n$-dimensional CW-complex. Let $X^{(n-1)}$ be the $(n-1)$-skeleton of $X$. Then $X$ sits in a cofiber sequence
    \begin{equation}
        X^{(n-1)} \to X \to \bigvee_{H \in \sF} \left(G/H \wedge\left(\bigvee_\alpha S^n\right)\right).
    \end{equation}
    Now $\exp_{\sF}X^{(n-1)} \leq n$ by induction, and the right hand side has $\sF$-exponent $\leq 1$ by \cref{prop:smashexp}, \cref{prop:wedgeexp} and \cref{prop:fingsetexp}.
\end{Proof}
\begin{Proposition}
    Let $X$ be a finite $G$-CW complex with isotropy in $\sF$. Then the equivariant Spanier-Whitehead dual $\sw(X_+)$ of $X_+$ is $\sF$-nilpotent, and $\exp_{\sF}\sw(X_+) = \exp_{\sF} X_+$.
    \label{prop:cwexp}
    \label{lem:cwexp}
\end{Proposition}
\begin{Proof}
    Write $n = \exp_{\sF} X_+$.
    Let $Y_+$ be an $(n-1)$ finite-dimensional $G$-CW complex with isotropy in $\sF$ such that there is a retraction
    \begin{equation}
        X_+ \to F(Y_+, X_+) \simeq \sw(Y_+) \wedge X_+ \to X_+.
    \end{equation}
    Applying $\sw(-)$ to this retraction exhibits $\sw(X_+)$ as a retraction of $Y_+ \wedge \sw(X_+)$, which has exponent $\leq n$ by \cref{prop:cwexp} and \cref{prop:smashexp}. Therefore $\exp_{\sF}\sw(X_+) \leq \exp_{\sF}X_+$, and replacing $X_+$ by $\sw(X_+)$ in this inequality shows equality.
\end{Proof}

\section{Comparison with the Lyndon-Hochschild-Serre spectral sequence}
\label{sec:lhsss}
Let $N$ be a normal subgroup of $G$, and let $\All_N$ be the family of subgroups of $N$. We then have the following identificaiton of spectral sequences.
\begin{Proposition}[{{\cite[Lem.\ A.3]{mnn}}}]
    The $\All_N$-homotopy limit spectral sequence
    \begin{equation}
        \sideset{}{^s}\lim_{\sO(G)_{\All_N}^{\op}}H^t(BH;\Ffield_2) \Rightarrow H^{s+t}(BG;\Ffield_2)
    \end{equation}
    is isomorphic to the Lyndon-Hochschild-Serre spectral sequence
    \begin{equation}
        H^s(BG/N; H^t(BN;\Ffield_2)) \Rightarrow H^{s+t}(BG;\Ffield_2).
    \end{equation}
    \label{prop:compserre}
\end{Proposition}
We will make use of this identification in all our computations. In some cases we extend the family $\sE_{(2)}$ to a family which is of the form $\All_N$, but this is not always possible because there is not always a proper normal subgroup containing all the elementary abelian subgroups.

\section{A decomposition of the $E_2$-page.}
\label{sec:absplit}
In \cref{sec:lhsss}, we saw that the $\sF$-homotopy limit spectral sequence reduces to the LHSSS in case the family $\sF$ is of the form $\sF = \All_N$, for some normal subgroup $N$ of $G$. In case the family $\sF$ is not of this form, we will compute the $E_2$-term of the $\sF$-homotopy limit spectral sequence following the strategy of \cite[App.\ B]{mnn}. In this section we describe this strategy in the abstract.
\subsection{Homotopy cofinality}
To compute the ordinary categorical limit of a diagram $I \to \CatC$ indexed by a category $I$, one can restrict along a left cofinal functor $I' \to I$. The corresponding notion for homotopy limits is homotopy left finality, which can be characterized as follows for homotopy limits indexed by ordinary categories. We have taken the dual of the cited Proposition. We denote by $N(-)$ the nerve of an ordinary category.
\begin{Proposition}[{{\cite[Prop.~4.1.1.8]{lurietop}}}]
   Let $\iota \colon I' \to I$ be a functor of ordinary categories. The following are equivalent.
   \begin{enumerate}[(1)]
       \item The functor $\iota$ is homotopy left cofinal.
       \item For $\CatC$ any $\infty$-category and $I \to \CatC$ any diagram, the induced map of homotopy limits by $\iota$ is a weak equivalence. That is, taking right Kan extensions of
           \begin{equation}
               \begin{tikzcd}
                   N(I) \arrow{r} \arrow{d} & \CatC \\
                   \pt \arrow[dashed]{ur} & 
               \end{tikzcd}
           \end{equation}
           and of 
           \begin{equation}
               \begin{tikzcd}
                   N(I') \arrow{r} \arrow{d} & \CatC \\
                   \pt \arrow[dashed]{ur} &
               \end{tikzcd}
           \end{equation}
           yields parallel dashed arrows
           \begin{equation}
               \begin{tikzcd}
                   N(I') \arrow{dr} \arrow{r} & N(I) \arrow{d} \arrow{r} & \CatC \\
                   & \pt \arrow[dashed, shift left=2pt]{ur} \arrow[dashed, shift right=2pt]{ur} & 
           \end{tikzcd}
           \end{equation}
           with a natural transformation between them. This natural transformation is an equivalence.
   \end{enumerate}
\end{Proposition}
\subsection{Homotopy cofinality in orbit categories.}
Let $\sF$ be a family of subgroups of $G$, and let $\{H_i\}$ be the maximal subgroups in $\sF$. Assume there is a subgroup $K$ of $G$ such that $K$ equals any intersection of at least two of the $H_i$:
\begin{equation}
    \bigcap_{i \in I} H_i = K,
\end{equation}
for all $I$ with $|I| \geq 2$. In particular, the set $\{H_i\} \cup \{K\}$ is closed under intersections. Because the set $\{H_i\}$ is closed under conjugations, so is $\{H_i\} \cup \{K\}$. Denote by $\sO_{\rmf}$ the subcategory of $\sO(G)_{\sF}$ spanned by the $G/J$ with $J$ in the set $\{H_i\} \cup \{K\}$. Now \cite[Prop.\ 6.31]{mnnnd} implies the following:
\begin{Proposition}
    The subcategory $\sO^{\op}_{\rmf} \to \sO(G)^{\op}_{\sF}$ is homotopy left cofinal. 
    \label{prop:intersecthofinal}
\end{Proposition}
The construction of the homotopy limit spectral sequence by Bousfield-Kan (\cite[Ch.~XI]{bousfieldkan}), applied to the diagram
\begin{align}
    \sO(G)^{\op}_{\sF} & \to \Sp, \\
    G/H_+ & \mapsto F(G/H_+,\underline{H\Ffield_2})^G,
\end{align}
yields the $\sF$-homotopy limit spectral sequence (\cite[Prop.~2.24]{mnn}).
\begin{Proposition}
    The inclusion functor $\iota \colon \sO_{\rmf} \to \sO(G)_{\sF}$ induces a morphism of homotopy limit spectral sequences which is an isomorphism from $E_2$ onward.    
\end{Proposition}
\begin{Proof}
    The functor $\iota$ induces a map of the cosimplicial spectra used in the construction of the homotopy limit spectral sequence, hence a map of homotopy limit spectral sequences. This map is an isomorphism on $E_2$ by the identification of $E_2$ in \cite[\S7.1]{bousfieldkan}.
\end{Proof}
In our computations with the $\sF$-homotopy limit spectral sequence we will make use of the previous two propositions by restricting to homotopy left cofinal subcategories of the orbit category when computing the $E_2$-page of the $\sF$-homotopy limit spectral sequence.
\subsection{A short exact sequence of coefficient systems}
We assume in addition that all the $H_i$ and the $K$ are normal subgroups of $G$.

Let $\underline{\Integers}$ be the constant $\Integers$-coefficient system on $\sO(G)_{\sF}$. For a subgroup $H$ of $G$, we denote by $\Integers[H]$ the coefficient system obtained from restricting $\underline{\Integers}$ to $\sO(G)_{\All_{H}}$ and left Kan extending back up to $\sO(G)_{\sF}$. Summing the counit maps $\Integers[H_i] \to \underline{\Integers}$ yields a short exact sequence of coefficient systems
\begin{equation}
    0 \to \bigoplus_{|I| - 1} \Integers[K] \xrightarrow{j} \bigoplus_{i \in I} \Integers[H_i] \to \underline{\Integers} \to 0
    \label{absplit:eq1}
\end{equation}
where the map $j$ is given by the matrix
\begin{equation}
    \begin{pmatrix}
        1 & 0 & \cdots & \ddots & \\
        -1 & 1 & 0 & \ddots & \\
        0 & -1 & 1 & \ddots &  \\
        \vdots & 0 & -1 & \ddots & \vdots \\
         & \vdots & 0 & \ddots & 0 \\
         & & \vdots & \ddots & 1 \\
         & & & \ddots & -1
    \end{pmatrix}
    \label{absplit:jmatrix}
\end{equation}
with $|I|$ rows and $|I| - 1$ columns. 
\begin{Lemma}
The short exact sequence (\ref{absplit:eq1}) lifts to a cofiber sequence of spectra
\begin{equation}
    \bigvee_{|I| - 1} {E\All_{K}}_+ \to \bigvee_{i \in I} {E\All_{H_i}}_+ \to E\sF_+.
    \label{absplit:eq2}
\end{equation}
\begin{Proof}
    The diagram of categories
    \begin{equation}
    \begin{tikzpicture}
    \matrix [matrix of math nodes, name=m, row sep=1cm, column sep=1cm]
    {
        \sO(G)_{\All_K} & \sO(G)_{\All_{H_i}} \\
        \sO(G)_{\All_{H_{i+1}}} & \\
    };
    \draw [{Hooks[right]}->] (m-1-1) -- (m-1-2);
    \draw [{Hooks[right]}->] (m-1-1) -- (m-2-1);
    \end{tikzpicture}
    \label{absplit:pushout}
\end{equation}
assemble together to a coequalizer diagram
\begin{equation}
    \begin{tikzcd}
        \coprod_{|I| - 1} \sO(G)_{\All_K} \arrow[shift left=4pt]{r} \arrow[shift right=4pt]{r} & \coprod_{i \in I} \sO(G)_{\All_{H_i}} \arrow[dashed]{r} & \sO(G)_{\sF}
        \label{absplit:coeqofcats}
    \end{tikzcd}
\end{equation}
with coequalizer the category $\sO(G)_{\sF}$. For any family $\sG$, denote by $\iota_{\sG} \colon \sO(G)_{\sG} \to S_G$ the inclusion into $G$-spaces. Then $E\sG \simeq \hocolim_{\sO(G)_{\sG}} \iota_{\sG}$ (\cite[p.~10]{mnn}). Applying this to (\ref{absplit:coeqofcats}) and the fact that homotopy colimits commute with coproducts yields a cofiber sequence of $G$-spaces
\begin{equation}
    \bigvee_{|I| - 1} E\All_{K} \to \bigvee_{i \in I} E\All_{H_i} \to E\sF,
    \label{absplit:cofibspaces}
\end{equation}
and, after taking suspension spectra, a cofiber sequence of $G$-spectra. After applying $H_*(-; \Integers)$ to this cofiber sequence of $G$-spectra, we get a long exact sequence that reduces to a short exact sequence because all the homology is concentrated in degree 0.  The left hand map is given by the matrix (\ref{absplit:eq1}) but with all $-1$ entries replaced by $1$. The bottom entry in every column comes from the vertical maps in the (\ref{absplit:pushout}), the top entry in every column entries comes from the horizontal mamps in the (\ref{absplit:pushout}), for all $i$. To change the bottom entries into $-1$, which is necessary to make a lift of (\ref{absplit:eq1}), we do the following: since suspension preserves homotopy colimits, we can replace at the level of $G$-spectra the vertical maps by the vertical maps smashed with $S \xrightarrow{-1} S$, which preserves cofibrations. We then get a strict map of homotopy colimit diagrams on objects given by either identity maps or the identity smashed with $S \xrightarrow{-1} S$, and therefore an equivalence of the homotopy colimits. We now have that
\begin{equation}
    0 \to \bigoplus_{|I| - 1} H_0(E\All^{(-)}; \Integers) \to \bigoplus_{i \in I} H_0(E\All_{H_i}^{(-)} ; \Integers) \to H_0(E\sF^{(-)}; \Integers) \to 0.
    \label{absplit:homgroupsses}
\end{equation}
coincides with (\ref{absplit:eq1}) on the level of objects and the left map.
For any family $\sG$ we have $\Integers[\sG](-) = H_0(E\sG^{(-)}; \Integers)$ (\cite[sec.~3]{mnn}), so the objects of the short exact sequences (\ref{absplit:eq1}) and (\ref{absplit:homgroupsses}) are equal. By construction, the left maps of (\ref{absplit:eq1}) and (\ref{absplit:homgroupsses}) are equal. The right hand maps possibly differ up to a sign, but postcomposing with the weak equivalence $E\sF_+ \xrightarrow{-1} E\sF_+$ makes sure that the right hand maps are equal on the level of $H_0(-)$.
\end{Proof}
\fxwarning{I need to explain what homotopy (co)finality is and why it is relevant.}
\end{Lemma}
For a normal subgroup $H$ of $G$, we denote by $E_*^{*,*}(H)$ the LHSSS associated to the group extension $H \to G \to G/H$. By \cref{prop:compserre}, the $\All_H$-homotopy limit spectral sequence is isomorphic to $E_*^{*,*}(H)$. The $\sF$-homotopy limit spectral sequence we denote by $E_*^{*,*}$. 
\begin{Lemma}
    Mapping (\ref{absplit:eq2}) into $\underline{H\Ffield_2}$ and taking Atiyah-Hirzebruch spectral sequences yields a long exact sequence of spectral sequences, which on the rows $E_2^{s,*}$ of the $E_2$-pages is given by
\begin{figure}[H]
    \centering
    \includegraphics[width=1.1\textwidth]{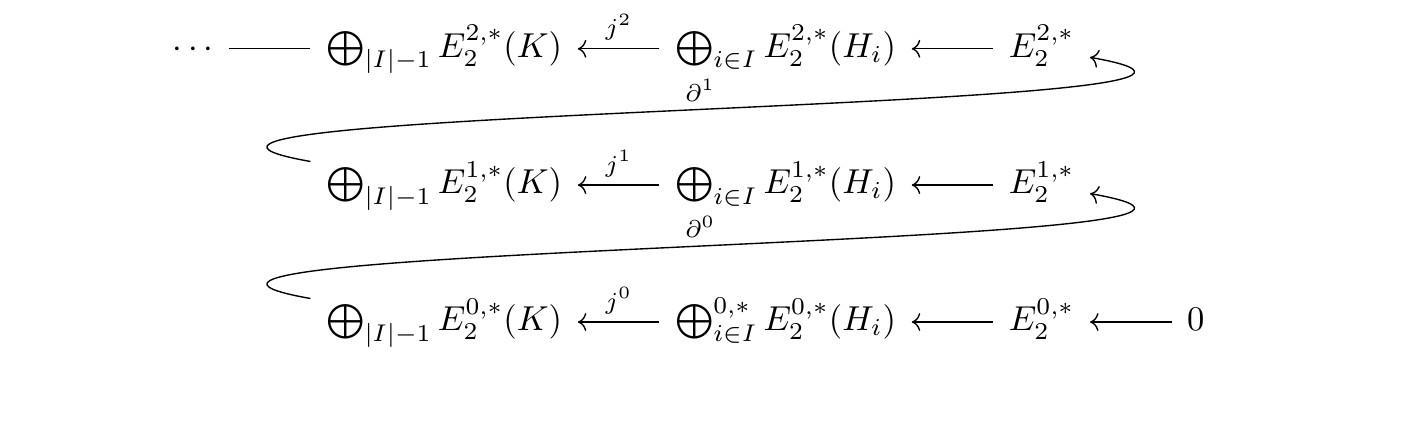} 
    \caption{A long exact sequence of rows $E_2^{s,*}$ of $E_2$-pages.}
    \label{absplit:eq3}
\end{figure}
where we have written $j^s$ for the map in degree $s$ in (\ref{absplit:eq3}) induced by $j$, and $\partial^s$ for the degree $s$-boundary morphism in (\ref{absplit:eq3}). The maps $j^s$ are given by the transpose of the matrix (\ref{absplit:jmatrix}). Summing the columns in (\ref{absplit:eq3}) gives the long exact sequence of spectral sequences.
\end{Lemma}
\begin{Proof}
    We use the following standard construction of the AHSS: take a Postnikov tower on $\underline{H\Ffield_2}$, and take homotopy fibers at each stage (i.e.\ the Postnikov sections of $\underline{H\Ffield_2}$). Apply $F(W,-)$ to this Postnikov tower for $W$ the three objects in the cofiber sequence (\ref{absplit:eq2}). Taking homotopy groups gives three exact couples, with associated spectral sequences the AHSS's $H^*(W; \pi_*^{(-)}\underline{H\Ffield_2})$. Moreover, since the objects $W$ sit in a cofiber sequence, we get a long exact sequence of exact couples, whence of spectral sequences. The identification of the AHSS's with $\sF$-homotopy limit spectral sequence from the $E_2$-page onward is \cite[Prop.~2.24]{mnn}), the identification of $\All_{K}$- and $\All_{H_i}$-homtopy limit spectral sequences with LHSSS's is \cref{prop:compserre}.

    The identification of the maps $j^s$ follows from the identifcations
    \begin{align}
        H^s(E\sG; \pi_*^{(-)} \underline{H\Ffield_2}) & \cong \Ext_{\Integers \sO(G)}^{s,t}(\Integers[\sG], \pi_*^{(-)} \underline{H\Ffield_2})
    \end{align}
    from \cite[Cor.~3.6]{mnn}, for every family $\sG$, and the fact that (\ref{absplit:cofibspaces}) lifts (\ref{absplit:eq1}).
\end{Proof}
 
We write
\begin{align}
    A^* & = \Im \partial^{*-1} = \coker j^{*-1}, \\
    B^* & = \coker \partial^{* - 1} = \ker j^{*},
    \label{absplit:eq4}
\end{align}
(where we have suppressed $t$-degrees)
and obtain the following proposition:
\begin{Proposition}
    The $E_2$-page of the homotopy limit spectral sequence fits in a short exact sequence
    \begin{equation}
        0 \to A^* \to E_2^* \to B^* \to 0
        \label{absplit:eq5}
    \end{equation}
    of spectral sequences  (where we suppressed $t$-degrees). This short exact sequence splits as a short exact sequence of bigraded $\Ffield_2$-modules.
    \label{prop:absplit}
\end{Proposition}
\begin{Proof}
    The short exact sequence is obtained from summing (\ref{absplit:eq4}) over the $s = *$-degree. Any short exact sequence splits as $\Ffield_2$-modules, and necessarily preserves the bigrading because the maps in the short exact sequence preserve the bigrading.
\end{Proof}
\fxwarning{I need to understand what happens at higher pages. Taking cohomology of a ses does not preserve exactness...}

\chapter{Exponent lemmas}
\label{ch:explems}
\section{Introduction}
We discuss some lemmas that will be of use in determining (bounds on) exponents. Some of these statements appear as exercises in \cite[sec.\ 4]{mnnnd}.

\section{Lemmas for $\sF$-exponents}
\label{sec:exponentlemmas}
We now give various lemmas which describe how $\sF$-exponents can change as the family $\sF$ varies.
\begin{Lemma}
    Let $\sF_1$, $\sF_2$ be two families of subgroups. Then $E(\sF_1 \cap \sF_2)_+ \simeq E\sF_{1_+} \wedge E\sF_{2_+}$.
    \label{lem:classspaceintersect}
\end{Lemma}
\begin{Proof}
    This follows from comparing fixed points.
\end{Proof}
\begin{Notation}
    For $\sF$ a family of subgroups of $G$, and $H$ a normal subgroup of $G$, denote by $\sF_H$ the subgroups of $\sF$ that are contained in $H$.  If  $\sF$ is a family that make sense for all groups $G$, such as the family of all subgroups, the family of elementary abelian $p$-groups, etc., we write $\sF(G)$ for this family of subgroups of $G$. For instance, we write $\sE_{(2)}(D_8)$ for the elementary abelian subgroups of the dihedral group of order 8.
\end{Notation}

\begin{Lemma}
    If $\sF$ is a family of subgroups of $G$, and $H$ is a normal subgroup of $G$, then $\Res^G_H E\sF \simeq E\sF_H$.
    \label{lem:restrictingclassifyingspaces}
\end{Lemma}
\begin{Proof}
    This follows from comparing fixed points.
\end{Proof}

\begin{Lemma}
    Let $M$ be an $\sF$-nilpotent $G$-spectrum, $H \subset G$ a subgroup. Then $\Res^G_H M$ is $\sF_H$-nilpotent, and $\exp_{\sF_H}\Res^G_H M \leq \exp_{\sF}M$. 
    \label{lem:expsubgp}
\end{Lemma}
\begin{Proof}
    This follows from \cite[Cor.\ 4.13]{mnnnd} and \cref{lem:restrictingclassifyingspaces}. 
\end{Proof}

Recall that for a group $G$, we denote the spectrum representing Borel $G$-equivariant $\Ffield_p$-cohomology by $\underline{H\Ffield_p}_{G}$ (\cref{not:borel}).
\begin{Corollary}
    Let $G$ be a group and $H \subset G$ a subgroup. Then
    \begin{equation}
        \exp_{\sE_{(p)}(H)} \underline{H\Ffield_p}_H \leq \exp_{\sE_{(p)}(G)} \underline{H\Ffield_p}_G.
    \end{equation}
    \label{cor:expofsubgroup}
\end{Corollary}
\begin{Lemma}
    Let $\sF_1$, $\sF_2$ be two families of subgroups, and let $M$ be a $G$-spectrum which is both $\sF_1$- and $\sF_2$-nilpotent, with exponents $m$, $n$ respectively. Then $M$ is $\sF_1 \cap \sF_2$-nilpotent, and
    \begin{equation}
        \exp_{\sF_1 \cap \sF_2} M \leq m+n - 1.
    \end{equation}
    \label{lem:expintersect}
\end{Lemma}
\begin{Proof}
    The fact that $M$ is $\sF_1 \cap \sF_2$-nilpotent is part of \cite[Prop.\ 6.39]{mnnnd}. The assumption on the exponents implies that both maps in 
    \begin{align}
        \sk_{m-1}E\sF_{1_+} \wedge \sk_{n-1}E\sF_{2_+} \wedge M & \to \sk_{m-1}E\sF_{1_+} \wedge E\sF_{2_+} \wedge M \\
        & \to E\sF_{1_+} \wedge E\sF_{2_+} \wedge M
        \label{explemmascap:eq1}
    \end{align}
    have a retraction, hence the composite has a retraction. The composite (\ref{explemmascap:eq1}) factors as
    \begin{equation}
        \begin{tikzcd}
            \sk_{m-1}E\sF_{1_+} \wedge \sk_{n-1}E\sF_{2_+}  \arrow{r}{\mathrm{(*)}}  \arrow{d}{\mathrm{(**)}} 
            & E\sF_{1_+} \wedge E\sF_{2_+} \wedge M \\
            \sk_{m+n-2}\left(E\sF_{1_+} \wedge E\sF_{2_+}\right) \wedge M  \arrow{ur}{\mathrm{(***)}} &
        \end{tikzcd}.
    \end{equation}
    Composing the retraction of $\mathrm{(*)}$ with $\mathrm{(**)}$ gives a retraction of $\mathrm{(***)}$. Since $E\sF_{1_+} \wedge E\sF_{2_+} \simeq E(\sF_1 \cap \sF_2)_+$ by \cref{lem:classspaceintersect}, the desired bound follows from \cref{expprop} (\ref{expprop:smash2}).
\end{Proof}
\begin{Notation}
    For a $G$-space $X$, we denote by $\sI(X)$ the minimal family containing the isotropy groups of $X$. 
\end{Notation}
\begin{Example}
    For an orthogonal $G$-representation $V$, the unit sphere $S(V)$ in $V$ inherits a $G$-action. Then $\sI(S(V))$ the smallest family containing the isotropy groups of $S(V)$.
\end{Example}
\begin{Lemma}[{{cf. Proof of \cite[Thm.\ 2.3]{mnn}}}]
    Let $R$ be a ring $G$-spectrum with multiplicative Thom classes (see \cite[Def.\ 5.1]{mnn}). Let $V$ be a $G$-representation with corresponding oriented Euler class
    \begin{equation}
        (\chi(V) \colon S^{-|V|} \to R) \in R^*
    \end{equation}
    Suppose $\chi(V)$ is nilpotent with $\chi(V)^n = 0$. Then $R$ is $\sI(S(V))$-nilpotent and 
    \begin{equation}
      \exp_{\sI(V)}R \leq n\dim_{\Reals} V.
    \end{equation}
\end{Lemma}
\begin{Proof}
    The fact that $\chi(V)^n = 0$ is equivalent to the oriented Euler class
    \begin{equation}
        R \xrightarrow{e_{nV}} S^{nV} \wedge R
        \label{lem:repstoexps}
    \end{equation}
    being nullhomotopic (see \cite[Lem.\ 5.3]{mnn} and \cite[Rem.\ 2.4]{mnn}). Hence the left map in the fiber sequence
    \begin{equation}
        S(nV)_+ \wedge R \to R \to S^{nV} \wedge R
    \end{equation}
    has a section: $R$ is a retract of $S(nV)_+ \wedge R$. But $S(nV)$ is a $(n \dim_{\Reals}V - 1)$-dimensional $G$-CW complex, hence by  \cref{lem:cwexp}, $\exp_{\sF}S(nV)_+ \leq n \dim_{\Reals}V$, and hence the same bound holds for $R$ by \cref{lem:expretracts}.
   \end{Proof}
\begin{Corollary}
    \label{cor:hfpeulerclass}
    Let $V$ be a $d$-dimensional $G$-representation with corresponding oriented Euler class (see \cite[Def.\ 5.1]{mnn})
    \begin{equation}
        (\chi(V) \colon S^{-d} \to \underline{H\Ffield_2}) \in H^{|V|}(BG;\Ffield_2).
    \end{equation}
    Suppose $\chi(V)$ is nilpotent with $\chi(V)^n = 0$. Then $\underline{H\Ffield_2}$ is $\sI(S(V))$-nilpotent with $\exp_{\sI(V)}\underline{H\Ffield_2} \leq nd$.
\end{Corollary}
\begin{Lemma}
    Let $f \colon G\to C_2 \cong O(1)$ be a 1-dimensional real representation with oriented Euler class $e \in H^1(BG;\Ffield_2)$. Suppose $e$ is nilpotent with $n$ the minimal integer $\geq 0$ such that $e^n = 0$. Then $\underline{H\Ffield_2}$ is nilpotent for the family $\All_{\ker f}$ of subgroups of $\ker f$ with $\exp_{\All_{\ker f}} \underline{H\Ffield_2} = n$.
    \label{cor:1dimrepexps}
\end{Lemma}
\begin{Proof}
    The upper bound for the exponent follows immediately from \cref{cor:hfpeulerclass}.

    For the lower bound on the exponent, denote a class representing $e$ on the $E_\infty$-page of the $\All_{\ker f}$-limit spectral sequence converging to $H^*(BG;\Ffield_2)$ by $\etilde$. Since $e$ restricts to 0 on $\ker f$, hence on all subgroups in $\All_{\ker f}$, we know that $\etilde$ lives in filtration degree $\geq 1$ on $E_\infty$. By assumption, $e^{n-1} \neq 0$. This does not imply $\etilde^{n-1} \neq 0$, however, if $\etilde^{n-1} = 0$, then there must be a class above it in higher filtration which detects $e^{n-1}$ on $E_{\infty}$, since $e^{n-1} \neq 0$. We see that on $E_{\infty}$ there are elements in filtration degree $\geq n-1$, hence $\exp_{\All_{\ker f}} R \geq n$.
\end{Proof}
\begin{Lemma}
    For a $G$-space $X$, $H \subset G$ a subgroup, the composite of maps of Borel cohomology rings
    \begin{equation}
        H_G^*(X; \Ffield_p) \xrightarrow{\Res^G_H} H_H^*(X;\Ffield_p) \xrightarrow{\Ind_H^G} H_G^*(X;\Ffield_p)
        \label{eq:borelresind}
    \end{equation}
    is multiplication by $[G:H]$.
    \label{lem:borelresind}
\end{Lemma}
\begin{Proof}
    This follows from \cite[Prop.\ III.2.4]{borel60} (note that the induction morphism also goes by the name of transfer morphism in this context, and that the morphism constructed in the proof of \cite[Prop.\ III.2.4]{borel60} is the transfer morphism).

    The result also follows by considering directly what happens on the cochain level.
\end{Proof}
    
\begin{Corollary}
    Let $G$ be a group and $P \subset G$ a $p$-Sylow. Then the maps of $G$-spectra
    \begin{equation}
        \underline{H\Ffield_p}_G \to G/P_+ \wedge \underline{H\Ffield_p}_G \to \underline{H\Ffield_p}_G
        \label{eq:psylowretract}
    \end{equation}
    representing the natural transformations $\Res^G_P$ and $\Ind_P^G$ exhibit $\underline{H\Ffield_p}_G$ as a retract of $G/P_+ \wedge \underline{H\Ffield_p}_G$.
    \label{cor:psylowretract}
\end{Corollary}
\begin{Proof}
    We need to show that for all subgroups $G' \subset G$, applying $\pi_*^{G'}$ to (\ref{eq:psylowretract}) yields an isomorphism. This composite is given precisely by \cref{lem:borelresind} for $X=G/G'$ and $H=P$. But multiplication by $[G:P]$ is an isomorphism because $p \nmid [G:P]$, whence $[G:P] \in \Ffield_p^{\times}$.
\end{Proof}
\begin{Lemma}
    Let $G$ be a group, $P$ a $p$-Sylow subgroup of $G$. Then 
    \begin{equation}
        \exp_{\sE_{(p)}(G)} \underline{H\Ffield_p}_G = \exp_{\sE_{(p)}(P)} \underline{H\Ffield_p}_P.
    \end{equation}
    \label{lem:psylowhfpexp}
\end{Lemma}
\begin{Proof}
    First,
    \begin{align}
        \Res^G_{P} \underline{H\Ffield_p}_G = \underline{H\Ffield_p}_{P}, \label{explemmas:eq3} \\
        \Res^G_{P} E\sE_{(p)}(G) = E\sE_{(p)}(P).       \label{explemmas:eq2}
    \end{align}
    Hence by \cref{cor:expofsubgroup} (cf. \cite[Cor.\ 4.13]{mnnnd}), 
    \begin{equation}
        \exp_{\sE_{(p)}(P)} \underline{H\Ffield_p}_P \leq \exp_{\sE_{(p)}(G)} \underline{H\Ffield_p}_G.
    \end{equation}

    For the upper bound, write $n = \exp_{\sE_{(p)}(P)} \underline{H\Ffield_p}_P$. We then have that 
    \begin{equation}
        \sk_{n-1}E \sE_{(p)}(P) \wedge \underline{H\Ffield_p}_P \to \underline{H\Ffield_p}_P
        \label{explemmas:eq1}
    \end{equation}
    admits a section.
    Note that for every $X \in \Sp_G$, we have $\Ind_H^G \Res^G_H X = G/H_+ \wedge X$. Furthermore we have that, if we use the model for $E\sE_{(p)}(P)$ from (\ref{explemmas:eq2}),
    \begin{align}
        \sk_{n-1}E\sE_{(p)}(P) & = \Res^G_P\sk_{n-1}E\sE_{(p)}(G), \\
    \end{align}
    Applying this, the fact that $\Res^G_P$ is monoidal, and (\ref{explemmas:eq3}) to (\ref{explemmas:eq1}) yields a section of
    \begin{equation}
        G/P_+ \wedge \sk_{n-1}E \sE_{(p)}(G) \wedge \underline{H\Ffield_p}_G \to G/P_+ \wedge \underline{H\Ffield_p}_G.
    \end{equation}
    Hence
    \begin{equation}
        \exp_{\sE_{(p)}(G)} G/P_+ \wedge \underline{H\Ffield_p}_G \leq n.
    \end{equation}
    Applying \cref{cor:psylowretract} and \cref{lem:expretracts} yields
    \begin{equation}
        \exp_{\sE_{(p)}} \underline{H\Ffield_p}_G \leq n.
    \end{equation}
\end{Proof}
\begin{Lemma}
\label{lem:twofamexp}
    Let $\sF$ and $\sG$ be families of subgroups, and let $M$ be both $\sF$- and $\sG$-nilpotent. Then for every $K \in \sF$, $\Res_K^G M$ is $\sG_K$-nilpotent by \cref{lem:expsubgp}. Write $n = \exp_{\sF} M$, $m_K = \exp_{\sG_K} \Res^G_K M$ for all $K \in \sF$, and $m = \max_K m_K$. Then
    \begin{equation}
        \exp_{\sG}M \leq mn.
    \end{equation}
\end{Lemma}
\begin{Proof}
    The $\sG_K$-nilpotence of $\Res^G_K M$ implies that, for all $K \in \sF$,
    \begin{equation}
        \sk_{m_K-1}E\sG_K \wedge \Res^G_K M \to \Res^G_K M
    \end{equation}
    admits a section. By inducing up to $G$, we see that $\exp_{\sG}(G/K_+ \wedge M) \leq m_K$. By taking the coproduct over all $K \in \sF$, we get that $\exp_{\sG} \sk_0 E\sF \wedge M \leq m$. We now argue by induction that
    \begin{equation}
        \exp_{\sG} \sk_{d-1} E\sF \wedge M \leq md,
        \label{eq:FGind}
    \end{equation}
    for all $d \geq 1$, of which the base case $d=1$ we have just established.

    Hence assume (\ref{eq:FGind}) has been established for some $d \geq 1$, and consider the cofiber sequence
    \begin{equation}
        \sk_{d-1}E\sF\wedge M \to \sk_d E\sF\wedge M \to \bigvee(S^d \wedge G/K_+)\wedge M
    \end{equation}
    that ends in a wedge of spheres with isotropy in $\sF$ smashed with $M$. By induction, the $\sG$-exponent of the left term is $\leq md$, and the $\sG$-exponent of the right hand side is $\leq m$ because we already saw that $\exp_{\sG}(G/K_+ \wedge M) \leq m$ for all $K \in \sF$. Hence by \cite[Prop.\ 4.9.2]{mnnnd} the $\sG$-exponent of the middle term is $\leq m(d+1) $, which completes the induction.
    
    We have by $\sF$-nilpotency of $M$ that $M$ is a retract of $\sk_{n-1}E\sF \wedge M$, hence taking $d=n$ in (\ref{eq:FGind}) yields the result.
\end{Proof}

\section{From representations to exponents via the projective bundle theorem}
In \cite[Ex.\ 5.16]{mnn} an upper bound on the $\sE_{(2)}$-exponent of $\underline{H\Ffield_2}_{Q_8}$ is determined using the projective bundle theorem. We will also repeatedly use this technique to give upper bounds on the $\sF$-exponent of $\underline{H\Ffield_2}_G$ for various 2-groups $G$ and families $\sF \supset \sE_{(2)}$. Therefore in this section we will describe this argument in some detail. We will need the following classical notion:
\begin{Definition}
    Let $V$ be a real or complex vector space. Then the \emph{projectivation} $\Proj(V)$ of $V$ is the space of all lines in $V$.
\end{Definition}
\begin{Remark}
    Note that if $V$ comes equipped with a linear $G$-action (i.e., is a $G$-representation), then $\Proj(V)$ inherits a natural $G$-action.
\end{Remark}
The goal of this section is to prove:
\begin{Proposition}
    Let $G$ be a finite group, $n \geq 0$ a natural number and suppose $G$ has a real $n$-dimensional representation $V$ such that the projectivation $\Proj(V)$ has isotropy groups contained in some family $\sF$. That is, for every real line $L \subset V$ we assume that the isotropy group $G_L \leq G$ of elements of $G$ fixing $L$ satisfies $G_L \in \sF$. Then $\underline{H\Ffield_2}$ is $\sF$-nilpotent and the exponent satisfies $\exp_{\sF} \underline{H\Ffield_2} \leq n$.
    \label{projbund:prop1}
\end{Proposition}
\begin{Proposition}
    Let $G$ be a finite group, $n \geq 0$ a natural number and suppose $G$ has a complex $n$-dimensional representation $V$ such that the projectivation $\Proj(V)$ has isotropy groups contained in some family $\sF$. That is, for every complex line $L \subset V$, we assume that the isotropy group $G_L \leq G$ of elements of $G$ fixing $L$ satisfies $G_L \in \sF$. Then $\underline{H\Integers}$ is $\sF$-nilpotent, and the exponent satisfies $\exp_{\sF} \underline{H\Integers} \leq 2n -1$.
    \label{cpxprojbund:prop1}
\end{Proposition}

\subsection{The projective bundle theorem}
To prove \cref{projbund:prop1}, we will use the projective bundle theorem, which can be used to develop the theory of Stiefel-Whitney classes and Chern classes. This is carried out for instance in \cite[Ch. 17]{husemoeller1994fibre}. We will follow the treatment and notation of \cite[\S 17.2]{husemoeller1994fibre}, but only discussing what we need. For an overview of all bundles and their relations that will appear we refer to diagram (\ref{projbund:eq1}) below, where the dashed arrows indicate pullbacks.

As in \cite[Ch.~17]{husemoeller1994fibre}, we consider real and complex vector bundles at the same time. For the case of real vector bundles, we write $c = 1$, we consider cohomology with coefficients in $K_c = \Integers/2$, and we let $F$ be the field $\Reals$ of real numbers. For the case of complex vector bundles, we write $c=2$,  we consider cohomology with coefficients in $K_c = \Integers$, and we let $F$ be the field $\Complex$ of complex numbers.

We will write $E(\eta)$ resp.\ $B(\eta)$ for the total, resp.\ base space, of a fiber (not necessarily vector) bundle $\eta$. Let $\xi \colon E \xrightarrow{p} B$ be an $n$-dimensional vector bundle. Let $E_{0}$ be the non-zero vectors in $E$. Let $E'$ be the quotient of $E_0$ where we identify non-zero vectors in the same line. This yields a factorization
\begin{equation}
    \begin{tikzcd}
        E_0 \arrow[rr] \arrow[dr] & & B \\
        & E' \arrow{ur}[dashed]{q}
    \end{tikzcd},
\end{equation}
and $E' \xrightarrow{q} B$ is a fibre bundle with fibre $F\Proj^{n-1}$, called the projectivation of $\xi$ and denoted $\Proj(\xi)$.

For every point $b \in B$ in the base, the inclusion of the fibre $F^n \to p^{-1}(b) \subset E$ defines a natural inclusion $j_b \colon F \Proj^{n-1} \to q^{-1}(b) \subset E(P(\xi))$. Pulling back the bundle $\xi$ along $q$ yields a bundle $q^*(\xi) \colon E(q^*(\xi)) \to E(P(\xi))$, which admits the canonical line bundle $\lambda_\xi$ as a subbundle, where a point in the total space $E(\lambda_\xi)$ is a pair $(L,x)$ where $L$ and $x$ lie over the same base point and $x \in L$.

\begin{Proposition}[{{see, e.g., \cite[Prop.\ 17.2.]{husemoeller1994fibre}}}]
    Let 
    \begin{equation}
    j_b \colon F\Proj^{n-1} \to E(P(\xi))
    \end{equation}
    be the inclusion of the fibre. Then $j_b^*(\lambda_\xi)$ is isomorphic to the tautological line bundle on $F \Proj^{n-1}$.
\end{Proposition}

All the bundles that appear and their relations are as follows, with dashed lines indicating pullbacks:
\begin{equation}
    \begin{tikzcd}
        E(j_b^*(\lambda_\xi))  \arrow[dashed]{r}  \arrow[dashed]{ddd} & E(\lambda_\xi) \arrow{d}[hook]{} & & \\
        & E(q^*(\xi))  \arrow[dashed]{dd} \arrow[dashed]{rr} & & E \arrow{dd}{p} \\
        & & E_0  \arrow{ur}[hook]{} \arrow{dl}  & \\
        F \Proj^{n-1} \arrow{r}{j_b} & E(\Proj(\xi)) \arrow{rr}{q} & & B
    \end{tikzcd}
    \label{projbund:eq1}
\end{equation}

Let $f \colon E(\Proj(\xi)) \to F \Proj^{\infty}$ be a classifying map for the line bundle $\lambda_\xi$, i.e.\ $f^*(\gamma) \cong \lambda_\xi$, where $\gamma$ is the tautological line bundle on $F\Proj^\infty$. Let $z \in H^c(F \Proj^{\infty} ; K_c)$ be the polynomial generator, and let $a_\xi = f^*(z)$. Since $f$ is unique up to homotopy, $a_\xi$ is well defined. We then have the following Theorem:
\begin{Theorem}[{{Projective bundle theorem, see \cite[Thm.\ 17.2.5]{husemoeller1994fibre}}}]
    \label{thm:projectivebundle}
    For an $n$-dimensional vector bundle $\xi$, the classes $1, a_\xi, \ldots, a_\xi^{n-1}$ form a basis of the free $H^*(B(\xi);K_c)$-module $H^*(E(\Proj(\xi)); K_c)$. Moreover, 
    \begin{equation} 
    q^* \colon H^*(B(\xi);K_c) \to H^*(E(\Proj(\xi)); K_c)
    \end{equation}
    is a monomorphism.
\end{Theorem}
\subsection{Free $R$-modules}
In this subsection we discuss some properties of free $R$-modules. The statements and notions are not new, but collected here for convenience. 
\begin{Definition}
    Let $R_*$ be a graded ring, $M_*$ a graded a graded $R_*$-module. We say that $M_*$ is a \textbf{free} $R_*$-module if there are integers $n_i$ and an isomorphism of $R_*$-modules
    \begin{equation}
        g \colon \bigoplus_i \Sigma^{n_i} R_* \xrightarrow{\cong} M_*.
    \end{equation}
    Here we write $(\Sigma N)_* = N_{*+1} $ for the suspension of graded abelian groups.
\end{Definition}
Likewise, we have:
\begin{Definition}
    Let $R$ be a (possibly equivariant) ring spectrum, and let $M$ be an $R$-module. Then $M$ is a \textbf{free} $R$-module if there are integers $n_i$ and a weak equivalence  of $R$-modules
    \begin{equation}
        g \colon \bigvee_I \Sigma^{n_i} R \xrightarrow{\simeq} M.
    \end{equation}
    This weak equivalence is understood to be an equivariant weak equivalence in the equivariant case.
\end{Definition}
A map $f \colon R \to T$ of homotopy ring spectra makes $T$ into an $R$-module. The next proposition says that one can check freeness of $T$ as an $R$-module on the level of homotopy groups. The next proposition considers the non-equivariant case, the equivariant case is treated later.
\begin{Proposition}
    \label{prop:freermod}
    Let $f \colon R \to T$ be a map of non-equivariant homotopy ring spectra. Then the following are equivalent:
    \begin{enumerate}[(1)]
        \item The map $f$ makes $T$ into a free $R$-module. \label{prop:freermod:1}
        \item The map $\pi_* f$ makes $\pi_* T$ into a free $\pi_* R$-module. \label{prop:freermod:2}
        \item There are integers $n_i$ and homotopy classes of maps $\alpha_i \colon S^{n_i} \to T$ such that the the $\alpha_i$ form a basis of $\pi_*T$ as a $\pi_* R$-module. \label{prop:freermod:3}
    \end{enumerate}
\end{Proposition}
\begin{Proof}
    (1) $\Rightarrow$ (2) This follows from applying $\pi_*$ to the equivalence $\bigvee_i \Sigma^{n_i} R \simeq T$ of $R$-modules.

    (2) $\Rightarrow$ (3) This follows from applying $\pi_* f$ to the canonical basis of $\bigoplus_i \Sigma^{n_i} \pi_* R$.

    (3) $\Rightarrow$ (1) The $\alpha_i$ define an isomorphism
    \begin{equation}
        \bigoplus_i \pi_*\Sigma^{n_i} R \to \pi_*T.
    \end{equation}
    The result follows by lifting this isomorphism to an equivalence of spectra, which one can do as follows. For every $i$, the composite
    \begin{equation}
        R \wedge S^{n_i} \xrightarrow{f \wedge \alpha_i} T \wedge T \xrightarrow{\mu_T} T
    \end{equation}
    is a map of $R$-modules. Taking the coproduct over $i$ yields a map of $R$-modules
    \begin{equation}
        \bigvee_i \Sigma^{n_i} R \to T
    \end{equation}
    which is the desired lift.
\end{Proof}
Combining this Proposition with \cref{thm:projectivebundle} yields the following Corollary:
\begin{Corollary}
    Let the notation be as in \cref{thm:projectivebundle}, and let $H\Ffield_2$ be the Eilenberg-MacLane spectrum representing cohomology with $\Ffield_2$-coefficients. Then the map $F(B(\xi)_+,H\Ffield_2) \to F(E(\Proj(\xi))_+, H\Ffield_2)$ induced by $q$ makes the mapping spectrum $F(E(\Proj(\xi))_+, H\Ffield_2)$ a free $F(B(\xi)_+, H\Ffield_2)$-module.
\end{Corollary}
\begin{Proof}
    According to \cref{thm:projectivebundle}, $q_*$ makes $H^*(E(\Proj(\xi)), \Ffield_2)$ into a free $H^*(B(\xi),\Ffield_2)$-module, which precisely means that the map 
    \begin{equation} 
    \pi_*F(B(\xi)_+,H\Ffield_2) \to \pi_*F(E(\Proj(\xi))_+,H\Ffield_2)
\end{equation}
    induced by $q$ makes $\pi_*F(E(\Proj(\xi)_+,H\Ffield_2))$ into a free $\pi_*F(B(\xi)_+,H\Ffield_2)$-module. Applying \cref{prop:freermod} yields the desired result.
\end{Proof}
\begin{Corollary}
    Some suspension of the mapping spectrum $F(B(\xi),H\Ffield_2)$ is a retract of the mapping spectrum $F(E(\Proj(\xi)), H\Ffield_2)$.
\end{Corollary}
\begin{Proof}
    Every free $R$-module admits some suspension of $R$ as a retract.
\end{Proof}
Similarly, equivariantly we have
\begin{Proposition}
    \label{prop:freegrmod}

    Let $R \to T$ be a map of homotopy $G$-ring spectra. Then the following are equivalent:
    \begin{enumerate}[(1)]
        \item The map $f$ makes $T$ into a free equivariant $R$-module,
        \item There are homotopy classes of maps $\alpha_i \colon S^{n_i} \to T$ such that for every subgroup $H \leq G$, the canonical map $\bigoplus_i \pi^H_{*} \Sigma^{n_i} R \to \pi_*^H T, x \mapsto x \cdot \Res^G_H \alpha_i$ is an isomorphism of $\pi_*^H R$-modules.
    \end{enumerate}
\end{Proposition}
\begin{Proof}
    (1) $\Rightarrow$ (2) follows from moving the canonical basis of $\bigoplus_i \pi_*^G \Sigma^{n_i} R$ through the equivalence of $R^G$-modules $\bigvee_i \Sigma^{n_i} R^G \simeq T^G$.

    (2) $\Rightarrow$ (1) follows from taking the coproduct over $i$ of the composite
    \begin{equation}
        R \wedge S^{n_i} \to T \wedge T \to T
    \end{equation}
    to get a map
    \begin{equation}
        \bigvee_i \Sigma^{n_i} R \to T.
        \label{eq:eqfreermod:1}
    \end{equation}
    The assumption on the restrictions says precisely that this is an equivariant equivalence.
\end{Proof}
\subsection{From representations to exponents}
In the proof of \cref{projbund:prop1}, we will take the projectivation of the Borel construction of a $G$-representation $V$. The proof will make use of the fact that the projectivation and the Borel construction commute, which is the content of the next lemma.
\begin{Lemma}
    Let $V$ be a real or complex $G$-representation, and let $H$ be a subgroup of $G$. Consider the Borel construction
    \begin{equation}
        V \to V \times_H EG \to BH.
    \end{equation}
    Then the projectivation of the Borel construction equals the Borel construction of the projectivation, that is,
    \begin{equation}
        \Proj(V \times _H EG) = \Proj(V) \times_H EG.
    \end{equation}
    \label{lem:projborcomm}
\end{Lemma}
\begin{Proof}
    Write $V_0$ for the non-zero vectors in $V$, write $P$ for the product $V \times EG$, an write $P_0 = V_0 \times EG$. Then $F^{*}$ acts on the left coordinate of $P_0$ by multiplication, and trivially on the right coordinate, and $H$ acts coordinatewise on $P_0$. Moreover, these actions commute, by linearity of the action of $G$ (hence of $H$) on $V$. Therefore $P_0$ has in fact an action of $F^* \times H$. Dividing out the action of $F^{*}$ first and then dividing out the action of $H$ gives the Borel construction of the projectivation, dividing out the action of $H$ first and then the action of $F^*$ gives the projectivation of the Borel construction. Both equal the space obtained by dividing out $F^* \times H$, and they are therefore equal.
\end{Proof}
We are now ready to prove \cref{projbund:prop1} and \cref{cpxprojbund:prop1}.
\begin{Proof}[{{of \cref{projbund:prop1} and \cref{cpxprojbund:prop1}}}]
    Consider the Borel construction on $V$ for an arbitrary subgroup $H \leq G$, and call the resulting bundle $\xi_H$:
    \begin{equation}
        \xi_H : V \to V \times_H EG\to BH.
    \end{equation}
    Note that this is natural with respect to inclusions of subgroups. The associated projective bundle is 
    \begin{equation}
        \Proj(\xi_H) \colon \Proj(V) \to \Proj(V \times_H EG) \to BH.
    \end{equation}
    Observe that $\Proj(V \times_H EG) = \Proj(V) \times_H EG$ by \cref{lem:projborcomm}. Hence by \cref{thm:projectivebundle}, 
    \begin{equation}
    F((\Proj(V) \times_H EG)_+, HK_c) \cong F(\Proj(V)_+, \underline{HK_c})^H
\end{equation}
is a free $F(BH,K_c) = \underline{HK_c}^H$-module. Recall that the \emph{real} dimension of $V$ is $cn$. A basis is given by the proof of \cite[Prop.\ 17.3.3]{husemoeller1994fibre}, which shows that there are classes $1, a_{\xi_H}, \ldots, a_{\xi_H}^{n-1}$ that form a basis of $\pi_*^H F(\Proj(V),\underline{HK_c})$ as a $\pi_*^H \underline{HK_c}$-module. Moreover, the element $a_{\xi_H}$ is natural with respect to restriction to subgroups. It follows that all basis elements, being powers of $a_{\xi}^H$, are natural with respect to restriction to subgroups. In particular,
    \begin{equation}
        \Res^G_H a_{\xi_G} = a_{\xi_H}
    \end{equation}
    for all $H \leq G$. Hence by \cref{prop:freegrmod}, $F(\Proj(V), \underline{HK_c})$ is a free $\underline{HK_c}$-module. Hence a suspension of $\underline{HK_c}$ is a retract of $F(\Proj(V),\underline{HK_c})$. Therefore by \cite[Prop.\ 4.9]{mnnnd} and \cref{cor:expsuspinv},
    \begin{align}
        \exp_{\sF}\underline{HK_c} & \leq \exp_{\sF}F(\Proj(V),\underline{HK_c}),
        \intertext{but $V$ is $cn$-dimensional, therefore $\Proj(V)$ admits the structure of an $(cn-c)$-dimensional $G$-CW-complex (\cite[Cor.~7.2]{illman83}) with isotropy (by assumption) contained in $\sF$, therefore by \cite[Prop.\ 2.26]{mnn} the exponent of the right hand side is}
         & \leq cn - c +1.
    \end{align}
    In the case of a real representation we have $c=1$ and this equals $n$, in the case of a complex representation we have $c=2$ and this equals $2n-1$.
\end{Proof}

\section{A bound on the exponent for the family of proper subgroups}
Let $G$ be a finite non-abelian $2$-group of size $2^k$. The goal of this section is to prove
\begin{equation}
    \exp_{\sP} \underline{H\Ffield_2}_G \leq 2\floor{\sqrt{|G| - k}}-1,
\end{equation}
see \cref{cor:propsubgpexp} below. The argument is an adaption of the ones found in \cite[Lem.\ 4.3]{pakianathan03} and \cite{symonds91}.
\begin{Lemma}
    Let $G$ be a finite non-abelian $p$-group of size $p^k$. Then $G$ has an irreducible complex representation $V$ with $\dim_{\Complex} V \geq 2$, and moreover all such $V$ satisfy
    \begin{equation}
        \dim_{\Complex} V \leq \floor{\sqrt{|G| - 1}}.
    \end{equation}
    \label{lem:irredrepspgps}
\end{Lemma}
\begin{Proof}
    Denote by $n_1,\ldots,n_h$ the $\Complex$-dimensions of the irreducible $\Complex$-representations of $G$. We have
    \begin{equation}
        n_1^2 + n_2^2 + \cdots + n_h^2 = |G|,
    \end{equation}
    (see, e.g., \cite[Cor.\ 2.4.2]{Serrelinrep}). Since $G$ is non-abelian, some $n_i \geq 2$ (\cite[Thm.\ 3.1.9]{Serrelinrep}. Assume without loss of generality that it is $n_1$. Then
    \begin{align}
        |G| & = n_1^2 + n_2^2 +  \cdots + n_h^2 \\
        & \geq n_1^2 + 1, \\
    \end{align}
    because there is also a trivial representation, hence
    \begin{equation}
        n_1 \leq \sqrt{|G| - 1}.
    \end{equation}
    Flooring this preserves the inequality and does not change the integer on the left hand side. This yields the result.
\end{Proof}
\begin{Corollary}
    Let $G$ be a finite non-abelian $2$-group of size $2^k$, and let $\sP$ be the family of proper subgroups of $G$. Then
    \begin{equation}
        \exp_{\sP} \underline{H\Ffield_2}_G \leq 2 \floor{\sqrt{|G|-1}} -1.
    \end{equation}
    \label{cor:propsubgpexp}
\end{Corollary}
\begin{Proof}
    Let $V$ be an irreducible complex representation of $G$ satisfying
    \begin{equation}
        2 \leq \dim_{\Complex} V \leq \floor{\sqrt{|G| - 1}},
    \end{equation}
    as furnished by \cref{lem:irredrepspgps}. Then the complex projectivation $\Proj(V)$ has isotropy in $\sP$ the family of proper subgroups of $G$, for if $L \in \Proj(V)$ were fixed by all of $G$, $V$ would not be irreducible, since $\dim_{\Complex} V \geq 2$. An application of \cref{cpxprojbund:prop1} yields the result.
\end{Proof}

\chapter{Mod-2 group cohomology computations of 2-groups using the $\sF$-homotopy limit spectral sequence}
\label{ch:gpcohfhtylimss}
\section{Introduction}
In this chapter we apply the $\sF$-homotopy limit spectral sequence to the case $X = \pt$, $\sF \supset \sE_{(2)}$ and $G$ all finite 2-groups of order $\leq 16$, to compute group cohomology. For the indexing of the $\sF$-homotopy limit spectral sequence, we use the following non-standard convention:

\begin{Convention}
    For a 2-group $G$ and family $\sF \supset \sE_{(2)}$, the $\sF$-homotopy limit spectral sequence is given by
    \begin{equation}
        E_2^{s,t} = \sideset{}{^s}\lim_{(-) \in \sO(G)^{\op}_{\sF}} H^t(B(-);\Ffield_2) \Rightarrow H^{s+t}(BG; \Ffield_2).
    \end{equation}
    We use Adams' grading convention for displaying the spectral sequences, which means that we put $s+t$ on the horizontal axis and $s$ on the vertical axis. This yields a spectral sequence drawn in the first quadrant, with $d_n$ going 1 over to the right and $n$ up.
\end{Convention}
The group cohomology rings that we study in this chapter are well known. A classical technique for computing group cohomology is to use the Lyndon-Hochschild-Serre spectral sequence, as was for example done by Quillen to compute the group cohomology of the extra-special 2-groups (\cite{quillen71mod2}). Rusin computed the cohomology of 2-groups up to order 32 in \cite{rusin89} using the Eilenberg-Moore spectral sequence. Using large initial segments of projective resolutions, and a completion criterion of Benson (\cite{benson04}), Carlson computed the group cohomology of all 2-groups up to order 64 (\cite{carlson2003cohomology}) using a computer. Green and King improved these methods by deriving a variant of Benson's completion criterion that allows for earlier detection of completion (\cite[Thm.\ 3.3]{green11}) and using this computed the cohomology of all 2-groups up to order 64 (\cite{green11} and \cite{green15}).

The $\sF$-homotopy limit spectral sequence provides a fourth method for computing group cohomology. The advantage of this spectral sequence is that, if one takes $\sF \supset \sE_{(2)}$, there will be a horizontal vanishing line at a finite page. Hence the computation is a finite one, and moreover the knowledge that a vanishing line must eventually appear can also help determine differentials. However, these two advantages come with a trade-off compared to using the Lyndon-Hochschild-Serre or Eilenberg-Moore spectral sequence: compared to the latter two, it is in general rather difficult to compute the $E_2$-page of the $\sF$-homotopy limit spectral sequence. Moreover, the product structure on the $E_2$-page of the $\sF$-homotopy limit spectral sequence is rather mysterious, which makes propagating differentials difficult. But the trade-off is at least good enough to still be able to compute the group cohomology of all the 2-groups up to order 16, at least additively, as shown in this chapter. 

\section{The 2-groups of order at most 16}
\label{sec:classif2gps}
The classification of 2-groups up to order 16 has been long known (see, e.g., \cite[\S 118]{burnside}). In this section we recall this classification, and provide forward references to the bounds on the exponents of $\underline{H\Ffield_2}_G$ for all these groups.

Each classification subsection is preceded by one or two subsections defining families of 2-groups needed in the classification.
\subsection{Abelian 2-groups}
For all finite abelian 2-groups $A$, we compute the $\sE_{(2)}$-homotopy limit spectral sequence converging to $H^*(BA;\Ffield_2)$ in \cref{sec:ab2gps}. Additionally, we obtain the following exponent:
\begin{nnProposition}[{{see \cref{prop:expelab}}}]
    Let $A$ be a finite abelian 2-group, and let $n$ be the minimal number of generators of the subgroup of 2-divisible elements of $A$, i.e., of the image of the multiplication-by-2 map $A \xrightarrow{\cdot 2} A$. (We consider the trivial group to be generated by 0 elements.) Then
    \begin{equation}
        \exp_{\sE_{(2)}} \underline{H\Ffield_2}_A = n + 1.
    \end{equation}
\end{nnProposition}
Because all 2-groups up to order 4 are abelian, \cref{prop:expelab} in particular determines all the $\sE_{(2)}$-exponents of all 2-groups up to order 4.
\subsection{The dihedral groups}
Denote by $R_{\theta} = \begin{pmatrix} \cos \theta & - \sin \theta \\ \sin \theta & \cos \theta \end{pmatrix}$ the rotation by $\theta$-matrix, and by $S = \begin{pmatrix} 1 & 0 \\ 0 & -1 \end{pmatrix}$ the reflection-in-the-$x$-axis matrix. Then $R_{2\pi/2^{n-1}}$ and $S$ generate a subgroup of order $2^n$ of $O(2)$, called the dihedral group of order $2^n$, and denoted $D_{2^n}$. It is the automorphism group of the regular $2^{n-1}$-gon centered at the origin and with two of its sides bisected by the $x$-axis.

The map given by $R_{2\pi /2^{n-1}}  \mapsto \rho$ and $S \mapsto \sigma$ gives an isomorphism of $D_{2^n}$ with the finitely presented group
\begin{equation}
    \langle \sigma, \rho \, | \, \sigma^2 = \rho^{2^{n-1}} = e, \sigma \rho \sigma^{-1} = \rho^{-1} \rangle.
\end{equation}
\subsection{The quaternion groups}
The generalized quaternion group $Q_{2^n}$ of order $2^n$ is the finite subgroup of quaternionic space $\HH$ generated multiplicatively by the elements of unit length $\{ e^{2\pi i / 2^{n-1}},\, j\}$ (see, e.g., \cite[XII.\S7]{cartaneilenberg}). Denoting these generators by $r$ and $s$, respectively, one gets the presentation
\begin{equation}
    Q_{2^n} = \langle r,\, s\, | \, r^{2^{n-2}} = s^2,\, rsr = s \rangle.
\end{equation}
\subsection{The groups of order 8}
Now that we have defined the families of dihedral and quaternion groups, we are able to state the classification of the groups of order 8.
\begin{Theorem}[{{see, e.g., \cite[Thm.\ 13.3]{armstrong}}}]
    There are, up to isomorphism, exactly 5 groups of order 8. They are given by $C_2^{\times 3}$, $C_2 \times C_4$, $C_8$, $D_8$ and $Q_8$.
\end{Theorem}
The $\sE_{(2)}$-homotopy limit spectral sequence converging to $H^*(BD_8; \Ffield_2)$ is computed in \cref{sec:dihgps}, and $\exp_{\sE_{(2)}} \underline{H\Ffield_2}_{D_{2^n}}$ is determined.

The $\sE_{(2)}$-homotopy limit spectral sequence converging to $H^*(BQ_{2^n};\Ffield_2)$ was already computed for $n=3$ in \cite[Ex.\ 5.18]{mnn}, as well as $\exp_{\sE_{(2)}} \underline{H\Ffield_2}_{Q_{2^n}}$ in that case. We give a straightforward generalization for both computations to the case $n > 3$ in \cref{sec:quatgps}. In summary, the exponents of the groups of order 8 are then given by
\begin{table}[H]
    \centering
    \begin{tabular}{r|r|r}
        $G$ & $\exp_{\sE_{(2)}} \underline{H\Ffield_2}_G$ & Reference \\
        \hline
        $C_2^{\times 3}$ & 1 & \cref{prop:expelab} \\
        $C_2 \times C_4$ & 2 & \cref{prop:expelab} \\
        $C_8$ & 2 & \cref{prop:expelab} \\
        $D_8$ & 2 & \cref{prop:d2nexp} \\
        $Q_8$ & 4 & \cite[Ex.\ 5.18]{mnn}
    \end{tabular}
    \caption{The $\sE_{(2)}$-exponents for the groups of order 8.}
    \label{table:groups8}
\end{table}
\subsection{The groups of order 16}
We state the classification of groups of order 16, following the treatment in \cite{wild}. We first need to introduce some notation as used in \cite{wild}.
\begin{Notation}
    Given groups $N$ and $H$ and an homomorphism $\phi \colon H \to \Aut(N)$, we denote the semidirect prodct of $N$ by $H$ with respect to $\phi$ by $N \overset{\phi}{\rtimes} H$.

    In the special case $H = C_2$, $N = C_8 =  \langle x \rangle$, the possible actions of $C_2$ on $C_8$ are given by $x \mapsto x^n$ with $n$ odd, and for all such actions we denote the induced semidirect product by $C_8 \overset{n}{\rtimes} C_2$. In particular, $C_8 \overset{1}{\rtimes} C_2 \cong C_8 \times C_2$ and $C_8 \overset{7}{\rtimes} C_2 \cong D_{16}$. These semidirect products for $n=3,5$ also have names: $C_8 \overset{3}{\rtimes} C_2$ is also known as the semidihedral group of order 16, denoted $SD_{16}$. The group $C_8 \overset{5}{\rtimes} C_2$ is also known as the modular group of order 16, denoted $M_{16}$.
\end{Notation}
\begin{Notation}
    If $N$ and $H$ are groups such that there is an up to isomorphism unique non-trivial semi-direct product of $N$ by $H$, we denote it by $N \rtimes H$.
\end{Notation}
\begin{Definition}
    \label{def:classact}
    Let $C_4 = \langle x \rangle$, and $C_2 = \langle y \rangle$. Following \cite[Fact 4]{wild}, we denote by $\psi_5, \psi_6$ the elements of $\Aut(C_4 \times C_2)$ given by $\psi_5(x) = xy$, $\psi_5(y) = y$, $\psi_6(x) = x^3$ and $\psi_6(y) = x^2y$.
\end{Definition}
We are now able to list all groups of order 16, and the bounds on the exponents that we obtained.
\begin{Theorem}
    The following table lists in its first column all groups of order 16, see, e.g., \cite[Thm.\ 2]{wild}. The second column either gives the $\sE_{(2)}$-exponent or an interval in which the exponent lies. The third column gives a forward reference for the proof of the contents of the second column. The fourth column gives a forward reference to the computation of an $\sF$-homotopy limit spectral sequence with $\sF \supset \sE_{(2)}$ converging to the group cohomology of the first column.
    \begin{table}[H]
        \centering
        \begin{tabular}{r|r|r|r}
            $G$ & $\exp_{\sE_{(2)}}\underline{H\Ffield_2}_G$ & Reference & holim SS \\
            \hline
            $C_2^{\times 4}$ & 1 & \cref{prop:expelab} & \cref{sec:ab2gps} \\
            $C_2^{\times 2} \times C_4$ & 2 & \cref{prop:expelab} & \cref{sec:ab2gps} \\
            $C_4 \times C_4$ & 3 & \cref{prop:expelab} & \cref{sec:ab2gps} \\
            $C_8 \times C_2$ & 2 & \cref{prop:expelab} & \cref{sec:ab2gps} \\
            $C_{16}$ & 2 & \cref{prop:expelab} & \cref{sec:ab2gps}  \\
            $D_{16}$ & 2 & \cref{prop:d2nexp} & \cref{sec:dihgps} \\
            $Q_{16}$ & 4 & \cref{prop:q2ne2exp} & \cref{sec:quatgps} \\
            $SD_{16} =  C_8 \overset{3}{\rtimes} C_2$ & $[3,4]$ & \cref{prop:sd16e2exp} & \cref{sec:sd16} \\
            $M_{16} = C_8 \overset{5}{\rtimes} C_2$ & $[3,4]$ & \cref{prop:m16e2exp} & \cref{sec:m16} \\
            $D_8 \ast C_4 = (C_4 \times C_2) \overset{\psi_6}{\rtimes} C_2$ & $4$ & \cref{d8cpc4:e2exp} & \cref{sec:d8cpc4} \\
            $C_4 \rtimes C_4$ & $[3,4]$ & \cref{prop:c4sdc4e2exp} & \cref{sec:c4sdc4} \\
            $(C_4 \times C_2) \overset{\psi_5}{\rtimes} C_2$ & 2 & \cref{prop:c4xc2sdc2exp} & \cref{sec:c4xc2sdc2} \\
            $Q_8 \times C_2$ & 4 & \cref{prop:q8xc2e2exp} & \cref{sec:q8xc2} \\
            $D_8 \times C_2$ & 2 & \cref{prop:d8xc2e2exp} & \cref{sec:d8xc2} 
        \end{tabular}
        \caption{Summary of results for groups of order 16.}
        \label{table:groups16}
    \end{table}
\end{Theorem}

\section{Abelian 2-groups}
\label{sec:ab2gps}
In this section we compute the $\sE_{(2)}$-homtopy limit spectral sequence for finite 2-groups. We first recall some basic results on the cohomology rings of cyclic 2-groups and the LHSSS's computing these. We will omit the $\Ffield_2$-coefficients from the notation for cohomology groups in this section.
\subsection{Cyclic 2-groups}
\label{subsec:cycl2gps}
Let $C_{2^n}$ the cyclic group of order $2^n$, generated by $t$. Denote by $\delta_t \in H^1(BC_{2^n})$ the cohomology class dual to $t$, where $C_{2^n}$ necessarily acts trivially on $\Ffield_2$ since $\Ffield_2$ has only the identity as automorhpism. Denote by $\beta_n(-)$ the $n$-th order Bockstein. Then $\beta_n\delta_t \in H^2(BC_{2^n})$ is the cohomology class corresponding to the unique non-trivial group extension of $C_{2^n}$ by $C_2$, and furthermore (see, e.g., \cite[sec. 3.2]{evens91}),
\begin{equation}
    H^*(BC_{2^n}) \cong \Ffield_2[\delta_t]/(\delta_t^2) \otimes \Ffield_2[\beta_n\delta_t],
\end{equation}
with, for $n = 1$, the multiplicative extension $\beta_1 \delta_t = \Sq^1 \delta_t = \delta_t^2$.

The group $C_{2^n}$ has a unique maximal elementary abelian subgroup, which is generated by $t^{2^{n-1}}$. By \cref{prop:compserre}, the $\sE_{(2)}$-homotopy limit spectral sequence for the group $C_{2^n}$ reduces to the LHSSS of the group extension
\begin{equation}
    C_2 \to C_{2^n} \to C_{2^{n-1}}.
\end{equation}
For $n=1$, the right hand side group is the trivial group, and therefore the LHSSS collapses to the $s=0$-line at $E_2$.

For $n > 1$, we write $a = \delta_{\tbar}$, $\beta_{n-1} a = \beta_{n-1} \delta_{\tbar}$ the $(n-1)$-th order Bockstein on $a$, and $x = \delta_{t^{2^{n-1}}}$. Because the the group $C_{2^n}$ is abelian, the local coefficient system is trivial, and the $E_2$-page is isormorhic to
\begin{equation}
    E_2^{*,*} \cong \Ffield_2[a, \beta_{n-1}a, x]/(a^2),
\end{equation}
with, for $n = 2$, the multiplicative extension $\beta_1 a = a^2$, and with $(s,t)$-degrees given by $|a| = (1,0)$, $|\beta_{n-1} a| = (2,0)$ and $|x| = (0,1)$.

Because the $\Ffield_2$-dimension of $\Hom(C_{2^n}, \Ffield_2)$ is 1, the dimension of $H^1(BC_{2^n})$ is 1. But $E_2$ has 2 classes in the $s+t = 1$-stem, so one of them has to die, which can only be $x$ for degree reasons, which can only possibly hit $\beta_{n-1} a$ with a $d_2$. The class $a$ is a permanent cycle for degree reasons, and hence
\begin{equation}
    E_3 \cong \Ffield_2[a, [x^2]]/(a^2).
\end{equation}
This has a vanishing line of height 2, and therefore the spectral sequence collapses at this stage.
\subsection{Abelian 2-groups}
Let $A$ be a finite abelian 2-group. Then $A$ is isomorphic to
\begin{equation}
    \prod_{j \in J} C_{2^{n_j}} \times \prod_{k \in K} C_2
\end{equation}
with $n_j \geq 2$, for some indexing sets $J$ and $K$.
We compute the $\sE_{(2)}$-homotopy limit spectral sequence converging to $H^*(BA)$. By the K\"unneth isomorphism, this ring is isomorphic to $\Ffield_2[x_j, \beta_{n_j}x_j, x_k]/(x_j^2)$, with degrees $|x_j| = |x_k| = 1$, and $|\beta_{n_j - 1}(x_j)| = 2$. 

Since $A$ has a unique maximal elementary abelian 2-subgroup $\prod_{j \in J} C_2 \times \prod_{k \in K} C_2$, this subgroup is normal and hence by \cref{prop:compserre} this reduces the $\sE_{(2)}$-homotopy limit spectral sequence to the LHSSS of the extension
\begin{equation}
    \prod_{j \in J} C_2 \times \prod_{k \in K} C_2 \to A \to \prod_{j \in J} C_{2^{n_j -1}}.
    \label{abgps:eq1}
\end{equation}

Write $t_j$ for the generator of $C_{2^{n_j}}$, and $t_k$ for the generator of $C_2$ indexed by $k$.

For all $j \in J$, the LHSSS corresponding to the extension
\begin{equation}
    C_2 \to C_{2^{n_j}} \to C_{2^{n_j - 1}}
    \label{abgps:eq3}
\end{equation}
has as $E_2 = \Ffield_2[a_j, \beta_{n_j-1}a_j, x_j]/(a_j^2)$ with $a_j = \delta_{\overline{t_j}}$, $\beta_{n_j-1}a_j = \beta_{n_j-1} \delta_{\overline{t_j}}$, $x_i = \delta_{t_i^{2^{n_i - 1}}}$, and for $n_j = 2$ with multiplicative extension $\beta_1 a_j = a_j^2$. The $(s,t)$-degrees are given by $|a_j| = (1,0)$, $|\beta_{n_j - 1} a_j| = (2,0)$ and $|x_j| = (0,1)$. Furthermore, $d_2 x_j = \beta_{n_j -1 } a_j$ and
\begin{equation}
    E_3 = E_\infty = \Ffield_2[a_j, [x_j^2]]/(a_j^2).
\end{equation}

Using the K\"unneth isomorphism, the $E_2$-page of the LHSSS of (\ref{abgps:eq1}) is given by
\begin{equation}
    E_2^{s,t} = \Ffield_2[a_j, \beta_{n_j - 1} a_j, x_j, x_k]/(a_j^2 \, | \, j \in J)
\end{equation}
with $(s,t)$-degrees as above, and for $n_j = 2$ with the multiplicative extensions $\beta_1 a_j = a_j^2$.

By naturality, the generators on $E_2$ supporting a $d_2$ are the  $x_j$, and
\begin{align}
    d_2 x_j = \beta_{n_{j-1}} a_j, \\
\end{align}
for all $j$. Hence
\begin{equation}
    E_3 = \Ffield_2[a_j, [x_j^2], x_k]/(a_j^2\, |\, j \in J).
\end{equation}
By the K\"unneth theorem, the stems of $E_3$ have Betti numbers equal to to the Betti numbers of the target, hence $E_3 = E_\infty$.

In particular, $\prod_{j \in J}a_j \in E^{\# J,0}_\infty$ is a non-zero element in filtration degree $\# J$, which shows that
\begin{equation}
    \exp_{\sE_{(2)}} \underline{H\Ffield_2} \geq \# J+1.
    \label{abgps:eq2}
\end{equation}

Conversely, for $j \in  J$, consider the projection of $A$ onto the $C_2$ in the $j$-th factor:
\begin{equation}
    p_j \colon A \to C_2 \cong O(1).
\end{equation}
The corresponding Euler class is $a_j \in H^1(BA;C_2)$. Since $a_j^2 = 0$, we get $\exp_{\All_{\ker p_j}} \underline{H\Ffield_2} = 2$. Now $\bigcap_{j  \in J}\ker p_j = \sE_{(2)}$, and hence by \cref{lem:expintersect}, $\exp_{\sE_{(2)}} \underline{H\Ffield_2} \leq \# J+1$. Combined with \cref{abgps:eq2} this shows that
\begin{Proposition}
    \label{prop:expelab}
    Let $A$ be a group isomorphic to
    \begin{equation}
      \prod_{j \in J} C_{2^{n_j}}  \times \prod_{k \in K} C_2
    \end{equation}
    with $n_j \geq 2$.
    Then
    \begin{equation}
        \exp_{\sE_{(2)}}\underline{H\Ffield_2}_A =  \#J + 1.
        \label{eq:expab2gps}
    \end{equation}
\end{Proposition}

\begin{Corollary}
    Let $G$ be a group with subgroup isomorphic to $\prod_{j \in J} C_{2^{n_j}}$ with all $n_j \geq 2$. Then
    \begin{equation}
        \exp_{\sE_{(2)}} \underline{H\Ffield_2} \geq \#J +1.
    \end{equation}
    \label{cor:p2rankexp}
\end{Corollary}
\begin{Proof}
    This follows from \cref{cor:expofsubgroup} and \cref{eq:expab2gps}.
\end{Proof}

\section{Dihedral groups}
\label{sec:dihgps}
\subsection{Introduction}
Denote by $D_{2^n}$ the dihedral group of order $2^n$. We apply the $\sE_{(2)}$-homotopy limit spectral sequence to compute the group cohomology of $D_{2^n}$. We first do computation for $D_{2^n} = D_8$, and in the last section we do the computation for $n \geq 4$.

The cohomology rings $H^*(BD_{2^n}; \Ffield_2)$ are all isomorphic to
\begin{equation}
    \Ffield_2[x,y,w]/(xy)
\end{equation}
with degrees $|x| = |y| = 1$, $|w| = 2$, see, e.g., \cite[Thm.\ 2.7]{adem}.
Since $\underline{H\Ffield_2}$ is $\sE_{(2)}$-nilpotent (\cite[Prop.\ 5.14]{mnn}), the spectral sequence will admit a vanishing line, which is bounded by the exponent $\exp_{\sE{(2)}} \underline{H\Ffield_2}$. In fact, the computation will show that the spectral sequence has a vanishing line of height 1, which is equivalent to the fact that $D_{2^n}$ has its cohomology detected on elementary abelian 2-subgroups, which is \cite[Lem.\ 4.6]{quillen71adams}.

Before proceeding with the calculation, we first identify the exponent of $D_{2n}$-Borel equivariant $\Ffield_2$-cohomology for all $n$.

\subsection{The exponent}
In this section we will show the following. 
\begin{Proposition}
\begin{equation}
    \exp_{\sE{(2)}}\underline{H\Ffield_2}_{D_{2n}} = 
    \begin{cases}
        2 & \text{for } n \equiv 0 \pmod{4}, \\
        1 & \text{otherwise}.
    \end{cases}
\end{equation}
\label{prop:d2nexp}
\end{Proposition}
\begin{Proof}
    Write $n = 2^v m$ with $2 \nmid m$. Then $D_{2^{v+1}} \subset D_{2n}$ is a 2-Sylow subgroup and by \cref{lem:psylowhfpexp} we are reduced to showing that
    \begin{equation}
        \exp_{\sE_{(2)}}\underline{H\Ffield_2}_{D_{2^k}}
        \begin{cases}
            = 2 & \text{ for } k \geq 3, \\
            = 1 & \text{ for } k = 2.
        \end{cases}
    \end{equation}
    For the case $k \geq 3$, we have $C_{2^{k-1}} \subset D_{2^k}$, and hence by \cref{cor:p2rankexp}, $\exp_{\sE_{(2)}} \underline{H\Ffield_2} \geq 2$. Combining this with the upper bound furnished by \cref{cor:expd8upperbound} below gives the claimed exponent.

    For the case $k = 2$, we have $D_4 \cong C_2 \times C_2$, and hence by \cref{lem:expallsubgroups} we get the claimed exponent.
\end{Proof}
\subsubsection{The upper bound}
In this section we will show using \cref{projbund:prop1} that for $D_{2^n}$ the dihedral group of order $2^n$ we have
\begin{equation}
    \exp_{\sE_{(2)}} \underline{H\Ffield_2}_{D_{2^n}} \leq 2,
\end{equation}
see \cref{cor:expd8upperbound} below.
To this end, consider the presentation
\begin{equation}
    D_{2^n} = \langle \sigma, \, \rho \, | \, \sigma^2 = \rho^{2^{n-1}} = e, \, \sigma \rho \sigma^{-1} = \rho^{-1} \rangle.
\end{equation}
For an angle $\theta$, denote the matrix representing rotation in $\Reals^2$ by $\theta$ by
\begin{equation}
    R_{\theta} = \begin{pmatrix}
        \cos\theta & - \sin\theta \\
        \sin\theta & \cos\theta
    \end{pmatrix}.
\end{equation}
Denote by $T$ the matrix representing reflection in the $x$-axis:
\begin{equation}
    T = \begin{pmatrix}
        1 & 0 \\
        0 & -1
    \end{pmatrix}.
\end{equation}
The group $D_{2^n}$ is the subgroup of $O(2)$ fixing the regular $2^{n-1}$-gon centered at the origin in $\Reals^2$ with one of its sides bisected in the right half-plane by the $x$-axis. The inclusion $D_{2^n} \to O(2)$ is on elements given by $\sigma \to T$, $\rho \to R_{2\pi/2^{n-1}}$. This yields a linear 2-dimension real representation of $D_{2^n}$, which we call $V$.  We apply \cref{projbund:prop1} to compute an upper bound on the exponent with respect to the minimal family $\sF$ containing the isotropy of $\Proj(V)$. This minimal family is determined in the next lemma.
\begin{Lemma}
    The minimal family containing the isotropy groups of the projectivation $\Proj(V)$ is $\sE_{(2)}$, the family of elementary abelian 2-groups in $D_{2^n}$.
    \label{lem:d8isotropy}
\end{Lemma}
\begin{Proof}
    The proof is by elementary linear algebra.
    Since $V$ is an orthogonal representation, a linear subspace $L$ spanned by a vector $v$ is fixed by an element $g$ of $D_{2^n}$ if and only if $g$ has eigenvalue $1$ or $-1$ (we say that $g$ has eigenvalue $\lambda$ if the matrix by which $g$ acts on $V$ has eigenvalue $\lambda$).

    First consider the rotations in $D_{2^n}$, that is the elements of the form $\rho^j$. We note that $R_{\theta}$ has characteristic polynomial
    \begin{equation}
        (\cos(\theta) - \lambda)^2 + \sin^2 \theta = 0,
    \end{equation}
    which for $\lambda = 1$ is 0 if and only if $\theta \equiv 0 \pmod{2\pi}$, and for $\lambda = -1$ is 0 if and only if $\theta = \pi \pmod{2\pi}$. Hence $\rho^j$ fixes a line if and only if $2 \pi j/2^{n-2} \equiv \pi \pmod{2\pi}$, that is if and only if $j \equiv 2^{n-2} \pmod{2^{n-1}}$. Hence $\rho^{2^{n-2}}$ is the only non-trivial rotation that fixes a line. Since $\rho^{2^{n-2}}$ acts via $R_{\pi}$, the elements $\rho^{2^{n-2}}$ fixes every line in $V$.

    For the lines fixed by a reflection $\sigma \rho^j$, we first note that, writing $e_1 = (\, 1 \quad 0\, )^t$, $TR_{\theta}$ has eigenvectors $R_{-\theta/2}e_1$ and $R_{-\theta/2 - \pi}e_1$, with eigenvalues $1$ and $-1$, respectively. Therefore these are, up to scalar, all eigenvectors of $TR_{\theta}$, since $V$ is 2-dimensional. Assume $\sigma \rho^j$ and $\sigma \rho^k$ fix a common line, spanned by a vector $v$. If $v$ had the same eigenvalue for $\sigma \rho^j$ and $\sigma \rho^k$, we would have
    \begin{equation}
        \frac{-2\pi j}{2^{n-1}} \equiv \frac{-2\pi k}{2^{n-1}} \pmod{2\pi},
    \end{equation}
    which is equivalent to $k \equiv j \pmod{2^{n-1}}$, in other words $\sigma \rho^j = \sigma \rho^k$.

    Alternatively, if $v$ has different eigenvalues for $\sigma \rho^j$ and $\sigma \rho^k$, then we have
    \begin{equation}
        \frac{-2\pi j}{2^{n-2}} \equiv \frac{- 2 \pi k}{2^{n-2}} - \pi \pmod{2 \pi},
    \end{equation}
    which is equivalent to $j - k \equiv 2^{n-2} \pmod{2^{n-1}}$. Hence if a line is fixed by a relfection $\sigma \rho^j$, it is fixed by exactly one more reflection, namely $\sigma \rho^{j+2^{n-2}}$. Note that
    \begin{equation}
        \sigma \rho^j \sigma \rho^{j+2^{n-2}} = \rho^{2^{n-2}} = \sigma \rho^{j+2^{n-2}} \sigma \rho^j.
    \end{equation}
    Therefore $\langle \sigma \rho^j, \sigma \rho^{j+2^{n-2}} \rangle$ is isomorphic to $C_2 \times C_2$.

    We conclude that a line in $V$ is fixed either by $\langle \rho^{2^{n-2}}\rangle$, or by $\langle \sigma \rho^{j}, \sigma \rho^{j + 2^{n-2}} \rangle$ for exactly one $j \pmod{2^{n-2}}$. All these subgroups are elementary abelian, proving the result.

\end{Proof}

\begin{Corollary}
    The $\sE_{(2)}$-exponent of $\underline{H\Ffield_2}_{D_{2^n}}$ satisfies
    \begin{equation}
        \exp_{\sE_{(2)}} \underline{H\Ffield_2}_{D_{2^n}} \leq 2.
    \end{equation}
    \label{cor:expd8upperbound}
\end{Corollary}
\begin{Proof}
    Immediate from \cref{projbund:prop1} and \cref{lem:d8isotropy} and the fact that $V$ is 2-dimensional.
\end{Proof}

\subsection{The spectral sequence}
\label{sec:dih:ss}
\subsubsection{Summary of the computation}
We will now compute the group cohomology ring $H^*(BD_8;\Ffield_2)$ using the $\sE_{(2)}$-homotopy limit spectral sequence converging to $H^*(BD_8;\Ffield_2)$. Since $\sE_{(2)}$ is the derived defect base of $\underline{H\Ffield_2}$ (\cite[Prop.\ 5.16]{mnn}), the spectral sequence will have a horizontal vanishing line at a finite page. In fact, by \cref{prop:d2nexp}, the spectral sequence will have a horizontal vanishing line of height 2 on the $E_3$-page. In particular, $E_3 = E_\infty$.

The computation, which follows the strategy of \cite[App.\ B]{mnn}, can be summarized as follows. The most work goes into computing the $E_2$-term. We find a homotopy final subcategory $\sO_{\rmf}$ of $\sO(D_8)_{\sE{(2)}}$ (using \cref{prop:intersecthofinal}). We then write a short exact sequence of coefficent systems ending in the constant coefficient system $\underline{\Integers}$ on $\sO_{\rmf}$, which lifts to a long exact sequence of spectral sequences, constituting of two Lyndon-Hochschild-Serre spectral sequences and the $\sE_{(2)}$-homotopy limit spectral sequence. We compute the $E_2$-pages of the LHSSS's, which allows us to determine the $E_2$-page of the $\sE_{(2)}$-homotopy limit spectral sequence.

The vanishing line of height 2 on $E_3$ implies that the only non-trivial differential is a $d_2$, which implies enough about the possible patterns of differentials on the $E_2$-page to determine $E_3 = E_\infty$ (see \cref{fig:D8HF2E2}).  
\subsubsection{The orbit category}
A generating set of morphism for the orbit category $\sO(D_8)_{\sE_{(2)}}$ is given by
\begin{figure}[H]
	\centering
	\includegraphics[width=\textwidth]{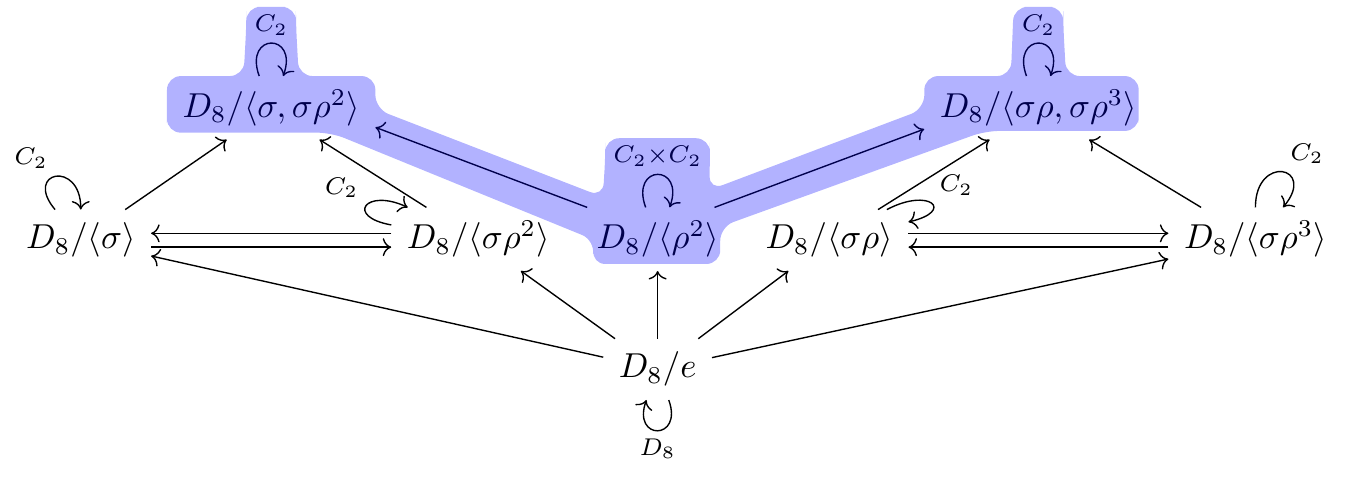}
	\caption{The orbit category $\sO(D_8)_{\sE_{(2)}}$. The highlighted subcategory is homotopy final.}
\end{figure}
Write
\begin{align}
	Z & = Z(D_8) = \langle \rho^2 \rangle, \\
	H_1 & = \langle \sigma, \sigma \rho^2 \rangle, \\
	H_2 & = \langle \sigma \rho, \sigma \rho^3 \rangle. \\
	\label{d8hf2:eq3}
\end{align}
With these $H_i$ and $Z$, we apply the docomposition of \cref{sec:absplit}. Therefore we next compute the LHSSS's obtained from the extensions $H_i \to D_8 \to C_2$ and $Z \to D_8 \to C_2 \times C_2$.
\subsubsection{The LHSSS of $H_i \to D_8 \to C_2$}
We compute the $E_2$-page of the LHSSS of the extension $H_1 \to D_8 \to C_2$, the case for $H_2$ follows from the automorphism $\sigma \leftrightarrow \sigma \rho$ of $D_8$.

The group $C_2 = \langle \rhobar \rangle$ acts on $H_1$ by
\begin{align}
    \rho (\sigma) \rho^{-1} = \sigma \rho^2, \\
    \rho (\sigma \rho) \rho^{-1} = \sigma.
\end{align}
Hence it acts on $H^*(H_1)$ by interchanging $\delta_{\sigma}$ and $\delta_{\sigma \rho^2}$. Write $u_{1,1} = \delta_{\sigma} + \delta_{\sigma \rho^2}$ and $u_{1,2} = \delta_{\sigma} \delta_{\sigma \rho^2}$ for the first and second elementary symmetic polynomial. It is well known that the invariants of a polynomial algebra of the action of the symmetric group on the generators is given by the polynomial algebra on the symmetric polynomials.

Taking the minimal periodic $C_2$-resolution of $\Integers$, applying the hom-functor of $C_2$-modules $\Hom_{C_2}(-, H^*(H_1))$, and taking cohomology, we see that the $s=0$-row of the LHSSS is given by the $C_2$-invariants
\begin{equation}
    \Ffield_2[u_{1,1}, u_{2,1}].
\end{equation}
Write $a_1 = \delta_{\rhobar}$. The image of the differential in the long exact sequence computing $E_2^{s,*}$ maps a polynomial to its $C_2$-symmetrization:
\begin{equation}
    P(\delta_{\sigma}, \delta_{\sigma \rho^2}) \mapsto P(\delta_{\sigma}, \delta_{\sigma \rho^2}) + P(\delta_{\sigma \rho^2}, \delta_{\sigma}).
\end{equation}
The image are the symmetrized  polynomials, which is $\Ffield_2[u_{1,1},u_{2,1}]\{u_{1,1}\}$. The rows $E_2^{s,*}$ for $s>1$ are then given by the $C_2$-invariants modulo the symmetrized polynomials, which are
\begin{equation}
    \Ffield_2[u_{2,1}]\{a_1^s\}.
\end{equation}
Assembling the rows together gives the $E_2$-page
\begin{equation}
    E_2 = \Ffield_2[a_1, u_{1,1}, u_{2,1}]/(a_1u_{1,1}),
\end{equation}
with $(s,t)$-degrees $|a_1| = (1, 0)$, $|u_{1,1}| = (0,1)$ and $|u_{2,1}| = (0,2)$.

For the extension $H_2 \to D_8 \to C_2$, we write $u_{1,2} = \delta_{\sigma \rho} + \delta_{\sigma \rho^3}$, $u_{2,2} = \delta_{\sigma \rho} \delta_{\sigma \rho^3}$, $a_2 = \delta_{\rhobar}$, and obtain as an $E_2$-page
\begin{equation}
    E_2 = \Ffield_2[a_2, u_{1,2}, u_{2,2}]/(a_2u_{1,2}).
\end{equation}
Because the number of classes in the $1$-stem and $2$-stem equals the Betti numbers of $H^1(BD_8)$ respectively $H^2(BD_8)$ (\cite[Thm.\ 2.7]{adem}), the LHSSS collapses at $E_2$.
\begin{figure}[H]
    \centering
    \includegraphics[width=0.8\textwidth]{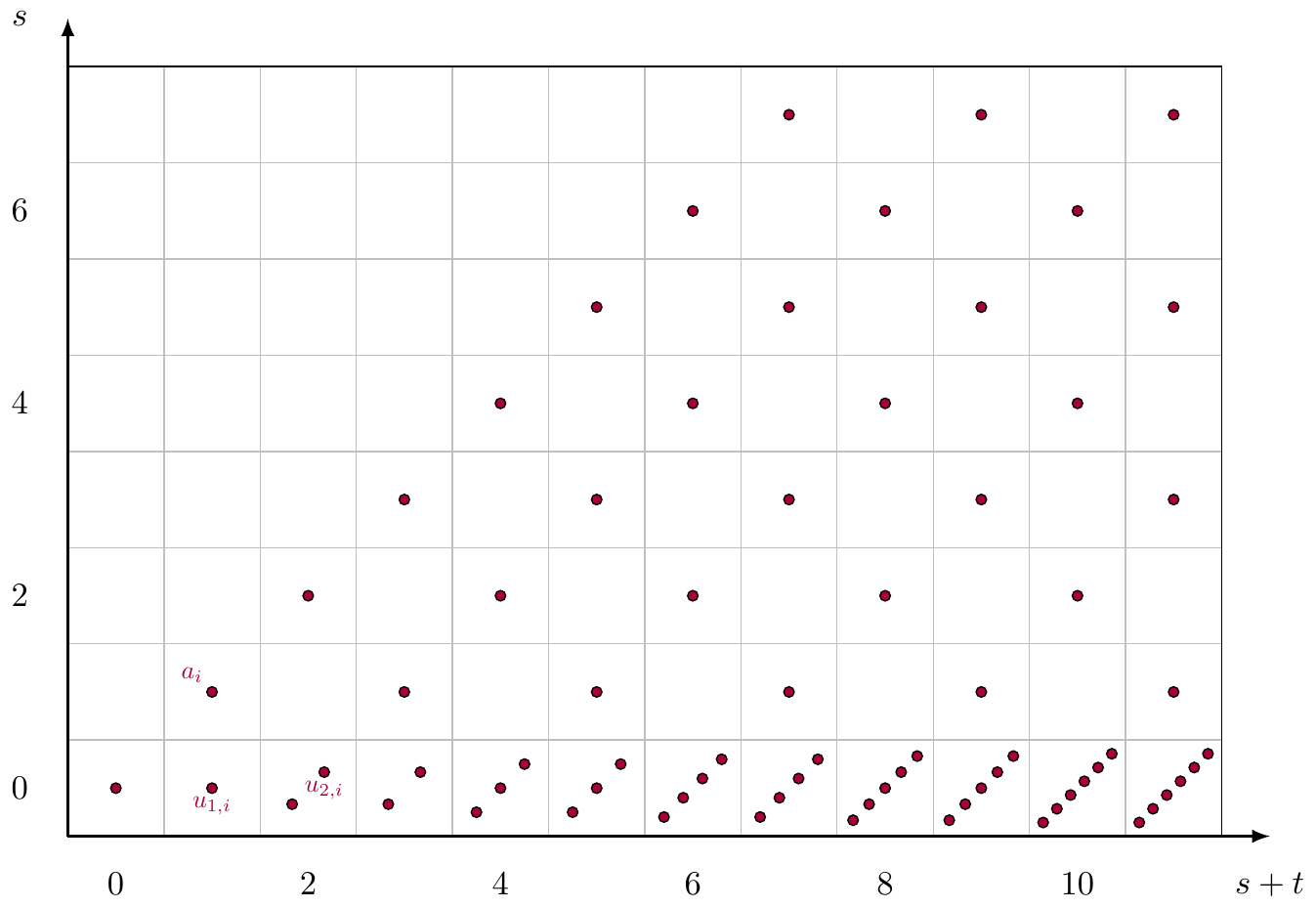}
    \caption{The $E_2 = E_\infty$-page of the LHSSS obtained from the extension $H_i \to D_8 \to C_2$.}
    \label{fig:vextss}
\end{figure}
\subsubsection{The LHSSS of $Z \to D_8 \to C_2 \times C_2$}
The extension
\begin{equation}
    \overset{\rho^2}{Z} \to D_8 \to \overset{\sigmabar}{C_2}\times \overset{\overline{\sigma \rho}}{C_2}
\end{equation}
is central, and therefore the associated LHSSS has trivial coefficient system. Writing $b = \delta_{\sigmabar}$, $c = \delta_{\overline{\sigma \rho}}$ and $x = \delta_{\rho^2}$, the $E_2$-page then is
\begin{equation}
    E_2 = \Ffield_2[b,c,x],
\end{equation}
with $(s,t)$-degrees $|b| = |c| = (1,0)$, $|x| = (0,1)$. The class $x$ supports a $d_2$ with $d_2(x) = bc$ (\cite[\S IV.2]{adem}), and the LHSSS collapses at $E_3$.
\begin{figure}[H]
    \centering
    \includegraphics[width=0.8\textwidth]{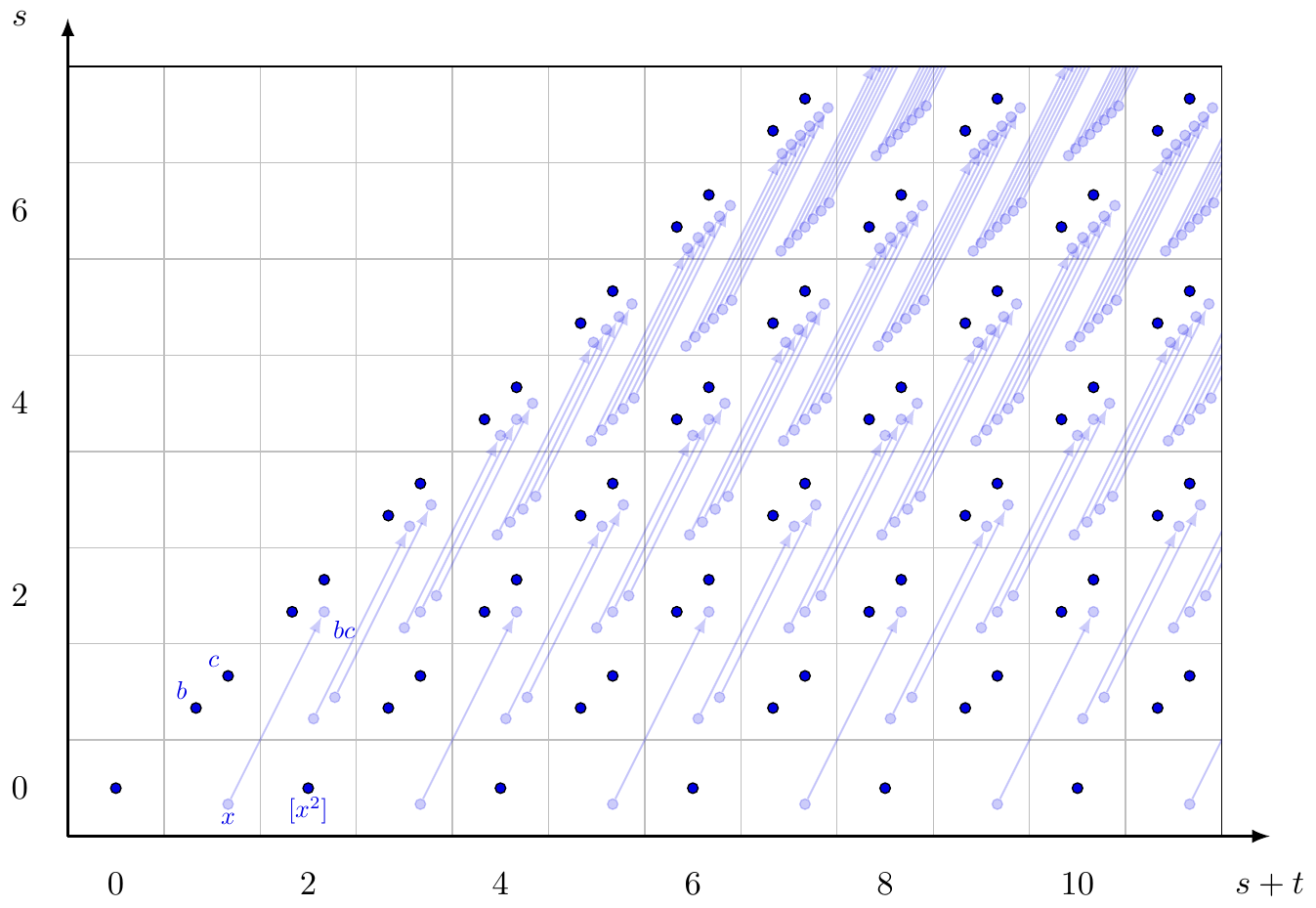}
    \caption{The $E_2$-page of the LHSSS obtained from the extension $Z \to D_8 \to C_2 \times C_2$.}
    \label{fig:zextss}
\end{figure}
\subsubsection{The map $j$}
In this subsection we prove the following proposition, which describes the map $j$ from (\ref{absplit:eq3}) from \cref{sec:absplit}.
\begin{Proposition}
    The map $j$ from (\ref{absplit:eq3}) is given on cohomology classes by
    \begin{align}
        u_{1,i} & \mapsto 0, \\
        u_{2,i} & \mapsto x^2, \\
        a_1 & \mapsto c, \\
        a_2 & \mapsto b.
    \end{align}
\label{d8:mapj}
\end{Proposition}
\begin{Proof}
 The map $j$ is in cohomology determined by the inclusions $Z \to H_i$ and the corresponding quotient maps $D_8/Z \to D_8 /H$.

 The map $Z \to H_1$ is on elements given by $\rho^2 \mapsto \sigma \cdot \sigma \rho^2$, therefore on cohomology classes by
 \begin{align}
     \delta_{\sigma} & \mapsto \delta_{\rho^2} = x, \\
     \delta_{\sigma \rho^2} & \mapsto \delta_{\rho^2} = x, \\
     \intertext{and therefore}
     u_{1,1} = \delta_{\sigma} + \delta_{\sigma \rho^2} & \mapsto 0, \\
     u_{2,1} & \mapsto x^2.
 \end{align}

The map $D_8/Z \cong C_2 \times C_2 \to D_8/H_1 \cong C_2$ is on elements given by $\sigmabar \mapsto \ebar$ and $\overline{\sigma \rho} \mapsto \rhobar$, and therefore on cohomology by
\begin{equation}
    a_1 = \delta_{\rhobar} \mapsto \delta_{\overline{\sigma \rho}} = c.
\end{equation}
This completes the formulas for $i=1$, the formulas for $i=2$ are obtained from the automorphism $\sigma \leftrightarrow \sigma \rho$ of $D_8$.
\end{Proof}
\subsubsection{The $E_2$-page}
We write $u_{1,1} = (u_{1,1}, 0)$, $u_{1,2} = (0, u_{1,2})$ and $u = (u_{2,1}, u_{2,2})$.  \cref{d8:mapj} shows that the kernel of $j$ is concentrated in the $s=0$-line, and equals
\begin{equation}
    \ker(j) = \Ffield_2[u_{1,1}, u_{1,2}, u]/(u_{1,1} u_{1,2}),
\end{equation}
with bidegrees given by $|u_{1,1}| = |u_{1,2}| = (0,1)$, and $|u| = (0, 2)$.
Moreover, it shows that the image of $j$ is given by
\begin{equation}
    \Im(j) = \Ffield_2[c, x^2]\{c\} \oplus \Ffield_2[b, x^2]\{b\} \oplus \Ffield_2[x^2].
\end{equation}
The split exact sequence (\cref{prop:absplit}) then implies that the $E_2$-page of the $\sE_{(2)}$-homotopy limit spectral sequence is given by
\begin{figure}[H]
	\centering
	\includegraphics[width=0.9\textwidth]{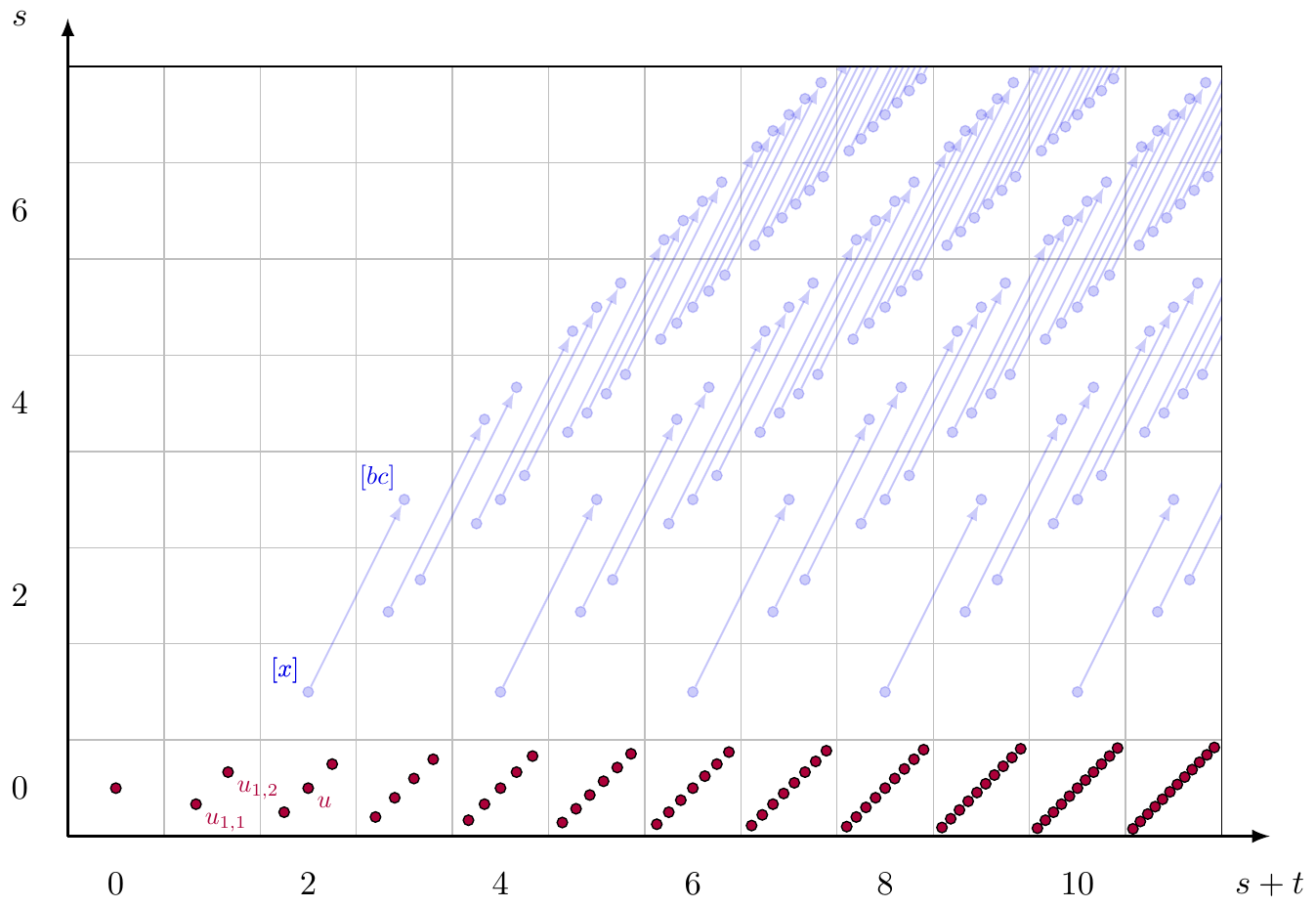}
    \caption{The $E_2$-page. Red classes come from the $B$-summand, blue classes from the $A$-summand. The labeled red classes generated the red summand. The entire page is not finitely genereted as a ring, but we labeled two blue classes for ease of comparison with \cref{fig:zextss}.}
    \label{fig:D8HF2E2}
\end{figure}
Because the $\sE_{(2)}$-exponent is 2, there can only be non-zero differentials on $E_2$, and on $E_3$ there must be a vanishing line of height 2. The differentials are coming from the LHSSS of the extension $Z \to D_8 \to C_2$, by naturality, see \cref{fig:zextss}.
So there is in fact a horizontal vanishing line of height 1 on the $E_3 = E_\infty$-page.
\subsubsection{Multiplicative extensions}
Since the $\sE_{(2)}$-homotopy limit spectral sequence is concentrated in the $s=0$-line on $E_\infty$, there are no multiplicative extension problems to solve, and
\begin{align}
    H^*(BD_8) & \cong E_\infty^{0,*} \\
    & = \Ffield_2[u_{1,1}, u_{1,2}, u]/(u_{1,1} u_{1,2}).
\end{align}
\subsubsection{The lower bound illustrated}
For $G= D_{2n}$ with $n \equiv 0 \pmod{4}$ and $n \geq 4$ we established in \cref{prop:d2nexp} that $\exp_{\sE_{(2)}} \underline{H\Ffield_2} = 2$. We illustrate this exponent with a $\sE_{(2)}$-homotopy limit spectral sequence which has a vanishing line of height 2 on the $E_\infty$-page.

Let $C_n = \langle \rho \rangle \subset D_{2n}$ act on $z \in S^1 \subset \Complex$ via
\begin{equation}
	\rho \cdot z = e^{2 \pi i /n} z.
\end{equation}
Let $X$ be the $D_{2n}$-space given by inducing the $C_n$-space $S^1$ up to $D_{2n}$:
\begin{equation}
	X := D_{2n} \times_{C_n} S^1
    \label{eq:inducedcircle}
\end{equation}
The space $X$ looks like
\begin{figure}[H]
	\centering
	\includegraphics[width=0.3\textwidth]{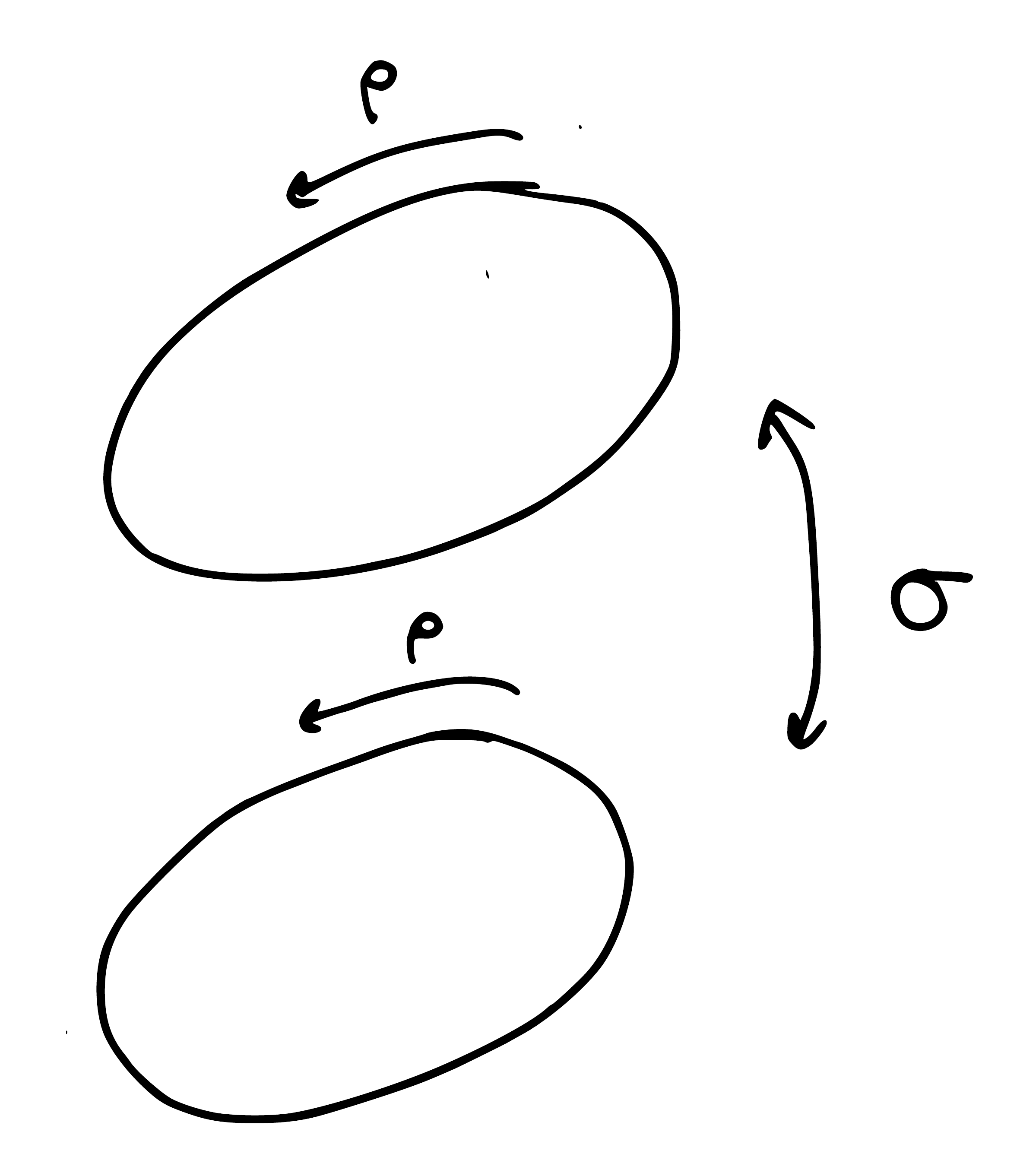}
	\caption{The free $D_{2n}$-space $X$.}
\end{figure}
The $D_{2n}$-action on $X$ is evidently free. As coset representatives of $D_{2n}/C_n$ we pick $\{e,\sigma\}$. The underlying non-equivariant space is the coproduct $S_e^1 \sqcup S_\sigma^1$, where the subscripts are decorations for the corresponding coset. Denoting the elements of $X$ by $(e,z_1)$ and $(\sigma,z_2)$ (i.e., by choosing the representatives of the elements of the quotient in \cref{eq:inducedcircle} which have as their first coordinate our chosen coset representatives), the $D_{2n}$-action on $X$ is explicitly given by
\begin{align}
    \rho \cdot (e,z) & = (e,\rho z), \\
    \rho \cdot (\sigma,z) & = (\sigma, \rho^{-1} z), \\
    \sigma \cdot (e,z) & = (\sigma,z), \\
    \sigma \cdot (\sigma,z) & = (e,z).
\end{align}
In particular, $\sigma \rho (\sigma,z) = (e,\rho^{-1} z)$, whereas $\rho \sigma (\sigma,z) = (e,\rho z)$, hence the actions of $\sigma$ and $\rho$ do not commute.

Since $X$ is a free $D_{2n}$-space, its homotopy orbits coincide with the actual orbits:
\begin{equation}
    X\sslash D_{2n} \simeq \frac{X \times E D_{2n}}{D_{2n}}  \simeq X/D_{2n} \simeq S^1/C_n \simeq S^1.
\end{equation}
Hence the Borel equivariant $\Ffield_2$-cohomology of $X$ is:
\begin{equation}
	H_{D_{2n}}^*(X;\Ffield_2) \cong H^*(S^1;\Ffield_2) \cong \Ffield_2[x]/(x^2).
	\label{d2nlowbound:eq1}
\end{equation}

Consider the $\sE_{(2)}$-homotopy limit spectral sequence of $X$:
\begin{equation}
	E_2^{s,t} = \sideset{}{^s}\lim_{\sO(G)^{\op}_{\sE_{(2)}}} H_{(-)}^t(X,\Ffield_2) \Rightarrow H_{D_{2n}}^{s+t}(X;\Ffield_2).
\end{equation}
We will show that its $E_\infty$-page has an element in filtration degree 1, which implies the desired lower bound on the exponent. Note that, given \cref{d2nlowbound:eq1} and the fact that the functors $H^t_{(-)}(X;\Ffield_2)$ are 0 for $t <0$, only two $E_\infty$-pages are possible:
\begin{figure}[H]
	\centering
	\begin{minipage}{0.4\textwidth}
		\centering
		\includegraphics[width=0.7\linewidth]{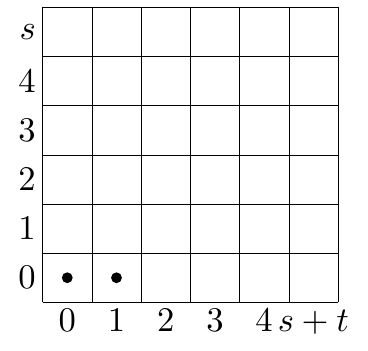}
		\caption{First possibility for $E_\infty$.}
	\end{minipage}%
	\begin{minipage}{0.2\textwidth}
		\centering
		and
	\end{minipage}%
	\begin{minipage}{0.4\textwidth}
		\centering
		\includegraphics[width=0.7\linewidth]{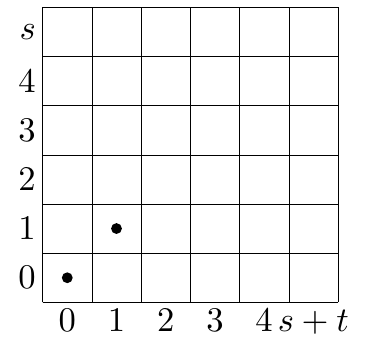}
		\caption{Second possibility for $E_\infty$.}
	\end{minipage}
\end{figure}
We will show that the $E_\infty$-page is the second one, by showing that $x \in H^1_{D_{2n}}(X;\Ffield_2)$ is in filtration $\geq 1$, which is equivalent to showing that $x$ is in the kernel of the edge map
\begin{equation}
	H_{D_{2n}}^1(X;\Ffield_2) \to \sideset{}{^0}\lim_{\sO(G)^{\op}_{\sE_{(2)}}} H_{(-)}^1(X;\Ffield_2).
\end{equation}
In other words, we must show that for all elementary abelian 2-groups $A \subset D_{2n}$ we have that the map induced by the quotient
\begin{equation}
	q \colon X\sslash A \to X\sslash D_{2n}
	\label{d2nlowbound:eq2}
\end{equation}
in cohomology
\begin{equation}
	q^* \colon H^1_{D_{2n}}(X;\Ffield_2) \to H^1_A(X;\Ffield_2)
\end{equation}
satisfies $q^*(x) = 0$. For this, we may restrict to the maximal elementary abelian 2-groups in $D_{2n}$, because all elementary abelian subgroups are a subgroup of a maximal elementary abelian subgroup. Hence if an element restricts to $0$ on all maximal elementary abelian subgroup, it also restricts further down to $0$ on all elementary abelian subgroups.

The elements Nf order $2$ of $D_{2n}$ are the reflections $\sigma \rho^i$, and the generator of the center $\rho^{n/2}$. Two reflections $\sigma \rho^i$ and $\sigma \rho^j$ commuting is equivalent to $\rho^{j-i} = \rho^{i-j}$, which is equivalent to $i \equiv j \pmod{n/2}$. Since $\sigma \rho^i \sigma \rho^{i+n/2} = \rho^{n/2}$, we see that all maximal elementary abelian subgroups of $D_{2n}$ are given by $\langle \sigma \rho^i,\, \sigma \rho^{i+n/2} \rangle$ with $0 \leq i <n/2$.

First dividing out $\sigma \rho^i$ identifies
\begin{equation}
    (e,z) \leftrightarrow (\sigma, \rho^i z).
\end{equation}
Hence every element of $X/\langle \sigma \rho^i \rangle$ is represented by uniquely by a $(e,z)$, and this space is isomorphic to $S^1$. We suppress the coset coordinate from now on.

Dividing this out by the central element $\rho^{n/2}$ identifies
\begin{equation}
    z \leftrightarrow \rho^{n/2}z,
\end{equation}
and every element becomes represented uniquely by a $e^{\theta i}$ with $0 \leq \theta < \pi$.

Finally, dividing out by the entire group $D_{2n}$ identifies
\begin{equation}
    z \leftrightarrow \rho z,
\end{equation}
hence every element of $X/D_{2n}$ is uniquely represented by a $e^{\theta i}$ with $0 \leq \theta < 2\pi/n$.

The map $S^1 \cong X/\langle \sigma \rho^i, \, \sigma \rho^{i+n/2}\rangle \to X/D_{2n} \cong S^1$ is given by
\begin{equation}
    \theta \pmod{\pi} \mapsto \theta \pmod{2\pi /n},
\end{equation}
which is a $(\frac{n}{2} : 1)$-cover. Hence the map $q^*$ is given by multiplication by $n/2$.

We assumed $n \equiv 0 \pmod{4}$, hence that $n/2 \equiv 0 \pmod{2}$, and therefore $q^*(x) = 0$, as desired.

\subsection{The dihedral groups $D_{2^n}$ with $n \geq 4$}
For the dihedral groups $D_{2^n}$ with $n \geq 4$ we will not compute the $\sE_{(2)}$-homotopy limit spectral sequence, but instead we will compute the $\sF$-homotopy limit spectral sequence for some family $\sF$ which strictly contains $\sE_{(2)}$, to be described next.

For $n \geq 4$, the dihedral group has exactly two subgroups isomorphic to $D_{2^{n-1}}$, which we denote by $H_1 = \langle \sigma, \sigma \rho^2 \rangle$ and $H_2 = \langle \sigma \rho, \sigma \rho^3 \rangle$. Additionally we write $J = H_1 \cap H_2 = \langle \rho^2 \rangle \cong C_{2^{n-2}}$. Write $\sF$ for the family $\sF = \All_{H_1} \cap \All_{H_2}$. Then $\sF \supset \sE_{(2)}$, and therefore $\underline{H\Ffield_2}_{D_{2^n}}$ is $\sF$-nilpotent. In this subsection we will compute the $\sF$-homotopy limit spectral sequence for $D_{2^n}$. This spectral sequence will turn out to be additively isomorphic to the $\sE_{(2)}$-homotopy limit spectral sequence for $D_8$. The computation will also be almost identical.
\subsubsection{The decomposition}
We apply the decomposition of \cref{sec:absplit} with the $H_i$ as in the previous section and $K = J$. Therefore we proceed by computing the LHSSS's for te extensions $H_i \to D_{2^n} \to C_2$ and $J \to D_{2^n} \to C_2 \times C_2$.
\subsubsection{The LHSSS of $D_{2^{n-1}} \to D_{2^n} \to C_2$}
We consider the case $D_{2^{n-1}} = H_1$, the other case being obtained by the automorphism $\sigma \leftrightarrow \sigma \rho$. Write $x_1 = \delta_{\overline{\sigma}}$ and $y_1 = \delta_{\overline{\sigma \rho^2}}$ for the classes in $H^1(BD_{2^{n-1}}) = \Hom(D_{2^{n-1}}, \Ffield_2)$. Let $w_1$ be the unique class in $H^2(BD_{2^{n-1}})$ that classifies $C_2 \to D_{2^n} \to D_{2^{n-1}}$ and that restricts to zero on both $\langle \sigma \rangle$ and $\langle \sigma \rho^2\rangle$ (see, e.g., \cite[Thm.\ 2.7]{adem}). Then $C_2 = \langle \rhobar \rangle$ acts on $H_1^{\ab} = H_1/\langle \rho^4 \rangle \cong C_2 \times C_2$ by
\begin{align}
    \rho \sigma \rho^{-1} & \equiv \sigma \rho^{-2} \equiv \sigma \rho^2 \mod{\langle \rho^4 \rangle}, \\
    \rho \sigma \rho^2 \rho^{-1} & \equiv \sigma \mod{\langle \rho^4 \rangle},
\end{align}
hence it interchanges $x_1$ and $y_1$. The class $w_1$ is uniquely determined by conditions which are preserved by the action of $C_2$, hence is fixed.

Write $u_{1,1} = x_1 + y_1$ and $a_1 = \delta_{\rhobar}$. Taking a projective $C_2$-resolution of $\Integers$, we see that the $s=0$ row of the $E_2$-page of the LHSSS is given by the $C_2$-invariants
\begin{equation}
    E_2^{0,*} \cong \Ffield_2[u_{1,1}, w],
\end{equation}
and the rows for $s \geq 1$ are given by
\begin{equation}
    E_2^{s,*} \cong \Ffield_2[w]\{a^s \}.
\end{equation}
Assembling all rows together gives the $E_2$-page
\begin{equation}
    E_2^{*,*}(H_1) = \Ffield_2[a, u_{1,1}, w_1]/(a_1 u_{1,1}),
\end{equation}
with $(s,t)$-degrees $|a| = (1,0)$, $|u_{1,1}| = (0,1)$, and $|w_1| = (0,2)$.

Consider now the extension $H_2 \to D_{2^n} \to C_2$. Write $x_2 = \delta_{\overline{\sigma \rho}}$, $y_2 = \delta_{\overline{\sigma \rho^3}}$, $w_2$ the unique class in $H^2(BH_2)$ classifying  $C_2 \to D_{2^n} \to D_{2^{n-1}}$ and restricting to zero on $\langle \sigma \rho \rangle$ and $\langle \sigma \rho^3 \rangle$, $a_2 = \delta_{\sigmabar}$ and $u_{1,2} = x_2 + y_2$. Then by applying the automorphism of $D_{2^n}$ determined by $\sigma \leftrightarrow \sigma \rho$ we get that the $E_2$-page of the LHSSS of the extension $H_2 \to D_{2^n} \to C_2$ is given by
\begin{equation}
    E_2^{*,*}(H_2) = \Ffield_2[a_2, u_{1,2}, w_1](a_2 u_{1,2}).
\end{equation}
\subsubsection{The LHSSS of $C_{2^{n-2}} \to D_{2^n} \to C_2 \times C_2$}
The extension
\begin{equation}
    J \cong C_{2^{n-2}} \to D_{2^n} \to C_2 \times C_2
\end{equation}
gives a LHSSS with trivial coefficient system, for $H^*(BC_{2^{n-2}})$ has only trivial automorphisms. Writing $b = \delta_{\sigmabar}$, $c= \delta_{\overline{\sigma \rho}}$, $z = \delta_{\rho^2}$ and $\beta_{n-2}(z) = \beta_{n-2}(\delta_{\rho^2})$ (the $(n-2)$-th order Bockstein) we get
\begin{equation}
    E_2^{*,*}(J) = \Ffield_2[b, c, z, \beta_{n-2}(z)]/(z^2).
\end{equation}
\subsubsection{The map $j$}
We now determine the effect of the map $j$ from (\ref{absplit:eq3}) from \cref{sec:absplit}.
\begin{Proposition}
    The map $j$ in the long exact sequence (\ref{absplit:eq3}) is given by
    \begin{align}
        u_{1,i} & \mapsto 0,  \\
        w_i & \mapsto \beta_{n-2}(z), \\
        a_1 & \mapsto c, \\
        a_2 & \mapsto b.
    \end{align}
    \label{d2n:mapj}
\end{Proposition}
\begin{Proof}
    The map $j$ is determined by the inclusions $J \to H_i$ and the corresponding quotient maps $D_{2^n}/J \to D_{2^n} /H_i$.

    For $i=1$, the inclusion $J = \langle \rho^2 \rangle \hookrightarrow H_1 = \langle \sigma, \sigma \rho^2 \rangle$ is on elements given by $\rho^2 \mapsto \sigma \cdot \sigma \rho^2$, hence on abelianizations by $\overline{\rho^2} \mapsto \sigmabar \cdot \overline{\sigma \rho^2}$, hence on degree-1 cohomology classes by
    \begin{align}
        x_1 = \delta_{\sigmabar} & \mapsto \delta_{\rho^2} = z, \\
        y_1 = \delta_{\overline{\sigma \rho^2}} & \mapsto \delta_{\rho^2} = z,
        \intertext{and hence}
        u_{1,1}  = x_1 + y_1 & \mapsto 0.
    \end{align}
    To see how the class $w_1$ restricts, we recall that it classifies the extension $C_2 \to D_{2^n} \to D_{2^{n-1}}$, which pulls back along $C_{2^{n-1}} \to D_{2^{n-1}}$ as
    \begin{equation}
        \begin{tikzcd}
            C_2 \arrow[d] \arrow[r] & C_{2^{n-1}} \arrow[d,hook] \arrow[r] & C_{2^{n-2}} \arrow[d, hook] \\
            C_2 \arrow[r] & D_{2^n} \arrow[r] & D_{2^{n-1}}
        \end{tikzcd},
    \end{equation}
    which is non split, and therefore the class $w_1$ restricts to a non zero class, the only possibility being
    \begin{equation}
        w_1 \mapsto \beta_{n-2}(z).
    \end{equation}

    The map $D_{2^n} / J \to D_{2^n}/H_1$ is on elements given by $\overline{\sigma \rho} \mapsto \rhobar$ and $\sigmabar \mapsto e$, and therefore on cohomology classes by
    \begin{equation}
        a_1 = \delta_{\rhobar} \mapsto \delta_{\overline{\sigma \rho}} = c.
    \end{equation}
    The formulas for $i=2$ are obtained by applying the automorphism $\sigma \leftrightarrow \sigma \rho$. Observe that in particular this automorphism interchanges $H_1$ and $H_2$, hence $b$ and $c$.
\end{Proof}
\subsubsection{The kernel of $j$}
We write $u_{1,1} = (u_{1,1},0)$, $u_{1,2} = (0, u_{1,2})$ $w = (w_1, w_2)$. 
Then \cref{d2n:mapj} shows that the kernel of $j$ is concentrated in the 0-line and equals 
\begin{equation}
    \ker(j) = \Ffield_2[u_{1,1}, u_{1,2}, w]/(u_{1,1} u_{1,2}).
\end{equation}
\subsubsection{The image of $j$}
\cref{d2n:mapj} shows that the image of $j$ is given by
\begin{equation}
    \Ffield_2[b, c, \beta_{n-2}(z)].
\end{equation}
\subsubsection{The $E_2$-page}
This gives an $E_2$-page which additively is isomorphic to the one in \cref{fig:D8HF2E2}. But also all the differentials in \cref{fig:D8HF2E2} need to happen in this case, because it is the only possibility for a vanishing line of height 2 on $E_3$, which must occur by \cref{cor:expd8upperbound}.
\subsubsection{Multiplicative structure}
Because the $E_3 = E_\infty$-page is concentrated in the $s=0$-line, there are no multiplicative extension problems, and we get
\begin{equation}
    H^*(BD_{2^{n}}; \Ffield_2) \cong \Ffield_2[u_{1,1}, u_{1,2}, w]/(u_{1,1}u_{1,2}).
\end{equation}

\section{Quaternion groups}
\label{sec:quatgps}
The generalized quaternion group $Q_{2^n}$ of order $2^n$ is the finite subgroup of quaternionic space $\HH$ generated multiplicatively by the elements of unit length $\{ e^{2\pi i / 2^{n-1}},\, j\}$ (see, e.g., \cite[XII.\S7]{cartaneilenberg}). Denoting these generators by $r$ and $s$, respectively, one gets the presentation
\begin{equation}
    Q_{2^n} = \langle r,\, s\, | \, r^{2^{n-2}} = s^2,\, rsr = s \rangle.
\end{equation}
for all $n \geq 3$. The subgroup generated by $\langle s^2 \rangle$ is central and of order 2, which makes $Q_{2^n}$ into a central extension (\cite[IV.2]{adem}). 
\begin{equation}
    C_2 \to Q_{2^n} \to D_{2^{n-1}}.
    \label{eq:Q2ncentral}
\end{equation}
See \cref{q8:lhsss} and \cref{lem:Q2ng3LHS} below for a recollection on this cohomology ring and the LHSSS computing it.
\subsection{An upper bound on the exponent}
In this subsection we will prove the following upper bound on the $\sE_{(2)}$-exponent.
\begin{Proposition}
    The $\sE_{(2)}$-exponent satisfies
    \begin{equation}
        \exp_{\sE_{(2)}} \underline{H\Ffield_2}_{Q_{2^n}} \leq 4.
    \end{equation}
    \label{q2n:expupper}
\end{Proposition}
\begin{Proof}
    The proof is a straightforward adaption of the argument in \cite[Ex.\ 5.18]{mnn}.

    Let $\HH \cong \Reals^4$ be the 4-dimensional real representation coming from the embedding $Q_{2^n} \hookrightarrow \HH$. This is a free action, and restricts to a free action on $S^3$. The subgroup $\langle \pm 1 \rangle$ is central, and therefore $Q_{2^n}/\langle \pm 1 \rangle$ acts on $S^3/\langle \pm 1 \rangle =  \Proj(\Reals^4)$ with isotropy in $\langle \pm 1 \rangle$. The result now follows from \cref{projbund:prop1}.
\end{Proof}

\subsection{The $\sE_{(2)}$-homotopy limit spectral sequence}
The group $Q_{2^n}$ has a unique maximal elementary abelian 2-subgroup, given by $\langle \pm 1 \rangle$. Therefore, by \cref{prop:compserre}, the $\sE_{(2)}$-homotopy limit spectral sequence reduces to the LHSSS corresponding to the central extension
\begin{equation}
    \langle -1 \rangle = Z(Q_{2^n}) \to Q_{2^n} \to D_{2^{n-1}}.
\end{equation}
This spectral sequence is well known, we recall it here for convenience.
\begin{Proposition}[{{see, e.g., \cite[Lem.\ 2.10]{adem}}}]
    For $n=3$, the $E_2$-page of the LHSSS associated to the extension (\ref{eq:Q2ncentral}) is
    \begin{equation}
        E_2^{*,*} = H^*(BC_2^{\times 2};\Ffield_2) \otimes_{\Ffield_2} H^*(BC_2;\Ffield_2) \cong \Ffield_2[x,y,e],
    \end{equation}
    with all the generators in degree $1$.
    The differentials are generated under the Leibniz rule by
    \begin{align}
        d_2(e) & = x^2 + xy + y^2, \\
        d_3(e^2) & = \Sq^1(x^2 + xy +y^2) = x^2y + xy^2.
    \end{align}
    The spectral sequence collapses at $E_4$ with a vanishing line of height 4, with resulting cohomology ring
    \begin{equation}
        H^*(BQ_8;\Ffield_2) \cong \Ffield_2[x,y,e_4]/(x^2 + xy + y^2,\, x^2y + xy^2).
    \end{equation}
    with degrees given by $|x| = |y| = 1$, and $|e_4| = 4$.
    \label{q8:lhsss}
\end{Proposition}
\begin{Proposition}[{{see, e.g., \cite[Lem.\ 2.11]{adem}}}]
    Denote the classes $x= \delta_{\sigmabar}$ and $y= \delta_{\overline{\sigma\rho}}$.

    For $n \geq 4$, the $E_2$-page of the LHSSS associated to the extension (\ref{eq:Q2ncentral}) is given by
    \begin{equation}
        E_2^{*,*} = H^*(BD_{2^{n-1}};\Ffield_2) \otimes_{\Ffield_2} H^*(BC_2;\Ffield_2) \cong \Ffield_2[x,y,w,e]/(xy). 
    \end{equation}
    The differentials are generated under the Leibniz rule by
    \begin{align}
        d_2(e) & = x^2 + y^2 + w,  \\
        d_3(e^2) & = \Sq^1(x^2 + y^2 + w) = (x+y)w.
    \end{align}
    The spectral sequence collapses at $E_4$ with a vanishing line of height 4, and the resulting cohomology ring is
    \begin{equation}
        H^*(BQ_{2^n};\Ffield_2) \cong \Ffield_2[x,y,e_4]/(xy,x^3+y^3).
    \end{equation}
    with $|x| = |y| = 1$, and $|e_4| = 4$.
    \label{lem:Q2ng3LHS}
\end{Proposition}
\subsection{The $\sE_{(2)}$-exponent}
Combining \cref{q2n:expupper} with the vanishing lines of height 4 in \cref{q8:lhsss} and \cref{lem:Q2ng3LHS} gives the following.
\begin{Proposition}
    \label{prop:q2ne2exp}
    The $\sE_{(2)}$-exponent satisfies
    \begin{equation}
        \exp_{\sE_{(2)}}\underline{H\Ffield_2}_{Q_{2^n}} = 4.
    \end{equation}
\end{Proposition}

\section{$SD_{16}$}
\label{sec:sd16}
Let
\begin{equation}
    SD_{16} = \langle s,r \, | \, s^2 = r^8 = e, \,  srs^{-1} = r^3 \rangle
\end{equation}
be the semidihedral group of order 16. This group has no unique maximal elementary abelian 2-subgroup, but the subgroup $\langle s, r^2 \rangle \cong D_8$ does contain all the elementary abelian 2-subgroups. We will compute the $\All_{D_8}$-homotopy limit spectral sequence converging to $H^*(BSD_{16};\Ffield_2)$, which is isomorphic to (\cite[App.\ C, \#13(16)]{carlson2003cohomology})
\begin{equation}
    \Ffield_2[z, y, x, w]/(zy, y^3, yx, z^2w + x^2),
\end{equation}
with degrees $|z| = |y| = 1$, $|x| = 3$ and $|w| = 4$.

By \cref{prop:compserre}, the $\All_{D_8}$-homotopy limit spectral sequence is isomorphic to the LHSSS corresponding to the group extension
\begin{equation}
    \overset{s, r^2}{D_8} \to SD_{16} \to \overset{\rbar}{C_2}
\end{equation}
\subsection{The $E_2$-page}
Since $rsr^{-1} = sr^2$, $C_2$ acts on $D_8$ by interchanging the generators $s$ and $sr^2$ of $D_8$. Write $H^*(BD_8;\Ffield_2) = \Ffield_2[x,y,w]/(xy)$, with $x=\delta_{s}$, $y = \delta_{sr^2}$, and where $w$ is the unique class in degree 2 that restricts to 0 on $\langle s \rangle$ and $\langle sr^2 \rangle$ and classifies an extension of $D_8$ isomorphic to $D_{16}$ (see, e.g., \cite[Thm.\ 2.7]{adem}). Then $C_2 = \langle \rbar \rangle$ acts by interchanging $x$ and $y$. Of course it sends an extension isomomorphic to $D_{16}$ to another such extension, and since $s$ and $sr^2$ get interchanged, we see that the $C_2$-action leaves $w$ fixed, by the restriction property mentioned before.

To determine $E_2^{s,t} = H^s(BC_2;H^t(BD_8;\Ffield_2))$, let
\begin{equation}
    \cdots \xrightarrow{1-\rbar} \Integers C_2 \xrightarrow{1+\rbar} \Integers C_2 \xrightarrow{1-\rbar} \Integers C_2 \xrightarrow{\epsilon} \Integers
\end{equation}
be the minimal free $C_2$-resolution of $\Integers$. Applying $\Hom_{C_2}(-,H^*(BD_8;\Ffield_2))$, we see that the differentials are on monomials given by
\begin{align}
    (x^i + y^i)w^j & \mapsfrom x^i w^j, \\
    (x^k + y^k) w^l & \mapsfrom y^k w^l.
\end{align}
Hence the kernel has as an $\Ffield_2$-basis the elements of the form $w^j P(x,y)$, where $P(x,y)$ is a symmetric polynomial in $x$ and $y$. Taking the cokernel divides out the symmetrized polynomials in $x$ and $y$. Hence the result has as an $\Ffield_2$-basis the polynomials of the form $x^i y^i w^j$, which is 0 unless $i=0$. Writing $a = \delta_{\rbar}$, and $\sigma_1 = x+y$, we get
\begin{equation}
    E_2^{s,t} = \Ffield_2[a,\sigma_1,w]/(a\sigma_1),
\end{equation}
with $(s,t)$-degrees given by $|a| = (1,0)$, $|\sigma_1| = (0,1)$ and $|w| = (0,2)$.
\subsection{Differentials and $E_4 = E_\infty$}
Since $SD_{16}^{\ab} = SD_{16}/\langle r^4 \rangle \cong C_2 \times C_4$, we have $\dim_{\Ffield_2} H^1(BSD_{16};\Ffield_2) = 2$. Hence $x+y$ is a permanent cycle. Since there is no vanishing line on $E_2$, $w$ will have to support a differential, which can only be $d_3 w = a^3$. Hence the $E_3$-page is
\begin{figure}[H]
    \centering
    \includegraphics[width=0.9\textwidth]{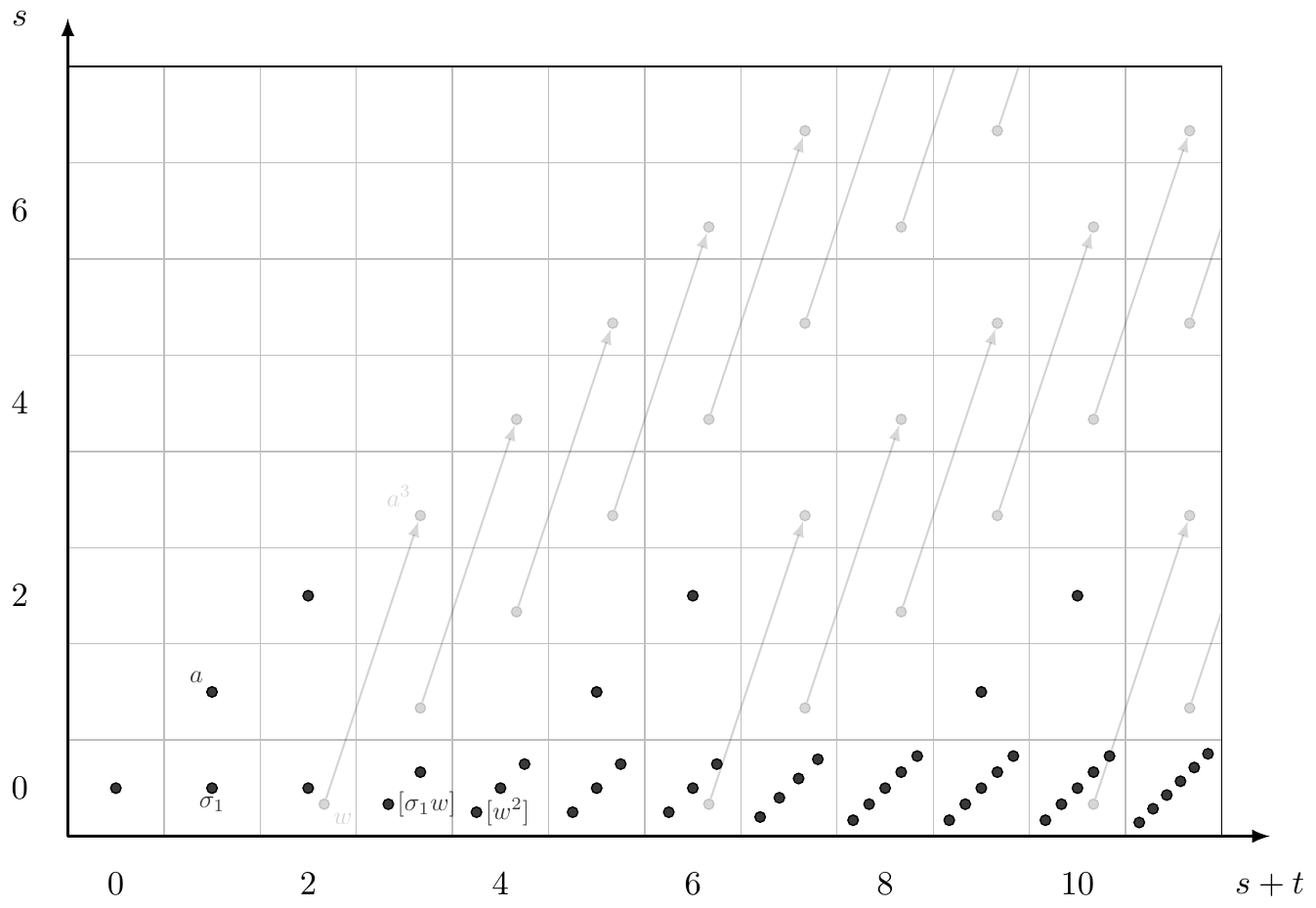}
    \caption{The $E_3$-page.}
\end{figure}
Therefore the $E_4$-page is
\begin{equation}
    \Ffield_2[a,\sigma_1,[\sigma_1w],[w^2]]/(a\sigma_1,a[\sigma_1 w], \sigma_1^2[w^2] - [\sigma_1 w]^2, a^3).
\end{equation}
Hence the $E_4$-page has a vanishing line of height 3, and therefore $E_4 = E_\infty$. There are several multiplicative extension problems, which we leave unresolved.
\fxwarning{Mention if I leave multiplicative extension problems unsolved or solve them.}
\subsection{Poincar\'e Series}
The $E_\infty$-page has as $s=0$-line
\begin{align}
    & \Ffield_2[\sigma_1, [\sigma_1 w], [w^2]]/(\sigma_1^2[w^2] - [\sigma_1 w]^2) \\
    & \cong \Ffield_2[\sigma_1, [w^2]]\{1, [\sigma_1 w]\},
\end{align}
which has Poincar\'e series
\begin{equation}
    \frac{1+t^3}{(1-t)(1-t^4)}.
\end{equation}
The $s=1$-line of the $E_\infty$-page is given by
\begin{equation}
    \Ffield_2[[w^2]]\{a\},
\end{equation}
which has Poincar\'e series
\begin{equation}
    \frac{t}{1-t^4}.
\end{equation}
The $s=2$-line of the $E_\infty$-page is given by
\begin{equation}
    \Ffield_2[[w^2]]\{a^2\},
\end{equation}
which has Poincar\'e series
\begin{equation}
    \frac{t^2}{1-t^4}.
\end{equation}
This shows that the Poincar\'e series of $H^*(BSD_{16}; \Ffield_2)$ is given by
\begin{align}
    \frac{1+t^3}{(1-t)(1-t^4)} + \frac{t + t^2}{1-t^4} & = \frac{(1+t^3) + t(1-t) + t^2(1 - t)}{(1-t)(1-t^4)} \\
    & = \frac{1+t}{(1-t)(1-t^4)} \\
    & = \frac{1}{(1-t)^2(1+t^2)},
\end{align}
cf.\ \cite[App.\ B and App.\ G]{carlson2003cohomology}.
\subsection{An upper bound on the $\sE_{(2)}$-exponent}
In this subsection we will prove the following upper bound on the $\sE_{(2)}$-exponent.
\begin{Proposition}
    The $\sE_{(2)}$-exponent satisfies
    \begin{equation}
        \exp_{\sE_{(2)}} \underline{H\Ffield_2}_{SD_{16}} \leq 4.
    \end{equation}
\end{Proposition}
\begin{Proof}
    We will use the notation $R_\theta$ and $T$ from the proof of \cref{m16:expupper}, and also use the fact proven there that $R_{\theta}$ fixes a line only if $\theta \in \Integers\{\pi\}$. Write $\alpha = 2 \pi /8$. Let $V$ be the real $SD_{16}$-representation given by
    \begin{align}
        s & \mapsto T, \\
        r & \mapsto R_\alpha \oplus R_\alpha
    \end{align}
    from \cite[Lem.\ 13.4]{totaro14}. The same argument as in the proof of \cref{m16:expupper} shows that the only powers of $r$ fixing a line are $e$ and $r^4$, and that the only elements of the form $sr^k$ that can possibly fix a line are $s$, $sr^2$, $sr^4$, $sr^6$.
    The elements
    \begin{equation}
    e, r^4, s, sr^2, sr^4, sr^6.
    \end{equation}
     do not generate an elementary abelian 2-group as in the proof of \cref{m16:expupper}, but all the isotropy groups of lines that occur consist of proper subsets of the above elements that do form elementary abelian 2-groups, as we will now show.

First, $r^4$ fixes every line, so is in every isotropy group.

Write $[x_1:y_1:x_2:y_2]$ for the real homogenous coordinates in $\Reals P^3$. Then
\begin{equation}
    s \cdot [x_1 : y_1 : x_2 : y_2 ] = [x_2 : y_2 : x_1 : y_1],
\end{equation}
which is fixed if and only if $(x_1,y_1) = (x_2,y_2)$ or $(x_1,y_1) = - (x_2,y_2)$. Computation of $sr^4 \cdot [x_1 : y_1 : x_2: y_2]$ shows that this element fixes exactly the same lines.

Computing
\begin{equation}
    sr^2 \cdot [x_1 : y_1 : x_2 : y_2] = [-y_2: x_2: y_1:-x_1]
\end{equation}
shows that $sr^2$ fixes precisely those lines with $(x_1,y_1) = (-y_2,x_2)$ or $(x_1,y_1) = (y_2,-x_2)$. Computing $sr^6 \cdot [x_1:y_1:x_2:y_2]$ shows that this element fixes precisely the same lines.

The sets of lines $\{x_1 = x_2, y_1 = y_2\} \cup \{x_1 = -x_2, y_1 = -y_2\}$ and $\{x_1 = -y_2, y_1 = x_2\} \cup \{x_1 = y_2, y_1 = -x_2\}$ have empty intersection, because a line in the intersection $\{x_1 = x_2, y_1 = y_2\}  \cap \{x_1 = -y_2, y_1 = x_2\}$ would satisfy
\begin{equation}
    y_2 = y_1 = x_2 = x_1 = -y_2,
\end{equation}
a line in the intersection $\{x_1 = x_2, y_1 = y_2\} \cap \{x_1 = y_2, y_1 = -x_2\}$ would satisfy
\begin{equation}
    y_1 = y_2 = x_1 = x_2 = -y_1,
\end{equation}
a line in the intersection $\{x_1 = -x_2, y_1 = -y_2\} \cap \{x_1 = -y_2, y_1 = x_2\}$ would satisfy
\begin{equation}
    y_2 = -y_1 = -x_2 = x_1 = -y_2,
\end{equation}
and finally a line in the interesection $ \{x_1 = -x_2, y_1 = -y_2\} \cap \{x_1 = y_2, y_1 = -x_2\}$ would satisfy
\begin{equation}
    y_2 = x_1 = -x_2 = y_1 = -y_2.
\end{equation}

Hence the isotropy groups that occur are
\begin{equation}
    \langle r^4 \rangle, \langle s, sr^4 \rangle, \langle sr^2, sr^6 \rangle,
\end{equation}
all of which are elementary abelian 2-groups. Applying \cref{projbund:prop1} gives the desired result.
\end{Proof}

\subsection{The exponent}
The 1-dimensional real representation
\begin{equation}
    SD_{16} \to SD_{16}/D_8 \cong C_2
\end{equation}
has Euler class $a$. This has isotropy in $\All_{D_8}$, and the vanishing line shows that $a$ is of nilpotence degree 3, hence
\begin{equation}
    \exp_{\All_{D_8}} \underline{H\Ffield_2} = 3.
\end{equation}
Since increasing the family only can make the exponent go down, we have as a consequence
\begin{Proposition}
    The exponent $\sE_{(2)}$-exponent satisfies 
    \begin{equation}
    3 \leq \exp_{\sE_{(2)}} \underline{H\Ffield_2}_{SD_{16}} \leq 4.
\end{equation}
\label{prop:sd16e2exp}
\end{Proposition}

\section{$M_{16}$}
\label{sec:m16}
\subsection{Introduction}
Let
\begin{equation}
    M_{16}= \langle r,\, f \, | \, r^8 = f^2 = e, \, frf^{-1}=r^5 \rangle
\end{equation}
be the modular group of order 16. The group has this name because its lattice of subgroups is modular (see, e.g., \cite[I.\S7]{birkhoff67}).

In this section we will compute the $\sE_{(2)}$-homotopy limit spectral sequence converging to $H^*(BM_{16}; \Ffield_2)$, which is isomorphic to (\cite[App.\ C, \#11(16)]{carlson2003cohomology})
\begin{equation}
    \Ffield_2[z, y, x, w]/(z^2, zy^2, zx, x^2)
\end{equation}
with degrees $|z| = |y| = 1$, $|x| = 3$, $|w| = 4$.

\subsection{An upper bound on the exponent}
In this subsecton we prove the following upper bound on the $\sE_{(2)}$-exponent.
\begin{Proposition}
    The $\sE_{(2)}$-exponent satisfies
    \begin{equation}
        \exp_{\sE_{(2)}} \underline{H\Ffield_2}_{M_{16}} \leq 4.
    \end{equation}
    \label{m16:expupper}
\end{Proposition}
\begin{Proof}
    For an angle $\theta$, let $R_{\theta}$ be the rotation-by-$\theta$ matrix $\begin{pmatrix} \cos \theta & - \sin\theta \\ \sin \theta & \cos \theta \end{pmatrix}$. Let $T$ be the matrix
    \begin{equation}
        \begin{pmatrix}
            0 & 0 & 1 & 0 \\
            0 & 0 & 0 & 1 \\
            1 & 0 & 0 & 0 \\
            0 & 1 & 0 & 0
        \end{pmatrix}
    \end{equation}
    which interchanges the summands of $\Reals^2 \oplus \Reals^2$. Write $\alpha = 2\pi/8$. Let $V$ be the 4-dimensional real $M_{16}$-representation given by
    \begin{align}
        f & \mapsto T, \\
        r & \mapsto R_\alpha \oplus R_{5\alpha}.
        \label{}
    \end{align}
    This is the representation from \cite[Lem.\ 13.3]{totaro14}.

    We will now determine the isotropy of the projectivation $\Proj(V)$.

    A power of $r$ fixes a line in $\Reals^2 \oplus \Reals^2$ if and only if it fixes a line in one of the summands. The matrix $R_{\theta}$ has characteristic polynomial $P(\lambda) = \lambda^2 - 2 \lambda \cos\theta + 1$, which has discriminant $\Delta = -4 \sin^2 \theta$. Hence $R_\theta$ has real eigenvalues only if $\theta \in \Integers\{\pi \}$. Since $r \mapsto R_{2\pi/8}$, $r^k$ fixes a line only if $\alpha k \in \Integers\{\pi\}$ or $5 \alpha k \in \Integers\{\pi\}$, that is only if $k \in \Integers\{4\}$. Therefore $e$ and $r^4$ are the only elements powers of $r$ that can possibly fix a line. The element $r^4$ acts by $-\Id$, and therefore fixes all lines. The element $e$ of course also fixes all lines.

    We now consider the remaining elements of $M_{16}$, which are of the form $fr^k$, and which act by $T(R_{\alpha k}\oplus R_{5\alpha k})$. This matrix certainly only fixes a line if its square does, which is $R_{6\alpha k} \oplus R_{6\alpha k}$. These summands only fix a line if $6 \alpha k \in \Integers\{\pi\}$, which implies that $k$ is even. Hence the only possible elements of the form $fr^k$ fixing a line are $f$, $fr^2$, $fr^4$, and $fr^6$. One verifies directly that $fr^2$ acts by the matrix
    \begin{equation}
        \begin{pmatrix}
            0 & 0 & 0 & -1 \\
            0 & 0 & 1 & 0 \\
            0 & -1 & 0 & 0 \\
            1 & 0 & 0 & 0
        \end{pmatrix}
    \end{equation}
    which fixes no lines. Consequently its inverse $fr^6$ neither fixes lines. One also easily checks that $fr^4$ acts by
    \begin{equation}
        \begin{pmatrix}
            0 & 0 & -1 & 0 \\
            0 & 0 & 0 & -1 \\
            -1 & 0 & 0 & 0 \\
            0 & -1 & 0 & 0
        \end{pmatrix}
    \end{equation}
    which fixes the lines $[x_1 : y_1 : x_2 : y_2]$ in $\{x_1 = -x_2, y_1 = -y_2\}$. Hence the only elements of $M_{16}$ fixing a line are $e, f, fr^4, r^4$, which together form a Klein four group. Applying \cref{projbund:prop1} gives the desired result.
\end{Proof}

\subsection{The orbit category}
The group $M_{16}$ has a unique maximal elementary abelian 2-subgroup, given by
\begin{equation}
    C_2 \times C_2 \cong \langle f, \, fr^4 \rangle.
\end{equation}
By \cref{prop:compserre}, this identifies the $\sE_{(2)}$-homotopy limit spectral sequence with the LHSSS associated to the group extension
\begin{equation}
    C_2 \times C_2 \to M_{16} \to \langle \rbar \rangle \cong C_4.
\end{equation}
\subsection{The $E_2$-page}
The action of $C_4 = \langle \rbar \rangle$ on $C_2 \times C_2 = \langle f, \, fr^4 \rangle$ is given by $rfr^{-1} = fr^4$, i.e. $\rbar$ interchanges the generators $f$ and $fr^4$. Denote the duals by $x \coloneqq \delta_f,\, y \coloneqq \delta_{fr^4} \in H^1(BC_2 \times C_2;\, \Ffield_2)$. Then $H^*(BC_2 \times C_2;\, \Ffield_2) = \Ffield_2[x,y]$, and $\rbar$ acts by interchanging $x$ and $y$.

A 2-periodic projective $C_4$-resolution of $\Integers$ is given by
\begin{equation}
    \cdots \to \Integers C_4 \xrightarrow{1-\rbar} \Integers C_4 \xrightarrow{1+\rbar+\rbar^2+\rbar^3} \Integers C_4 \xrightarrow{1-\rbar} \Integers C_4 \xrightarrow{\epsilon} \Integers.
\end{equation}
Applying $\Hom_{\Integers C_4}(-,\, \Ffield_2[x,y])$ (where $\Ffield_2[x,y]$ carries the $C_4$-action described above) yields
\begin{equation}
    \cdots \xleftarrow{0} \Ffield_2[x,y] \leftarrow \Ffield_2[x,y] \xleftarrow{0} \Ffield_2[x,y] \leftarrow \Ffield_2[x,y],
\end{equation}
where the maps going from odd to even degrees are 0, and the maps from even to odd degrees are on monomials given by
\begin{equation}
    x^i y^j + x^j y^i \mapsfrom x^i y^j.
\end{equation}
This shows that the even rows of the $E_2$-page are given by the $C_4$-invariants of $\Ffield_2[x,y]$, whereas the odd rows are given by the coinvariants.

We denote $c_0 \coloneqq \delta_{\rbar} \in H^1(BC_4;\, \Ffield_2)$, and let $\beta_2(c_0) \in H^2(BC_4;\, \Ffield_2)$ be the second order Bockstein of $c_0$. Then $H^*(BC_4;\, \Ffield_2) = \Ffield_2[c_0, \beta_2(c_0)]/(c_0^2)$. In addition we write $\sigma_1 = x+y$ and $\sigma_2 = xy$ for the first and second elementary symmetric polynomials in $x$ and $y$.

The above discussion shows that $E^{0,*}_2 = \Ffield_2[\sigma_1,\sigma_2]$, and more generally that the the even rows of the $E_2$-page equals
\begin{equation}
    E_2^{0\!\!\!\!\pmod{2}, *} \cong \Ffield_2[\sigma_1, \sigma_2, \beta_2(c_0)].
    \label{M16:evenrows}
\end{equation}

The $s \equiv 1 \pmod{2}$-rows of $E_2$ are isomorphic to the coinvariants, and are as such modules over the $s=0$-row of $E_2$, which is the invariants. In the following lemma we describe this module structure, because we will use it to organize the calculation.
\begin{Lemma} 
    The $\Ffield_2[x,y]^{C_2}$-module $\Ffield_2[x,y]_{C_2}$ is (non-freely) generated by $\{[1], [x]\}$. (Recall that we denote by $[P]$ the coinvariance class of $P \in \Ffield_2[x,y]$.)
    \label{lem:coinvinvfg}
\end{Lemma}
\begin{Proof}
    Writing $\sigma_1 = x+y$, $\sigma_2 = xy$ for the first and second elementary symmetric polynomials, we prove that the degree $n$-part of $\Ffield_2[x,y]_{C_2}$ is generated by $\{[1], [x]\}$ by induction on $n$. Since the degree 0 and 1 part of $\Ffield_2[x,y]_{C_2}$ are given by $\Ffield_2\{[1]\}$ and $\Ffield\{[x]\}$ respectively, the lemma holds for $n=0,1$.

    Assume now $n \geq 2$. The degree $n$-part of $\Ffield_2[x,y]_{C_2}$ is given by
    \begin{equation}
        \Ffield_2\{[x^n], [x^{n-1}y], \ldots, [x^{\ceil{n/2}}y^{\floor{n/2}}]\}.
    \end{equation}
    All elements $[x^i y^j] \neq [x^n]$ in this basis have strictly positive $x$- and $y$-exponent, and therefore
    \begin{equation}
        [x^i y^j] = xy[x^{i-1} y^{i-1}],
    \end{equation}
    and $xy$ is an invariant and $[x^{i-1}y^{i-1}]$ is a coinvariant of degree strictly smaller then $n$, hence by induction in the $\Ffield_2[x,y]^{C_2}$-module generated by $\{[1], [x]\}$.

    The element $[x^n]$, finally, equals
    \begin{equation}
        [x^n] = (x+y)[x^{i-1}] + xy[x^{i-1}],
    \end{equation}
    where the coinvariance classes on the right hand side are also of degree strictly smaller than $n$. 
    
    Notice that we used the induction hypothesis for $n-1$ and $n-2$, but not smaller, so that we do not use more than the base case provided.
\end{Proof}
Let $\Ffield_2[\sigma_1, \sigma_2]$ be the free graded ring with $|\sigma_1| = 1$, $|\sigma_2| = 2$, and let 
\begin{equation}
\Ffield_2[\sigma_1, \sigma_2]\{c_0, c_1\}
\end{equation}
be the free graded module with $|c_0| = 0$ and $|c_1| = 1$. By \cref{lem:coinvinvfg}, the morphism of $\Ffield_2[\sigma_1, \sigma_2]$-modules given by
\begin{align}
    \Ffield_2[\sigma_1, \sigma_2]\{c_0, c_1\} & \to \Ffield_2[x,y]_{C_2} \\
    c_0 & \mapsto [1] \\
    c_1 & \mapsto [x]
    \label{M16:eq3}
\end{align}
is surjective.
\begin{Lemma}
    The kernel of the map (\ref{M16:eq3}) is $\Ffield_2[\sigma_1, \sigma_2]\{\sigma_1 c_0\}$.
    \label{lem:keroddrows}
\end{Lemma}
\begin{Proof}
    This follows from a Poincar\'e series calculation. The Poincar\'e series of the right hand side of (\ref{M16:eq3}) is 
    \begin{align}
        \frac{1+t}{(1-t^2)^2} & = \frac{1}{(1-t)(1-t^2)}.
    \end{align}
    The Poincar\'e series of the left hand side of (\ref{M16:eq3}) is 
    \begin{equation}
        \frac{1+t}{(1-t)(1-t^2)}.
    \end{equation}
    
    Certainly, $\Ffield_2[\sigma_1, \sigma_2]\{\sigma_1 c_0\}$ is contained in the kernel, for $\sigma_1 c_0 \mapsto (x+y)[1] = 0$. The module $\Ffield_2[\sigma_1, \sigma_2]\{\sigma_1 c_0\}$ has Poincar\'e series 
    \begin{equation}
        \frac{t}{(1-t)(1-t^2)}.
    \end{equation}
    Because 
    \begin{align}
        \frac{1+t}{(1-t)(1-t^2)} - \frac{t}{(1-t)(1-t^2)} = \frac{1}{(1-t)(1-t^2)},
    \end{align}
    we see that in fact this is all of the kernel.
\end{Proof}

Returning to our description of the $E_2$-page, \cref{lem:keroddrows} shows that the odd rows are given by
\begin{align}
    E_2^{1\!\!\!\!\pmod{2},*} & \cong\frac{\Ffield_2[\sigma_1, \sigma_2]\{c_0, c_1\}}{\Ffield_2[\sigma_1, \sigma_2]\{\sigma_1 c_0\}} \otimes \Ffield_2[\beta_2(c_0)]  \\
    & \cong \frac{E_2^{0,*}\{c_0, c_1\}}{E_2^{0,*}\{\sigma_1 c_0\}} \otimes \Ffield_2[\beta_2(c_0)].
    \label{M16:oddrows}
\end{align}
Combining (\ref{M16:oddrows}) with (\ref{M16:evenrows}) gives the $E_2$-page:
\begin{equation}
    E_2^{*,*} \cong \Ffield_2[\sigma_1, \sigma_2, c_0, \beta_2(c_0), c_1]/(c_0^2, c_1^2, c_0c_1, \sigma_1 c_0)
\end{equation}
with $(s,t)$-degrees given by $|\sigma_1| = (0,1)$, $|\sigma_2| = (0,2)$, $|c_0|= (1,0)$, $|\beta_2(c_0)| = (2,0)$, $|c_1| = (1,1)$. The $E_2$-page is depicted (without differentials) in the following figure:
\begin{figure}[H]
    \centering
    \includegraphics[width=0.9\textwidth]{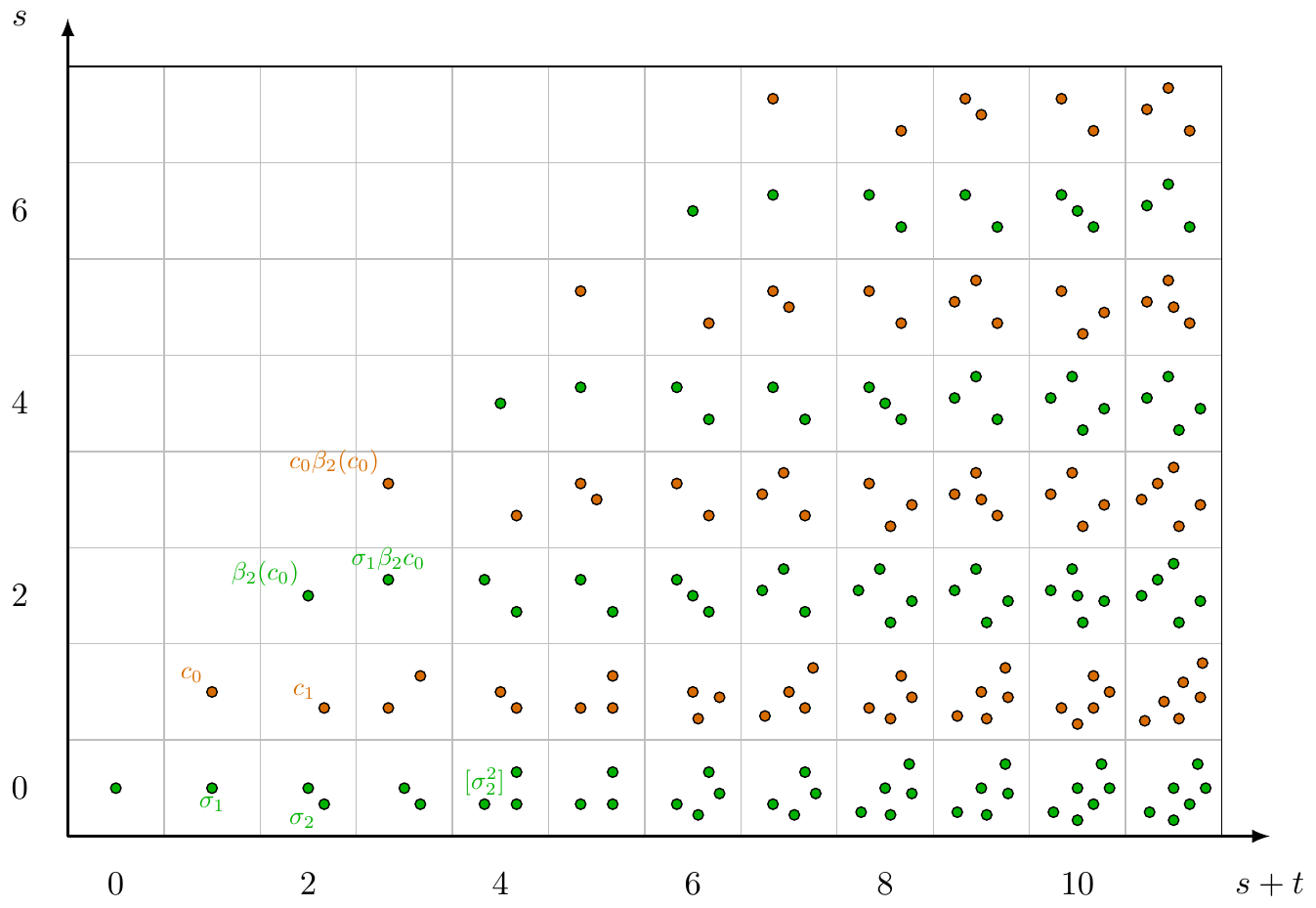}
    \caption{The $E_2$-page, without differentials. The even rows, which are the invariants, are colored green, the odd rows, which are the coinvariants, are colored orange.}
\end{figure}
\subsection{Differentials on the $E_2$-page}
\subsubsection{Generating differentials}
Denote the dimension of $H^i(BM_{16};\Ffield_2)$, i.e.\ the $i$-th Betti number, by $h^i$. We assume known that $h^1 = 2$ and $h^2 = 2$ in determining the differentials. The first follows from $M_{16}^{\ab} \cong C_2 \times C_2$, the second can be verified either by hand, or using a computer algebra system such as GAP.

The classes $c_0$ and $\beta_2(c_0)$ are permanent cycles for degree reasons, and $\sigma_1$ is a permanent cycle because $h^1 = 1$ and there are 2 classes in the 1-stem. There are 4 classes in the 2-stem, but $h^2 = 2$, so two will have to support a differential, and by the previous discussion, these will have to be $\sigma_2$ and $c_1$. The only possibility for $c_1$ to support a differential is
\begin{equation}
    d_2(c_1) = c_0 \beta_2(c_0).
    \label{M16:d2-1}
\end{equation}
Given this, the only remaining possibility for $\sigma_2$ to support a differential is
\begin{equation}
    d_2(\sigma_2) = \sigma_1 \beta_2(c_0).
    \label{M16:d2-2}
\end{equation}
We will study the propagation of these differentials using the Leibniz rule for the even and odd rows seperately.
\subsubsection{Differentials in the even rows}
The even rows of the $E_2$-page are given by
\begin{align}
    E_2^{0\!\!\!\!\pmod{2},*} & \cong \Ffield_2[\sigma_1, \sigma_2, \beta_2(c_0)] \\
     & \cong \Ffield_2[\sigma_1, [\sigma_2^2], \beta_2(c_0)] \oplus \Ffield_2[\sigma_1, [\sigma_2^2], \beta_2(c_0)]\{\sigma_2\}.
\end{align}
Propagating the differential (\ref{M16:d2-2}) gives the following pattern of differentials:
\begin{figure}[H]
    \centering
    \includegraphics[width=0.9\textwidth]{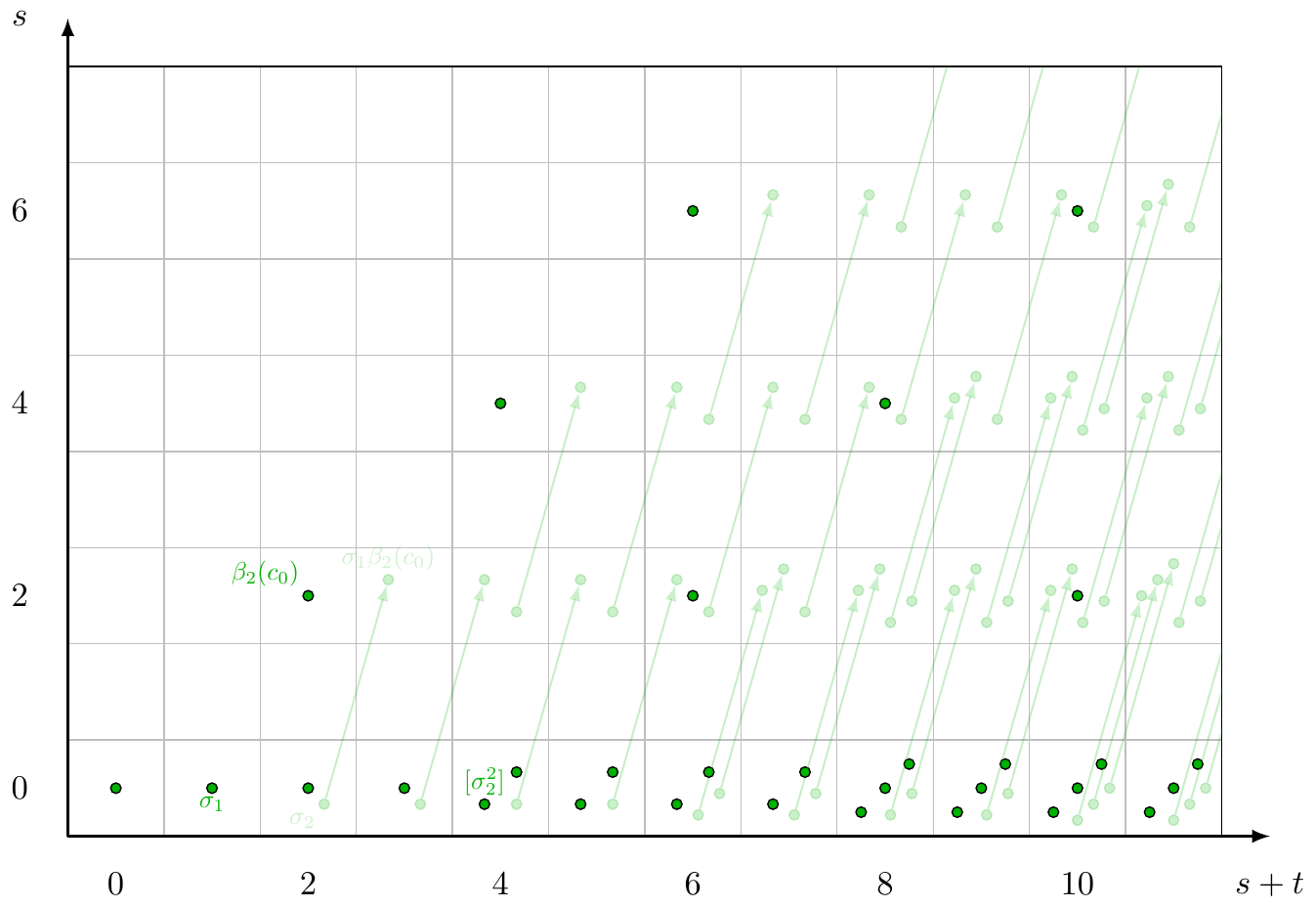}
    \caption{The differentials in the even rows on the $E_2$-page.}
\end{figure}
This also shows that there is no possibility for more $d_2$'s in the even rows on the $E_2$-page, and that the even rows of the $E_3$-page are given by
\begin{equation}
    E_3^{0\!\!\!\!\pmod{2},*} \cong \Ffield_2[\sigma_1, [\sigma_2^2], \beta_2(c_0)]/(\sigma_1 \beta_2(c_0)).
\end{equation}
\subsubsection{Differentials in the odd rows}
In the odd rows, propagating the differentials (\ref{M16:d2-1}) and (\ref{M16:d2-2}) using the Leibniz rule, we get
\begin{align}
    d_2(\sigma_2 c_1) & = \sigma_1 c_1 \beta_2(c_0) + \sigma_2 c_0 \beta_2(c_0), \\
    d_2(c_1) & = c_0 \beta_2(c_0), \\
    d_2(\sigma_1 \sigma_2 c_1) & = \sigma_1^2 c_1 \beta_2(c_0).
\end{align}
This implies that the classes in the odd rows supporting a differential generate the submodule 
\begin{align}
    & \left( \Ffield_2[\sigma_2^2]\{\sigma_2c_1\} \oplus \Ffield_2[\sigma_2^2]\{c_1\}\right. \\
    & \left. \oplus \Ffield_2[\sigma_1, \sigma_2^2]\{\sigma_1 \sigma_2 c_1\}\right) \otimes \Ffield_2[\beta_2(c_0)],
\end{align}
and that the classes that get hit are given by
\begin{align}
    & \left( \Ffield_2[\sigma_2^2]\{\sigma_1 c_1 + \sigma_2 c_0\} \oplus \Ffield_2[\sigma_2^2]\{c_0\} \right. \\
    & \left. \oplus \Ffield_2[\sigma_1, \sigma_2^2]\{\sigma^2_1 c_1\} \right) \otimes \Ffield_2[\beta_2(c_0)]\{\beta_2(c_0)\}.
\end{align}
This pattern of differentials is depicted in the following figure.
\begin{figure}[H]
    \centering
    \includegraphics[width=0.9\textwidth]{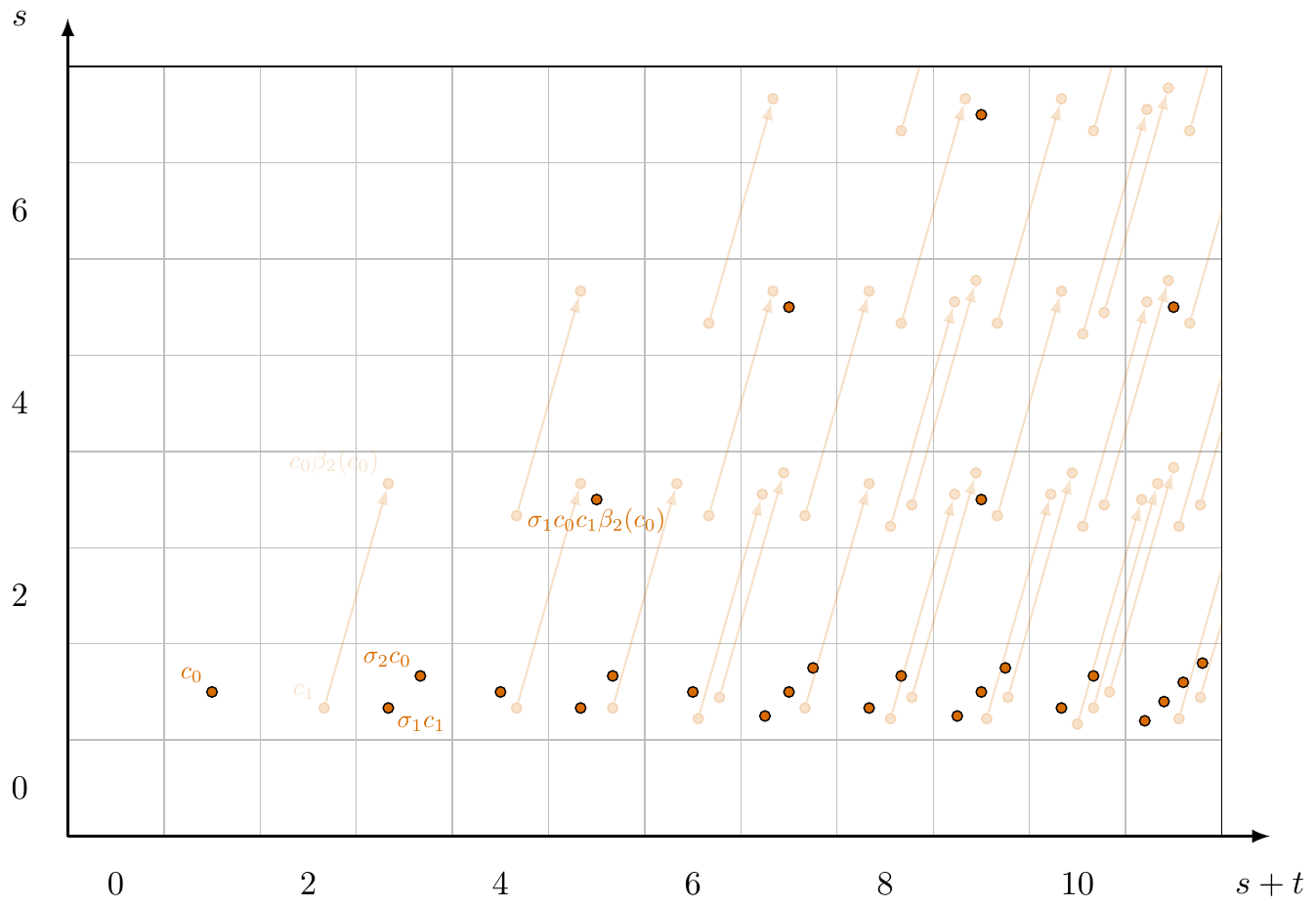}
    \caption{The differentials in the odd rows of $E_2$.}
\end{figure}
\fxwarning{Make the picture more faithful by showing that some classes hit a sum of classes.}
This implies that the classes in the $s=1$-line surviving to $E_3$ are
\begin{equation}
    E_3^{1,*} \cong \Ffield_2[\sigma_1, [\sigma_2^2]]\{\sigma_1 c_1\} \oplus \Ffield_2[[\sigma^2_2]]\{c_0, [\sigma_2 c_0]\}
\end{equation}
and that the odd $s$-rows with $s \geq 3$ surviving to $E_3$ are
\begin{align}
    E_3^{s,*} & \cong \frac{\Ffield_2[\sigma_1, [\sigma_2^2]]\{[\sigma_1 c_1]\}}{\Ffield_2[\sigma_1, [\sigma_2^2]]\{\sigma_1[\sigma_1 c_1]\}} \\
    & \cong \Ffield_2[[\sigma_2^2]]\{[\sigma_1 c_1]\}.
\end{align}
Putting the differentials in the even and odd rows together gives the following $E_2$-page:
\begin{figure}[H]
    \centering
    \includegraphics[width=0.9\textwidth]{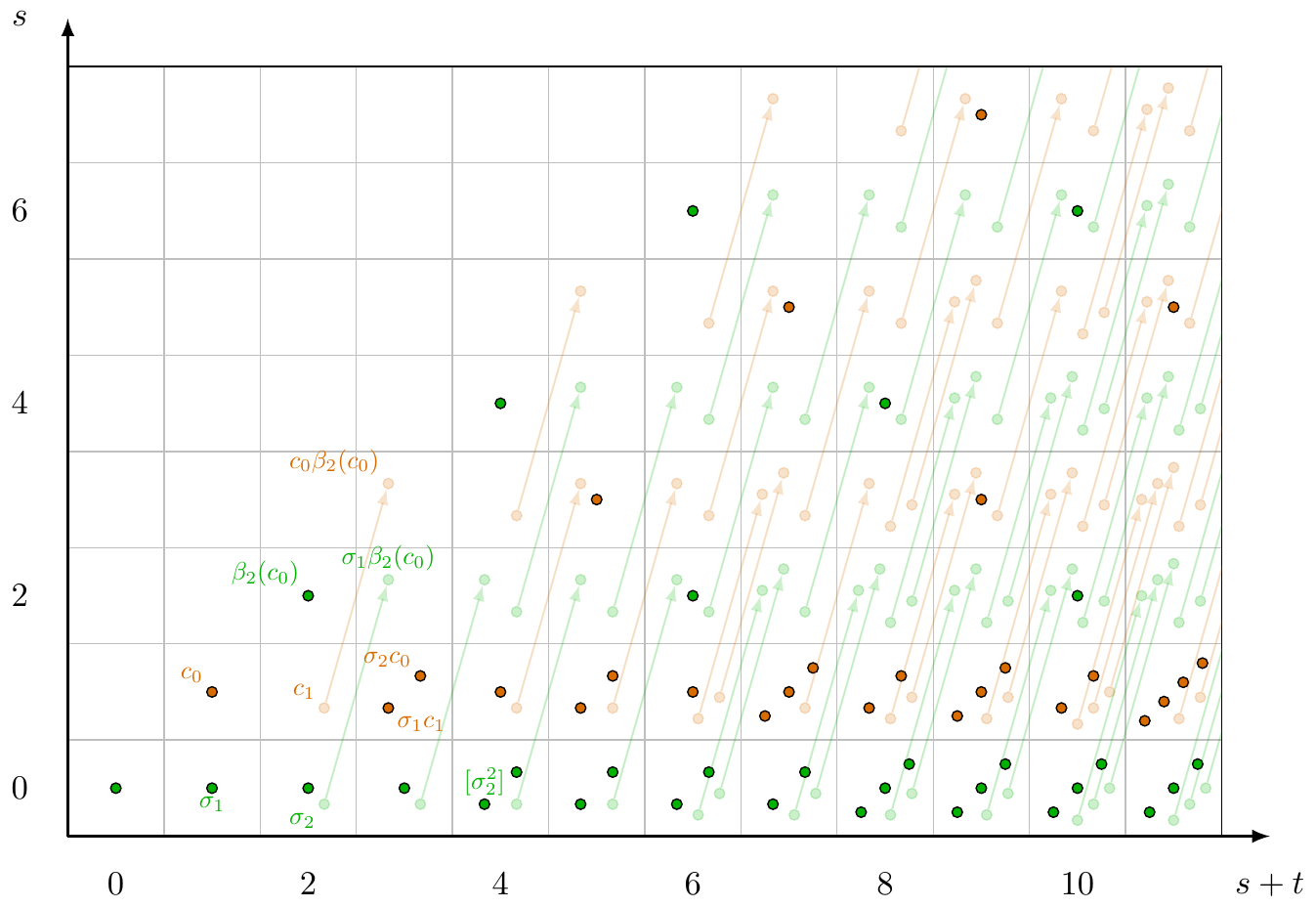}
    \caption{The $E_2$-page with differentials.}
\end{figure}
\subsection{The $E_3$-page}
For a horizontal vanishing line of height 4 on $E_5$ there must be a non-zero
\begin{equation}
    d_3 \colon E_3^{1,2} \to E_3^{4,0}.
\end{equation}
We have $E_3^{1,2} = \Ffield_2\{\sigma_1 c_1, \sigma_2 c_0\}$, hence at least one of these two classes has to support a differential, the only possible target for both is to hit $\beta_2(c_0)^2$ with a $d_3$. Because $\beta_2(c_0)^2$ is the only class in $E_3^{4,0}$, we have that $d_3([\sigma_1 c_1]) \neq 0$ if and only if $d_3([\sigma_1 c_1]) = \beta_2(c_0)^2$. Similarly, we have that $d_3([\sigma_1 c_1]\beta_2(c_0)) \neq 0$ if and only if $d_3([\sigma_1 c_1]\beta_2(c_0)) = \beta_2(c_0)^3$. By the Leibniz rule, $d_3([\sigma_1 c_1]) = \beta_2(c_0)^2$ implies that $d_3([\sigma_1 c_1]\beta_2(c_0)) = \beta_2(c_0)^3$. Furthermore, if $d_3([\sigma_1 c_1]) = 0$ then, also by the Leibniz rule, we have $d_3([\sigma_1 c_1]\beta_2(c_0)) = 0$. Therefore $d_3([\sigma_1 c_1]) = \beta_2(c_0)^2 $ if and only if $d_3([\sigma_1 c_1]\beta_2(c_0)) = \beta_2(c_0)^3$.

Similarly, $d_3(\sigma_2 c_0) = \beta_2(c_0)^2$ if and only if $d_3(\sigma_2 c_0 \beta_2(c_0)) = \beta_2(c_0)^3$.
But since
\begin{equation}
    [\sigma_1 c_1 \beta_2(c_0)] = [\sigma_2 c_0 \beta_2(c_0)],
\end{equation}
as a consequence of
\begin{equation}
    d_2(\sigma_2 c_1) = (\sigma_1 c_1 +  \sigma_2 c_0)\beta_2(c_0)
\end{equation}
we have:
\begin{equation}
    d_3(\sigma_1 c_1) = d_3(\sigma_2 c_0) = \beta_2(c_0)^2.
\end{equation}
\fxwarning{Make the picture more faithful by showing that the differential is really supported on the sum of two classes.}
This is depicted in
\begin{figure}[H]
    \centering
    \includegraphics[width=0.9\textwidth]{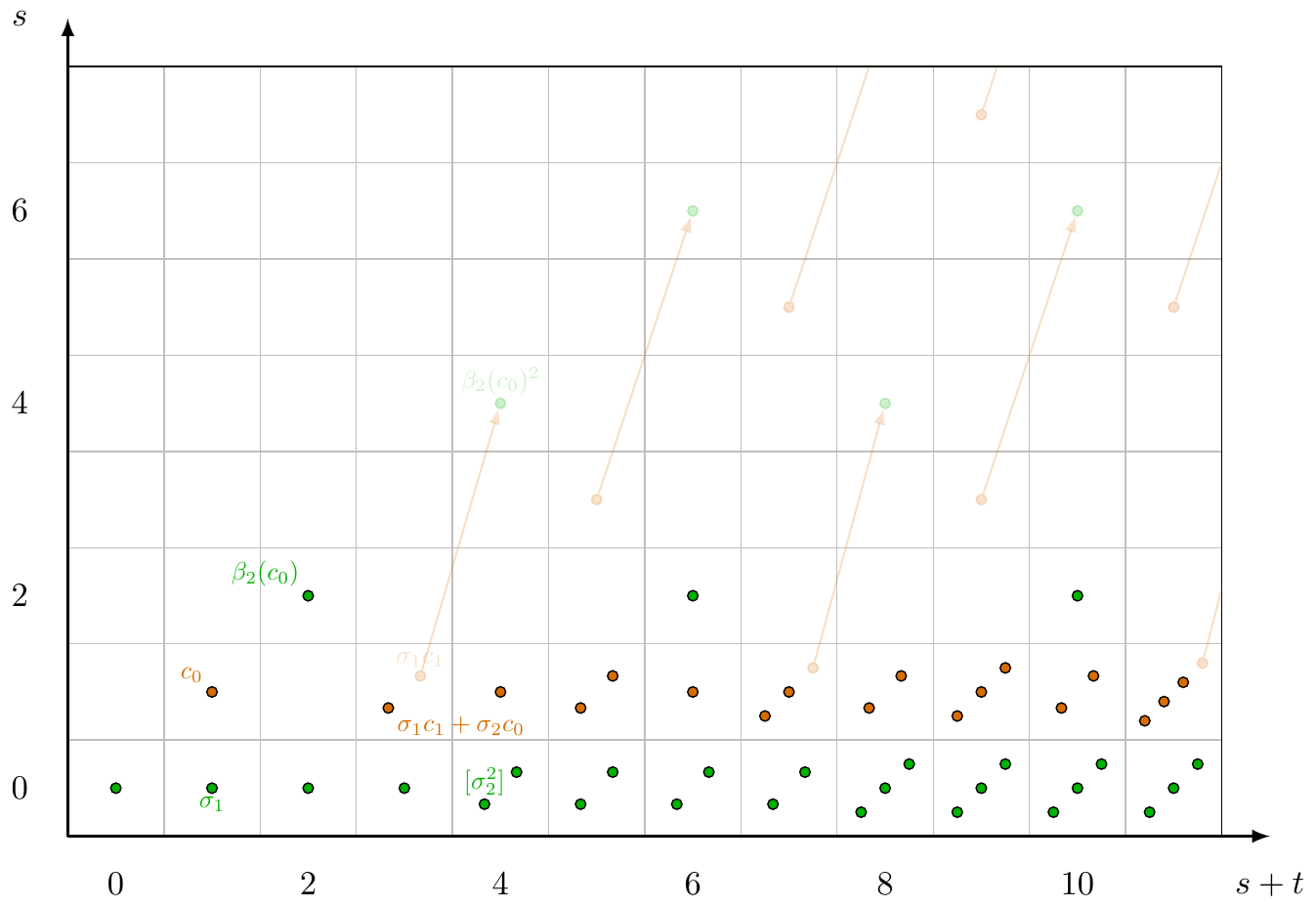}
    \caption{The $E_3$-page with differentials.}
\end{figure}
This shows that $E_4$ has a horizontal vanishing line of height 3, and is therefore the $E_\infty$-page.
\subsection{Poincar\'e Series}
The $0$-line $E_{\infty}^{0,*}$ is given by $\Ffield_2[\sigma_1, [\sigma_2^2]]$, which has Poincar\'e series $1/((1-t)(1-t^4))$. The $2$-line $E_{\infty}^{2,*}$ is given by $\Ffield_2[[\sigma_2^2]]\{\beta_2(c_0)\}$, which has Poincar\'e series $t^2/(1-t^4)$.

The $1$-line of $E_3$ is given by
\begin{equation}
    E_3^{1,*} \cong \Ffield_2[\sigma_1, [\sigma_2^2]]\{\sigma_1 c_1\} \oplus \Ffield_2[[\sigma_2^2]]\{c_0, [\sigma_2 c_0]\},
\end{equation}
which has Poincar\'e series
\begin{equation}
    \frac{t^3}{(1-t)(1-t^4)} + \frac{t}{1-t^2}.
\end{equation}
The $d_3$-differentials originating from the $s=1$-line kill 1 class in each $s+t$-stem with $s+t \equiv 3 \pmod{4}$, hence $E_4^{1,*} = E_{\infty}^{1,*}$ has Poincar\'e series
\begin{equation}
    \frac{t^3}{(1-t)(1-t^4)} + \frac{t}{1-t^2} - \frac{t^3}{1-t^4}.
\end{equation}
Taking everything together, and using
\begin{equation}
    1+(1-t)t^2 + t^3 + t(1+t^2)(1-t) - t^3(1-t) = 1+t,
\end{equation}
we get that $H^*(BM_{16}; \Ffield_2)$ has Poincar\'e series
\begin{equation}
    \frac{1+t}{(1-t)(1-t^4)} = \frac{1}{(1-t)^2(1+t^2)},
\end{equation}
\subsection{The exponent}
The above computation shows that $\exp_{\sE_{(2)}} \underline{H\Ffield_2} \geq 3$, which together with \cref{m16:expupper} shows that 
\begin{Proposition}
    \label{prop:m16e2exp}
    The $\sE_{(2)}$-exponent satisfies
\begin{equation}
    3 \leq \exp_{\sE_{(2)}} \underline{H\Ffield_2} \leq 4.
\end{equation}
\end{Proposition}

\section{$D_8 \ast C_4$}
\label{sec:d8cpc4}
\subsection{Introduction}
Let $D_8 = \langle \sigma, \rho \rangle$ be the dihedral group of order 8 and let $C_4 = \langle \gamma \rangle$ be the cyclic group of order 4. Both these groups have central cyclic subgroups of order 2, for $D_8$ this is $\langle \rho^2 \rangle$ and for $C_4$ this is $\langle \gamma ^2 \rangle$. The \emph{central product} of $D_8$ and $C_4$ is defined to be the direct product with these central subgroups identified:
\begin{equation}
    D_8 \ast C_4 := D_8 \times C_4 / \langle \rho^2 \gamma^{-2} \rangle.
\end{equation}
In this section we will evaluate the $\sF$-homotopy limit spectral sequence converging to $H^*(BD_8 \ast C_4;\Ffield_2)$ for the family $\sF = \sA$ of abelian subgroups of $D_8 \ast C_4$. This ring is isomorphic to (\cite[App.\ C, \#8(16)]{carlson2003cohomology})
\begin{equation}
    \Ffield_2[z, y, x, w]/(zx + y^2 + x^2, z^2x +  zx^2).
\end{equation}
with degrees $|z| = |y| = |x| = 1$, and $|w| = 4$.
\subsection{Summary of the computation}
Because the computation of the $\sA$-homotopy limit spectral sequence of $D_8 \ast C_4$ is rather long, we summarize it here.

The cohomology of $D_8 \ast C_4$ is detected on abelian subgroups (see \cref{d8cpc4:lem1} below). This lemma is part of the input of the calculation, and the proof of this lemma uses the depth of the cohomology ring of $D_8 \ast C_4$ as input. In particular, the computation of the $\sA$-homotopy limit spectral sequence carried out here does not provide a new computation of the cohomology of $D_8 \ast C_4$.

The strategy is the same as the one followed for the calculation of the $\sE_{(2)}$-homotopy limit spectral sequence of $D_8$ (\cref{sec:dihgps}), that is, we apply the decomposition from \cref{sec:absplit}.

Therefore we compute the LHSSS's needed for this decomposition. We need 4 LHSSS's, but 3 of them we can identify using automorphisms of $D_8 \ast C_4$, which we construct first.

After the computation of the LHSSS's we compute the map $j$ from the LES in (\ref{absplit:eq3}). This then allows us to calculate the $A$- and $B$- summand from (\ref{absplit:eq4}), and hence the $E_2$-page.

Next we compute an upper bound of 3 on the $\sA$-exponent of $\underline{H\Ffield_2}_{D_8 \ast C_4}$, which implies that the $\sA$-homotopy limit spectral sequence will collapse at $E_4$ with a vanishing line of height 4. However, because $D_8 \ast C_4$ has its cohomology detected on abelian subgroups, there will be even a vanishing line of height 1, i.e., the spectral sequence will collapse to the $s=0$-line.

Using naturality and the fact that the spectral sequence will collapse to the $s=0$-line at $E_4$, we can compute the rank of each differential in the $\sA$-homotopy limit spectral sequence. We are not able to determine precisely which class hits which, but we are able to determine which classes in the $s=0$-line support a differential.
\subsection{The orbit category}
The family $\sF = \sA$ of abelian subgroups gives the following subcategory of the orbit category:
\begin{figure}[H]
    \centering
    \includegraphics[angle=90,height=\textheight]{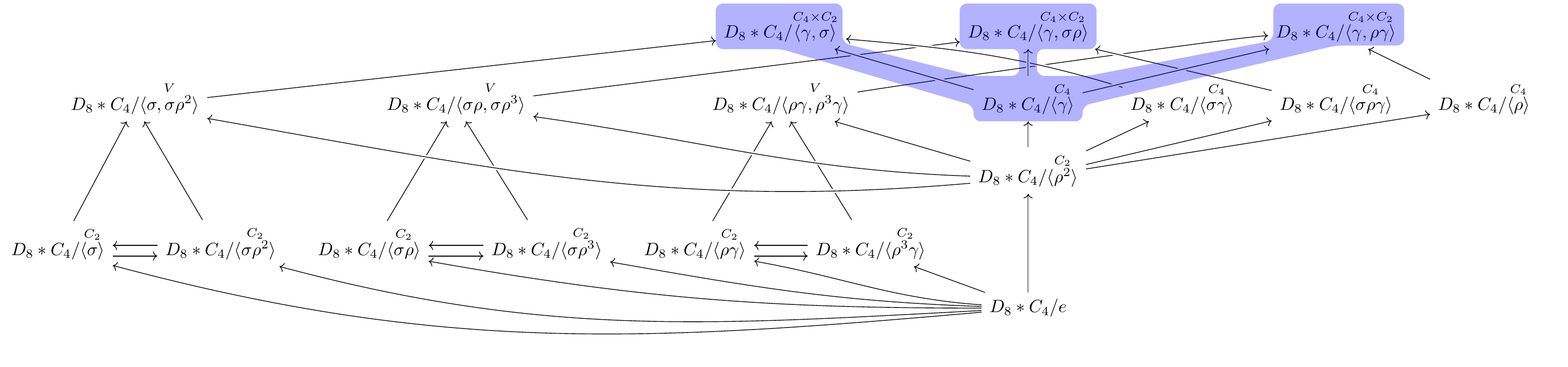}
    \caption{The category $\sO(D_8 \ast C_4)_{\sA}$. All Weyl groups are elementary abelian. The subcategory spanned by the highlighted objects is homotopy final.}
    \label{d8cpc4:fig1}
\end{figure}
The Weyl groups are not displayed in \cref{d8cpc4:fig1}. They all turn out to be elementary abelian.

We introduce the following notation for the maximal abelian subgroups of $D_8 \ast C_4$, all of which are isomorphic to $C_4 \times C_2$
\begin{align}
    A_1 & = \langle \gamma, \rho \gamma \rangle, \\
    A_2 & = \langle \gamma, \sigma \rangle, \\
    A_3 & = \langle \gamma, \sigma \rho \rangle.
\end{align}
In addition, we denote the cyclic group of order 4 generated by $\gamma$ by 
\begin{equation}
    C = \langle \gamma \rangle.
\end{equation}
The group $C$ equals any non-redundant intersection of the $A_i$'s:
\begin{equation}
    C = A_i \cap A_j = A_1 \cap A_2 \cap A_3
\end{equation}
for all $i \neq j$. We apply the decomposition of \cref{sec:absplit} with the $H_i = A_i$, and with $K = C$. Therefore we will compute the LHSSS's obtained from the extensions $C \to D_8 \ast C_4 \to C_2 \times C_2$ and $A_i \to D_8 \ast C_4  \to C_2$. Before doing so, we determine an upper bound on the $\sA$-exponent of $\underline{H\Ffield_2}_{D_8 \ast C_4}$.

\subsection{An upper bound on the $\sA$-exponent}
To give an upper bound on $\exp_{\sA} \underline{H\Ffield_2}_{D_8 \ast C_4}$, we will use the following lemma.
\begin{Lemma}
    The $\Ffield_2$-cohomology of $D_8 \ast C_4$ is detected on abelian subgroups.
    \label{d8cpc4:lem1}
\end{Lemma}
\begin{Proof}
    The centralizers of the elementary abelian 2-subgroups of rank 2 of $D_8 \ast C_4$ are precisely the $A_i$. If there would be a cohomology class which is not detected on the $A_i$, then this would imply by \cite[Thm.\ 2.3]{carlson95} that the depth of $H^*(BD_8 \ast C_4; \Ffield_2)$ is $< 2$. But  the ring $H^*(BD_8 \ast C_4; \Ffield_2)$ has depth 2 (see, e.g., \cite[App.\ C, \#8(16)]{carlson2003cohomology}), which gives the desired result.
\end{Proof}
\begin{Corollary}
    The $\sA$-homotopy limit spectral sequence converging to \\ $H^*(BD_8 \ast C_4)$ collapses at a finite page to the $s=0$-line.
    \label{d8cpc4:cor1}
\end{Corollary}
In addition, \cref{d8cpc4:lem1} allows us to prove the following.
\begin{Proposition}
    The $\sA$-exponent of $\underline{H\Ffield_2}_{D_8 \ast C_4}$ is $\leq 3$.
    \label{d8cpc4:prop2}
\end{Proposition}
\begin{Proof}
    For each $i$, let $\tau_i$ denote the pullback of the sign representation along the quotient map $D_8 \ast C_4 \to D_8 \ast C_4/A_i \cong C_2$. Then $T = \bigoplus_i \tau_i$ is a real representation of dimension 3 with isotropy in $\sA$. Moreover, $T$ has a trivial summand when restricted to any of the $A_i$, hence the Euler class $e_T$ is 0 by \cref{d8cpc4:lem1}. By \cref{cor:hfpeulerclass} the result follows.
\fxwarning{Need to fix overfull hboxes.}

\end{Proof}
\subsection{The decomposition}
We apply the decomposition of \cref{sec:absplit} with $K = C$ and the $H_i = A_i$.

We therefore proceed by computing the LHSSS's $E^{*,*}_*(A_i)$ and $E^{*,*}_*(C)$, and subsequently $A^*$ and $B^*$, which will then give us the $E_2$-page of the $\sA$-homotopy limit spectral sequence. But first we discuss some automorphisms of the group of $D_8 \ast C_4$ and some low-dimensional Betti numbers that will go into these computations.
\subsection{Automorphisms of $D_8 \ast C_4$}
The elements $\rho \gamma$, $\sigma$ and $\sigma \rho$ form a generating set of $D_8 \ast C_4$. In this section we construct an automorphism of $D_8 \ast C_4$ that allows us to identify the LHSSS's of $A_i \to D_8 \ast C_4 \to C_2$ for all $i$. 
\begin{Proposition}
    Every permutation of the generators $\rho \gamma$, $\sigma$ and $\sigma \rho$ extends to an automorphism of $D_8 \ast C_4$. 
    \label{prop:d8cpc4aut}
\end{Proposition}
\begin{Proof}
    Number the 16 elements of $D_8 \ast C_4$ as follows:
    \begin{table}[H]
        \centering
    \begin{tabular}[H]{rrcrr}
        1: & $e$ & & 9: & $\gamma$ \\
        2: & $\rho$ & & 10: & $\rho \gamma$ \\
        3: & $\rho^2$ & & 11: & $\rho^2 \gamma$ \\
        4: & $\rho^3$ & & 12: & $\rho^3 \gamma$ \\
        5: & $\sigma$ & & 13: & $\sigma \gamma$ \\
        6: & $\sigma\rho$ & & 14: & $\sigma \rho \gamma$ \\
        7: & $\sigma \rho^2$ & & 15: & $\sigma \rho^2 \gamma$ \\
        8: & $\sigma \rho^3$ & & 16: & $\sigma \rho^3 \gamma$
    \end{tabular}
    \centering
    \end{table}
    This numbering induces an inclusion $i \colon D_8 \ast C_4 \to \Sigma_{16}$ into the symmetric group on 16 elements, which on the generators given above is given by
    \begin{align}
        \rho \gamma & \mapsto (1, 10)(2, 11)(3, 12)(4, 9)(5, 16)(6, 13)(7, 14)(8, 15), \\
        \sigma & \mapsto (1, 5)(2, 6)(3, 7)(4, 8)(9, 13)(10, 14)(11, 15)(12, 16), \\
        \sigma \rho & \mapsto (1, 6)(2, 7)(3, 8)(4, 5)(9, 14)(10, 15)(11, 16)(12, 13).
    \end{align}
    One easily checks that the element of $\Sigma_{16}$ given by
    \begin{equation}
        (5, 6)(2, 4)(8, 7)(14, 15)(11, 9)(13, 16)
    \end{equation}
    induces an inner automorphism $\phi$ of $\Sigma_{16}$ which interchanges $i(\sigma)$ and $i(\sigma \rho)$ and leaves $i(\rho \gamma)$ fixed. Identifying $D_8 \ast C_4$ with its image in $\Sigma_{16}$ and restricting $\phi|_{\Im(i)}$ yields an automorphism of $D_8 \ast C_4$ which interchanges $\sigma$ and $\sigma \rho$, and leaves $\rho \gamma$ fixed.

    Likewise, the element of $\Sigma_{16}$ given by
    \begin{equation}
        (5, 10)(14, 16)(12, 7)(4, 15)(9, 11)(13, 2)
    \end{equation}
    yields an automorphism of $D_8 \ast C_4$ that interchanges $\sigma$ and $\rho \gamma$, and leaves $\sigma\rho$ fixed. 

    The two automorphisms we constructed restrict to two generating transpositions on the set $\{\rho \gamma, \sigma, \sigma \rho\}$, and hence we are done.
\end{Proof}
Using the automorphisms constructed in the proof of \cref{prop:d8cpc4aut}, we can show the following.
\begin{Proposition}
    There is an automorphism of $D_8 \ast C_4$ that cyclically permutes the elements of the set $\{\rho \gamma, \sigma, \sigma \rho\}$ and fixes $\gamma$. In particular this automorphism permutes the subgroups in $\{A_1, A_2, A_3\}$ cyclically.
    \label{prop:d8cpc4aut2}
\end{Proposition}
\begin{Proof}
    The existence of an automorphism that cyclically permutes the elements of $\{\rho \gamma, \sigma, \sigma \rho \}$ immediately follows from \cref{prop:d8cpc4aut}. For the statement about $\gamma$, we note that $\gamma = (\sigma \cdot \sigma \rho)^3 \cdot \rho \gamma$, and a direct computation shows that the effect of any transposition on $\gamma$ is given by
    \begin{align}
        (\sigma\rho \cdot \sigma )^3 \cdot \rho \gamma = \gamma^{-1}, \\
        (\rho \gamma \cdot \sigma \rho)^3 \cdot \sigma = \gamma^{-1}, \\
        (\sigma \cdot \rho \gamma)^3 \cdot \sigma \rho = \gamma^{-1}.
    \end{align}
    Hence the cyclic permutation of 3 elements, being the composite of two transpositions, fixes $\gamma$.
\end{Proof}

\subsection{Low dimensional Betti numbers}

We will denote the dimension of $H^i(BD_8 \ast C_4; \Ffield_2)$, i.e.\ the $i$-th Betti number, by $h^i$.
In the subsequent computations, we will freely use the following.
\begin{Proposition}
    The Betti numbers $h^i$ for $0 \leq i \leq 4$ are given by the following table:
    \begin{table}[H]
        \centering
        \begin{tabular}{r|rrrrr}
            $i$ & 0 & 1 & 2 & 3 & 4 \\
            \hline
            $h^i$ & 1 & 3 & 5 & 6 & 7 
        \end{tabular}
        \caption{The low dimensional Betti numbers $h^i$.}
    \end{table}
    \label{d8cpc4:prop1}
\end{Proposition}
\begin{Proof}
    For $i=0$ this is clear, and for $i=1$ we observe that $(D_8 \ast C_4)^{\ab} = D_8 \ast C_4/\langle \rho^2 \rangle \cong C_2^{\times 3}$.

    The statements for $i \geq 2$ can, for instance, be verified using computer algebra software such as GAP.
    
\end{Proof}
\subsection{The LHSSS for $C \to D_8 \ast C_4 \to C_2 \times C_2$}
\begin{figure}[H]
    \centering
    \includegraphics[width=0.9\textwidth]{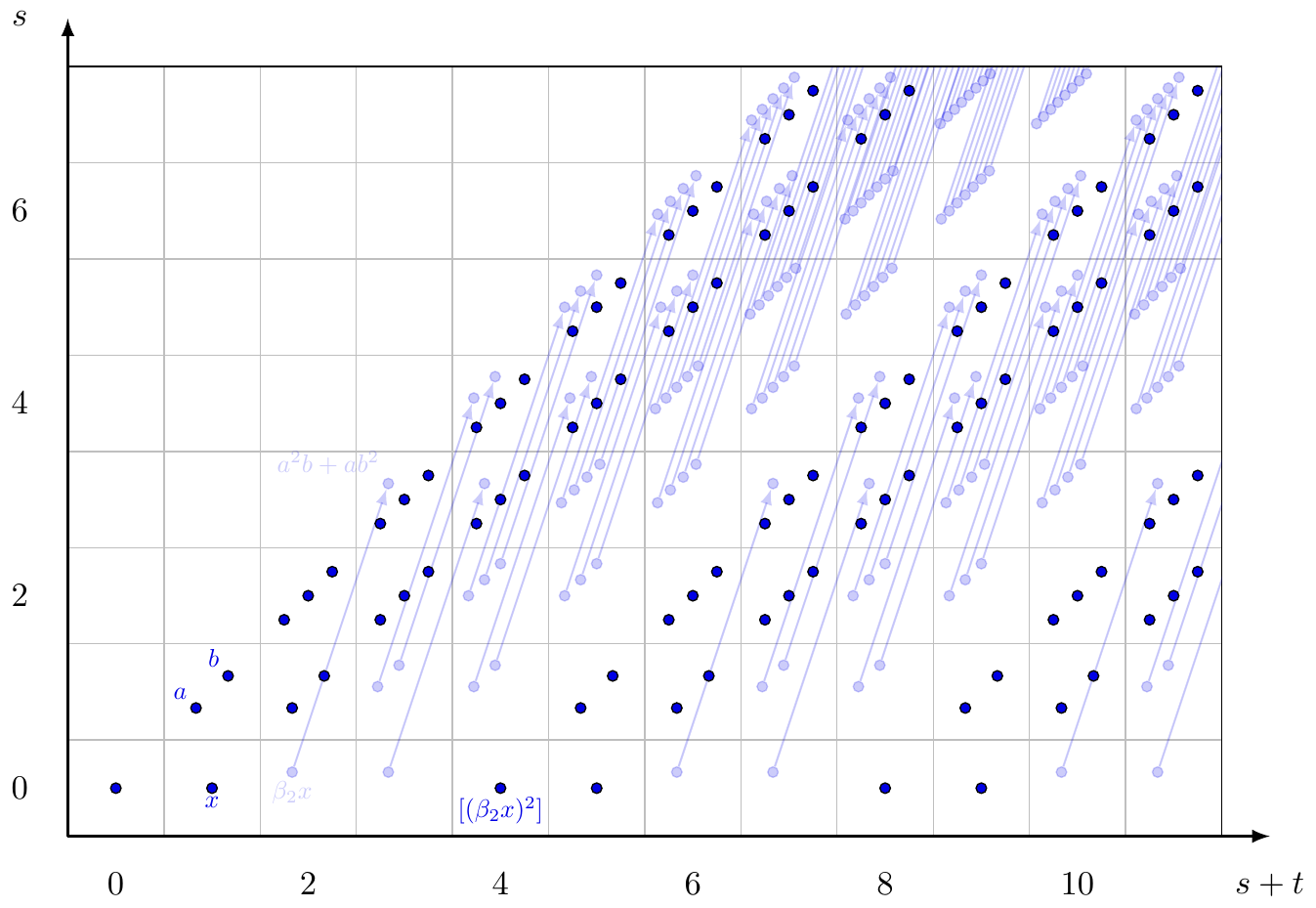}
    \caption{The $E_3$-page.}
    \label{d8cpc4:fig3}
\end{figure}

We now compute the LHSSS for the central extension
\begin{equation}
    C_4 \cong \overset{\gamma}{C} \to D_8 \ast C_4 \to \overset{\sigmabar}{C_2} \times \overset{\overline{\sigma \rho}}{C_2},
    \label{d8cpc4:eq2}
\end{equation}
where the elements above the groups are generators.
We introduce the following notation for the elements of $H^*(BC_2 \times C_2)$:
\begin{align}
    a & = \delta_{\sigmabar}, \\
    b &= \delta_{\overline{\sigma \rho}}, \\
    \intertext{and of $H^*(BC)$:}
    x & = \delta_{\gamma}, \\
    \beta_2 x & = \beta_2 \delta_{\gamma},
\end{align}
where $\beta_2$ denotes the second order Bockstein. Since the extension (\ref{d8cpc4:eq2}) is central, the local coefficient system of the associated LHSSS is trivial, and we get
\begin{equation}
    E_2^{*,*} = \Ffield_2[a,b,x,\beta_2 x]/(x^2),
\end{equation}
with $(s,t)$-degrees given by $|a| = |b| = (1,0)$, $|x| = (0,1)$ and $|\beta_2 x| = (0,2)$. See \cref{d8cpc4:fig3} for a depiction of the $E_3$-page, which up to differentials will turn out to be isomorphic to the $E_2$-page, because there are no non-zero differentials on $E_2$, as we will see in the next section.
\subsubsection{Differentials on $E_2$ and $E_3$}
The classes $a$ and $b$ are permanent cycles for degree reasons, and $x$ is a permanent cycle because $h^1 = 3$. Since the 2-stem of the $E_2$-page has 6 classes, but $h^2 = 5$, $\beta_2 x$ will have to support a differential, which for degree reasons will be either a $d_2$ or a $d_3$. There are finitely many possibilities, namely:
\begin{align}
    d_2(\beta_2 x) & \in \Ffield_2\{a^2x, abx, b^2x\}, \\
    d_3(\beta_2 x) & \in \Ffield_2\{a^3, a^2b, ab^2, b^3 \}.
\end{align}
We first apply naturality to the map of central extensions
\begin{equation}
    \begin{tikzcd}
        \overset{\gamma}{C} \arrow{r} \arrow[equals]{d} & \overset{\gamma}{C_4} \times \overset{\rho}{C_2} \arrow{r} \arrow[hook]{d} &  \overset{\rhobar}{C_2} \arrow[hook]{d} \\
        \overset{\gamma}{C} \arrow{r} & D_8 \ast C_4 \arrow{r} & \overset{\sigmabar}{C_2} \times \overset{\overline{\sigma \rho}}{C_2}
    \end{tikzcd}
\end{equation}
The rightmost vertical map is given by $\rhobar  \mapsto \sigmabar \cdot \overline{\sigma \rho}$, and therefore on cohomology by $a = \delta_{\sigmabar} \mapsto \delta_{\rhobar}$ and $b = \delta_{\overline{\sigma \rho}} \mapsto \delta_{\rhobar}$. The LHSSS of the top row collapses at $E_2$ by the K\"unneth theorem. Therefore $d_2(\beta_2 x)$ and $d_3(\beta_2 x)$, whatever they are, will have to map to 0 in the LHSSS of the top row. This cuts down the possibilities to
\begin{align}
    d_2(\beta_2 x) & \in \Ffield_2\{(a^2+b^2)x\}, \\
    d_3(\beta_2 x) & \in \Ffield_2\{ a^3 + b^3, a^2b + ab^2 \}.
    \label{d8cpc4:eq3}
\end{align}
Applying naturality again, this time to the map of central extensions
\begin{equation}
    \begin{tikzcd}
        \overset{\gamma}{C} \arrow{r} \arrow[equals]{d} & \overset{\gamma}{C_4} \times \overset{\sigma}{C_2} \arrow{r} \arrow[hook]{d} &  \overset{\sigmabar}{C_2} \arrow[hook]{d} \\
        \overset{\gamma}{C} \arrow{r} & D_8 \ast C_4 \arrow{r} & \overset{\sigmabar}{C_2} \times \overset{\overline{\sigma \rho}}{C_2}
    \end{tikzcd},
\end{equation}
for which the right hand vertical map is given by $\sigmabar \mapsto \sigmabar$, hence on cohomology classes by $a = \delta_{\sigmabar}  \mapsto \delta_{\sigmabar}$ and $b = \delta_{\overline{\sigma \rho}} \mapsto 0$. Again, the LHSSS of the top row collapses at $E_2$. Combined with (\ref{d8cpc4:eq3}) this shows that $d_2(\beta_2 x) = 0$, and that
\begin{equation}
    d_3(\beta_2 x) \in \Ffield_2\{a^2 b + a b^2\}.
\end{equation}
Since $\beta_2 x$ has to support some differential, and we have determined the only possible non-zero differential, we arrive at
\begin{equation}
    d_3(\beta_2 x) = a^2b + ab^2.
    \label{d8cpc4:eq5}
\end{equation}
\subsubsection{The $E_4 = E_\infty$-page}
The fact that $d_2 = 0$ implies $E_3 \cong E_2$, and the differential (\ref{d8cpc4:eq5}) shows that the $E_3$-page is given by
Therefore the $E_4$-page is given by
\begin{equation}
    E_4 \cong \Ffield_2[a, b, x, [(\beta_2x)^2]]/(x^2, a^2b + ab^2).
\end{equation}
The $E_4$-page has 7 classes in the 4-stem, and since $h^4 = 7$ (\cref{d8cpc4:prop1}), the class $[(\beta_2 x)^2]$ is a permanent cycle. Alternatively, one can apply Kudo's transgression theorem: the class $\beta_2 x$ transgresses with $d_3(\beta_2 x) = a^2 b + ab^2$. This implies that the class $[(\beta_2 x)^2] = [\Sq^2(\beta_2 x)]$ also transgresses, i.e., $d_3([(\beta_2 x)^2] = 0$ and $d_4([(\beta_2 x)^2]) = 0$, and  $d_5([(\beta_2 x)^2]) = [\Sq^2(a^2 b + ab^2)]$. Applying the Cartan formula several times shows 
\begin{align}
    \Sq^2(a^2 b) & = \Sq^2(a^2)b + \Sq^1(a^2) \Sq^1(b) + a \Sq^2(b) \\
    & = a^4 b,
\end{align}
and hence
\begin{equation}
    \Sq^2(a^2 b + a b^2) = a^4 b + a b^4.
\end{equation}
But since $a$ and $b$ are permanent cycles, we have
\begin{align}
    d_3((a^2 + ab + b^2)\beta_2 x) & = (a^2 + ab + b^2)(a^2 b + a b^2) \\
    & = a^4 b + a b^4,
\end{align}
so $d_5([(\beta_2 x)^2])= 0$ also. In any case, $[(\beta_2 x)^2]$ is a permanent cycle and the spectral sequence collapses at $E_4$.
\subsection{The LHSSS for $A_i \to D_8 \ast C_4 \to C_2$}
We now compute the LHSSS for the extensions
\begin{equation}
    A_i \to D_8 \ast C_4 \to C_2.
\end{equation}
By \cref{prop:d8cpc4aut2}, we can restrict to the case $i=1$, which is the LHSSS of the extension
\begin{equation}
    \overset{\gamma}{C_4} \times \overset{\rho\gamma}{C_2} \to D_8 \ast C_4 \to \overset{\sigmabar}{C_2}.
\end{equation}
We introduce the following notation for the following element of $H^*(BC_2)$:
\begin{align}
    c_1 & = \delta_{\sigmabar}, \\
    \intertext{and for the elements of $H^*(BC_4 \times C_2)$:}
    y_1 & = \delta_\gamma, \\
    \beta_2 y_1 & = \beta_2 \delta_\gamma, \\
    z_1 & = \delta_{\rho \gamma},
\end{align}
where $\beta_2$ denotes the second order Bockstein.
\subsubsection{The local coefficient system}
Since
\begin{align}
    \sigma \gamma \sigma^{-1} & = \gamma, \\
    \sigma (\rho \gamma) \sigma^{-1} & = \gamma^2 \cdot \rho \gamma, \\
    \label{d8cpc4:loccoeff:eq1}
\end{align}
we have that the action of $C_2$ on $H^1(BC_4 \times C_4; \Ffield_2)$ is determined by
\begin{align}
    y_1 & \mapsto y_1, \\
    z_1 & \mapsto z_1.
\end{align}
To calculate the action on $\beta_2 y_1$, we note that the second order Bockstein is natural, but only well-defined up to the image of the first order Bockstein, which in degree 2 is given by the squares of the elements in degree 1, i.e., $\Ffield_2\{z_1^2\}$. Therefore $\sigma \cdot \beta_2 y_1 \in \{\beta_2 y_1, \beta_2 y_1 + z_1^2\}$. To determine which of the two possibilities we have, we apply naturality to the inclusion of groups
\begin{align}
    \overset{\gamma^2}{C_2} \times \overset{\rho \gamma}{C_2} & \hookrightarrow  \overset{\gamma}{C_4} \times \overset{\rho \gamma}{C_2}
    \label{d8cpc4:loccoeff:eq2}
\end{align}
which on cohomology is given by
\begin{align}
    y_1 = \delta_{\gamma} & \mapsto 0, \\
    z_1 = \delta_{\rho \gamma} & \mapsto \delta_{\rho \gamma}, \\
    \beta_2 y_1 & \mapsto \delta_{\gamma^2}^2.
\end{align}
The automorphism (\ref{d8cpc4:loccoeff:eq1}) on $C_4 \times C_2$ restricts along the inclusion (\ref{d8cpc4:loccoeff:eq2}) to the automorphism of $C_2 \times C_2$ which on elements is given by
\begin{align}
    \gamma & \mapsto \gamma, \\
    \rho \gamma & \mapsto \gamma^2 \cdot \rho \gamma,
\end{align}
and therefore on the cohomology classes by
\begin{align}
    \delta^2_{\gamma^2} & \mapsto \delta^2_{\gamma^2} + \delta_{\rho \gamma}^2, \\
    \delta_{\rho \gamma} & \mapsto \delta_{\rho \gamma}.
    \label{d8cpc4:eq8} 
\end{align}
This implies, by naturality, that $\sigma \cdot \beta_2 y_1$ must be class of $H^*(BC_4 \times C_2;\Ffield_2)$ that restricts via (\ref{d8cpc4:loccoeff:eq2}) to a class with a non-zero $\delta^2_{\rho \gamma}$-coefficient in $H^*(BC_2 \times C_2; \Ffield_2)$. But there is only one such class in $\{\beta_2 y_1, \beta_2 y_1 + z_1^2\}$ by (\ref{d8cpc4:eq8}), namely $\beta_2 y_1 + z_1^2$ and therefore
\begin{equation}
    \sigma \cdot \beta_2 y_1 = \beta_2 y_1 + z_1^2.
\end{equation}
\subsubsection{The $E_2$-page}
We have $H^*(BC_4 \times C_2) \cong \Ffield_2[y_1, z_1, \beta_2y_1]/(y_1^2)$, which has a $C_2 \cong \langle \sigmabar\rangle$-action which we determined in the previous subsection: $y_1$ and $z_1$ are fixed and $\beta_2 y_1  \mapsto \beta_2 y_1 + z_1^2$. The $E_2$-page of the LHSSS is given by the group cohomology of $C_2$ with these coefficients.

The $C_2$-module $H^*(BC_4 \times C_2)$ splits as
\begin{equation}
    \Ffield_2[y_1, z_1, \beta_2 y_1]/(y_1^2)  = \Ffield_2[z_1, \beta_2 y_1] \oplus \Ffield_2[z_1, \beta_2 y_1]\{y_1\}
\end{equation}
(as $C_2$-modules). The $C_2$-module $\Ffield_2[z_1, \beta_2 y_1]$ splits further as
\begin{equation}
    \Ffield_2[z_1^2,\beta_2 y_1] \oplus \Ffield_2[z_1^2,\beta_2 y_1]\{z_1\}
    \label{d8cpc4:eq9}
\end{equation}
(as $C_2$-modules). This splitting enables the computation of the $C_2$-cohomology with coefficients in $\Ffield_2[y_1, z_1, \beta_2 y_1])/(y_1^2)$ by computing the $C_2$-cohomology with coefficients in the $C_2$-module (\ref{d8cpc4:eq9}). 

The $C_2$-cohomology $H^*(BC_2; M)$ with coefficients in any $C_2$-module $M$ is in degree $* = 0$ given by the $C_2$-invariants of $M$, and in degrees $* \geq 1$ by the $C_2$-invariants of $M$ modulo the $C_2$-coinvariants of $M$, as can be seen from mapping the standard minimal free $\Ffield_2[C_2]$-resolution of $\Ffield_2$ into $M$.

The $C_2$-module $\Ffield_2[z_1^2, \beta_2 y_1]$, equals $\Ffield_2[z_1^2 + \beta_2 y_1, \beta_2 y_1]$. On this latter polynomial algebra, $C_2$ acts by interchanging the polynomial generators. The invariants of this action are well-known and are given by
\begin{equation}
    \Ffield_2[z_1^2, (z_1^2 + \beta_2 y_1)\beta_2 y_1].
\end{equation}
We introduce the notation
\begin{equation}
    q_1 = (z_1^2 + \beta_2 y_1)\beta_2 y_1,
\end{equation}
so that the invariants are given by $\Ffield_2[z_1^2,q_1]$, and the invariants modulo the coinvariants are given by $\Ffield_2[q_1]$.

Assembling back the splitting gives the invariants
\begin{align}
    \Ffield_2[y_1, z_1, q_1]/(y_1^2),
    \intertext{and the invariants modulo coinvariants}
    \Ffield_2[y_1, z_1, q_1]/(y_1^2, z_1^2).
\end{align}
Assembling all the lines together then gives as the $E_2$-page (recall that $c_1 = \delta_{\sigmabar}$):
\begin{equation}
    E_2 = \Ffield_2[c_1, y_1, z_1, q_1]/(c_1z_1^2, y_1^2)
\end{equation}
with $(s,t)$-degrees given by $|c_1| = (1,0)$, $|y_1| = |z_1| = (0,1)$ and $|q_1| = (0,4)$. A depiction of the $E_2$-page is
\begin{figure}[H]
    \centering
    \includegraphics[width=0.9\textwidth]{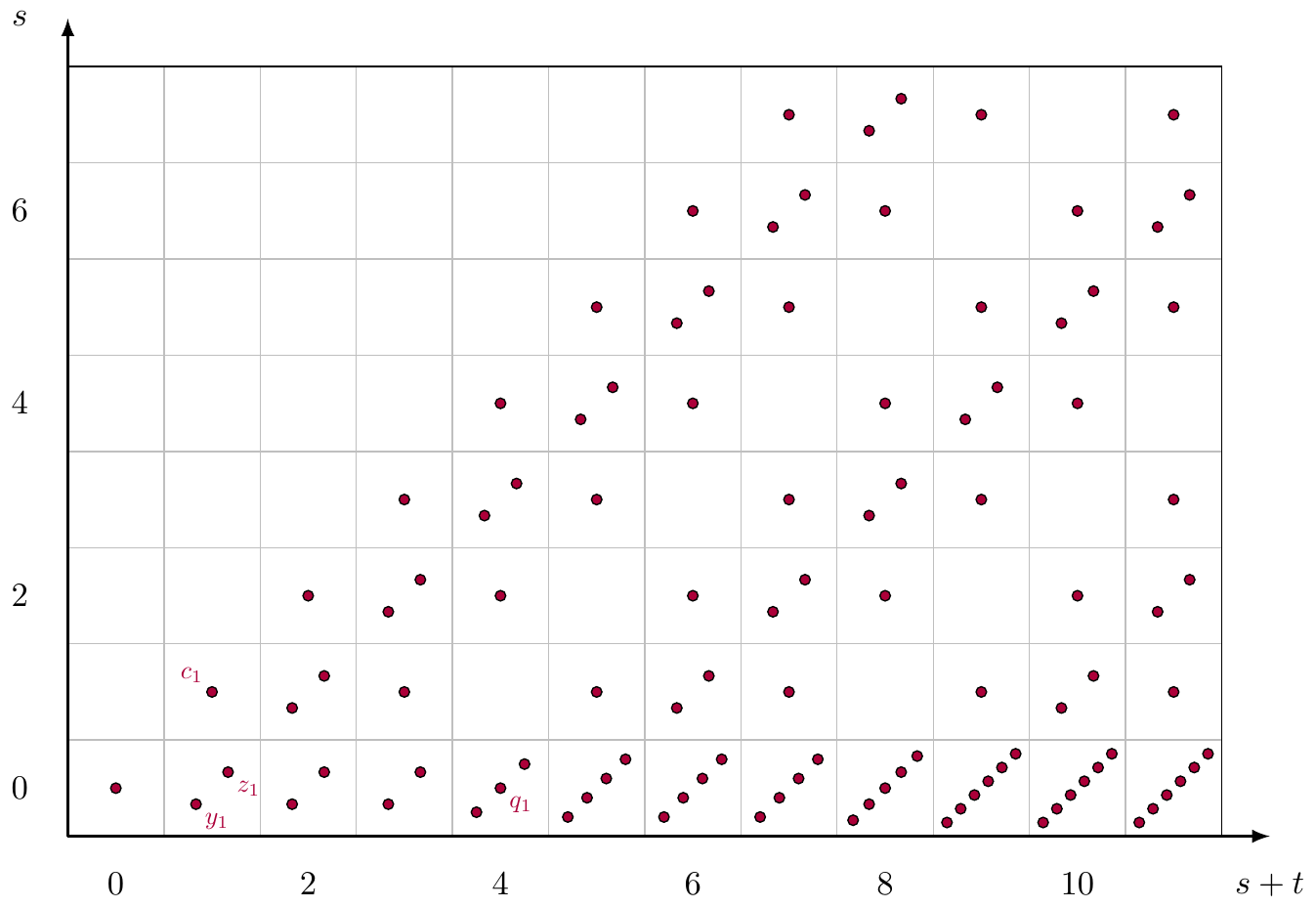}
    \caption{The $E_2 = E_\infty$-page.}
    \label{d8cpc4:chart1}
\end{figure}
The class $c_1$ is a permanent cycle for degree reasons, and $y_1$, $z_1$ and $q_1$ are permanent cycles because $h^1$ and $h^4$ equal to the number of classes in the 1-stem and 4-stem respectively (\cref{d8cpc4:prop1}). Hence the spectral sequence collapses at $E_2$.
\begin{Remark}
    To get a better understanding of the multiplicative structure depicted in \cref{d8cpc4:chart1} it can be helpful to observe that $E_2 \cong \Ffield_2[y_1]/(y_1^2) \otimes \Ffield_2[c_1, z_1, q_1]/(c_1z_1^2)$ and to display only the second tensor factor:
    \begin{figure}[H]
        \centering
        \includegraphics[width=0.9\textwidth]{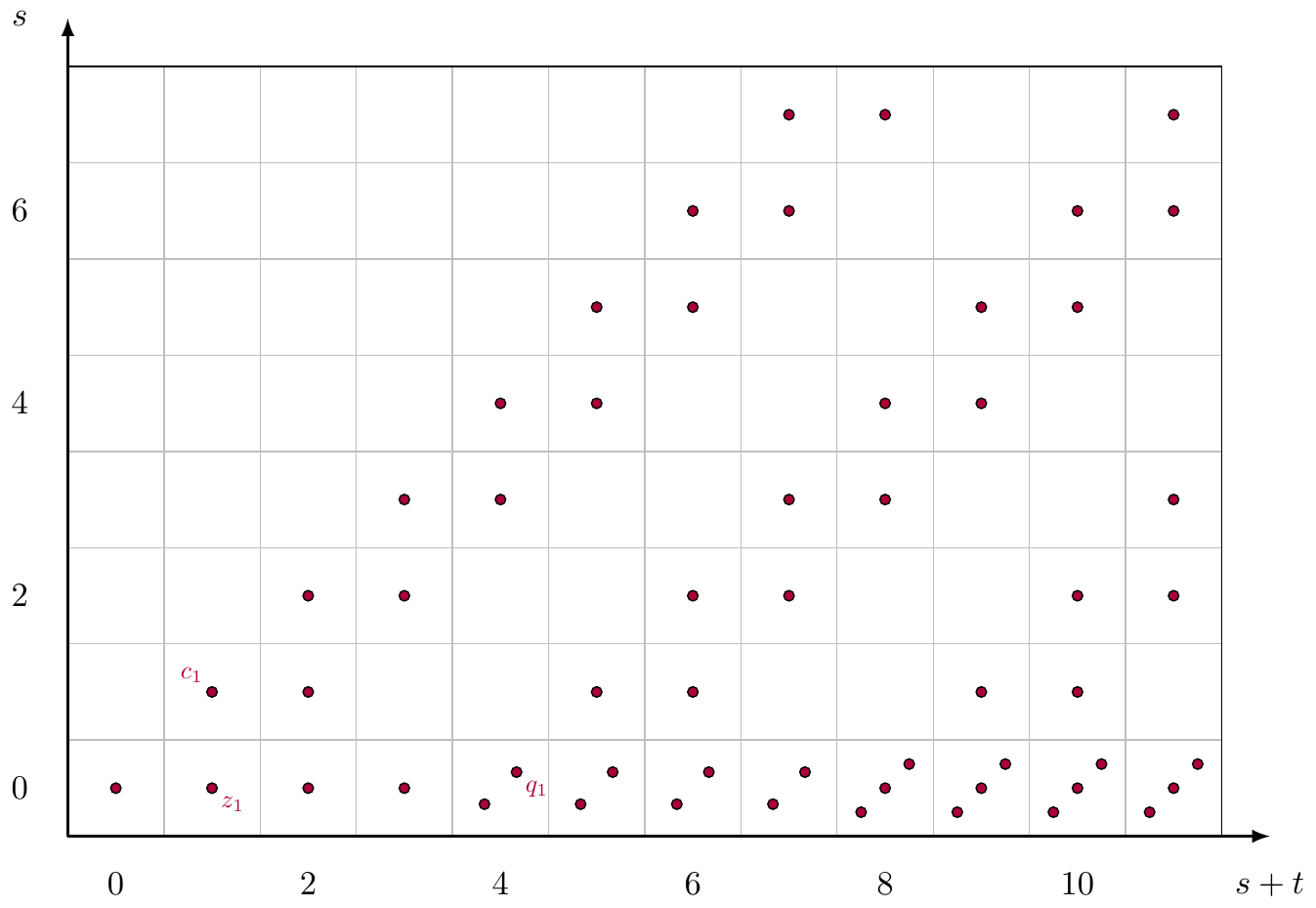}
        \caption{The tensor factor $\Ffield_2[c_1, z_1, q_1]/(c_1 z_2^2)$ of the $E_2$-page.}
    \end{figure}
\end{Remark}
\subsection{Summary of LHSSS computations}
We briefly summarize the LHSSS computations that we have carried out. The $E_2 = E_\infty$-page of the LHSSS of the extension $A_1 \to D_8 \ast C_4 \to C_2$ is given by
\begin{equation}
    E_2 = E_\infty \cong \Ffield_2[c_1, y_1, z_1, q_1]/(c_1 z_1^2, y_1^2).
\end{equation}

For the subgroup $A_2$, we introduce the notation
\begin{align}
    c_2 & = \delta_{\rhobar}\\
    \intertext{for the element of $H^*(BD_8 \ast C_4/A_2)$ and}
    y_2 & =\delta_\gamma, \\
    z_2 & =  \delta_\sigma, \\
    q_2 & = (\delta_\sigma^2 + \beta_2 \delta_\gamma) \beta \delta_\gamma \\
    \intertext{for the elements of $H^*(BA_2)$.}
\end{align}
For the subgroup $A_3$ we introduce the notation
\begin{align}
    c_3 & = \delta_{\sigmabar}, \\
    y_3 & = \delta_\gamma, \\
    z_3 &  = \delta_{\sigma \rho}, \\
    q_3 & = (\delta_{\sigma \rho}^2 + \beta_2 \delta_{\gamma})\beta_2 \delta_\gamma.
\end{align}
Then, by the automorphism of \cref{prop:d8cpc4aut2} the $E_2 = E_\infty$-page of the LHSSS of the extensions $A_i \to D_8 \ast C_4 \to C_2$ is given by
\begin{equation}
    E_2 = E_\infty \cong \Ffield_2[c_i, y_i, z_i, q_i]/(c_iz_i^2, y_i^2)
\end{equation}
with $(s,t)$-degrees given by $|c_i| = (1,0)$, $|y_i| = |z_i| = (0,1)$ and $|q_i| = (0,4)$.

The $E_2$-page of the LHSSS of the central extension $C \to D_8 \ast C_4 \to C_2 \times C_2$ is given by
\begin{equation}
    E_2 \cong \Ffield_2[a, b, x, \beta_2 x]/(x^2).
\end{equation}
The differential $d_2 = 0$, and $d_3$ is generated by $d_3(\beta_2 x) = a^2b +a b^2$. We introduce the notation $r = [(\beta_2 x)^2]$. Then the $E_4 = E_\infty$-page is given by
\begin{equation}
    E_4 = E_\infty \cong \Ffield_2[a,b,x,r]/(x^2, a^2b + a b^2)
\end{equation}
with $(s,t)$-degrees given by $|a| = |b| = (1,0)$, $|x| = (0,1)$ and $|r| = (0,4)$.
\subsection{The map $j^*$}
In this subsection we determine the effect of the map $j^*$ in the long exact sequence (\ref{absplit:eq3}), which in terms of group cohomology is given by
\begin{equation}
    \bigoplus_{i=1}^3 H^*(BC_2; H^*(BA_i; \Ffield_2)) \xrightarrow{j^*} \bigoplus_{k=1}^2 H^*(BC_2 \times C_2; H^*(BC; \Ffield_2)),
    \label{d8cpc4:eq10}
\end{equation}
given by the matrix
\begin{equation}
    \begin{pmatrix}
        1  & -1  &  0 \\
        0  &  1  & -1
    \end{pmatrix},
    \label{d8cpc4:eq11}
\end{equation}
and the inclusions $C \to A_i$ and the quotient maps $D_8 \ast C_4/C \to D_8 \ast C_4/A_i$. In (\ref{d8cpc4:eq10}) we have used that finite products and coproducts coincide in graded $\Ffield_2$-algebras.
\begin{Notation}
    For the elements of $\bigoplus_{i=1}^3 H^*(BC_2; H^*(BA_i; \Ffield_2))$ we will write $c_1 = (c_1, 0, 0)$, $y_1 = (y_1, 0, 0)$, $c_2 = (0, c_2, 0)$, etc. 

    For the elements of $\bigoplus_{k=1}^2 H^*(BC_2 \times C_2; H^*(BC; \Ffield_2)$ we write $a_1 = (a, 0)$, $b_1 = (b, 0)$, $a_2 = (0, a)$, etc. For example, we have $a_1 + a_2 = (a, a)$.
\end{Notation}
\begin{Proposition}
    The map $j$ from the long exact sequence (\ref{absplit:eq3}) in terms of (\ref{d8cpc4:eq10}) is given by
    \begin{align}
        c_1 & \mapsto a_1 + b_1 = (a + b, 0), \\
        c_2 & \mapsto b_1 + b_2 = (b, b), \\
        c_3 & \mapsto a_2 = (0, a), \\
        y_1 & \mapsto x_1 = (x, 0), \\
        y_2 & \mapsto x_1 + x_2 = (x, x), \\
        y_3 & \mapsto x_2 = (0, x), \\
        z_i & \mapsto 0, \\
        \beta_2 y_1 & \mapsto \beta_2 x_1 = (\beta_2 x, 0), \\
        \beta_2 y_2 & \mapsto \beta_2 x_1 +  \beta_2 x_2 = (\beta_2 x, \beta_2 x), \\
        \beta_2 y_3 & \mapsto \beta_2 x_2 = (0, \beta_2 x), \\
        q_1 & \mapsto (\beta_2 x_1)^2 = ((\beta_2 x)^2, 0), \\
        q_2 & \mapsto (\beta_2 x_1)^2 + (\beta_2 x_2)^2 = ((\beta_2 x)^2, (\beta_2 x)^2), \\
        q_3 & \mapsto (\beta_2 x_2)^2 = (0, (\beta_2 x)^2).
        \label{d8cpc4:eq4}
    \end{align}
\end{Proposition}
\begin{Proof}
    The map $j^*$ in the long exact sequence (\ref{absplit:eq3}) is determined by the inclusions $C \to A_i$, which is given by $\gamma \mapsto \gamma$, and the induced maps on Weyl groups $D_8 \ast C_4/C \to D_8 \ast C_4 /A_i$, which for $i = 1$ is given by
\begin{align}
    \sigmabar & \mapsto \sigmabar, \\
    \overline{\sigma \rho} & \mapsto \sigmabar,\\
    \intertext{for $i=2$ is given by}
    \sigmabar & \mapsto \ebar, \\
    \overline{\sigma \rho} & \mapsto \rhobar, \\
    \intertext{ and for $i = 3$ is given by}
    \sigmabar & \mapsto \sigmabar, \\
    \overline{\sigma \rho} & \mapsto \ebar.
\end{align}
This and the matrix (\ref{d8cpc4:eq11}) shows that on cohomology classes we have
\begin{align}
    c_1 = \delta_{\sigmabar} & \mapsto {\delta_{\sigmabar}}_1 + {\delta_{\overline{\sigma \rho}}}_1 = a_1 + b_1, \\
    c_2 = \delta_{\rhobar} & \mapsto {\delta_{\overline{\sigma \rho}}}_1 + {\delta_{\overline{\sigma \rho}}}_2  = b_1 + b_2, \\
    c_3 = \delta_{\sigmabar} & \mapsto {\delta_{\sigmabar}}_2 = a_2,
\end{align}
which takes care of the first three assignments of (\ref{d8cpc4:eq4}).

For the remaining four assignments of (\ref{d8cpc4:eq4}), we observe that the restrictions
\begin{equation}
   H^*(BA_i; \Ffield_2) \to H^*(BC; \Ffield_2) 
\end{equation}
are for all $i$ on low dimensional cohomology classes given by
\begin{align}
    y_i = \delta_\gamma &\mapsto \delta_{\gamma} = x, \\
    z_i  & \mapsto 0, \\
    \beta_2 y_i = \beta_2 \delta_\gamma & \mapsto \beta_2 x, \\
    q_i = (z_i^2 + \beta_2 y_i)\beta_2 y_i & \mapsto (\beta_2 x)^2.
\end{align}
This together with the matrix (\ref{d8cpc4:eq11}) completes the proof of the remaining assignments.
\end{Proof}
\subsubsection{The $B$-summand}
Having determined $j^*$, we now compute the $B=\ker(j^*)$-summand of the short exact sequence (\ref{absplit:eq5}). The degree $\geq 1$ part of the domain of $j^*$ is given by
\begin{equation}
    \bigoplus_{i=1}^3 \Ffield_2[c_i, y_i, z_i, q_i]/(c_iz_i^2, y_i^2) \{c_i\}.
\end{equation}
We have that the map $j^*$ maps
\begin{align}
    c_1 & \mapsto a_1 + b_1, \\
    c_2 & \mapsto b_1 + b_2, \\
    c_3 & \mapsto a_2,
\end{align}
in particular, the $c_i$ map to algebraically independent polynomial generators of the codomain. Therefore the kernel of $j^*$ in degree $\geq 1$ is given by the direct sum of the kernels of the 3 summands.

The classes $q_i$ map to a polynomial generator in the codomain. The classes $y_i$ map to a non-zero element that squares to 0, but the classes $y_i$ themselves also square to 0. Finally, the classes $z_i$ map to 0. Therefore the kernel of each summand is given by
\begin{align}
    \Ffield_2[c_i, y_i, z_i, q_i]/(y_i^2, c_i z_i^2) \{c_i z_i\} & =  \Ffield_2[c_i, y_i, q_i]/(y_i^2) \{c_iz_i\},\\
\end{align}
and hence the kernel of $j^*$ is in degrees $\geq 1$ given by
\begin{equation}
    \ker(j^*)^{s \geq 1} = \bigoplus_{i=1}^3 \Ffield_2[c_i, y_i, q_i]/(y_i^2)\{c_i z_i\}.
\end{equation}
For the computation of the degree $=0$ part of $\ker(j^*)$, we introduce the notation
\begin{align}
    y & = (y_1, y_2, y_3), \\
    q & = (q_1, q_2, q_3),
\end{align}
and recall the notation
\begin{align}
        z_1 & = (z_1, 0, 0), \\
        z_2 & = (0, z_2, 0), \\
        z_3 & = (0, 0, z_3).
\end{align}
The domain of $j^*$ in degree $=0$ is given by
\begin{equation}
    \bigoplus_{i=1}^3 \Ffield_2[y_i, z_i, q_i]/(y_i^2).
\end{equation}
The classes $y_i$ map to a non-zero class that squares to 0, but the $y_i$ also square to 0. The classes $q_i$ map to polynomial generators, and the classes $z_i$ map to 0. The images of the $y_i$ satisfy the single relation
\begin{equation}
    \sum_{i=1}^3 j^*(y_i) = 0,
\end{equation}
and the images of the $q_i$ satisfy also the single relation
\begin{equation}
    \sum_{i=1}^3 j^*(q_i) = 0.
\end{equation}
Therefore, for a class to be in the kernel of $j^*$, it needs to be a multiple of at least one of the following:
\begin{equation}
    z_1, z_2, z_3, y, q,
\end{equation}
which shows that the kernel in degree $= 0$ is given by
\begin{equation}
    \ker(j^*)^{s=0} = \Ffield_2[z_1, z_2, z_3, y, q]/(z_1z_2, z_1z_3, z_2z_3, y^2).
\end{equation}
The entire kernel is depicted in the following diagram:
\begin{figure}[H]
    \centering
    \includegraphics[width=\textwidth]{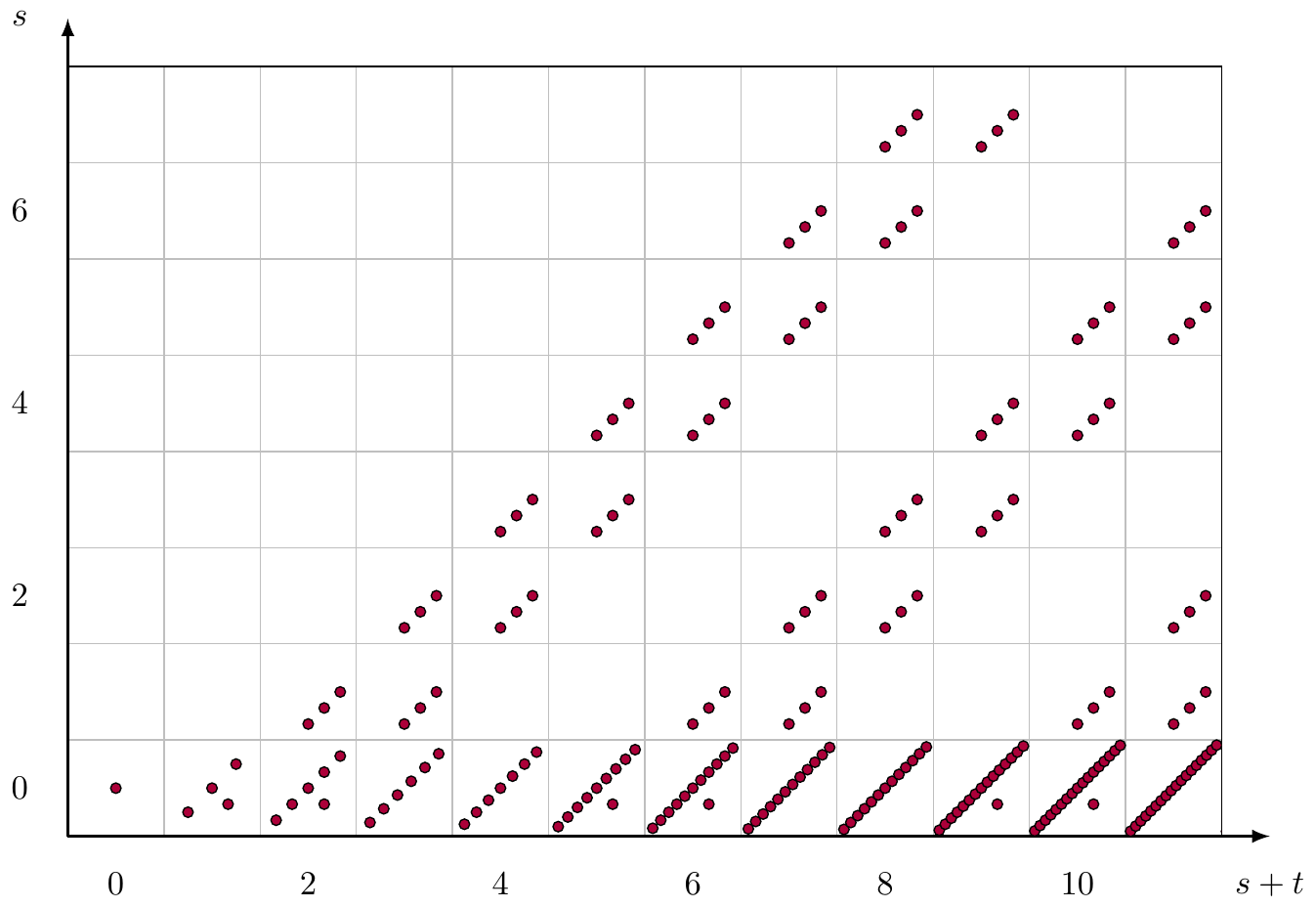}
    \caption{The $B$-summand, which equals $\ker(j^*)$.}
\end{figure}
\subsubsection{The $A$-summand}
The determination of $j^*$ allows us to write down $\Im(j^*)$. For $s \geq 1$ the image splits because the map $j^*$ itself splits:
\begin{align}
    \Im(j^*)^{s \geq 1} = & \Ffield_2[(a + b,0), (x,0), ((\beta_2 x)^2, 0)]/((x,0)^2) \\
    & \oplus \Ffield_2[(b, b), (x,x), ((\beta_2 x)^2, (\beta_2 x)^2)]/((x,x)^2) \\
    & \oplus \Ffield_2[(0, a), (0, x) (0,(\beta_2 x)^2)]/((0,x)^2).
\end{align}
For $s=0$ the image is
\begin{equation}
    \Im(j^*)^{s=0} = \Ffield_2[(x,0), ((\beta_2 x)^2,0), (0, x), (0, (\beta_2 x)^2)].
\end{equation}
The cokernel $\coker(j^*)$ is then given by
\begin{figure}[H]
    \centering
    \includegraphics[width=0.9\textwidth]{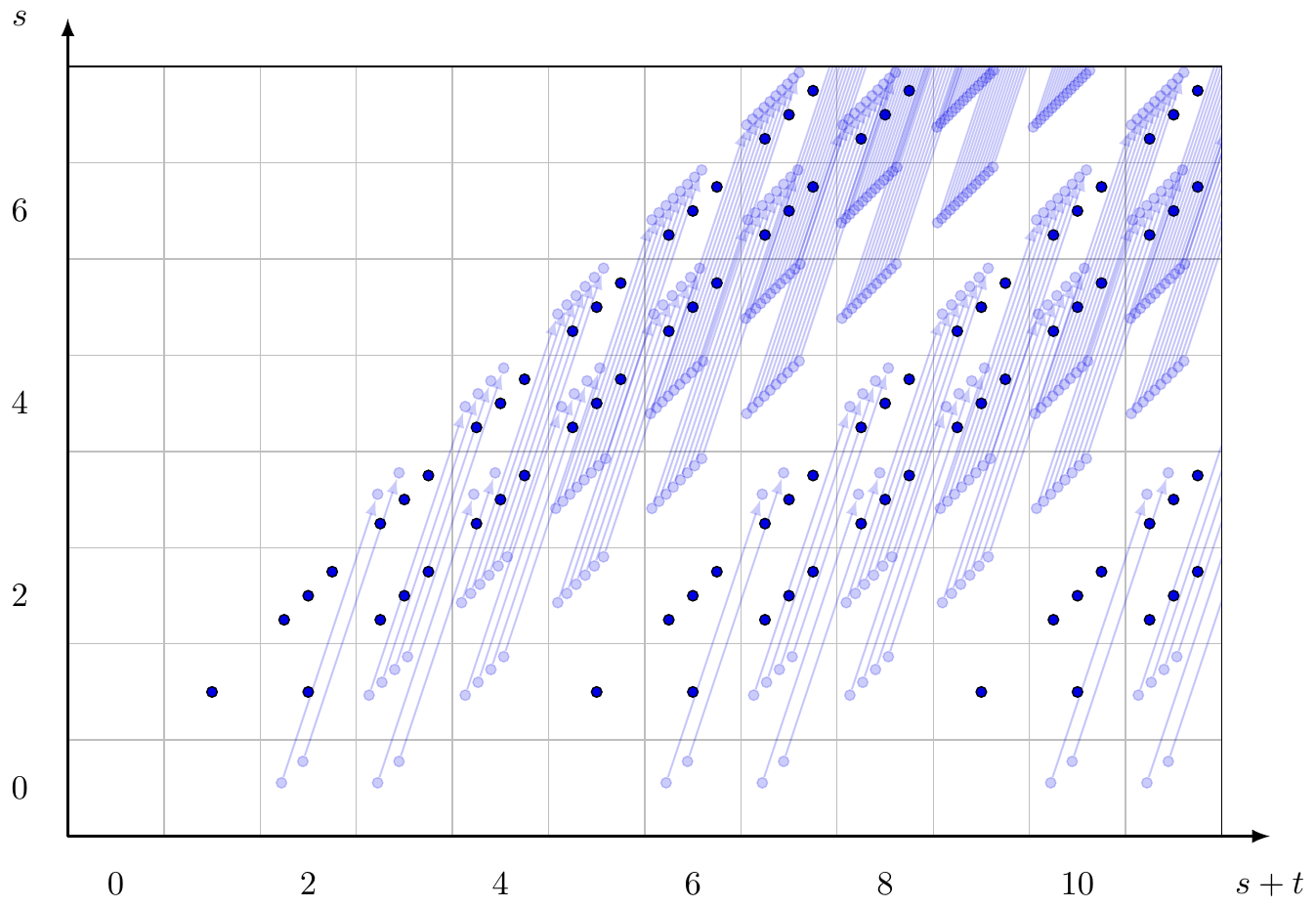}
    \caption{The cokernel of $j^*$, with the differentials induced from the map of spectral sequences.}
\end{figure}
\subsubsection{The $E_2$-page}
Since $E_2 \cong A \oplus B$ as bigraded $\Ffield_2$-modules, the $E_2$-page, without differentials, is given by
\begin{figure}[H]
    \centering
    \includegraphics[width=0.9\textwidth]{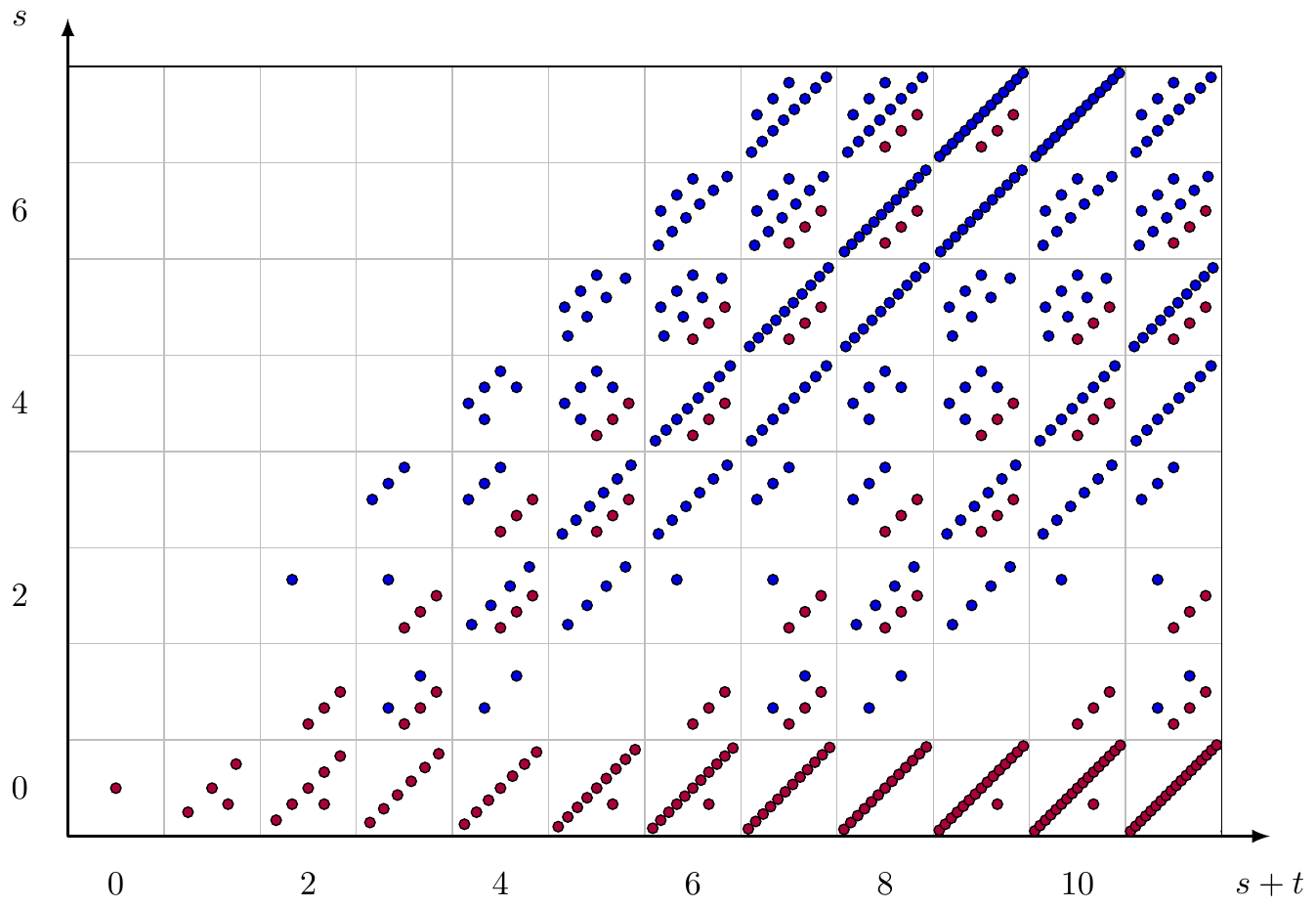}
    \caption{The $E_2$-page of the $\sA$-homotopy limit spectral sequence converging to the group cohomology of $D_8 \ast C_4$, without differentials.}
\end{figure}
\subsection{Differentials}
\label{subsec:d8cpc4diffs}
By \cref{d8cpc4:prop2}, the $\sA$-exponent of $\underline{H\Ffield_2}_{D_8 \ast C_4}$ is $\leq 3$, and hence the $\sA$-homotopy limit spectral sequence will collapse at $E_4$ with a vanishing line of height 3. Furthermore, by \cref{d8cpc4:cor1}, the spectral sequence will even collapse to the $s=0$-line.

We will see that the horizontal vanishing line implies that the $d_2$-differentials supported at the $E_2^{s,t} \cap B$ with $s>0$ and $t \equiv 1, 2\pmod{4}$ map injectively to  $E_2^{s+t, t-1} \cap A$ (cf.\ \cref{d8cpc4:d2d3abovezero}). In particular, all these classes in the $B$-summand die at $E_2$. We are however not able to determine precisely which class in the $B$-summand hits which class in the $A$-summand.

We also show that the vanishing line implies that the $d_2$-differentials supported at $E_2^{0,t}$ with $t \equiv 1,2 \pmod{4}$ map surjectively to $E_2^{2,t-1} \cap A$. Using this statement about surjectivity, we are after the fact able to determine which class gets in $E_2^{0,t}$ supports a differential. We postpone this to \cref{subsec:diffsfroms0}.

The long exact sequence (\ref{absplit:eq3}) respects differentials. This implies that no (red) B-summand class will hit another B-summand class, and that we have the $d_3$ coming from the LHSSS:
\begin{figure}[H]
    \centering
    \includegraphics[width=0.9\textwidth]{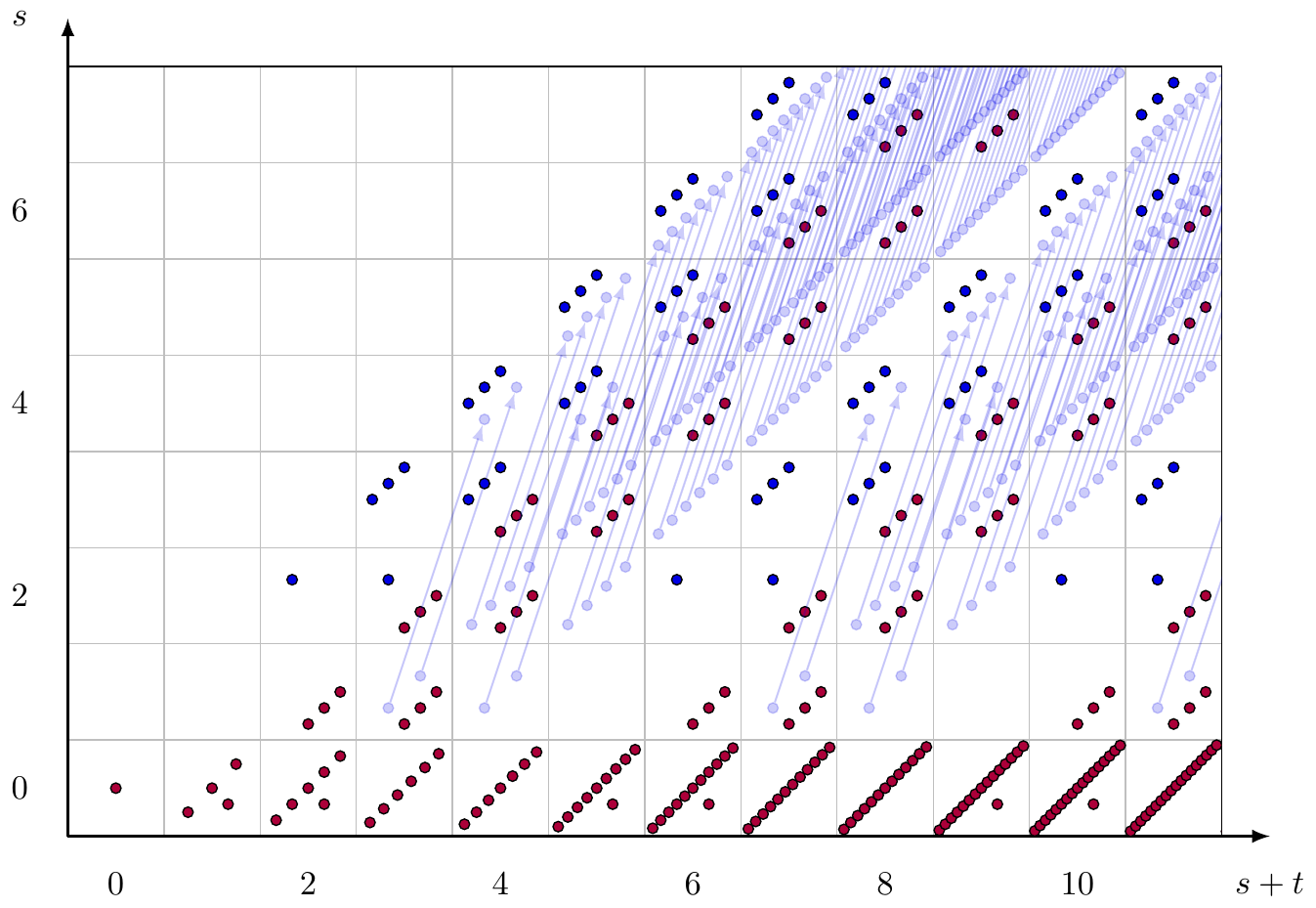}
    \caption{A $d_3$ coming from naturality.}
    \label{d8cpc4:fig2}
\end{figure}
Note however that, a priori, some of the classes involved supporting or getting hit by the $d_3$ in \cref{d8cpc4:fig2} could already have died on the $E_2$-page. We will see however that this is not the case: all the differentials in \cref{d8cpc4:fig2} actually occur. To see this, observe that the the (blue) A-summand classes in the $t \equiv 0, 1 \pmod{4}$ diagonals that do not get hit by the $d_3$ from \cref{d8cpc4:fig2} need to die (since the spectral sequence will collapse to the $s=0$-line), and the only possibility for them to do so is to get hit by (red) B-summand classes (by naturality they cannot get hit by (blue) A-summand classes). 

The classes in the $A$-summand with degree $(s,t)$ with $s+t \geq 5$ and $t \equiv 1 \pmod{4}$ that do not get hit by a $d_3$ from the $A$-summand can only get hit by classes from the $B$-summand with degrees $(s-2, t + 1)$ or $(s-3, t + 2)$. But there are no classes of degree $(s-3, t+2)$ in the $B$-summand satisfying this inequality and congruence, therefore these classes in the $A$-summand will have to get hit by a $d_2$ originating from the $B$-summand.

Now the classes in the $B$-summand with degree $(s,t)$ with $s=2$ and $t \equiv 1 \pmod{4}$ also will die. By naturality, they will not get hit by other classes from the $B$-summand, hence they will hit a class from the $A$-summand, which for degree reasons can only be a class in degree $(4, t-1)$ with a $d_2$. If such a class in the $B$-summand would hit the target of a $d_3$ in the $A$-summand, then the pre-image of this target in the $A$-summand, which lives in degree $(1, t+ 1)$, which would then be a permanent cycle, contrary to the fact that the spectral sequence will collapse to the $s=0$-line. 

For the same reason, classes in the $B$-summand with degree $(s, t)$ with $s=1$ and $t \equiv 2 \pmod{4}$ cannot hit a class of the $A$-summand of degree $(4, t-2)$, for there would remain permanent cycles in the $s=1$-line. The only remaining possibility for these classes in the $B$-summand to die is to hit the three classes of the $A$-summand with degree $(3, t-1)$ that do not get hit by a $d_3$ by a class in the $A$-summand.

Finally, the only remaining possibility for the classes in the $B$-summand with degree $(s,t)$ with $s=1$ and $t \equiv 1 \pmod{4}$ to die is to hit the three classes in the $A$-summand with degree $(3, t-1)$.

This pattern of differentials is expanded in \cref{d8cpc4:d2d3abovezero}.
\begin{figure}[H]
    \centering
    \includegraphics[width=0.9\textwidth]{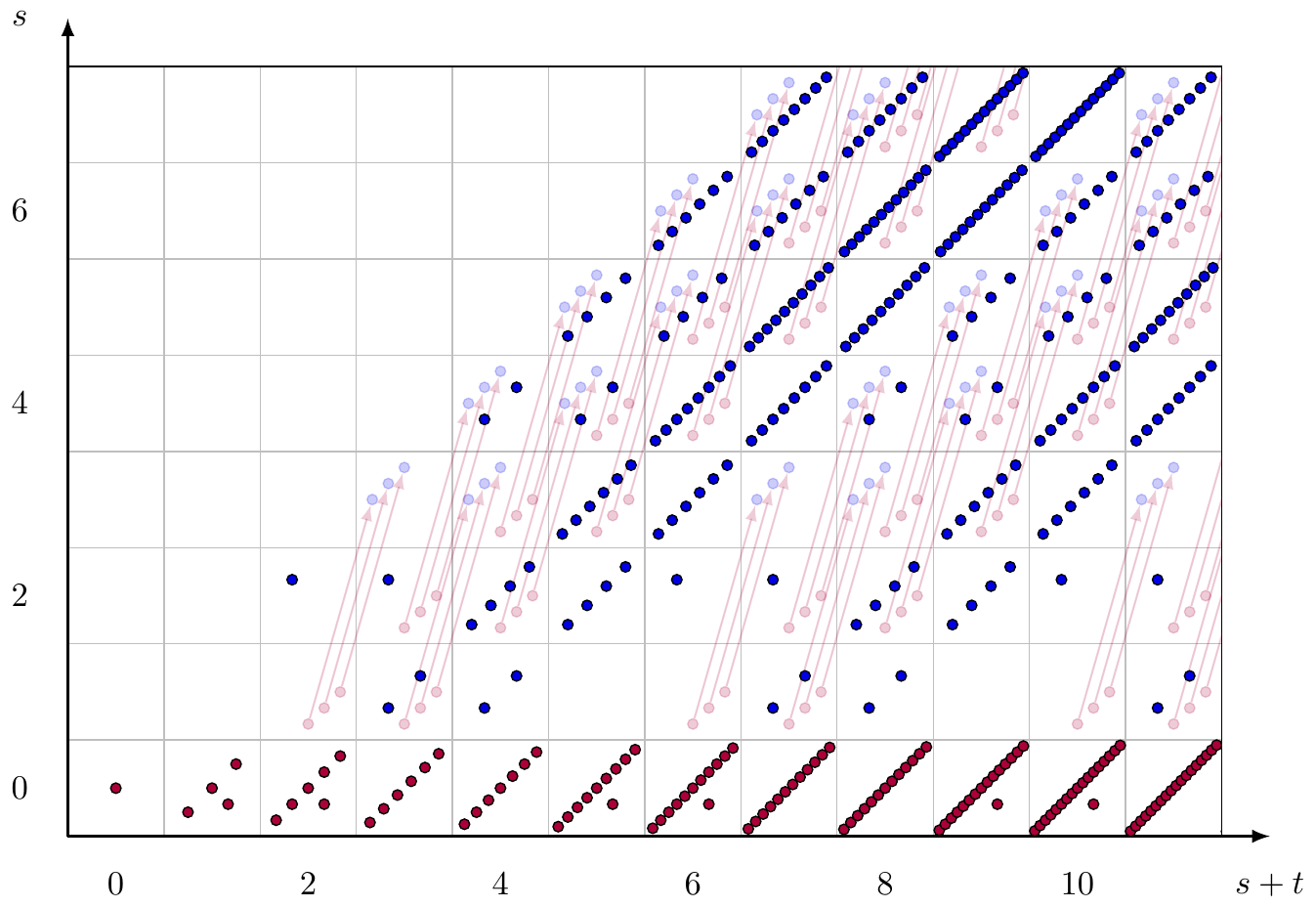}
    \caption{The only possible pattern of $d_2$'s and $d_3$'s originating from the classes in the $B$-summand with $s$-degree $\geq 1$. The pattern only involves non-zero $d_2$-differentials.}
    \label{d8cpc4:d2d3abovezero}
\end{figure}
With these $d_2$'s determined, we also see that all the $d_3$'s from \cref{d8cpc4:fig2} need to happen.

Putting the $d_2$ and $d_3$ differentials that we've deduced so far together gives
\begin{figure}[H]
    \centering
    \includegraphics[width=0.9\textwidth]{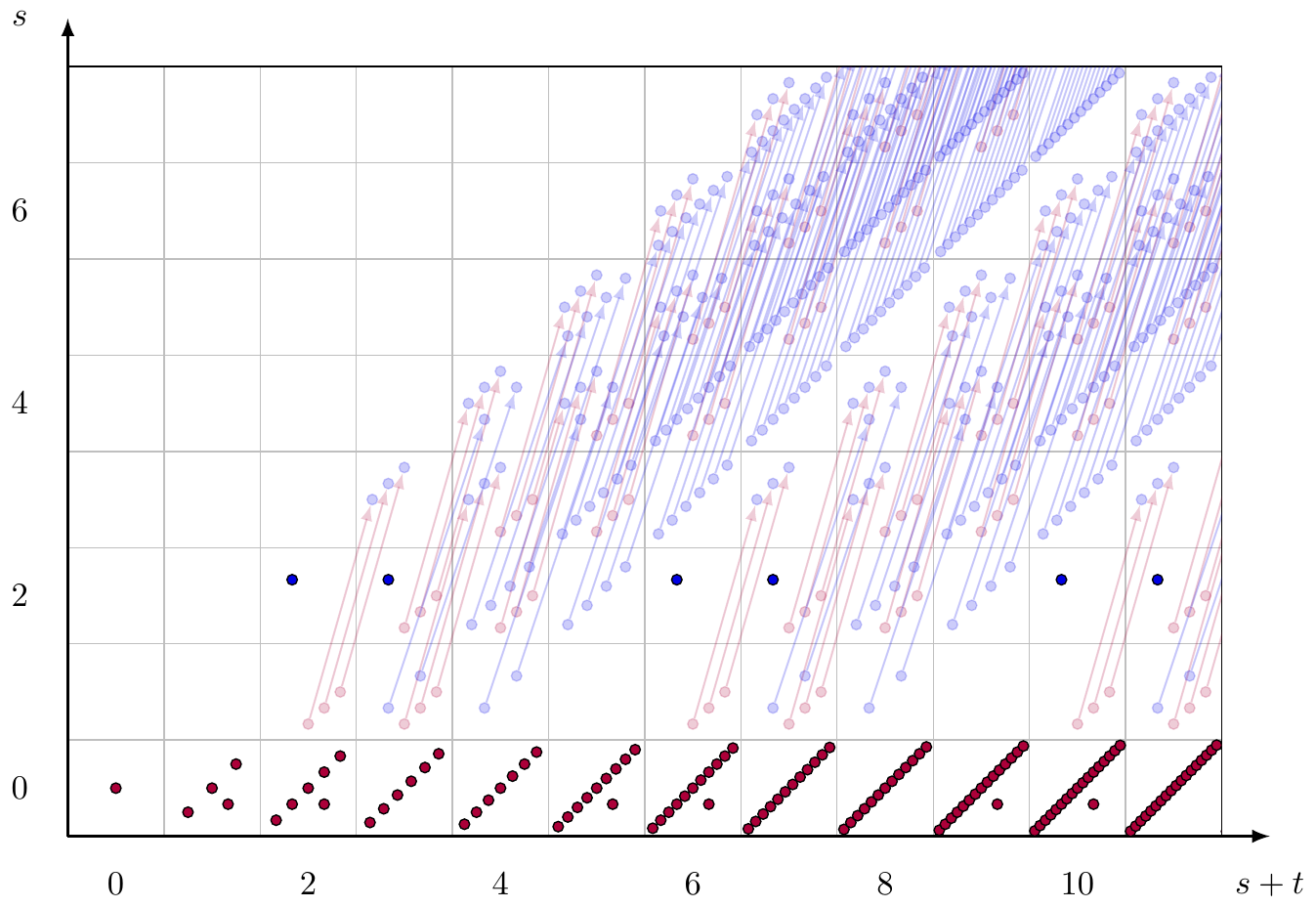}
    \caption{The $d_2$ and $d_3$ differentials deduced so far.}
\end{figure}
Of course the (blue) A-summand classes with $s=2$ also need to die, and they can only be killed by (red) B-summand classes on the $s=0$-line. It is an interesting question which classes in the $s=0$-line precisely support a differential, for this information is required if we want to have any hope of deducing the ring structure on $H^*(BD_8 \ast C_4)$ from this spectral sequence. We will leave this question for the next section, and for now only depict the $\sA$-homotopy limit spectral sequence with all the differentials drawn:
\begin{figure}[H]
    \centering
    \includegraphics[width=0.9\textwidth]{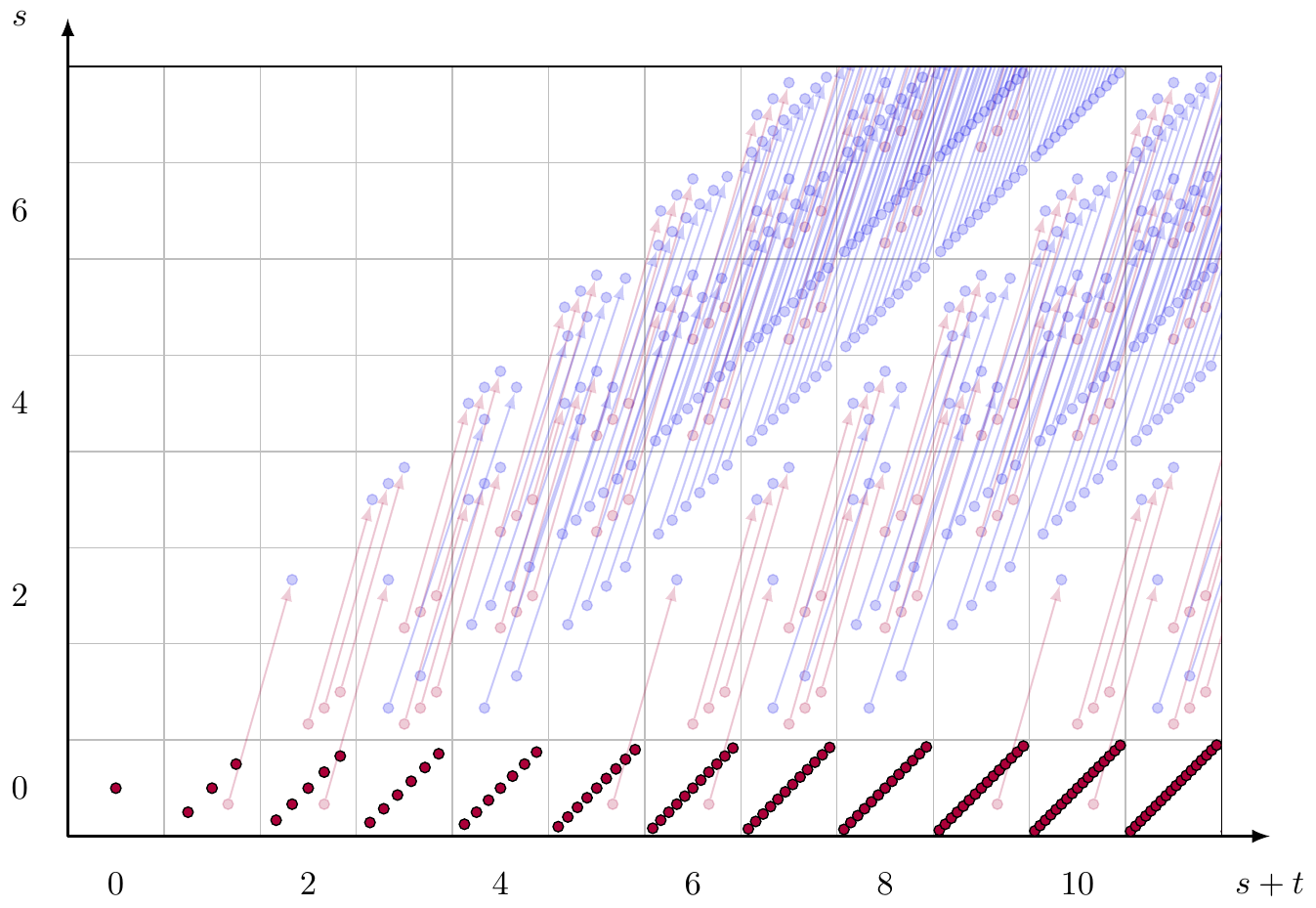}
    \caption{The $\sA$-homotopy limit spectral sequence converging to the $\Ffield_2$-group cohomology of $D_8 \ast C_4$.}
    \label{d8cpc4:completeholimss}
\end{figure}
\fxwarning{Maybe I need to add the differentials originating from s degree 0 t degree 1.}
Since there are non-zero $d_3$-differentials, we have $\exp_{\sA} \underline{H\Ffield_2}_{D_8 \ast C_4} \geq 3$. Combined with \cref{d8cpc4:prop2}, this gives
\begin{Proposition}
    \label{d8cpc4:abexp}
    The $\sA$-exponent satisfies
    \begin{equation}
        \exp_{\sA} \underline{H\Ffield_2}_{D_8 \ast C_4} = 3.
    \end{equation}
\end{Proposition}
\subsection{Poincar\'e series}
To compute the Poincar\'e series, we first recall that the $s=0$-line of the $E_2$-page was given by
\begin{equation}
    E_2^{0,*} = \Ffield_2[z_1, z_2, z_3, y, q]/(z_1 z_2, z_1 z_3, z_2 z_3, y^2).
\end{equation}
Write $R = \Ffield_2[y, q]/(y^2)$ for the subring of this 0-line. Then, as $\Ffield_2[q]$-modules, we have
\begin{align}
    R \cong \Ffield_2[q]\{1, y\},
\end{align}
and therefore $R$ has Poincar\'e series $(1+t)/(1-t^4)$. Furthermore, the $s=0$-line is, as an $R$-module, isomorphic to
\begin{align}
    E_2^{0, *} \cong R \oplus \bigoplus_{i = 1}^3 R[z_i]\{z_i\}.
\end{align}
Therefore $E_2^{0,*}$ has Poincar\'e series
\begin{align}
    \frac{1+t}{1-t^4} + 3 \cdot \frac{t(1+t)}{(1-t)(1-t^4)} & = \frac{(1+t)(1-t) + 3t(1+t)}{(1-t)(1-t^4)} \\
    & = \frac{1+3t + 2t^2}{(1-t)(1-t^4)}.
\end{align}
There is exactly one differential originating from $E_*^{0,k}$ for every $k \equiv 1, 2 \pmod{4}$, and therefore $E_{\infty}^{0,*}$ has Poincar\'e series
\begin{align}
    \frac{1+3t+2t^2}{(1-t)(1-t^4)} - \frac{t + t^2}{1-t^4} & = \frac{1+2t+2t^2 + t^3}{(1-t)(1-t^4)} \\
    & = \frac{(1+t)(1+t+t^2)}{(1-t)^2(1+t)(1+t^2)} \\
    & = \frac{1+t+t^2}{(1-t)^2(1+t^2)},
\end{align}
cf.\ \cite[App.\ C, \#8(16)]{carlson2003cohomology}.
\subsection{Differentials originating from the $s=0$-line}
\label{subsec:diffsfroms0}
We now determine the $d_2$-differentials originating from the $s=0$-line. These classes kill the (blue) A-summand classes with $s=2$, which are given by $[(r^kx^{\epsilon}a,0)]$ with $k \geq 0$ and $\epsilon = 0,1$.

\begin{Proposition}
    The equations
    \begin{align}
        d_2(y) = d_2(z_1) = d_2(z_2) = d_2(z_3) = [(a,0)]
    \end{align}
    completely describe the differential $d_2 \colon E_2^{0,1} \to E_2^{2,0}$.
    \label{d8cpc4:prop3}
\end{Proposition}
\begin{Proof}
We are looking for the classes that hit $[(a,0)] \in E^{2,0}$. If the element $y = (y_1, y_2, y_3)$ were to survive to $E_\infty$, this would imply that $H^*(BD_8 \ast C_4)$ had an element that restricted to $y_i$ on all abelian subgroups. However, the classes that do restrict to $y_i$ for the various $i$ are given by:
\begin{table}[H]
    \centering
    \begin{tabular}{rrr}
        $A_1$: & & $\{\delta_{\sigmabar}, \delta_{\overline{\sigma \rho}}\}$ \\
        $A_2$: & & $\{\delta_{\overline{\sigma\rho}}, \delta_{\overline{\rho \gamma}}\}$ \\
        $A_3$: & & $\{\delta_{\sigmabar}, \delta_{\overline{\rho \gamma}}\}$
    \end{tabular}
\end{table}
Since the intersection of these sets is empty, such a class does not exist, and $y$ will have to die, hence $d_2(y) = [(a,0)]$.

Suppose $(z_1, 0, 0)$ did \emph{not} support a $d_2$. Then $H^*(BD_8 \ast C_4)$ would have a class that restricted to $\delta_{\rho \gamma}$ on $A_1$ and to 0 on $A_2$ and $A_3$. But such a class does not exist. Therefore $d_2((z_1, 0, 0)) = [(a,0)]$, and by symmetry (\cref{prop:d8cpc4aut2} and the fact that there is a unique (blue) A-summand class in $E^{2,0}$) we also get $d_2((0,z_2,0)) = d_2((0,0,z_3)) = [(a,0)]$. This completes the proof.
\end{Proof}
\begin{Proposition}
    The equations
    \begin{equation}
        d_2(z_1y) = d_2(z_2 y) = d_2(z_3 y) = [(xa,0)]
    \end{equation}
    completely describe differential $d_2 \colon E_2^{0,2} \to E^{2,1}$.
    \label{d8cpc4:prop5}
\end{Proposition}
\begin{Proof}
    We know that the class $[(xa,0)]$ dies. Squares cannot support $d_2$-differentials, hence some $z_i y$, and therefore all $z_i y$ by \cref{prop:d8cpc4aut2}, must hit $[(xa,0)]$.
\end{Proof}
\begin{Remark}
    The Leibniz rule together with the proof of \cref{d8cpc4:prop5} shows that we have at $E_2$ the relation
    \begin{equation}
        y[(a,0)] + z_i[(a,0)] = [(xa, 0)]
    \end{equation}
    for all $i$.
\end{Remark}
\begin{Lemma}
    For all $k \geq 0$, $l \geq 3$, $\epsilon = 0, 1$ and $i=1,2,3$ we have
    \begin{equation}
        d_2(q^k z_i^l y^\epsilon) =0
    \end{equation}
    \label{d8cpc4:lem2}
\end{Lemma}
\begin{Proof}
    The class $q$ is a permanent cycle, since there are no differentials originating from $E_*^{0,4}$, by our computations in \cref{subsec:d8cpc4diffs}.
    First observe that
    \begin{equation}
        d_2(z_i^3) = z_i^2 dz_i = 0,
    \end{equation}
    because there is no differential originating from $E_*^{0,3}$, by our computations in \cref{subsec:d8cpc4diffs}. This then implies for $l \geq 3$ that
    \begin{equation}
        d_2(q^k z_i^l) = q^k z_i^{l-1} dz = 0.
    \end{equation}

    Likewise, because there are no differentials originating from $E_2^{0,3}$, by our computations in \cref{subsec:d8cpc4diffs}, we have
    \begin{equation}
        d_2(z_i^2 y) = z_i^2 dy = 0,
    \end{equation}
    and because there are no differentials originating from $E_2^{0,4}$ either, we have
    \begin{equation}
        d_2(z_i^3 y) = yz_i^2 dz_i + z_i^3 dy = 0.
    \end{equation}
    This then shows for $k \geq 0$, $l \geq 3$ that
    \begin{equation}
        d_2(q^k z_i^l y) = 0.
    \end{equation}
\end{Proof}
\begin{Proposition}
    The equations
    \begin{equation}
        d_2(q^k z_i) = d_2(q^k y) = [(r^k a,0)]
    \end{equation}
    completely describe the differentials $d_2 \colon E_2^{0, 4k+1} \to E_2^{2,4k}$.
    \label{d8cpc4:prop4}
\end{Proposition}
\begin{Proof}
    By \cref{d8cpc4:lem2}, the only classes that possibly can support a non-zero differential are the $q^k z_i$ and $q^k y$. We also know that there must be \emph{some} non-zero differential. By \cref{d8cpc4:prop3}, we have
    \begin{align}
        d_2(q^k z_i) & = q^k [(a,0)], \\
        d_2(q^k y) & = q^k [(a,0)].
    \end{align}
    Since at least one of these differentials must be non-zero, they both are (since they hit the same class).
\end{Proof}
\begin{Remark}
    The proof of \cref{d8cpc4:prop4} shows that we have at $E_2$ the relation 
    \begin{equation}
        q^k[(a,0)] = [(r^k a,0)].
    \end{equation}
\end{Remark}
In a completely similar manner we have
\begin{Proposition}
    The equations
    \begin{equation}
        d_2(q^k z_i y) = [(r^k x a,0)]
    \end{equation}
    completely describe the differentials $d_2 \colon E_2^{0,4k+2} \to E_2^{2,4k+1}$.
    \label{d8cpc4:prop6}
\end{Proposition}
\begin{Proof}
    There is at least some class in $E_2^{0,4k+2}$ that kills $[(r^k x a,0)]$. By \\ \cref{d8cpc4:lem2}, the only candidates are the $q^k z_i y$. Since for all $i$ we have
    \begin{equation}
        d_2(q^k z_i y) = q^k[(xa,0)],
    \end{equation}
    we must have that $q^k[(xa,0)] = [(r^k x a, 0)]$, since that is the only possibility for the latter class to die.
\end{Proof}
\begin{Remark}
    The proof of \cref{d8cpc4:prop6} shows that we have at $E_2$ the relation
    \begin{equation}
        q^k[(xa,0)] = [(r^k x a, 0)].
    \end{equation}
\end{Remark}
\subsection{Exponents for various families}
In this subsection we will determine upper bounds on the $\sF$-exponent of \\ $\underline{H\Ffield_2}_{D_8 \ast C_4}$ for two different families, one of which is $\sE_{(2)}$.

The group $D_8 \ast C_4$ has three subgroups isomorphic to $D_8$: $\langle \sigma, \sigma \rho\rangle$, $\langle \sigma \rho, \rho \gamma \rangle$, and $\langle \rho \gamma, \sigma \rangle$, which are all normal. Denote them by $H_1$, $H_2$ and $H_3$. Denote by $\tau_i$ the real 1-dimensional representation of $D_8 \ast C_4$ obtained from pulling back the sign representation along the quotient map $D_8 \ast C_4 \to D_8 \ast C_4/H_i \cong C_2$, and let $e_{\tau_i}$ be the corresponding Euler classes. Let $e$ be the Euler class of $\bigoplus_{i=1}^3\tau_i$. 
\begin{Lemma}
    The $\bigcup_{i=1}^3 \All_{H_i}$-exponent satisfies
    \begin{equation}
        \exp_{\bigcup_{i=1}^3 \All_{H_i}} \underline{H\Ffield_2} \leq 3.
    \end{equation}
    \label{d8cpc4:dihexp}
\end{Lemma}
\begin{Proof}
    The Euler classes $e_{\tau_i}$ are all distinct, because $\res^{D_8 \ast C_4}_{H_i} e_{\tau_j} = 0$ if and only if $i = j$. They are one dimensional cohomology classes, hence detected on the $E_\infty$-page of the $\sA$-homotopy limit spectral sequence considered in the previous sections. Since the Euler classes $e_{\tau_i}$ are of cohomological degree 1, and since the $\sA$-homotopy limit spectral sequence collapsed to the 0-line, they are detected in
    \begin{equation}
        E_\infty^{0,1} = \Ffield\{[z_1 + y], [z_2 + y], [z_3 + y]\}.
        \label{d8cpc4:eq7}
    \end{equation}
    The automorphism from \cref{prop:d8cpc4aut2} that cyclically permutes the generators $\{\sigma, \sigma\rho, \rho \gamma\}$ also cyclically permutes the groups $\{H_1, H_2, H_3\}$, since it cyclically permutes the generators of these dihedral groups. This automorphism also cyclically permutes the basis elements of (\ref{d8cpc4:eq7}).

    We now consider three cases. The first case is that $e_{\tau_1}$ is detected by $[z_i + y]$ for some $i$. Then $e_{\tau_2}$ is detected by $[z_{i+1} + y]$ and $e_{\tau_3}$ is detected $[z_{i+2} + y]$, where the indices of the $z$'s are understood to be classes modulo 3. Because
    \begin{align}
        \prod_{i=1}^3 [z_i + y] & = [z_1y + z_2y][z_3 + y] = 0,
    \end{align}
    we have that $e$ is 0 up to higher filtration, of which there is none. Hence in this case $e = 0$. 
    
    The second case is that $e_{\tau_1}$ is detected by some $[z_i + y] + [z_j + y]$ with $i \neq j$. Again the automorphism from \cref{prop:d8cpc4aut2} then shows that $e_{\tau_2}$ is detected by $[z_{i+1} + y] + [z_{j+1} + y]$ and $e_{\tau_3}$ is detected by $[z_{i+2} + y] + [z_{j+2} + y]$, where the indices of the $z$'s are understood to be classes modulo 3. Because
    \begin{align}
        \prod_{i=1}^3 ([z_i + y] + [z_{i+1} + y]) & = \prod_{i=1}^3 [z_i + z_{i+1}] \\
        & = [z_2^2][z_3 + z_1] \\
        & = 0,
    \end{align}
    we have again that $e$ is 0 up to higher filtration, of which there is none. Hence also in this case $e=0$.

    The third and final case is that $e_{\tau_1}$ is detected by $\sum_{i=1}^3 [z_i + y]$. But the latter is fixed by the automorphism of \cref{prop:d8cpc4aut2}, hence this implies that $e_{\tau_1}$ is fixed up to higher filtration, of which there is none. This contradicts the fact that the automorphism maps $e_{\tau_1} \mapsto e_{\tau_2}$ and that the $e_{\tau_i}$ are distinct, so in fact this case does not occur.

    We conclude that in any case, $e = 0$, and by \cref{cor:hfpeulerclass}, the result follows.
\end{Proof}
\begin{Corollary}
    The $\sE_{(2)}$-exponent satisfies $\exp_{\sE_{(2)}} \underline{H\Ffield_2}_{D_8 \ast C_4} \leq 5$.
    \label{d8cpc4:elexpupperbound}
\end{Corollary}
\begin{Proof}
    This follows from \cref{d8cpc4:prop2}, \cref{d8cpc4:dihexp}, $\sA \cap \left(\bigcup_{i=1}^3 \All_{H_i}\right) = \sE_{(2)}$, \cref{lem:expintersect}, and $3 + 3 - 1 = 5$.
\end{Proof}
We can slightly improve the bound in \cref{d8cpc4:elexpupperbound} as follows.
\begin{Proposition}
    The $\sE_{(2)}$-exponent satisfies $\exp_{\sE_{(2)}} \underline{H\Ffield_2}_{D_8 \ast C_4} \leq 4$.
    \label{d8cpc4:elexpupperbound2}
\end{Proposition}
\begin{Proof}
    Let
    \begin{align}
        \sigma \mapsto \begin{pmatrix} 1 & 0 & 0 & 0 \\ 0 & 1 & 0 & 0 \\ 0 & 0 & -1 & 0 \\ 0 & 0 & 0 & -1 \end{pmatrix},
        \rho \mapsto \begin{pmatrix} 0 & 0 & -1 & 0 \\ 0 & 0 & 0 & -1 \\ 1 & 0 & 0 & 0 \\ 0 & 1 & 0 & 0 \end{pmatrix},
        \gamma \mapsto \begin{pmatrix} 0 & -1 & 0 & 0 \\ 1 & 0 & 0 & 0 \\ 0 & 0 & 0 & -1 \\ 0 & 0 & 1 & 0  \end{pmatrix}
    \end{align}
    be the underlying 4-dimensional real representation from \cite[Lem.\ 13.5]{totaro14}. The following table gives, up to inverses, all elements of $D_8 \ast C_4$, the effect of this representation on a line $[x_1 : y_1 : x_2 : y_2]$, and the set of lines each element fixes.
    \begin{table}[H]
        \centering
        \begin{tabular}{r|r|r}
            $g$ & $g[x_1 : y_1 : x_2 : y_2]$ & fixes the lines in the set(s) \\
            \hline
            $\rho$ & $[-x_2 : -y_2 : x_1 : y_1 ]$ & $\varnothing$ \\
            $\rho^2$ & $[-x_1 : -y_1 : -x_2 : -y_2]$ & $\varnothing$ \\
            $\gamma$ & $[-y_2 : x_1 : -y_2 : x_2 ]$ & $\varnothing$ \\
            $\sigma \rho \gamma$ & $[y_2 : -x_2 : y_1 : -x_1 ]$ & $\varnothing$ \\
            $\sigma \gamma$ & $[-y_1 : x_1 : y_2 : -x_2]$ & $\varnothing$ \\
            $\sigma$ & $[x_1 : y_1 : -x_2 : -y_2]$ & $\{x_1 = y_1 = 0\}$, $\{x_2 = y_2 = 0\}$ \\
            $\sigma \rho^2$ & $[-x_1 : -y_1 : x_2 : y_2]$ & $\{x_1 = y_1 = 0\}$, $\{x_2 = y_2 = 0\}$ \\
            $\sigma \rho$ & $[-x_2 : -y_2 : -x_1 : -y_1]$ & $\{x_1 = x_2, y_1=y_2\}$, $\{x_1 = -x_2, y_1 = -y_2\}$\\
            $\sigma \rho^3$ & $[x_2 : y_2 : x_1 : y_1]$ & $\{x_1 = x_2, y_1 = y_2\}$, $\{x_1 = -x_2, y_1 = -y_2\}$ \\
            $\rho \gamma$ & $[y_2 : -x_2 : -y_1 : x_1]$ &  $\{x_1 = y_2, y_1 = -x_2\}$, $\{x_1 = -y_2, y_1 = x_2\}$\\
            $\rho^3 \gamma$ & $[-y_2 : x_1 : y_2 : -x_2 ]$ & $\{x_1 = y_2, y_1 = -x_2\}$, $\{x_1 = -y_2, y_1 = x_2\}$ \\
        \end{tabular}
        \caption{The elements of $D_8 \ast C_4$, their effect on a line in the above representation, and the lines fixed by each element.}
        \label{d8cpc4:tabfixlines}
    \end{table}
    Write
    \begin{align}
        L_1 & = \{x_1 = y_1 = 0\} \cup \{x_2 = y_2 = 0\}, \\
        L_2 & = \{x_1 = -x_2, y_2 = -y_1\} \cup  \{x_1 = x_2, y_1 = y_2 \}, \\
        L_3 & = \{x_1 = y_2 = y_1 = -x_2\} \cup  \{x_1 = -y_2 = y_1 = x_2 \}.
    \end{align}
    Then the $L_i$, together with the set of all lines, are the sets of lines that occur as the lines fixed by an element of $D_8 \ast C_4$. We immediately see from \cref{d8cpc4:tabfixlines} that the lines in $L_1$ are fixed by $\langle \sigma, \sigma \rho^2 \rangle$, the lines in $L_2$ are fixed by $\langle \sigma \rho, \sigma \rho^3 \rangle$, and the lines in $L_3$ are fixed by $\langle \rho \gamma, \rho \gamma^3 \rangle$. Moreover, the $L_i$ are pairwise disjoint. To see this, note that a line in $L_1 \cap L_2$ would be fixed by both $\sigma$ and $\sigma \rho$, hence by their product $\sigma \cdot \sigma \rho = \rho$. But $\rho$ does not fix any line, by \cref{d8cpc4:tabfixlines}. Similarly, a line in $L_2 \cap L_3$ would be fixed by both $\sigma \rho$ and $\rho^3 \gamma$, but $\sigma \rho \cdot \rho^3 \gamma = \sigma \gamma$ does not fix any line. Lastly, a line in $L_3 \cap L_1$ would be fixed by $\sigma \cdot \rho \gamma$, which does not fix any line. Therefore, no line is fixed by a strict supergroup of $\langle \sigma, \sigma \rho^2 \rangle$, $\langle \sigma \rho, \sigma \rho^3 \rangle$, $\langle \rho \gamma, \rho^3 \gamma \rangle$, which are therefore the maximal isotropy groups of a line. All these three groups are elementary abelian, so the result follows from \cref{projbund:prop1}.
\end{Proof}
We can use the previous proposition to determine the $\sE_{(2)}$-exponent. 
\begin{Proposition}
    The $\sE_{(2)}$-exponent satisfies
    \begin{equation}
       \exp_{\sE_{(2)}} \underline{H\Ffield_2} = 4.
    \end{equation}
    \label{d8cpc4:e2exp}
\end{Proposition}
\begin{Proof}
    The upper bound is \cref{d8cpc4:elexpupperbound2}, the lower bound follows from the fact that $D_8 \ast C_4$ has a subgroup isomorphic to the quaternion group $Q_8$ of order 8, generated by $\langle \rho, \sigma \gamma \rangle$. Now \cref{prop:q2ne2exp} and \cref{lem:expsubgp} imply the desired lower bound.
\end{Proof}

\section{$C_4 \rtimes C_4$}
\label{sec:c4sdc4}
Let $r$ and $s$ generate two copies of $C_4$: $\langle r \, | \, r^4 = e \rangle,\, \langle s \, | \, s^4 = e \rangle \cong C_4$. Let $\langle s \rangle$ act on $\langle r \rangle$ by $s \cdot r = r^{-1}$. Then $C_4 \rtimes C_4$ is defined to be the semi-direct product $\langle r \rangle \rtimes \langle s \rangle$. A presentation is $C_4 \rtimes C_4 = \langle r,s \, | \, r^4 = s^4 = e, srs^{-1} = r^3 \rangle$. We will compute the $\sE_{(2)}$-homotopy limit spectral sequence converging to $H^*(BC_4 \rtimes C_4;\Ffield_2)$ which is isomorphic to (\cite[App.\ C, \#10(16)]{carlson2003cohomology})
\begin{equation}
    \Ffield_2[z, y, x, w]/(z^2 + y^2, zy).   
\end{equation}
with degrees $|z| = |y| = 1$, $|x| = |w| = 2$.
\subsection{An upper bound on the exponent}
In this section we show the following upper bound on the $\sE_{(2)}$-exponent:
\begin{equation}
    \exp_{\sE_{(2)}} \underline{H\Ffield_2}_{C_4 \rtimes C_4} \leq 8.
\end{equation}
\begin{Lemma}
    For $\sP$ the family of proper subgroups of $C_4 \rtimes C_4$, we have the upper bound:
    \begin{equation}
        \exp_{\sP} \underline{H\Ffield_2}_{C_4\rtimes C_4} \leq 4.
    \end{equation}
\end{Lemma}
\begin{Proof}
    Consider the linear action of $C_4 \rtimes C_4$ on $\Reals^4$ from \cite[Lem.\ 13.7]{totaro14} given by
\begin{equation}
    s \mapsto
    \begin{pmatrix}
        1 & 0 & 0 & 0 \\
        0 & -1 & 0 & 0 \\
        0 & 0 & 0 & -1 \\
        0 & 0 & 1 & 0
    \end{pmatrix}, \qquad
    r  \mapsto
    \begin{pmatrix}
        0 & -1 & 0 & 0 \\
        1 & 0 & 0 & 0 \\
        0 & 0 & 1 & 0 \\
        0 & 0 & 0 & 1
    \end{pmatrix}.
\end{equation}
We compute the lines fixed by this action. A group element $g$ fixes a line $L$ if and only if $g^{-1}$ fixes $L$. Furthermore, the central element $r^2 s^2$ fixes all lines. Hence the table below suffices to determine the isotropy groups of all lines:
\begin{table}[H]
    \centering
    \begin{tabular}{r|r|r}
        $g$ & $g[x_1 : y_1 : x_2 : y_2]$ & fixes the lines in the set(s) \\
        \hline
        $r$ & $[-y_1 : x_1 : x_2 : y_2]$ & $\{x_1 = y_1 = 0\}$ \\
        $r^2$ & $[-x_1 : -y_1 : x_2 : y_2]$ & $\{x_1 = y_1 = 0\}$, $\{x_2 = y_2 = 0\}$ \\
        $s$ & $[x_1 : -y_1 : -y_2 : x_2]$ & $\{x_2 = y_2 = x_1 = 0\}$, $\{x_2 = y_2 = y_1 = 0\}$ \\
        $rs$ & $[y_1 : x_1 : -y_2 : x_2]$ & $\{x_2 = y_2 = x_1 - y_1=0\}$, $\{x_2 = y_2 = x_1 + y_1=0\}$ \\
        $r^2 s$ & $[-x_1 : y_1 : -y_2 : x_2]$ & $\{x_2 = y_2 = x_1 = 0\}$, $\{x_2 = y_2 = y_1 = 0\}$ \\
        $r^3 s$ & $[-y_1 : -x_1 : -y_2 : x_2]$ & $\{x_2 = y_2 = x_1 - y_1 = 0\}$, $\{x_2 = y_2 = x_1 + y_1 = 0\}$ \\
        $s^2$ & $[x_1 : y_1 : -x_2 : -y_2]$ & $\{x_1 = y_1 = 0\}$, $\{x_2 = y_2 = 0\}$ \\
        $rs^2$ & $[-y_1 : x_1 : -x_2 : -y_2]$ & $\{x_1 = y_1 = 0\}$
    \end{tabular}
\end{table}
From this table we read of that the maximal isotropy groups are:
\begin{enumerate}
    \item The group $\langle r , s^2 \rangle \cong C_4 \times C_2$, fixing the lines in the set $\{x_1 = y_1 = 0\}$.
    \item The group $\langle r^2, s \rangle \cong C_2 \times C_4$ fixing the lines in the sets $\{x_2 = y_2 = x_1 = 0\}$ and $\{x_2 = y_2 = y_1 = 0\}$.
    \item The group $\langle r^2, rs \rangle \cong C_2 \times C_4$, fixing the lines in the sets $\{x_2 = y_2 = x_1 - y_1 = 0\}$ and $\{x_2 = y_2 = x_1 + y_1 = 0\}$.
\end{enumerate}
These are precisely the maximal proper subgroups of $C_4 \rtimes C_4$, hence the minimal family containing all the isotropy of $\Proj(\Reals^4)$ is $\sP$. An application of \cref{projbund:prop1} yields the desired result.
\end{Proof}
\begin{Proposition}
    \label{prop:c4sdc4expupbound}
    For $\sE_{(2)}$ the family elementary abelian subgroups, we have
    \begin{equation}
        \exp_{\sE_{(2)}} \underline{H\Ffield_2}_{C_4 \rtimes C_4} \leq 8.
    \end{equation}
\end{Proposition}
\begin{Proof}
    We apply \cref{lem:twofamexp} with $\sF = \sP$ and $\sG = \sE_{(2)}$. All maximal subgroups $K$ in $\sP$ are isomorphic to $C_4 \times C_2$, and for those we know (\cref{prop:expelab}) that $\exp_{\sE_{(2)}}\Res^G_K \underline{H\Ffield_2} = 2$. Hence, in the notation of \cref{lem:twofamexp}, $m := \max_{K \in \sP} m_K = 2$. Combined with $n = \exp_{\sP} \underline{H\Ffield_2} = 4$ this yields by \cref{lem:twofamexp}
    \begin{equation}
        \exp_{\sE_{(2)}} \underline{H\Ffield_2} \leq 4 \cdot 2 = 8.
    \end{equation}
\end{Proof}

\subsection{The spectral sequence}
We will compute the $\sE_{(2)}$-homotopy limit spectral sequence of $\underline{H\Ffield_2}$.

The elements of order precisely 2  in $G$ are $r^2$, $s^2$, $r^2s^2$.  Hence all elementary abelian 2-subgroups are contained in $\langle r^2, s^2 \rangle =: V_4 \cong C_2 \times C_2$, a copy of the Klein four-group, which is also the center of $G$. Hence, by \cref{prop:compserre}, the $\sE_{(2)}$-homotopy limit spectral sequence reduces to the LHSSS associated to the central extension
\begin{equation}
	 V_4  \to G \to \langle \bar{r}, \bar{s} \rangle,
\end{equation}
We observe that $G/V_4$ is also isomorphic to the Klein four-group, and call this group $V_4'$. Hence the spectral sequence takes the form.
\begin{equation}
	E_2^{s,t} = H^s(BV_4'; H^t(BV_4; \Ffield_2)) \Rightarrow H^{s+t}(BG; \Ffield_2).
\end{equation}
Since the action of $V_4'$ on the center is trivial the local coefficient system is also trivial. Hence the $E_2$-page is isomorphic to
\begin{equation}
	\Ffield_2[\delta_{\bar{r}},\delta_{\bar{s}}] \otimes \Ffield_2[\delta_{r^2},\delta_{s^2}],
\end{equation}
with $(s,t)$-degrees $|\delta_{\rbar}| = |\delta_{\sbar}| = (1,0)$, and $|\delta_{r^2}| = |\delta_{s^2}| = (0,1)$.
\begin{figure}[H]
	\centering
	\includegraphics[width=0.4\textwidth]{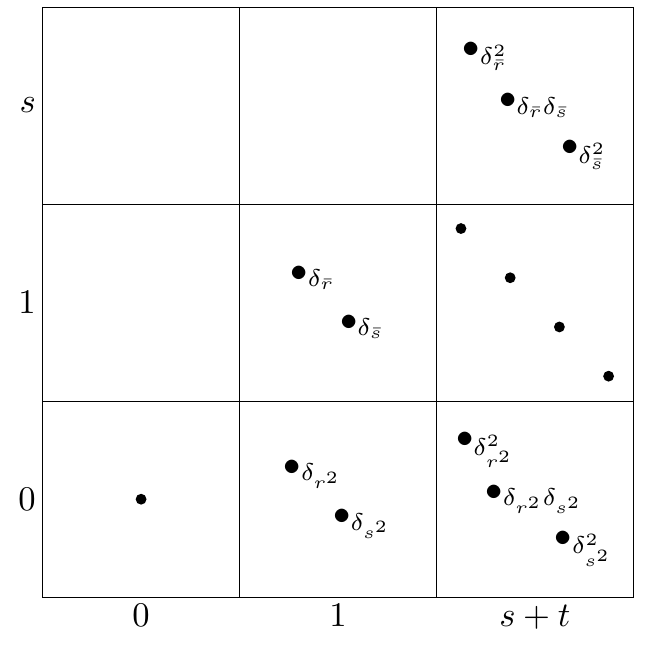}
	\caption{$E_2$-page for $C_4 \rtimes C_4$, with no differentials drawn.}
\end{figure}
To compute $d_2$ we consider the map of spectral sequences given by restriction to various subgroups.

First consider the restriction $\langle r \rangle \cong C_4 \hookrightarrow G$. This has one elementary 2-group, $\langle r^2 \rangle \cong C_2$, and we get
\begin{equation}
	\begin{tikzcd}
	1\to	\underbrace{C_2}_{=\langle r^2 \rangle} \arrow{d} \arrow{r} & \underbrace{C_4}_{=\langle r \rangle} \arrow{d} \arrow{r} & \underbrace{C_2}_{=\langle \bar{r} \rangle} \arrow{d} \to 1 \\
	1\to	\underbrace{V_4}_{=\langle r^2,s^2 \rangle} \arrow{r} & G \arrow{r} & \underbrace{C_2 \times C_2}_{= \langle\bar{r},\bar{s}\rangle} \to 1
	\end{tikzcd}
\end{equation}
We already saw that the LHSSS has a $d_2$ given by $d_2(\delta_{r^2}) = \delta_{\bar{r}}^2$ (\cref{subsec:cycl2gps}), where we use the same notation for dual elements in the cohomology of the subgroups. This gives the following two conclusions for $d_2$ of the LHSSS $G$:
\begin{enumerate}[1.]
	\item $d_2 \delta_{r^2}$ has a non-zero $\delta_{\bar{r}}^2$-coefficient;
	\item $d_2 \delta_{s^2}$ has a zero $\delta_{\bar{r}}^2$-coefficient.
\end{enumerate}

Second we consider the restriction $\langle s \rangle \cong C_4 \hookrightarrow G$. This gives a similar diagram
\begin{equation}
	\begin{tikzcd}
1\to		\underbrace{C_2}_{=\langle s^2 \rangle} \arrow{d} \arrow{r} & \underbrace{C_4}_{=\langle s \rangle} \arrow{d} \arrow{r} & \underbrace{C_2}_{=\langle \bar{s} \rangle} \arrow{d} \to 1 \\
1\to		\underbrace{V_4}_{=\langle r^2,s^2 \rangle} \arrow{r} & G \arrow{r} & \underbrace{C_2 \times C_2}_{= \langle\bar{r},\bar{s}\rangle} \to 1
	\end{tikzcd}
\end{equation}
Again by dimension and degree reasons we infer that
\begin{enumerate}
		\setcounter{enumi}{2}
	\item $d_2\delta_{s^2}$ has a non-zero $\delta_{\bar{s}}^2$-coefficient;
	\item $d_2 \delta_{r^2}$ has a zero $\delta_{\bar{s}}^2$-coefficient.
\end{enumerate}
Finally, we consider the restriction to $\langle rs \rangle \cong C_4 \hookrightarrow G$. Then $(rs)^2 = s^2$, and we get the diagram
\begin{equation}
	\begin{tikzcd}
1 \to		\underbrace{C_2}_{=\langle s^2 \rangle} \arrow{d} \arrow{r} & \underbrace{C_4}_{=\langle rs \rangle} \arrow{d} \arrow{r} & \underbrace{C_2}_{= \langle \overline{rs} \rangle} \arrow{d} \to 1 \\
1 \to		\underbrace{V_4}_{=\langle r^2,s^2 \rangle} \arrow{r} & G \arrow{r} & \underbrace{C_2 \times C_2}_{ = \langle\bar{r},\bar{s}\rangle} \to 1
	\end{tikzcd}
\end{equation}
Now the inclusion $C_2 \hookrightarrow C_2 \times C_2\colon \overline{rs} \mapsto \bar{r} \bar{s}$ gives the restriction map in cohomology, which is given by
\begin{align}
	H^*(BC_2 \times C_2;\Ffield_2) & \to H^*(BC_2;\Ffield_2) \\
	\delta_{\bar{r}} & \mapsto \delta_{\overline{rs}} \\
	\delta_{\bar{s}} & \mapsto \delta_{\overline{rs}}.
\end{align}
Again, the LHSSS supports a differential at $\delta_{s^2}$ with image $\delta_{\overline{rs}}^2$, which implies that
\begin{enumerate}
		\setcounter{enumi}{4}
	\item $d_2 \delta_{s^2}$ has an \emph{odd} number of non-zero $\{\delta_{\bar{r}}^2, \delta_{\bar{r}}\delta_{\bar{s}}, \delta_{\bar{s}}^2\}$-coefficients.
	\item $d_2 \delta_{r^2}$ has an \emph{even} number of non-zero $\{\delta_{\bar{r}}^2, \delta_{\bar{r}}\delta_{\bar{s}}, \delta_{\bar{s}}^2\}$-coefficients.
\end{enumerate}
Taking all these things together we obtain
\begin{align}
    d_2 \delta_{r^2} & = \delta_{\bar{r}}^2 + \delta_{\rbar}\delta_{\sbar}, \\
	d_2 \delta_{s^2} & = \delta_{\bar{s}}^2.
    \label{c4sdpc4:diffs1}
\end{align}
\begin{figure}[H]
	\centering
	\includegraphics[width=0.4\textwidth]{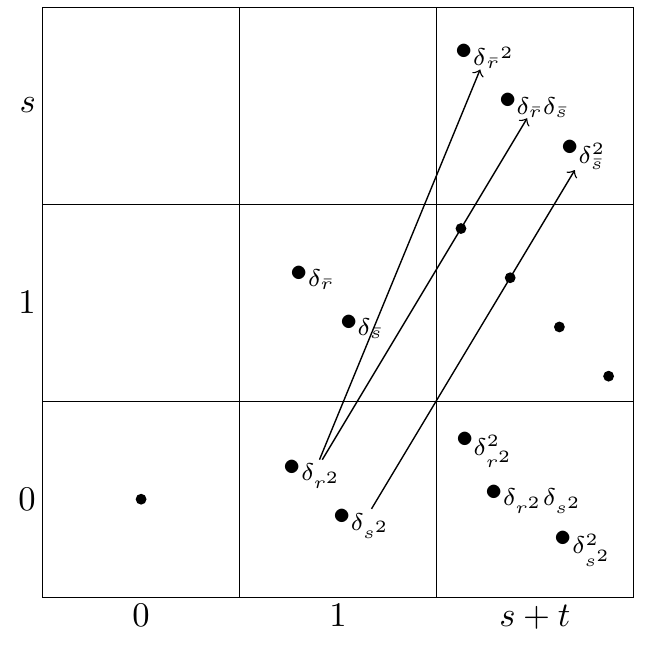}
	\caption{$E_2$-page for $C_4 \rtimes C_4$, with the generating differentials drawn.}
\end{figure}
Therefore we get that the $E_3$-page is
\begin{equation}
	E_3 \cong \Ffield_2[\delta_{\bar{r}}, \delta_{\bar{s}}]/(\delta_{\bar{r}}^2, \delta_{\bar{r}}\delta_{\bar{s}} + \delta_{\bar{r}}^2) \otimes_{\Ffield_2} \Ffield_2[[\delta_{r^2}^2],[\delta_{s^2}^2]].
\end{equation}
\begin{figure}[H]
	\centering
	\includegraphics[width=0.9\textwidth]{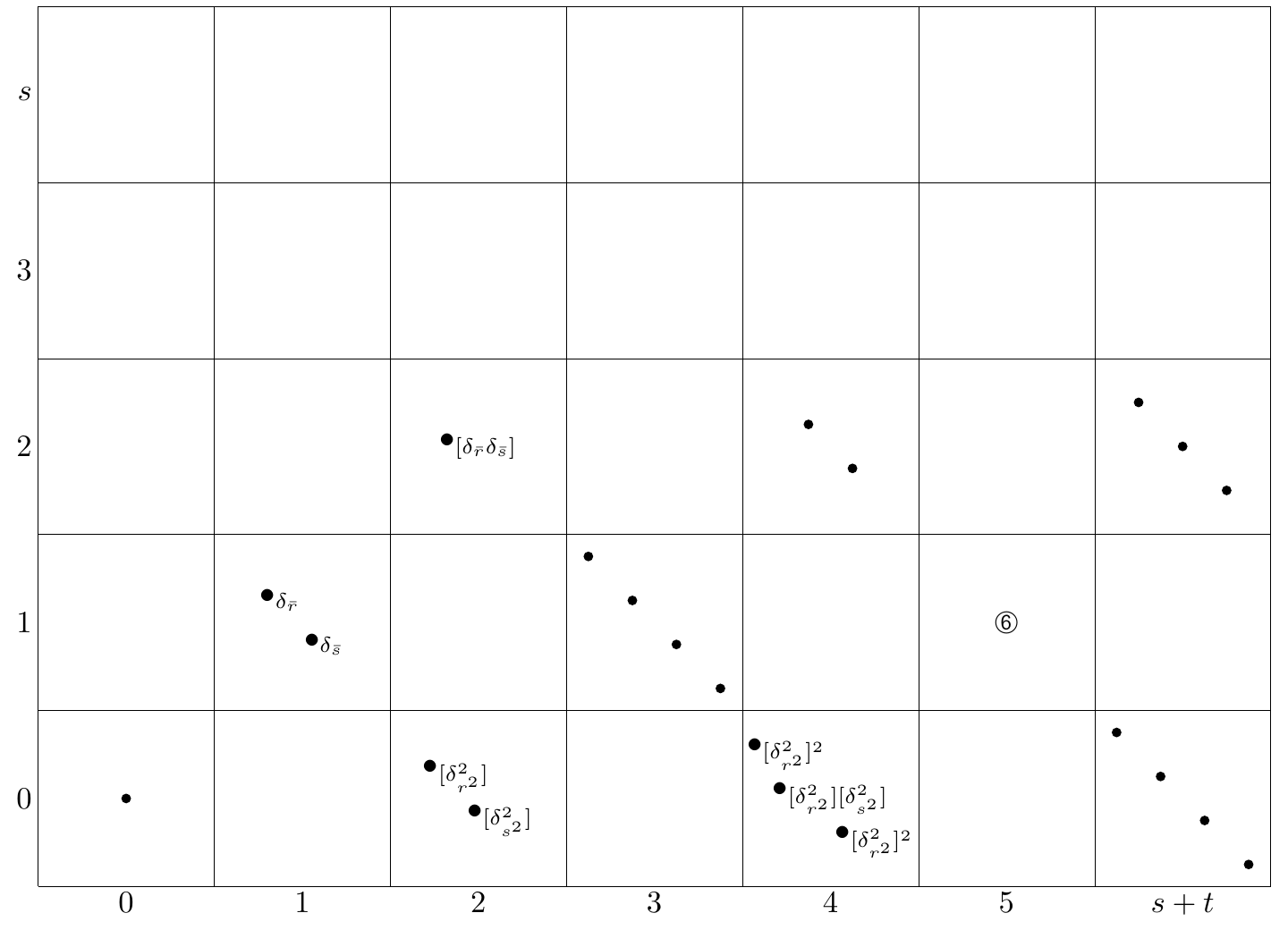}
	\caption{$E_3$-page for $C_4 \rtimes C_4$. There is a vanishing line of height 3.}
\end{figure}
This page has a vanishing line of height 3, and therefore
\begin{align}
    E_\infty & = E_3 \\
	& = \Ffield_2[\delta_{\bar{r}}, \delta_{\bar{s}}]/(\delta_{\bar{r}}^2, \delta_{\bar{r}}\delta_{\bar{s}} + \delta_{\bar{r}}^2) \otimes_{\Ffield_2} \Ffield_2[[\delta_{r^2}^2],[\delta_{s^2}^2]].
\end{align}
\subsection{Poincar\'e series}
The resulting Poincar\'e series of $H^*(BC_4 \rtimes C_4;\Ffield_2)$ is
\begin{equation}
    \frac{1+2t+t^2}{(1-t^2)^2} = \frac{1}{(1-t)^2},
\end{equation}
(cf.\ \cite[App.\ B and App.\ G, \#10(16)]{carlson2003cohomology}).
\subsection{The exponent}
Using the computation of the $\sE_{(2)}$-homotopy limit spectral sequence, we improve the upper bound on the $\sE_{(2)}$-exponent from \cref{prop:c4sdc4expupbound} by showing that
\begin{equation}
    \exp_{\sE_{(2)}} \underline{H\Ffield_2}_{C_4 \rtimes C_4} \leq 4.
\end{equation}
\begin{Lemma}
    The Euler class $e$ of the 1-dimensional real representation given by pulling back the sign representation along the quotient map
    \begin{equation}
        C_4 \rtimes C_4 \to C_4 \rtimes C_4/\langle r, s^2 \rangle \cong C_2:
    \end{equation}
    satisfies $e^2 = 0$.
    \label{c4sdpc4:rep1}
\end{Lemma}
\begin{Proof}
    We have the composite of quotient maps
    \begin{equation}
        \overset{r}{C_4} \rtimes \overset{s}{C_4} \to \overset{\sbar}{C_2} \times \overset{\rbar}{C_2} \to \overset{\sbar}{C_2},
    \end{equation}
    the composite being the quotient-by-$\langle r, s^2 \rangle$-map. The Euler class of the sign representation of $C_2$ is $\delta_{\sbar} \in H^1(BC_2)$. Euler classes are natural, so $e$ is the pullback of $\delta_{\sbar} \in H^1(BC_2 \times C_2)$. In other words, the edge map $E_2^{1,0} \to H^1(BC_4 \rtimes C_4)$ maps $\delta_{\sbar} \to e$. Hence $\delta_{\sbar}$ detects $e$ in the LHSSS. Since $\delta_{\sbar}^2= 0$ on $E_3$ on account of the differentials in (\ref{c4sdpc4:diffs1}), $e^2$ is zero up to higher filtration, of which there is none.
\end{Proof}
\begin{Lemma}
    The Euler class $e$ of the 1-dimensional real representation given by pulling back the sign representation along the quotient map
    \begin{equation}
        C_4 \rtimes C_4 \to C_4 \rtimes C_4/\langle r^2, s \rangle \cong C_2
    \end{equation}
    satisfies $e^3 = 0$.
    \label{c4sdpc4:rep2}
\end{Lemma}
\begin{Proof}
    The proof is the same as the proof of \cref{c4sdpc4:rep1}, but now $e$ is detected by $\delta_{\rbar}$ instead of $\delta_{\sbar}$. The class $\delta_{\rbar}^2$ is non-zero, but $\delta_{\rbar}^3 = 0$, hence $e^3$ is zero up to higher filtration, of which there is none. 
\end{Proof}
\begin{Proposition}
    \label{c4sdpc4:goodupper}
    The $\sE_{(2)}$-exponent satisfies
    \begin{equation}
        \underline{H\Ffield_2}_{C_4 \rtimes C_4} \leq 4.
    \end{equation}
\end{Proposition}
\begin{Proof}
    We just have to remark that $\langle r^2, s^2 \rangle$, which is the unique maximal elementary abelian 2-subgroup of $C_4 \rtimes C_4$, is the intersection of $\langle r^2, s \rangle$ and $\langle r, s^2 \rangle$, the groups we divided out by in \cref{c4sdpc4:rep1} and \cref{c4sdpc4:rep2}. Hence by \cref{cor:1dimrepexps} and \cref{lem:expintersect}, $\exp_{\sE_{(2)}} \underline{H\Ffield_2} \leq 2 + 3 - 1 = 4$.
\end{Proof}
\cref{c4sdpc4:goodupper} together with the vanishing line of height 3 on $E_\infty$ implies that
\begin{Proposition}
    \label{prop:c4sdc4e2exp}
    The $\sE_{(2)}$-exponent satisfies
\begin{equation}
    3 \leq \exp_{\sE_{(2)}} \underline{H\Ffield_2} \leq 4.
\end{equation}
\end{Proposition}

\section{$(C_{4} \times C_{2}) \protect\overset{\psi_{5}}{\rtimes} C_{2}$}
\label{sec:c4xc2sdc2}
Let $C_4 = \langle x \rangle$, $C_2 = \langle y \rangle$, and recall from \cref{def:classact} that $\psi_5$ was defined to be the element of $\Aut(C_4 \times C_2)$ defined by $x \mapsto xy$, $y \mapsto y$.
In this section we will compute the $\sE_{(2)}$-homotopy limit spectral sequence for the group $G = (C_4 \times C_2) \overset{\psi_5}{\rtimes} C_2$ converging to its cohomology ring, which is isomorphic to (\cite[App.\ C, \#9(16)]{carlson2003cohomology})
\begin{equation}
    \Ffield_2[u, v, w, r, t]/(u^2, uv, uw, v^2r + w^2).
\end{equation}
with degrees $|u| = |v| = 1$, $|w| = |r| = |t| = 2$.

Under the correspondence $a \leftrightarrow x$, $b\leftrightarrow y$ and $c \leftrightarrow z$, this group admits a presentation given by
\begin{equation}
	G = \langle a,b,c \, | \, a^4 = b^2 = c^2 = e, \, ab = ba, \, bc = cb, \, cac^{-1} = ab \rangle.
	\label{id16-3:eq1}
\end{equation}

Let $A$ be the normal subgroup $\langle a^2,bc,c\rangle$ of $G$. The group $A$ is the unique maximal elementary abelian 2-subgroup of $G$, hence by \cref{prop:compserre} the $\sE_{(2)}$-homotopy limit spectral sequence reduces to the LHSSS obtained from the extension $A \to G \to C_2$.

The Weyl group $W_G(A) = \langle \abar \rangle$ acts on the generators of $A$ by
\begin{align}
    a(a^2)a^{-1} & = a^2, \\
    a(bc) a^{-1} & = a(bc)a^3 = (ab)^4c= c, \\
    a(c) a^{-1} & = bc.
\end{align}
Writing $x = \delta_{a^2}$, $y = \delta_{bc}$, and $z = \delta_{c}$, we get that $C_2$ interchanges $y$ and $z$ and fixes $x$.

Write $\sigma_1 = y+z$, $\sigma_2 = yz$, and $b = \delta_{\abar}$. Taking a $C_2$-projective resolution of $\Integers$ and applying $\Hom(-,H^*(BA))$, we get that the $s=0$-line of the $E_2$-page is given by the invariants
\begin{equation}
    E_2^{0,*} = \Ffield[\sigma_1, \sigma_2, x],
\end{equation}
and the $s$-lines for $s\geq 1$ are given by the invariants modulo the symmetrized classed, which is
\begin{equation}
    E_2^{s,*} = \Ffield_2[\sigma_2, x]\{b^s\},
\end{equation}
which assembles to
\begin{equation}
    E_2 = \Ffield_2[b, \sigma_1, \sigma_2, x]/(b \sigma_1)
\end{equation}
with $(s, t)$-degrees given by $|b| = (1,0)$, $|\sigma_1| = (0,1)$, $|\sigma_2| = (0,2)$, $|x| = (0,1)$.\subsection{Differentials}
The subgroup $C_4 = \langle a \rangle$ is normal because $b$ commutes with $a$ and $ca^2c^{-1} = (ab)^2 = a^2$.

The map of extensions
\begin{equation}
    \begin{tikzcd}
        \overset{a^2}{C_2} \arrow[hook]{d} \arrow[r] & \overset{a}{C_4} \arrow[hook]{d} \arrow{r} & \overset{\abar}{C_2} \arrow[equal]{d} \\
        A \arrow{r} & G \arrow{r}  & \overset{\abar}{C_2}
    \end{tikzcd}
    \label{id16-3:eq2}
\end{equation}
gives a map of LHSSS's. This is induced by the maps in cohomology of the left and right vertical maps, which are the inclusion $C_2 \to A$ and the corresponding quotient $C_2 \to C_2$ (which is the identity). The inclusion $C_2 \hookrightarrow A$ is on elements given by $a^2 \mapsto a^2$, hence on cohomology by
\begin{align}
    x = \delta_{a^2} & \mapsto \delta_{a^2}, \\
    y = \delta_{bc} & \mapsto 0, \\
    z = \delta_{c} & \mapsto 0,
    \intertext{and therefore}
    \sigma_1 = y + z & \mapsto 0, \\
    \sigma_2 = yz & \mapsto 0.
\end{align}
The right hand vertical map $C_2 \to C_2$ of (\ref{id16-3:eq2}) is the identity, and hence in cohomology also given by the identity map $b = \delta_{\abar} \mapsto \delta_{\abar}$.

Because the differentials of the LHSSS of the top row are given by $d_2(\delta_{a^2}) = \delta_{\abar}^2$ (\cref{subsec:cycl2gps}), we must have
\begin{align}
    d_2(x) & = b^2, \\
    d_2(\sigma_1) & = 0, \\
    d_2(\sigma_2) & = 0,
\end{align}
Observe that $\sigma_1 x$ survives to $E_3$ because
\begin{align}
    d_2(\sigma_1 x) & = \sigma_1 b^2 = 0.
\end{align}
The $E_2$-page with differentials is depicted in the next picture.
\begin{figure}[H]
	\centering
	\includegraphics[width=0.9\textwidth]{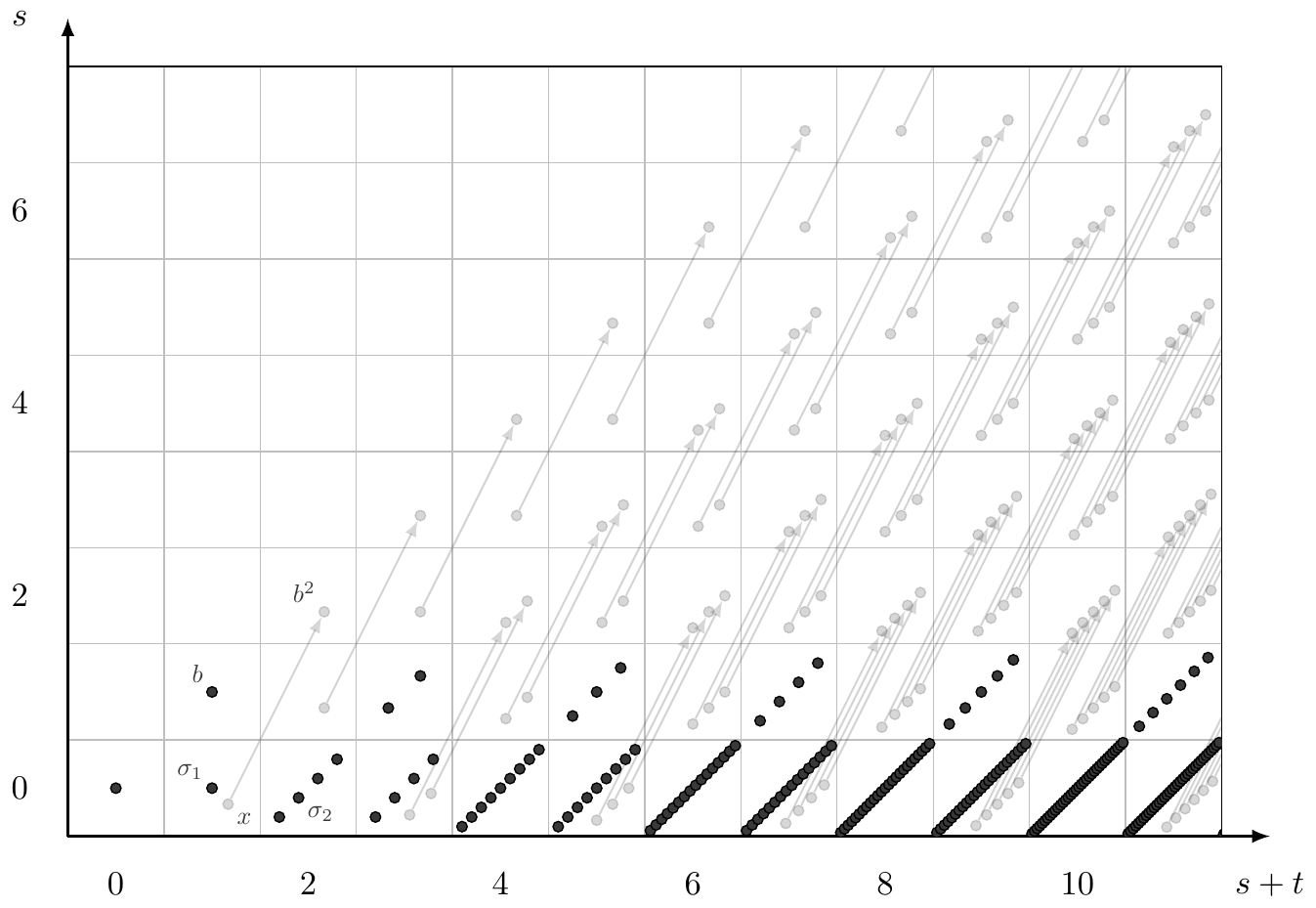}
	\caption{The $E_2$-page with differentials.}
\end{figure}
The $E_3$-page is given by
\begin{equation}
    E_3 =  \Ffield_2[b, \sigma_1, \sigma_2, [x^2], [\sigma_1 x]]/(b\sigma_1, b^2, [\sigma_1 x]^2 - \sigma_1^2[x^2], b[\sigma_1 x]),
\end{equation}
which has a vanishing line of height 2, and therefore the spectral sequence collapses as $E_3$.
\subsection{Poincar\'e series}
The $E_2$-page of the LHSSS is $E_2 = \Ffield_2[b, \sigma_1, \sigma_2, x]/(b \sigma_1)$. The polynomial ring $\Ffield_2[b, \sigma_1, \sigma_2, x]$ has Poincar\'e series $1/((1-t)^3(1-t^2))$, the ideal $(b\sigma_1)$ in this ring, which is a free module over the ring, has therefore Poincar\'e series $t^2/((1-t)^3(1-t^2))$. Again because the ideal is a free module, the quotient has Poincar\'e series $(1-t^2)/((1-t)^3(1-t^2))$.

The classes supporting the $d_2$-differential generate the sub-$\Ffield_2$-module of $E_2$ given by
\begin{equation}
    \Ffield_2[b, \sigma_2, [x^2]] \{x\}, 
\end{equation}
and their image is
\begin{equation}
    \Ffield_2[b, \sigma_2, [x^2]]\{b^2\}.
\end{equation}
The first has Poincar\'e series $t/((1-t)(1-t^2)^2)$, the second has Poincar\'e series $t^2/((1-t)(1-t^2)^2)$, and therefore $H^*(BG; \Ffield_2)$ has Poincar\'e series
\subsection{The exponent}
Consider the quotient map in the group extension
\begin{equation}
    A \to G \to C_2 \cong \langle \abar \rangle
\end{equation}
of which we just computed the LHSSS. This gives a 1-dimensional real representation of $G$ with oriented Euler class $e$. Because $e$ restricts to $0$ on $A$, we see that $e$ lives in filtration $>0$ on $E_\infty$, in other words $e = \delta_{\abar}$. The computation showed that $\delta_{\abar}^2 = 0$, since there is nothing in higher filtration. Hence by \cref{cor:1dimrepexps}, $\exp_{\sE_{(2)}}\underline{H\Ffield_2} \leq 2$, which together with the vanishing line of height 2 on the $E_\infty$-page implies
\begin{Proposition}
    The $\sE_{(2)}$-exponent satisfies
    \label{prop:c4xc2sdc2exp}
\begin{equation}
    \exp_{\sE_{(2)}} \underline{H\Ffield_2} = 2.
\end{equation}
\end{Proposition}

\section{$Q_8 \times C_2$}
\label{sec:q8xc2}
In this section we will compute the $\sF$-homotopy limit spectral sequence converging to the group cohomology of $Q_8\times C_2$ and using the family $\sF = \sE_{(2)}$ of elementary abelian subgroups. The computation suggests that there might be a better way to do this using a K\"unneth-type theorem for the $\sE_{(2)}$-homotopy limit spectral sequence for $\underline{H\Ffield_2}$, but we have not been able to prove such a theorem. The computation presented here is just a straightforward adaption of the computation of the $\sE_{(2)}$-homotopy limit spectral sequence for $Q_8$.

Let $i,j$ be the usual elements of $Q_8$, and let $\sigma$ be the generator of $C_2$. The group $Q_8 \times C_2$ has a maximal elementary abelian subgroup, $C_2 \times C_2 \cong \langle -1, \sigma \rangle$, which is furthermore central. Hence, by \cref{prop:compserre}, the $\sE_{(2)}$-homotopy limit spectral sequence reduces to the LHSSS for the central extension
\begin{equation}
    \langle -1, \sigma  \rangle \cong C_2 \times C_2 \to Q_8 \times C_2 \to C_2 \times C_2 \cong \langle \ibar, \jbar \rangle.
    \label{Q8xC2:eq1}
\end{equation}

Since $\langle -1, \sigma \rangle$ is central, the coefficient system of $B(-)$ applied to (\ref{Q8xC2:eq1}) is trivial. Writing $a = \delta_{\ibar}$, $b = \delta_{\jbar}$, $x = \delta_{-1}$ and $y=\delta_{\sigma}$, the $E_2$-page of the $\sE_{(2)}$ homotopy limit spectral sequence is therefore
\begin{align}
    E_2^{s,t} & \cong H^s(BC_2 \times C_2; H^t(BC_2 \times C_2)) \\
    & \cong \Ffield_2[a,b,x,y]
\end{align}
with $(s,t)$-degrees given by $|a| = |b| = (1,0)$, $|x| = |y| = (0,1)$.

Applying naturality to the map of extensions
\begin{equation}
    \begin{tikzcd}
        \langle -1, \sigma \rangle \arrow{d} \arrow{r} & Q_8 \times C_2 \arrow{d} \arrow{d} \arrow{r} & \langle \ibar, \jbar \rangle \arrow{d} \\
        \langle \sigma \rangle \arrow{r} & \langle \sigma \rangle \arrow{r} & \{e\}
    \end{tikzcd}.
\end{equation}
shows that $\delta_{\sigma} = y$ is a permanent cycle.

The differentials of the LHSSS for the extension
\begin{equation}
    \langle -1 \rangle \cong C_2 \to Q_8 \to C_2 \times C_2 \cong \langle \ibar, \jbar \rangle
\end{equation}
are well known (see \cref{q8:lhsss}). Using this, naturality, and the map of extensions
\begin{equation}
    \begin{tikzcd}
        \langle  -1, \sigma \rangle \arrow{r} \arrow{d} & Q_8 \times C_2\arrow{d} \arrow{r} & \langle \ibar, \jbar \rangle \arrow[equal]{d} \\
        \langle -1\rangle \arrow{r} & Q_8  \arrow{r} & \langle \ibar , \jbar \rangle
    \end{tikzcd}
\end{equation}
shows that
\begin{align}
    d_2(x) & = a^2 + ab + b^2 \\
    d_3([x^2]) & = a^2 b + ab^2.
\end{align}
Hence
\begin{equation}
    E_4 \cong \Ffield_2[a,b,[x^4],y]/(a^2 + ab + b^2, a^2b + ab^2)
    \label{q8xc2:eq2} 
\end{equation}
which has a vanishing line of height 4. Therefore the spectral sequence collapses and $E_4 = E_\infty$.
\subsection{The exponent}
The above computation allows us to prove
\begin{Proposition}
    \label{prop:q8xc2e2exp}
    The $\sE_{(2)}$-exponent is
    \begin{equation}
        \exp_{\sE_{(2)}} \underline{H\Ffield_2}_{Q_8 \times C_2} = 4.
    \end{equation}
\end{Proposition}
\begin{Proof}
    Elementary abelian subgroups of $Q_8$ pull back to elementary abelian subgroups of $Q_8 \times C_2$ along the projection map $Q_8 \times C_2 \to Q_8$. Therefore, pulling back the 4-dimensional real representation of $Q_8$ whose projectivation has isotropy in $\sE_{(2)}$ considered in \cite[Ex.\ 5.18]{mnn} along the projection map gives a 4-dimensional real representation of $Q_8 \times C_2$ with projectivation with isotropy in $\sE_{(2)}(Q_8 \times C_2)$. Hence, by \cref{projbund:prop1}, the exponent satisfies $\exp_{\sE_{(2)}} \underline{H\Ffield_2}_{Q_8 \times C_2} \leq 4$. 

    Conversely, the vanishing line of height 4 in (\ref{q8xc2:eq2}) shows the reverse inclusion, giving the desired result.
\end{Proof}

\section{$D_8 \times C_2$}
\label{sec:d8xc2}
\subsection{Introduction}
In this section we analyze the behaviour of the $\sE_{(2)}$-homotopy limit spectral sequence converging to $H^*(BD_8 \times C_2;\Ffield_2)$. The computation suggest, as the computation for $Q_8 \times C_2$ did, that there might be a better way to do this computation using a K\"unneth-type theorem for the $\sE_{(2)}$-homotopy limit spectral sequence for $\underline{H\Ffield_2}$, but we have not been able to prove such a theorem. The computation presented here is just a straightforward adaption of the computation for $D_8$.

The elementary abelian subgroups of $D_8$ pull back to elementary abelian subgroups of $D_8 \times C_2$ along the projection
\begin{equation}
    D_8 \times C_2 \to D_8,
\end{equation}
and since $D_8$ has a 2-dimensional real representation $V$ such that the projectivation $\Proj(V)$ has istropy in $\sE_{(2)}$, the same is true for $D_8 \times C_2$. Hence
\begin{equation}
    \exp_{\sE_{(2)}} \underline{H\Ffield_2} \leq 2.
\end{equation}
Therefore the $\sE_{(2)}$-homotopy limit spectral sequence will collapse on the $E_3$-page with a vanishing line of height $\leq 2$. In fact, since the cohomology of both $D_8$ and $C_2$ is detected on elementary abelian subgroups (see, e.g., \cite[Lem.\ 4.6]{quillen71adams}), the same is true for $D_8 \times C_2$, by the K\"unneth Theorem. Therefore the $\sE_{(2)}$-homotopy limit spectral sequence will be concentrated in the $s=0$-line on $E_3 = E_\infty$. 
\subsection{The decomposition}
Let $D_8 = \langle \sigma, \rho \, | \, \sigma^2 = \rho^4 = e, \sigma \rho \sigma^{-1} = \rho^{-1} \rangle$, $C_2 = \langle u \, | \, u^2 = e \rangle$. The maximal elementary abelian subgroups of $D_8 \times C_2$ are $F_1 := \langle \sigma, \sigma\rho^2, u \rangle \cong C_2^{\times 3}$ and $F_2 := \langle \sigma \rho, \sigma \rho^3, u \rangle \cong C_2^{\times 3}$. Closing this set of subgroups under intersections adds the additional group
\begin{equation}
    F_1 \cap F_2 = \langle \rho^2, u \rangle =: Z \cong C_2 \times C_2
\end{equation}
We apply the decomposition of \cref{sec:absplit} with $K = Z$, and $H_i = F_i$. We therefore proceed by computing the LHSSS's $E_*^{*,*}(F_i)$ and $E_*^{*,*}(Z)$.
\subsection{The Serre spectral sequence for the subgroups $F_i$}
We do the case $i = 1$, the case $i=2$ being completely similar. We consider the extension
\begin{equation}
    \langle \sigma, \sigma \rho^2, u \rangle \cong F_1 \to D_8 \times C_2 \to C_2 \cong \langle \rhobar \rangle.
    \label{D8xC2:eq2}
\end{equation}
Write $x_1 = \delta_{\sigma}$, $y_1 = \delta_{\sigma \rho^2}$, $z_1 = \delta_u$, and $a_1 = \delta_{\rhobar}$. The extension (\ref{D8xC2:eq2}) induces an action of $C_2$ on $F_1$ given by interchanging $\sigma$ and $\sigma \rho^2$ and fixing $u$. Hence the $s=0$-row of the Serre spectral sequence is isomorphic to the invariants in $x_1$ and $y_1$ tensored with $\Ffield_2[z_1]$. Writing $\sigma_{1,1} = x_1+y_1$ and $\sigma_{2,1} = x_1y_1$ for the first and second elementary symmetric polynomial in $x_1$ and $y_1$, we have 
\begin{equation}
    E_2^{0,*} = \Ffield_2[\sigma_{1,1}, \sigma_{2,1}] \otimes_{\Ffield_2} \Ffield_2[z_1].
\end{equation}
The rows with $s>0$ are given by the symmetric polynomials in $x_1$ and $y_1$ divided out by the symmetrized polynomials, i.e., for $s>0$,
\begin{equation}
    E_2^{s,*} = a_1^s \Ffield_2[\sigma_{2,1}] \otimes_{\Ffield_2} \Ffield_2[z_1].
\end{equation}
Combining this gives
\begin{equation}
    E_2^{*,*} = \Ffield_2[a_1,\sigma_{1,1},\sigma_{2,1}]/(a_1\sigma_{1,1}) \otimes_{\Ffield_2} \Ffield_2[z_1].
\end{equation}
with $(s,t)$ degrees given by $|a_1| = (1,0)$, $|\sigma_{1,1}| = |z_1| = (0,1)$, $|\sigma_{2,1}| = (0,2)$.

For every $k \geq 0$, the $k$-column $\bigoplus_s E_2^{s,k-s}$ of the $E_2$-page has the dimension of $H^k(BD_8 \times C_2;\Ffield_2)$, as can be seen from the K\"unneth-formula, hence the spectral sequence collapses at $E_2$.

We don't repeat the calculation for the case $i=2$, but we do define for later use for $i=2$ the classes $x_2 = \delta_{\sigma \rho}$, $y_2 = \delta_{\sigma \rho^3}$, $z_2 = \delta_u$, $a_2 = \delta_{\rhobar}$, $\sigma_{1,2} = x_2 + y_2$ and $\sigma_{2,2} = x_2y_2$.
\subsection{The Serre spectral sequence for the subgroup $Z$}
We consider the central extension
\begin{equation}
    \langle \rho^2, u \rangle = Z \to D_8 \times C_2 \to C_2 \times C_2 = \langle \rhobar \rangle.
\end{equation}
Writing $a_3 = \delta_{\sigmabar}$, $b_3 = \delta_{\overline{\sigma \rho}}$, $x_3 = \delta_{\rho^2}$, $z_3 = \delta_u$, we have as the $E_2$-page:
\begin{equation}
    E_2 = \Ffield_2[a_3, b_3, x_3, z_3] 
\end{equation}
with $(s,t)$-degrees given by $|a_3| = |b_3| = (1,0)$, $|x_3| = |z_3| = (0,1)$.

Applying naturality to the map of extensions
\begin{equation}
    \begin{tikzcd}
        \langle \rho^2 \rangle \arrow{r} \arrow{d} & D_8 \arrow{d} \arrow{r} & C_2 \times C_2 \arrow{d} \\
        \langle \rho^2, u \rangle \arrow{r} & D_8 \times C_2 \arrow{r} & C_2 \times C_2
    \end{tikzcd}
\end{equation}
shows that
\begin{equation}
    d_2(x_3) = a_3 b_3.
\end{equation}
Applying naturality to the map of extensions
\begin{equation}
    \begin{tikzcd}
        \langle \rho^2, u \rangle \arrow{r} \arrow{d} & D_8 \times C_2 \arrow{r} \arrow{d} & C_2 \times C_2 \arrow{d} \\
        \langle u \rangle \arrow{r} & \langle u \rangle \arrow{r} & \{e\}
    \end{tikzcd}
\end{equation}
shows that $z_3$ is a permanent cycle. The classes $a_3$ and $b_3$ are permanent cycles for degree reasons. Hence
\begin{equation}
    E_3 = \Ffield_2[a_3,b_3,[x_3^2],z_3]/(a_3b_3).
\end{equation}
Again by naturality, we see that $[x_3^2]$ is a permanent cycle, hence $E_3 = E_\infty$.
\subsection{The map $j^*$}
We now compute the map $j^*$ in the long exact sequence (\ref{absplit:eq3}) as a map of spectral sequences coming from the maps of extensions
\begin{equation}
    \begin{tikzcd}
        Z \arrow{d} \arrow{r} & D_8 \times C_2 \arrow{d} \arrow{r} & C_2 \times C_2 \arrow{d} \\
        F_i \arrow{r} & D_8 \times C_2 \arrow{r} & C_2.
    \end{tikzcd}
\end{equation}
Hence we determine the maps $Z \to F_i$ and $C_2 \times C_2 \to C_2$ and their effect on cohomology classes.
\subsubsection{Inclusion of the center $Z$ in $F_1$}
The inclusion
\begin{align}
    \langle \rho^2, u \rangle & \to \langle \sigma, \sigma \rho^2, u \rangle \\
    \intertext{is on elements given by}
    \rho^2 & \mapsto \sigma \cdot \sigma \rho^2, \\
    u & \mapsto u \\
    \intertext{and therefore is in cohomology given by}
    x_3 = \delta_{\rho^2} & \mapsfrom \delta_\sigma = x_1 \\
    x_3 = \delta_{\rho^2} & \mapsfrom \delta_{\sigma \rho^2} = y_1 \\
    z_3 = \delta_{u} & \mapsfrom \delta_u = z_1.
\end{align}
\subsubsection{Inclusion of the center $Z$ in $F_2$}
The inclusion
\begin{align}
    \langle \rho^2, u \rangle & \to \langle \sigma \rho, \sigma \rho^3, u \rangle \\
    \intertext{is on elements given by}
    \rho^2 & \mapsto \sigma \rho \cdot \sigma \rho^3, \\
    u & \mapsto u \\
    \intertext{and therefore is in cohomology given by}
    x_3 = \delta_{\rho^2} & \mapsfrom \delta_{\sigma \rho} = x_2 \\
    x_3 = \delta_{\rho^2} & \mapsfrom \delta_{\sigma \rho^3} = y_2 \\
    z_3 = \delta_u & \mapsfrom \delta_u = z_1.
\end{align}
\subsubsection{The quotient of $D_8\times C_2/F_1$ by $C_2$}
The quotient
\begin{align}
    D_8 \times C_2 /\langle \rho^2, u \rangle & \to D_8 \times C_2 /\langle \sigma, \sigma \rho^2, u \rangle \\
    \intertext{is on elements given by}
    \sigmabar & \mapsto \ebar \\
    \sigmabar \rhobar & \mapsto \rhobar \\
    \intertext{and therefore is in cohomology given by}
    b_3 = \delta_{\sigmabar\rhobar} & \mapsfrom \delta_{\rhobar} = a_1.
\end{align}
\subsubsection{The quotient of $D_8 \times C_2/F_2$ by $C_2$}
The quotient
\begin{align}
    D_8 \times C_2/ \langle \rho^2,u \rangle & \to D_8 \times C_2/\langle \sigma \rho, \sigma \rho^3, u \rangle \\
    \intertext{is on elements given by}
    \sigmabar & \mapsto \rhobar \\
    \sigmabar \rhobar & \mapsto \ebar \\
    \intertext{and therefore is in cohomology given by}
    a_3  = \delta_{\sigmabar}  & \mapsfrom \delta_{\rhobar} = a_2.
\end{align}
\subsubsection{Summary}
Summarizing, in the long exact sequence (\ref{absplit:eq3}), we have
\begin{align}
    E_2(Z) & = \Ffield_2[a_3,b_3,x_3,z_3], \\
    E_2(F_1) & = \Ffield_2[a_1,\sigma_{1,1},\sigma_{2,1},z_1]/(a_1 \sigma_{1,1}), \\
    E_2(F_2) & = \Ffield_2[a_2, \sigma_{1,2}, \sigma_{2,2},z_2]/(a_2 \sigma_{1,2}).
\end{align}
The map $j^*$ is given by
\begin{align}
    a_1 & \mapsto b_3, \\
    \sigma_{1,2} = x_1 + y_1 & \mapsto x_3 + x_3 = 0, \\
    \sigma_{2,1} & \mapsto x_3^2, \\
    z_1 & \mapsto z_3, \\
    a_2 & \mapsto a_3, \\
    \sigma_{1,2} & \mapsto 0, \\
    \sigma_{2,2} & \mapsto x_3^2, \\
    z_2 & \mapsto z_3.
    \label{D8xC2:eq3}
\end{align}
\subsubsection{The kernel of $j^*$}
The effect (\ref{D8xC2:eq3}) of $j^*$ on cohomology classes shows that in $s$-degree $s > 0$, $j^*$ is injective because $a_1$ and $a_2$ map to different non-nilpotent elements of $E_2(Z)$. Hence $\ker j^*$ is concentrated in $s$-degree equal to $0$, and from (\ref{D8xC2:eq3}) we deduce
\begin{equation}
    \ker j^* = \Ffield_2[(\sigma_{1,1},0), (0,\sigma_{1,2}), (\sigma_{2,1},\sigma_{2,2}), (z_1,z_2)]/((\sigma_{1,1},0)(0,\sigma_{1,2})).
    \label{D8xC2:eq4}
\end{equation}
\begin{Remark}
    As already mentioned, $D_8 \times C_2$ has its $\Ffield_2$-cohomology detected on the elementary abelian 2-subgroups. Hence the edge map to the $s=0$-line of the $\sE_{(2)}$-homotopy limit spectral sequence, which is precisely (\ref{D8xC2:eq4}), is injective. Moreover, the Poincar\'e series of (\ref{D8xC2:eq4}) is precisely equal to the one of $H^*(BD_8\times C_2;\Ffield_2)$, hence the edge morhpism is also surjective. In particular, there will be no differentials emitting from the $s=0$-line. We argued before that $\exp_{\sE_{(2)}} \underline{H\Ffield_2} \leq 2$, hence all classes above the $s=0$-line will either support a differential or be hit by one on $E_2$.
\end{Remark}
\subsubsection{The image of $j^*$}
The effect (\ref{D8xC2:eq3}) of $j^*$ on cohomology classes shows that
\begin{equation}
    \Im j^* = \Ffield_2\{a_3^{i_1} x_3^{2l_1} z_3^{k_1}\} \cup \{b_3^{i_2}x_3^{2l_2}z_3^{k_2}\}.
\end{equation}
In other words, a class $a_3^i b_3^{i'} x_3^l z_3^k$ \emph{not} being in the image of $j^*$ is exactly equivalent to at least one of the following conditions being true:
\begin{enumerate}[1.]
    \item The exponents $i$ and $i'$ are both $> 0$.
    \item The exponent $l \equiv 1 \pmod{2}$.
\end{enumerate}
\subsection{The $E_2$-page}
Having computed $\ker j^*$ and $\Im j^*$, we know $A = \coker j^*[1]$ and $B = \ker j^*$ in (\ref{absplit:eq5}). Moreover, since $\partial$ is a map of spectral sequences, and the fact that on $E_2(Z)$ we have
\begin{equation}
    d_3(x_3) = a_3 b_3,
    \label{D8xC2:eq7}
\end{equation}
shows that we also have this differential on the $\sE_{(2)}$-homotopy limit spectral sequence under the identification (\ref{absplit:eq5}).

Therefore all classes $a_3^i b_3^{i'} x_3^l z_3^k$ in $\coker j^*$ with $l \equiv 1 \pmod{2}$ will support a differential:
\begin{equation}
    d_2(a_3^i b_3^{i'} x_3^l z_3^k) = a_3^{i+1} b_3^{i'+1} x_3^{l-1}z_3^k.
    \label{D8xC2:diffs}
\end{equation}
Hence all classes with even $x_3$-exponent and strictly positive $a_3$ and $b_3$-exponent get hit, which accounts for all classes in $\coker j^*$: they all either hit or get hit by a $d_2$. Therefore there cannot be any $d_2$-differentials emitting from the $B$-summand, and $E_3 = B$. Since this is concentrated in the $s=0$-row, $E_3 = E_\infty$.
\subsection{Multiplicative extensions}
Since on $E_\infty$ everything is concentrated in the $s=0$ row, we have
\begin{align}
    B  & \cong E_\infty^{0,*}  \\
    & \cong F^0H^*(BD_8 \times C_2)/F^1H^*(BD_8 \times C_2) \\
    & = H^*(BD_8 \times C_2).
\end{align}

\subsection{The exponent}
We now prove
\begin{Proposition}
    \label{prop:d8xc2e2exp}
    The $\sE_{(2)}$-exponent satisfies
    \begin{equation}
        \exp_{\sE_{(2)}} \underline{H\Ffield_2} = 2.
    \end{equation}
\end{Proposition}
\begin{Proof}
    Pulling back transitive $D_8$-orbits with isotropy in elementary abelian groups along the projection map $D_8 \times C_2$ gives transitive $D_8 \times C_2$-orbits with isotropy in elementary abelian groups. Therefore pulling back the representation of $D_8$ considered in the proof of \cref{cor:expd8upperbound} along the projection map $D_8 \times C_2 \to D_8$ gives a 2-dimensional real representation of $D_8 \times C_2$ with isotropy in $\sE_{(2)}$. Hence by \cref{projbund:prop1} we have $\exp_{\sE_{(2)}} \underline{H\Ffield_2}_{D_8 \times C_2} \leq 2$.

    Conversely, in the above computation of the $\sE_{(2)}$-homotopy limit spectral sequence converging to $H^*(BD_8 \times C_2;\Ffield_2)$, we saw that $E_2 \neq E_\infty$ (\cref{D8xC2:diffs}). Hence we also have the reverse inequality, and the result follows.
\end{Proof}

\chapter{Properties of the $\sF$-homotopy limit spectral sequence for equivariant $K$-theory}
\label{ch:kthy}
\section{Elementary abelian 2-groups}
\subsection{Introduction}
In \cite[App.\ B]{mnn}, the $\sC$-homotopy limit spectral sequence of $KU$ is computed for $G = C_2^{\times 2}$. The target of this spectral sequence is $\pi_* KU \cong R(G)[\beta^{\pm}]$, the polynomial ring on the Bott periodicity generator $\beta$, with $|\beta| = 2$, and its inverse $\beta^{-1}$, with coefficients in the representation ring of $G$. The $E_3$-page is the last page with differentials, hence $\exp_{\sC} KU \geq 3$.

In this section, we will establish lower bound on the $\sC$-exponent of $KU$ for all elementary abelian 2-groups. More precisely, we will show that for $G = C_2^{\times n}$ we have $\exp_{\sC}KU \geq n$.

The proof can be summarized as follows. Most of the work goes into describing the edge map
\begin{equation}
    R(C_2^{\times n}) \to \sideset{}{^0}\lim_{\sO(C_2^{\times n})^{\op}_{\sC}} R(-) \cong E_2^{0,0}
    \label{eq1:el2gps}
\end{equation}
of the spectral sequence. The upshot is that the exponent (as a group) of the cokernel of (\ref{eq1:el2gps}) is $2^{n-1}$ (\cref{cor:expelab2gps}). A qualitative analysis of the $E_2$-page shows that the classes in positive filtration are (at most) 4-torsion, and Bott periodicity implies $E_2^{s,t} = 0$ if $s+t$ is odd. This then implies that there are non-zero differentials on at least $\ceil{(n-1)/2}$ odd pages, whence $\exp_{\sC} KU \geq n$ (\cref{prop:expkulowerbound}).
\subsection{The exponent of the cokernel of the edge map}
\label{subsec:expcok}
We will determine the exponent of the cokernel of the edge map
\begin{equation}
    R(C_2^{\times n}) \to \lim_{\sO(C_2^{\times n})_{\sC}^{\op}} R(-) = E_2^{0,0}
    \label{expcok:eq1}
\end{equation}
of the $\sC$-homotopy limit spectral sequence of $KU$. We will do this by finding bases of the left and right hand side that diagonalize the edge map.
\subsubsection{Splitting of the augmentation ideal}
For $e$ the trivial subgroup, identify $R(e) \cong \Integers$ by sending the trivial representation to $1$. For every group $H$, the restriction map then gives a map $R(H) \to \Integers$, the kernel of which is called the augmentation ideal of $R(H)$ and is denoted $I(H)$. The map $R(H) \to \Integers$ is naturally split by sending $1 \in \Integers$ to the trivial representation. 

For the coefficient system $R(-)$ on the orbit category $\sO(G)_{\sF}$, for every $H \in \sF$, the restriction map $R(H) \to R(e) \cong \Integers$ gives a map $R(-) \to \underline{\Integers}$, where $\underline{\Integers}$ is the constant coefficient system at $\Integers$. This map is naturally split, by applying the natural splitting to every individual group $H \in \sF$. 

Combining the two splittings splits the map (\ref{expcok:eq1}) as a direct sum of maps
\begin{equation}
    \Integers \oplus I(C_2^{\times n}) \to \lim_{\sO(C_2^{\times n})_{\sC}^{\op}}\underline{\Integers} \oplus \lim_{\sO(C_2^{\times n})_{\sC}^{\op}} I(-) \cong \Integers \oplus \lim_{\sO(C_2^{\times n})_{\sC}^{\op}} I(-).
\end{equation}

This direct sum decomposition reduces describing (\ref{expcok:eq1}) to describing
\begin{equation}
    I(C_2^{\times n}) \to \lim_{\sO(C_2^{\times n})_{\sC}^{\op}} I(-).
    \label{expcok:eq2}
\end{equation}
\subsubsection{A basis for $I(C_2^{\times n})$}
Denote
\begin{equation}
    C_2^{\times n} = \langle s_1,\, \ldots,\, s_n \, | \, s_i^2 = e, \, s_j s_k = s_ks_j \, \forall i,j,k \rangle.
\end{equation}
Denote by $\sigma_j$ the sign representation of $\langle s_j \rangle \cong C_2$, and also of the representation of $C_2^{\times n}$ given by
\begin{equation}
    \sigma_j(s_i) =
    \begin{cases}
        -1 & \text{if }j=i, \\
        1 & \text{otherwise.}
    \end{cases}
\end{equation}
For $J \subset \{1,\ldots,n\}$, denote
\begin{equation}
    a_J = \prod_{j \in J} (1- \sigma_j) \in R(C_2^{\times n}).
\end{equation}
Observe that $a_J \in I(C_2^{\times n})$ if and only if $J \neq \varnothing$. In fact, more is true.
\begin{Proposition}
    The set $A \coloneqq \{a_J\}_{\varnothing \neq J \subset \{1,\ldots,n\}}$ is a basis of $I(C_2^{\times n})$ as a $\Integers$-module.
    \label{prop:basisA}
\end{Proposition}
\begin{Proof}
    Since $\{1,\sigma_j\}$ is a basis of $R(\langle s_j \rangle)$ (\cite[Sec.\ 5.1]{Serrelinrep}) as a $\Integers$-module, the set $\{1,\, 1- \sigma_j\}$ is also basis of $R(\langle s_i \rangle)$. Since $R(C_2^{\times n}) = \bigotimes_i R(\langle s_i \rangle)$ (\cite[Thm.\ 3.2.10]{Serrelinrep}) the set $\{a_J\}_{J \subset \{1,\ldots,n\}}$ is a basis of $R(C_2^{\times n})$. Since $a_{\varnothing}$ generates the $\Integers$-summand of $R(C_2^{\times n}) \cong \Integers \oplus I(C_2^{\times n})$ and all $a_J \in I(C_2^{\times n})$ for $J \neq \varnothing$, the result follows.
\end{Proof}
\subsubsection{A basis for $\lim_{\sO(C_2^{\times n})_{\sC}^{\op}} I(-)$}
We will now describe a $\Integers$-basis of $\lim_{\sO(C_2^{\times n})_{\sC}^{\op}} I(-)$. This limit is given by the equalizer
\begin{equation}
    \eq \left( \prod_{C \in \sC} I(C) \to I(e) \right),
\end{equation}
where the map from the product to $I(e)$ comes from the universal property of the product and the restriction maps $I(C) \to I(e)$ for every $C$. But since $I(e) = 0$ this is equal to $\prod_{C \in \sC'} I(C)$, where $\sC'$ denotes the the set (not family) of non-trivial cyclic subgroups of $C_2^{\times n}$.
\begin{Notation}
For $J \subset \{1,\ldots,n\}$, denote $s_J = \prod_{j \in J} s_j \in C_2^{\times n}$. Then the elements of $\sC'$ are in $(1:1)$-correspondence with the non-empty $J$ via $J \leftrightarrow \langle s_J \rangle$.

Denote the sign representation of $\langle s_J \rangle$ by $\sigma_J$, and denote $1 - \sigma_J \in I(\langle \sigma_J \rangle)$ by $\overline{\sigma_J}$. Observe that $I(\langle s_J \rangle)$ is $\Integers$-module of rank 1, generated by $\overline{\sigma_J}$.

Denote the projection $\prod_{C \in \sC'} I(C) \to I(\langle s_J \rangle)$ by $\pr_J$, and let $b'_J$ be the element of $\prod_{C \in \sC'} I(C)$ satisfying $\pr_K(b'_J) = \overline{\sigma_J}$ if $K=J$ and 0 otherwise.
\end{Notation}

This discussion shows
\begin{Proposition}
    The set $B' \coloneqq \{b'_J \, | \, \varnothing \neq J \subset \{1,\ldots,n\}\}$ is a basis of $\lim_{\sO(C_2^{\times n})_{\sC}^{\op}} I(-)$.
    \label{prop:basisBprime}
\end{Proposition}

The map (\ref{expcok:eq2}) is not diagonal with respect to the bases $A$ and $B'$. To get a diagonal matrix, denote for non-empty $J \subset \{1,\ldots,n\}$
\begin{equation}
    b_J = \sum_{J \subset K \subset \{1,\ldots,n\}} b'_K.
\end{equation}
\begin{Proposition}
    The set $B = \{ b_J \}_{\varnothing \neq J \subset \{1,\ldots,n\}}$ is a basis of $\lim_{\sO(C_2^{\times n})_{\sC}^{\op}} I(-)$ as a $\Integers$-module.
    \label{prop:basisB}
\end{Proposition}
\begin{Proof}
    First note that $B \subset \Integers B'$. Conversely, by the principle of inclusion-exclusion,
    \begin{equation}
        b_J' = b_J - \sum_{\substack{J \subset K \\ \# (K-J) = 1}} b_K + \sum_{\substack{J \subset K \\ \#(K-J) = 2}} b_K - \cdots,
    \end{equation}
    hence also $B' \subset \Integers B$. \cref{prop:basisBprime} and $\# B' = \# B$ now imply the desired result.
\end{Proof}
\subsubsection{The edge map}
We now describe the edge map (\ref{eq1:el2gps}), which we reduced to describing (\ref{expcok:eq2}), which we will do in terms of the bases $A$ and $B$.
\begin{Warning}
    Although the next proposition is phrased in terms of the bases $A$ (\cref{prop:basisA}) and $B$ (\cref{prop:basisB}), the proof will also make use of the basis $B'$ (\cref{prop:basisBprime}).
\end{Warning}
\begin{Proposition}
    The map (\ref{expcok:eq2}) is given on the bases $A$ (\cref{prop:basisA}) and $B$ (\cref{prop:basisB}) by
    \begin{equation}
        a_J  \mapsto 2^{\# J -1} b_J.
        \label{edgeonbasis:eq1}
    \end{equation}
\end{Proposition}
\begin{Proof}
    For $\sigma$ the sign representation of a cyclic group of order 2, we have $(1-\sigma)^2 = 2\sigma$. Hence for basis elements $b'_J, b'_{K} \in B'$, we have
    \begin{equation}
        b'_J b'_{K} = 
        \begin{cases}
            2b'_J & \text{if } J = K, \\
            0 & \text{otherwise.}
        \end{cases}
        \label{prodstrL:eq1}
    \end{equation}
    By induction therefore, we also have $(b'_J)^m = 2^{m-1}b'_J$.

    First we verify (\ref{edgeonbasis:eq1}) on basis elements indexed by singletons:
    \begin{align}
        a_{\{j\}} & = 1 - \sigma_j \\
        & \mapsto \sum_{j \in J \subset\{1,\ldots,n\}} b'_J \\
        & =: b_{\{j\}}.
    \end{align}

    Then, by multiplicativity of the edge map, we get for all non-empty $J \subset \{1,\ldots,n\}$
    \begin{align}
        a_J & = \prod_{j \in J} (1-\sigma_j)  \\
        & \mapsto \prod_{j \in J} \left(\sum_{j \in K \subset\{1,\ldots,n\}} b'_K\right) \\
        \intertext{By the product structure on $\lim_{\sO(C_2^{\times n})_{\sC}^{\op}} I(-)$, all terms of the sum with $K \not\supset J$ will go to zero in the product. Hence we can discard those terms:}
        & = \prod_{j \in J} \left(\sum_{J \subset K \subset \{1,\ldots,n\}}b'_K \right)\\
        \intertext{The factors of this product do not depend on the index $j$ anymore, so we get:}
        & = \left(\sum_{J \subset K \subset \{1,\ldots,n\}} b'_K\right)^{\# J} \\
        \intertext{which by \cref{prodstrL:eq1} equals}
        & = \sum_{J \subset K \subset \{1,\ldots,n\}} (b'_K)^{\# J} \\
        & = 2^{\#J -1} \sum_{J \subset K \subset \{1,\ldots,n\}} b_K' \\
        & =: 2^{\# J -1 } b_J.
    \end{align}
\end{Proof}
\begin{Corollary}
    The cokernel of (\ref{expcok:eq2}) (which equals the cokernel of (\ref{eq1:el2gps})), is 
    \begin{equation}
        \bigoplus_{\varnothing \neq J \subset \{1,\ldots,n\}} \Integers/2^{\#J -1}\Integers \{b_J\}. 
    \end{equation}
    \label{cor2:expelab2gps}
\end{Corollary}
\begin{Corollary}
    The exponent of the cokernel of (\ref{eq1:el2gps}) is $2^{n-1}$.
    \label{cor:expelab2gps}
\end{Corollary}
\subsection{A lower bound on $\exp_{\sC} KU$}
Using \cref{cor:expelab2gps}, we now prove a lower bound on $\exp_{\sC} KU$. We do this by analyzing the $\sC$-homotopy limit spectral sequence
\begin{equation}
    E_2^{s,t} = \sideset{}{^s}\lim_{\sO(C_2^{\times n})^\op_{\sC}}R(-)[\beta^{\pm}] \Rightarrow R(C_2^{\times n})[\beta^{\pm}].
        \label{eqkth:eq2}
\end{equation}
where $\beta$ is the Bott periodicity generator with $(s,t)$-degree $|\beta| = (0,2)$.
We need the following lemma, which gives an upper bound on the torsion on the $E_2$-page for $s \geq 2$.
\begin{Lemma}
    \label{eqkth:lem1}
    In the $\sC$-homotopy limit spectral sequence (\ref{eqkth:eq2}), $E_2^{s,t}$ is annihilated by $4$ for $s \geq 1$.
\end{Lemma}
\begin{Proof}
    For a family $\sF \subset \sC$, write $\underline{\Integers}[\sF]$ for the constant coefficient system $\underline{\Integers}$ restricted to $\sO(C_2^{\times n})_{\sF}$ and then left Kan extended back up to $\sO(C_2^{\times n})_{\sC}$, as in \cite[App.\ B]{mnn}. As in \cite[B.3]{mnn}, the sum of the counit maps
    \begin{equation}
        \bigoplus_{C \in \sC} \underline{\Integers}[\{e,C\}] \to \underline{\Integers}
    \end{equation}
    is surjective and yields a short exact sequence of coefficient systems
    \begin{equation}
        0 \to \bigoplus_{|\sC| - 1} \underline{\Integers}[\{e\}] \to \bigoplus_{C \in \sC} \underline{\Integers}[\{e,C\}] \to \underline{\Integers} \to 0. \label{eqkth:eq3}
    \end{equation}
    This induces a long exact sequence of $\Ext$-groups as in \cite[eq.\ B.6]{mnn} that we will not spell out, but we describe what we need.

    The $\Ext$-groups corresponding to the left term in (\ref{eqkth:eq3}) are isomorphic to a direct sum of groups of the form $H^s(BC_2^{\times n}; \Integers)$, which are 2-torsion for $s \geq 1$ by the K\"unneth theorem.

    The $\Ext$-groups corresponding to the middle term in (\ref{eqkth:eq3}) are isomorphic to a direct sum of groups of the form $H^s(BC_2^{\times (n - 1)};R(C_2))$. Here $R(C_2)$ is a free $\Integers$-module with trivial $C_2^{\times(n-1)}$-action, hence $H^s(BC_2^{\times(n-1)}; R(C_2))$ is also 2-torsion for $s \geq 1$ by the K\"unneth theorem.

    As in the proof of \cite[Thm.\ B.9]{mnn}, we obtain a short exact sequence
    \begin{equation}
        0 \to A^{s,t} \to E_2^{s,t} \to B^{s,t} \to 0
    \end{equation}
    where $A$ and $B$ are subquotients of the $\Ext$-groups corresponding to the left and middle term of the short exact sequence (\ref{eqkth:eq3}), but with $A$ shifted by 1 in the positive $s$-direction. In particular, $A$ and $B$ are 2-torsion for $s \geq 2$, and hence $E_2$ is 4-torsion for $s \geq 2$.
\end{Proof}
\begin{Proposition}
    For $G = C_2^{\times n}$, we have $\exp_{\sC} KU \geq n$.
    \label{prop:expkulowerbound}
\end{Proposition}
\begin{Proof}
   We consider the $\sC$-homotopy limit spectral sequence
\begin{equation}
    E_2^{s,t} = \sideset{}{^s}\lim_{\sO(C_2^{\times n})^\op_{\sC}}R(-)[\beta^{\pm}] \Rightarrow R(C_2^{\times n})[\beta^{\pm}].
\end{equation}
where $\beta$ is the Bott periodicity generator which has $(s,t)$-degree $(0,2)$.
First, by \cref{eqkth:lem1}, the $E_2$-page is qualitatively (made precise after the figure) given by
\begin{figure}[H]
	\centering
	\includegraphics[width=0.5\textwidth]{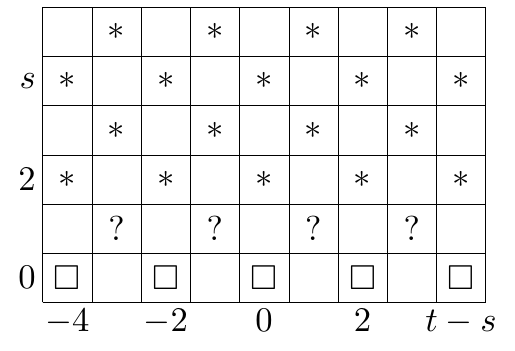}
	\caption{The $E_2$-page qualitatively.}
\end{figure}
Here qualitatively means that every $\square$ signifies a free $\Integers$-module (not necessarily of rank 1), every $\ast$ signifies a not necessarily free $\Integers/4\Integers$-module (not necessarily of rank 1), every empty square signifies the 0 group and every question mark $?$ signifies that it could be anything.

Second, by \cref{cor:expelab2gps}, we have that the edge map
\begin{equation}
    R(C_2^{\times n}) \xrightarrow{\widetilde{\res}} \sideset{}{^0}\lim_{\sO(C_2^{\times n})_\sC^{\mathrm{op}}} R(-)
\end{equation}
is injective and has cokernel with exponent $2^{n-1}$. That is, there is an $x$ in $E_2^{0,0}$ such that $2^{n-1}x \in \Im(\widetilde{\res})$ but $2^{n-2}x \notin \Im(\widetilde{\res})$.

Combining these two points, we see that $E_2^{0,0}$ will have to support at least differentials on $\ceil{\frac{n-1}{2}}$ pages, since all the groups in positive filtration $\geq 2$ are 4-torsion, and there is no possibility for a class in the $s=0$-line to hit a question mark in the $s=1$-line, so taking a kernel gives at least a subgroup of index 4. The even differentials are all 0, and so we get that $d_M \neq 0$ for some $M \geq 2\ceil{\frac{n-1}{2}} + 1 \geq n$, by the pigeonhole principle. Since $\exp_{\sC} KU \leq N$ implies that the spectral sequence degenerates at $E_{N+1}$ (\cite[Prop.\ 2.26]{mnn}), we get that $\exp_{\sC} KU \geq n$, as desired.
\end{Proof}
\begin{Remark}
    By distinguishing between $n$ even and odd, the proof in fact shows
    \begin{equation}
        \exp_{\sC} KU \geq
        \begin{cases}
            n + 1 & \text{ if $n \equiv 0 \pmod{2}$,} \\
            n & \text{ if $n \equiv 1 \pmod{2}$.} 
        \end{cases}
    \end{equation}
    \label{eqkth:rem1}
\end{Remark}
\subsection{Odd primes}
\label{subsec:cokoddprimes}
We do not have a generalization of the structure result of the cokernel like \cref{cor2:expelab2gps} to odd primes $p$ yet. However, experimental evidence does suggest a generalization. To describe it, we need the following combinatorial notion.
\begin{Definition}[{{\cite{andre1876}, \cite[Ex.\ I.16]{comtet}, \cite{fahssi12}}}]
    Define the \textbf{polynomial coefficient} $\binom{x,q}{k}$ by the coefficients of
	\begin{equation}
		(1+t+ \cdots + t^{q-1})^x = \sum_{k \geq 0} \binom{x,q}{k} t^k.
		\label{cokerrcpn:eq2}
	\end{equation}
	The coefficients for $q=2$ are also known as binomial coefficients. Polynomial coefficients for $q=3$ we call trinomial coefficients. For arbitrary $q$ we call them the $q$-nomial coefficients.
\end{Definition}
\begin{Warning}
    Polynomial coefficients are not the same thing as multinomial coefficients.
\end{Warning}
\begin{Conjecture}
    Let $p$ be a prime. Set
    \begin{equation}
        e_{n,k,p} = \sum_{j=0}^{p-1}\binom{n,p}{(p-1)(k+1) - j}
    \end{equation}
    (cf.\ for $p=3$ \cite{A104029}).
    Then the cokernel of the edge map $R(C_p^{\times n}) \to E_2^{0,0}$ of the $\sC$-homotopy limit spectral sequence converging to $\pi_*KU$ is isomorphic to
    \begin{equation}
        \bigoplus_{k=0}^{n-1} (\Integers/p^k)^{\oplus e_{n,k,p}}.
    \end{equation}
    \label{conj:exponents}
\end{Conjecture}
\begin{Remark}
    For $p=2$, \cref{conj:exponents} is \cref{cor2:expelab2gps}. 
\end{Remark}
Experimental evidence for \cref{conj:exponents} is given in \cref{ch:sage}.
\begin{Corollary}[{{of \cref{conj:exponents}}}]
    For $G = C_p^{\times n}$ the elementary abelian $p$-group of rank $n$, we have $\exp_{\sC} KU \geq n$.
\end{Corollary}
\begin{Proof}
    The same as the proof of \cref{prop:expkulowerbound}.
\end{Proof}

\subsection{The real case for $p=2$}
For $G$ an elementary abelian 2-group, the complex and real representation rings of $G$ coincide, which we use to give a lower bound on $\exp_{\sC} KO$ using real Bott periodicity. The argument is completely analogous to the complex case.
\begin{Proposition}
    For $G = C_2^{\times n}$ the $\sC$-exponent of $KO$ satisfies
    \begin{equation}
        \exp_{\sC} KO \geq n.
    \end{equation}
\end{Proposition}
\begin{Proof}
    As in the complex case, there need to be at least $\ceil{\frac{n-1}{2}}$ non-zero differentials, as again the $E_2$-page is $4$-torsion for $s \geq 2$. However, real Bott periodicity is 8-fold periodic with 4 non-zero groups in each period, which makes computing a lower bound on the $r$ for which $E_r = E_\infty$ a bit harder as it comes down to making a case distinction on the class of $n$ modulo 8 rather than $n$ modulo 2 as in the complex case. Nevertheless, one can do the counting, and the result is:
    \begin{equation}
        \exp_{\sC} KO \geq
        \begin{cases}
            2 + 2\ceil{\frac{n-1}{2}} & \text{ if $n \equiv 0,1 \pmod{8}$}, \\
            2 + 2(\ceil{\frac{n-1}{2}} - 1) + 1 & \text{ if $n \equiv 2,3 \pmod{8}$}, \\
            2 + 2(\ceil{\frac{n-1}{2}} - 2) + 3 & \text{ if $n \equiv 4,5 \pmod{8}$}, \\
            2 + 2(\ceil{\frac{n-1}{2}} - 3) + 7 &\text{ if $n \equiv 6,7 \pmod{8}$}.
        \end{cases}
    \end{equation}
    All these 4 expressions are $\geq n$.
\end{Proof}
\begin{Remark}
    The proof in fact shows
    \begin{equation}
        \exp_{\sC}KO \geq
        \begin{cases}
            n+2 & \text{ if $n \equiv 0 \pmod{8}$}, \\
            n+1 & \text{ if $n \equiv 1 \pmod{8}$}, \\
            n+1 & \text{ if $n \equiv 2 \pmod{8}$}, \\
            n & \text{ if $n \equiv 3 \pmod{8}$}, \\
            n+1 & \text{ if $n \equiv 4 \pmod{8}$}, \\
            n & \text{ if $n \equiv 5 \pmod{8}$}, \\
            n+3 & \text{ if $n \equiv 6 \pmod{8}$}, \\
            n+2 & \text{ if $n \equiv 7 \pmod{8}$}.
        \end{cases}
        \label{reqkth:eq1}
    \end{equation}
    The fact that $KU$ is a $KO$-module implies that $KU$ is a retract of $KO$ and therefore $\exp_{\sC} KU \leq \exp_{\sC}KO$. Hence the lower bounds on $\exp_{\sC} KU$ of \cref{eqkth:rem1} give lower bounds on $\exp_{\sC}KO$, but the ones in (\ref{reqkth:eq1}) are in all cases at least as good or better.
\end{Remark}

\appendix
\chapter{Sage computations}
\label{ch:sage}
\section{Introduction}
    Let $p$ be a prime, and consider the edge map $R(C_p^{\times n}) \to E_2^{0,0}$ of the $\sC$-homotopy limit spectral sequnce (cf.\ \cref{subsec:expcok} and \cref{subsec:cokoddprimes}). Denote the cokernel of this map by $Q_{p,n}$.
    Using Sage (\cite{sage}), we have computed the isomorphism type of $Q_{p,n}$ for the primes $p=2,3,5,7,11$ and various small values of $n$. The results are summarized in the following five tables. Of course, for $p=2$ we also have \cref{cor2:expelab2gps}. We also give the Sage code used for the computation. The computations provide experimental evidence for \cref{conj:exponents}.
\section{Summary of computations}
\subsection{$p=2$}
The following table gives $Q_{2,n}$ for $1 \leq n \leq 8$. It turns out that for those values of $n$, $Q_{2,n}$ is the direct sum of groups of the form $(\Integers/2^k \Integers)^{\oplus e_k}$,  for some exponents $e_k$. The entry in the $k$-th row and $n$-th column of the table gives the exponent $e_k$ of $\Integers/2^k \Integers$ in $Q_{2,n}$. An empty entry means that the exponent is zero.
\begin{table}[H]
  \centering
	\begin{tabular}{r|crrrrrrrr}
	  & $n$ & 1 & 2 & 3 & 4 & 5 & 6 & 7 & 8 \\
	  \hline
	  $\Integers/2^k$ & &  &  &  &  &  &  &  &  \\
	  $\Integers/1$ & & 1 & 2 & 3 & 4 & 5 & 6 & 7 & 8 \\
	  $\Integers/2$ & &   & 1 & 3 & 6 &10 &15 &21 &28 \\
	  $\Integers/4$ & &   &   & 1 & 4 &10 &20 &35 &56 \\
	  $\Integers/8$ & &   &   &   & 1 & 5 &15 &35 &70 \\
	  $\Integers/16$& &   &   &   &   & 1 & 6 &21 &56 \\
	  $\Integers/32$& &   &   &   &   &   & 1 & 7 &28 \\
	  $\Integers/64$& &   &   &   &   &   &   & 1 & 8 \\
	  $\Integers/128$&&   &   &   &   &   &   &   & 1
	\end{tabular}
	\caption{Exponent $e_k$ in decomposition of $Q_{2,n}$ as direct sum of groups of the form $(\Integers/2^k)^{\oplus e_k}$.}
\end{table}
For instance, from the table we read off  that
\begin{equation}
    Q_{2,4} \cong \Integers/8 \oplus (\Integers/4)^{\oplus 4} \oplus (\Integers/2)^{\oplus 6} \oplus (\Integers/1)^{\oplus 4}.
\end{equation}
The reason that we write down those trivial groups is that $Q_{2,3}$, for instance, is the quotient of a free abelian group of rank 7. Hence in the quotient ``there really are'' these trivial groups. Moreover, they fit nicely in the pattern of half of Pascal's triangle that this table has.
\subsection{$p=3$}
\label{coker:q3n}
Also for low values of $n$, $Q_{3,n}$ is isomorphic to a direct sum of groups of the form $(\Integers/3^k)^{\oplus e_k}$ for some exponents $e_k$. The following table gives these exponents for $1 \leq n \leq 6$ (we go less high because the algorithm takes too long to terminate for higher $n$).
\begin{table}[H]
  \centering
  \begin{tabular}{r|rrrrrrr}
                    & $n$   &   1 &   2 &   3 &   4 &   5 &   6 \\
		    \hline
    $\Integers/3^k$ &       &     &     &     &     &     &     \\
    $\Integers/1$   &       &   2 &   5 &   9 &  14 &  20 &  27 \\
    $\Integers/3$   &       &     &   3 &  13 &  35 &  75 & 140 \\
    $\Integers/9$   &       &     &     &   4 &  26 &  96 & 267 \\
    $\Integers/27$  &       &     &     &     &   5 &  45 & 216 \\
    $\Integers/81$  &       &     &     &     &     &   6 &  71 \\
    $\Integers/243$ &       &     &     &     &     &     &   7
  \end{tabular}
  \caption{Exponent $e_k$ in decomposition of $Q_{3,n}$ as direct sum of groups of the form $(\Integers/3^k)^{\oplus e_k}$.}
  \label{table:q3n}
\end{table}
For example,
\begin{equation}
    Q_{3,4} \cong (\Integers/27)^{\oplus 5} \oplus (\Integers/9)^{\oplus 26} \oplus (\Integers/3)^{\oplus 35} \oplus (\Integers/1)^{\oplus 14}
\end{equation}
The pattern in the table seems to be described by pairwise sums of trinomial coefficients, cf.\ \cite{A104029}.
\subsection{$p=5$}
We write down the same table as in the previous two sections for the prime 5:
\begin{table}[H]
  \centering
  \begin{tabular}{r|rrrrr}
     	              & $n$ &   1 &   2 &   3 &   4 \\
		      \hline
    $\Integers/5^k$   &     &     &     &     &     \\
    $\Integers/1$     &     &   4 &  14 &  34 &  69 \\
    $\Integers/5$     &     &     &  10 &  70 & 285 \\
    $\Integers/25$    &     &     &     &  20 & 235 \\
    $\Integers/125$   &     &     &     &     &  35
  \end{tabular}
  \caption{Exponent $e_k$ in decomposition of $Q_{5,n}$ as a direct sum of groups of the form $(\Integers/5^k)^{\oplus e_k}$.}
  \label{table:q5n}
\end{table}
The pattern of the table seems to be described by 4-term sums of 5-nomial coefficients.
\subsection{$p=7$}
For the prime 7 the table is
\begin{table}[H]
  \centering
  \begin{tabular}{r|rrrr}
    		        & $n$ &   1 &   2 &   3 \\
   \hline
   $ \Integers/7^k$     &     &     &     &     \\
   $ \Integers/1$       &     &   6 &  27 &  83 \\
   $ \Integers/7$       &     &     &  21 & 203 \\
   $ \Integers/49$      &     &     &     &  56
 \end{tabular}
   \caption{Exponent $e_k$ in decomposition of $Q_{7,n}$ as a direct sum of groups of the form $(\Integers/7^k)^{\oplus e_k}$.}
   \label{table:q7n}
\end{table}
The pattern of the table seems to be described by 6-term sums of 7-nomial coefficients.
\subsection{$p=11$}
For the prime 11 the table is
\begin{table}[H]
  \centering
  \begin{tabular}{r|rrr}
                        & $n$ &   1 &   2 \\
			\hline
     $ \Integers/11^k$  &     &     &     \\
     $ \Integers/1$     &     &  10 &  65 \\
     $ \Integers/11$    &     &     &  55 
   \end{tabular}
   \caption{Exponent $e_k$ in decomposition of $Q_{11,n}$ as direct sum of groups of the form $(\Integers/11^k)^{\oplus e_k}$.}
\end{table}
The pattern of the table seems to be described by 10-term sums of 11-nomial coefficients.
\section{Sage code}
The Sage code implementing the cokernel calculation is

\begin{lstlisting}
sage: def first_non_zero_entry(v):
...       r"""
...       Return the first non-zero entry of the row
...       vector `v`.
...       
...       The input is a row vector `v`. The output
...       is the first non-zero entry of `v`,
...       starting from the right, if it exists.
...       Otherwise the output is None.
...       """
...       for i in range(v.parent().rank()):
...           if v[i] != 0:
...               return v[i]
...       return None
sage: def next_p_vector(v,p):
...       r""" 
...       Return the next vector modulo p in lexicographic
...       ordering.
...       
...       Input is a prime `p` and a row vector `v` with entries
...       between 0 and p-1. 
...       Output is the next vector of this type in 
...       lexicographic ordering, or None if no such 
...       vector exists.
...       """
...       # Start at the right hand side of the row vector `v`.
...       n = v.parent().rank() - 1
...       # Move from the right to the left.
...       for i in range(v.parent().rank()):
...           # If we can increase the digit ...
...           if v[n-i] < p-1:
...               # ... do it, and return the resulting vector,
...               v[n-i] = v[n-i] + 1
...               return v
...           # otherwise the entry equals `p-1`, hence we need
...           # set the digit to zero and move one digit to
...           # the left
...           elif v[n-i] == p-1:
...               v[n-i] = 0
...       return None
sage: def matA(p,n):
...       r"""
...       Return a matrix with generators for all the
...       subgroups of `C_p^n` isomorphic to `C_p`.
...       
...       Input is a prime `p` and a natural number `n`. 
...       Output is a matrix `A` whose columns are 
...       precisely generators of all the subgroups 
...       of `C_p^n` isomorphic to `C_p`, namely the 
...       vectors that have the shape
...       `(0,0,...,1,*,*,...,*)^T,` where `*` means any 
...       entry and `^T` transpose.
...       """
...       # We first do some initialization.
...       # Make `v` the vector `(0,0,...,0,1)`.
...       v = (ZZ**n)(0)
...       v[n-1] = 1
...       # Make the matrix `A` the same thing,
...       # but then as a matrix:
...       A = MatrixSpace(ZZ,n,1)(0)
...       A[n-1,0] = 1
...       # Now `A` already has the current `v` as
...       # its first column, so we move on to the
...       # next `v`:
...       v = next_p_vector(v,p)
...       # While we haven't gone through all the
...       # vectors ...
...       while v != None:
...           # ... if the current vector is one of our
...           # chosen generators ...
...           if first_non_zero_entry(v) == 1:
...               # ... augment it to A,
...               A = A.augment(v)
...           # otherwise, skip this vector.
...           v = next_p_vector(v,p)
...       # Once we have augmented all the
...       # generators to `A`, return `A`.
...       return A
sage: def matB(p,n):
...       r"""
...       Return a matrix with columns the
...       non-identity elements of `C_p^n`.
...       
...       Input is a prime p and a natural number n. 
...       Output is a matrix whose columns are 
...       precisely all the non-identity elements of 
...       `C_p^n`.
...       """
...       # First we do some initialization
...       # Set `v` to be the first non-identity
...       # element in lexicographic ordering.
...       # `v = (0,0,...,0,1)`.
...       v = (ZZ**n)(0)
...       v[n-1] = 1
...       # Do the same for `B`, but then as a 
...       # matrix.
...       B = MatrixSpace(ZZ,n,1)(0)
...       B[n-1,0] = 1
...       # `B` now already has the current `v` as
...       # a column, so we move to the next `v`.
...       v = next_p_vector(v,p)
...       # While we haven't added all elements...
...       while v != None:
...           # add the current one.
...           B = B.augment(v)
...           v = next_p_vector(v,p)
...       return B
sage: def redmatmodp(Matr,p):
...       r"""
...       Return the integer matrix `Matr` with all
...       its entries reduced `\mod p`.
...       
...       Input is a matrix `Matr` with integer values
...       and a prime `p`. Output is the matrix `Matr` 
...       with integer entries, reduced `\mod p`,
...       hence the output matrix has entries between
...       `0` and `p-1`.
...       """
...       Matroutput = Matr.parent()(0)
...       for row in range(Matr.parent().nrows()):
...           for col in range(Matr.parent().ncols()):
...               Matroutput[row,col] = mod(Matr[row,col],p)
...       return Matroutput
sage: def matD(C,p,n):
...       r"""
...       Input is a prime p, a natural number n, 
...       and a matrix C with entries between 0 
...       and p-1. The idea is that all the 
...       entries between 1 and p-1 correspond 
...       to different basis elements of the 
...       representation ring. Therefore, each 
...       row of entries between 0 and p-1 will 
...       give p-1 rows in the output matrix. If 
...       the j-th entry of the row is a number 
...       1 <= k <= p-1, then this will give a 1 
...       in the k-th row and j-th column (of the 
...       p-1 columns corresponding to this 
...       column).
...       """
...       dim = p**n - 1
...       D = MatrixSpace(ZZ,dim,dim)(0)
...       for row in range((p**n - 1)/(p-1)):
...           for col in range(p**n-1):
...               if C[row,col] != 0:
...                   D[C[row,col] + (p-1)*row - 1,col] = 1
...       return D
sage: def count_list_elements(L):
...       r"""
...       Count the elements in a list.
...       
...       Input is a list. Output is a dictionary with
...       keys the elements of the list, values the
...       number of occurences.
...       """
...       counts = {}
...       for i in range(len(L)):
...           if L[i] in counts:
...               counts[L[i]] = counts[L[i]] + 1
...           else:
...               counts[L[i]] = 1
...       return counts
sage: def calculation_2(p,n):
...       """
...       For the n-th elementary abelian p-group, 
...       one can look at the orbit category 
...       given by the family of cyclic subgroups. 
...       The representation ring functor gives a 
...       contravariant functor on this category, 
...       and we can compute the limit. The 
...       representation ring of the total group 
...       maps injectively into this limit. This 
...       is in general not surjective, and this 
...       function computes the cokernel. The 
...       output is a dictionary, whose keys are 
...       moduli of the integers, and whose values 
...       are the exponents with which the 
...       integers modulus the moduli appear in 
...       the quotient. For example, 
...       calculation_1(2,3) yields 
...       {1: 3, 2: 3, 4: 1}, which means that 
...       the quotient is isomorphic to 
...       {e}^3 x C_2^3 x C_4.
...       """
...       A = matA(p,n)
...       print "Done A"
...       B = matB(p,n)
...       print "Done B"
...       C = A.transpose()*B
...       print "Done C"
...       C = redmatmodp(C,p)
...       print "Done reduction modulo p"
...       D = matD(C,p,n)
...       print "Done D"
...       # The elementary divisors of `D` 
...       # are the elements on the diagonal
...       # of a Smith Normal form of `D`.
...       # An elementary divisor `d`
...       # corresponds to an occurens of
...       # `\ZZ/d \ZZ` in the quotient.  
...       divisors = D.elementary_divisors()
...       print "Done Elementary Divisors"
...       counted = count_list_elements(divisors)
...       return counted
sage: calculation_2(2,3)
Done A
Done B
Done C
Done reduction modulo p
Done D
Done Elementary Divisors
{1: 3, 2: 3, 4: 1}
\end{lstlisting}
This last calculation shows, for example, that
\begin{equation}
    Q_{2,3} \cong (\Integers/1)^{\oplus 3} \oplus (\Integers/2)^{\oplus 3} \oplus \Integers/4.
\end{equation}

\bibliography{../bibliography}{}

\newcommand{\etalchar}[1]{$^{#1}$}
\begin{thebibliography}{CTVEZ03}

\bibitem[AE81]{alperinevens81}
J.~L. Alperin and L.~Evens.
\newblock Representations, resolutions and {Q}uillen's dimension theorem.
\newblock {\em J. Pure Appl. Algebra}, 22(1):1--9, 1981.

\bibitem[AM04]{adem}
Alejandro Adem and R.~James Milgram.
\newblock {\em Cohomology of finite groups}, volume 309 of {\em Grundlehren der
  Mathematischen Wissenschaften [Fundamental Principles of Mathematical
  Sciences]}.
\newblock Springer-Verlag, Berlin, second edition, 2004.

\bibitem[And76]{andre1876}
D{\'e}sir{\'e} Andr{\'e}.
\newblock M\'emoire sur les combinaisons r\'eguli\`eres et leurs applications.
\newblock {\em Ann. Sci. \'Ecole Norm. Sup. (2)}, 5:155--198, 1876.

\bibitem[Arm88]{armstrong}
M.~A. Armstrong.
\newblock {\em Groups and symmetry}.
\newblock Undergraduate Texts in Mathematics. Springer-Verlag, New York, 1988.

\bibitem[Ben93]{benson93}
D.~J. Benson.
\newblock The image of the transfer map.
\newblock {\em Arch. Math. (Basel)}, 61(1):7--11, 1993.

\bibitem[Ben04]{benson04}
Dave Benson.
\newblock Dickson invariants, regularity and computation in group cohomology.
\newblock {\em Illinois J. Math.}, 48(1):171--197, 2004.

\bibitem[Bir67]{birkhoff67}
Garrett Birkhoff.
\newblock {\em Lattice theory}.
\newblock Third edition. American Mathematical Society Colloquium Publications,
  Vol. XXV. American Mathematical Society, Providence, R.I., 1967.

\bibitem[BK72]{bousfieldkan}
A.~K. Bousfield and D.~M. Kan.
\newblock {\em Homotopy limits, completions and localizations}.
\newblock Lecture Notes in Mathematics, Vol. 304. Springer-Verlag, Berlin-New
  York, 1972.

\bibitem[Bor60]{borel60}
Armand Borel.
\newblock {\em Seminar on transformation groups}.
\newblock With contributions by G. Bredon, E. E. Floyd, D. Montgomery, R.
  Palais. Annals of Mathematics Studies, No. 46. Princeton University Press,
  Princeton, N.J., 1960.

\bibitem[Bur55]{burnside}
W.~Burnside.
\newblock {\em Theory of groups of finite order}.
\newblock Cambridge University Press, 1955.
\newblock 2d ed.

\bibitem[Car95]{carlson95}
Jon~F. Carlson.
\newblock Depth and transfer maps in the cohomology of groups.
\newblock {\em Math. Z.}, 218(3):461--468, 1995.

\bibitem[CE56]{cartaneilenberg}
Henri Cartan and Samuel Eilenberg.
\newblock {\em Homological Algebra}.
\newblock Princeton University Press, 1956.

\bibitem[Com74]{comtet}
Louis Comtet.
\newblock {\em Advanced Combinatorics}.
\newblock D.\ Reidel Publishing Company, Dordrecht-Holland/Boston-U.S.A, 1974.

\bibitem[CTVEZ03]{carlson2003cohomology}
Jon~F. Carlson, Lisa Townsley, Luis Valeri-Elizondo, and Mucheng Zhang.
\newblock {\em Cohomology rings of finite groups}, volume~3 of {\em Algebras
  and Applications}.
\newblock Kluwer Academic Publishers, Dordrecht, 2003.
\newblock With an appendix: Calculations of cohomology rings of groups of order
  dividing 64 by Carlson, Valeri-Elizondo and Zhang.

\bibitem[Duf81]{duflot81}
J.~Duflot.
\newblock Depth and equivariant cohomology.
\newblock {\em Comment. Math. Helv.}, 56(4):627--637, 1981.

\bibitem[Eve91]{evens91}
Leonard Evens.
\newblock {\em The cohomology of groups}.
\newblock Oxford Mathematical Monographs. The Clarendon Press, Oxford
  University Press, New York, 1991.
\newblock Oxford Science Publications.

\bibitem[{Fah}12]{fahssi12}
N.-E. {Fahssi}.
\newblock {Polynomial Triangles Revisited}.
\newblock {\em ArXiv e-prints}, February 2012.
\newblock \url{https://arxiv.org/abs/1202.0228}.

\bibitem[GK11]{green11}
David~J. Green and Simon~A. King.
\newblock The computation of the cohomology rings of all groups of order 128.
\newblock {\em J. Algebra}, 325:352--363, 2011.

\bibitem[GK15]{green15}
David~J. Green and Simon~A. King.
\newblock The cohomology of finite {$p$}-groups, 2015.
\newblock Website. \url{http://users.minet.uni-jena.de/cohomology/}.

\bibitem[Hus94]{husemoeller1994fibre}
Dale Husemoller.
\newblock {\em Fibre bundles}, volume~20 of {\em Graduate Texts in
  Mathematics}.
\newblock Springer-Verlag, New York, third edition, 1994.

\bibitem[Ill83]{illman83}
S\"oren Illman.
\newblock The equivariant triangulation theorem for actions of compact {L}ie
  groups.
\newblock {\em Math. Ann.}, 262(4):487--501, 1983.

\bibitem[Inc05]{A104029}
OEIS~Foundation Inc.
\newblock The {On-Line} {Encyclopedia} of {Integer} {Sequences}, 2005.
\newblock \url{http://oeis.org/A123456}.

\bibitem[LMSM86]{lms}
L.~G. Lewis, Jr., J.~P. May, M.~Steinberger, and J.~E. McClure.
\newblock {\em Equivariant stable homotopy theory}, volume 1213 of {\em Lecture
  Notes in Mathematics}.
\newblock Springer-Verlag, Berlin, 1986.
\newblock With contributions by J. E. McClure.

\bibitem[Lur09]{lurietop}
Jacob Lurie.
\newblock {\em Higher topos theory}, volume 170 of {\em Annals of Mathematics
  Studies}.
\newblock Princeton University Press, Princeton, NJ, 2009.

\bibitem[MM02]{mandellmay}
M.~A. Mandell and J.~P. May.
\newblock Equivariant orthogonal spectra and {$S$}-modules.
\newblock {\em Mem. Amer. Math. Soc.}, 159(755):x+108, 2002.

\bibitem[MNN15]{mnn}
Akhil Mathew, Niko Naumann, and Justin Noel.
\newblock {Derived induction and restriction theory}.
\newblock {\em ArXiv e-prints}, July 2015.
\newblock \url{http://arxiv.org/abs/1507.06867}.

\bibitem[MNN17]{mnnnd}
Akhil Mathew, Niko Naumann, and Justin Noel.
\newblock Nilpotence and descent in equivariant stable homotopy theory.
\newblock {\em Adv. Math.}, 305:994--1084, 2017.

\bibitem[PY03]{pakianathan03}
Jonathan Pakianathan and Erg\"un Yal\c{c}{\i}n.
\newblock On nilpotent ideals in the cohomology ring of a finite group.
\newblock {\em Topology}, 42(5):1155--1183, 2003.

\bibitem[Qui71a]{quillen71adams}
Daniel Quillen.
\newblock The {A}dams conjecture.
\newblock {\em Topology}, 10:67--80, 1971.

\bibitem[Qui71b]{quillen71mod2}
Daniel Quillen.
\newblock The mod 2 cohomology rings of extra-special 2-groups and the spinor
  groups.
\newblock {\em Mathematische Annalen}, 194:197--212, 1971.

\bibitem[Qui71c]{quillen71spectrum}
Daniel Quillen.
\newblock The spectrum of an equivariant cohomology ring. {I}, {II}.
\newblock {\em Ann. of Math. (2)}, 94:549--572; ibid. (2) 94 (1971), 573--602,
  1971.

\bibitem[Rus89]{rusin89}
David~J. Rusin.
\newblock The cohomology of the groups of order {$32$}.
\newblock {\em Math. Comp.}, 53(187):359--385, 1989.

\bibitem[S{\etalchar{+}}15]{sage}
W.\thinspace{}A. Stein et~al.
\newblock {\em {S}age {M}athematics {S}oftware}.
\newblock The Sage Development Team, 2015.
\newblock \url{http://www.sagemath.org}.

\bibitem[Ser77]{Serrelinrep}
Jean-Pierre Serre.
\newblock {\em Linear representations of finite groups}.
\newblock Springer-Verlag, New York-Heidelberg, 1977.
\newblock Translated from the second French edition by Leonard L. Scott,
  Graduate Texts in Mathematics, Vol. 42.

\bibitem[Sym91]{symonds91}
Peter Symonds.
\newblock The complexity of a module and elementary abelian subgroups: a
  geometric approach.
\newblock {\em Proc. Amer. Math. Soc.}, 113(1):27--29, 1991.

\bibitem[Tot14]{totaro14}
Burt Totaro.
\newblock {\em Group cohomology and algebraic cycles}, volume 204 of {\em
  Cambridge Tracts in Mathematics}.
\newblock Cambridge University Press, Cambridge, 2014.

\bibitem[Wil05]{wild}
Marcel Wild.
\newblock The groups of order sixteen made easy.
\newblock {\em Amer. Math. Monthly}, 112(1):20--31, 2005.

\end{thebibliography}
\bibliographystyle{alpha}
\end{document}